\documentclass[a4paper,12pt,reqno]{amsart}
\pdfoutput=1
\usepackage{amsfonts, amsmath, amssymb, amstext} 
\usepackage{tikz}
\usepackage[normalem]{ulem}
\usepackage{mathrsfs}
\usepackage{hyperref}
\usepackage{enumitem}
\usepackage{tikz-cd, float}
\usetikzlibrary{decorations.markings}

\usepackage{caption, subcaption}
\allowdisplaybreaks[2]

\numberwithin{equation}{section}

\theoremstyle{plain}
\newtheorem{thm}{Theorem}[section] 
\newtheorem{lem}[thm]{Lemma}
\newtheorem{prop}[thm]{Proposition}
\newtheorem{cor}[thm]{Corollary}

\theoremstyle{definition}

\newtheorem{exam}[thm]{Example}
\newtheorem*{hypoth}{Standing Hypotheses}

\theoremstyle{remark}
\newtheorem{rem}[thm]{Remark}

\newcommand{\norm}[1]{\left\|#1\right\|}
\newcommand{\cspan}{\overline{\mathrm{span}}}
\def\C{{\mathbb C}}
\def\N{{\mathbb N}}

\def\Z{{\mathbb Z}}
\def\T{{\mathbb T}}

\begin{document}

\title[Iterating the Cuntz--Nica--Pimsner construction]{Iterating the Cuntz--Nica--Pimsner construction for compactly aligned  product systems}

\author{James Fletcher}
\address{School of Mathematics and Statistics \\ Victoria University of Wellington \\ Wellington, New Zealand}
\email{james.fletcher@vuw.ac.nz}
\thanks{This research was supported by an Australian Government Research Training Program (RTP) Scholarship and by the Marsden grant 15-UOO-071 from the Royal Society of New Zealand.} 



\begin{abstract}  
We study how decompositions of a quasi-lattice ordered group $(G,P)$ relate to decompositions of the Nica--Toeplitz algebra and Cuntz--Nica--Pimsner algebra of a compactly aligned product system $\mathbf{X}$ over $P$. In particular, we are interested in the situation where $(G,P)$ may be realised as the semidirect product of quasi-lattice ordered groups. Our results generalise Deaconu's work on iterated Toeplitz and Cuntz--Pimsner algebras --- we show that the Nica--Toeplitz algebra and Cuntz--Nica--Pimsner algebra of a compactly aligned product system over $\N^k$ may be realised as $k$-times iterated Toeplitz and Cuntz--Pimsner algebras respectively. 
\end{abstract} 

\maketitle
\setcounter{tocdepth}{1} 
\tableofcontents

\section{Introduction}

In \cite{fowler}, Fowler introduced product systems of Hilbert bimodules as a generalisation of the continuous product systems of Hilbert spaces studied by Arveson \cite{arveson} and the discrete product systems studied by Dinh \cite{dinh}. Loosely speaking, a product system of Hilbert $A$-bimodules over a unital semigroup $P$ is a semigroup $\mathbf{X}=\bigsqcup_{p\in P} \mathbf{X}_p$, such that each $\mathbf{X}_p$ is a Hilbert $A$-bimodule, and the map $x\otimes_A y\mapsto xy$ extends to an isomorphism from $\mathbf{X}_p\otimes_A \mathbf{X}_q$ to $\mathbf{X}_{pq}$ for each $p,q\in P\setminus \{e\}$. 

Motivated by the work of Nica \cite{nica} and Laca and Raeburn on Toeplitz algebras associated to non-abelian groups \cite{laca--raeburn}, Fowler focused on product systems over quasi-lattice ordered groups that satisfied a condition he called compact alignment. Fowler then studied representations of such product systems satisfying a constraint called Nica covariance. Results of \cite{fowler} show how to associate to each compactly aligned product system $\mathbf{X}$  a $C^*$-algebra $\mathcal{NT}_\mathbf{X}$, generated by a universal Nica covariant representation, which we call the Nica--Toeplitz algebra of $\mathbf{X}$. In the same article, the faithful representations of Nica--Toeplitz algebras are characterised \cite[Theorem~7.2]{fowler}, generalising both Laca and Raeburn's uniqueness theorem for Toeplitz algebras of quasi-lattice ordered groups \cite[Theorem~3.7]{laca--raeburn}, and Fowler and Raeburn's uniqueness theorem for Toeplitz algebras of Hilbert bimodules \cite[Theorem~2.1]{fowler--raeburn}. Subsequently, it was shown in \cite[Theorem~3.2]{NTuniqueness}, that the hypothesis present in \cite[Theorem~7.2]{fowler} of each fibre being (left) nondegenerate is not necessary. 

Fowler also proposed a notion of Cuntz--Pimsner covariance for representations of a compactly aligned product system over a quasi-lattice ordered group $(G,P)$, and an associated universal $C^*$-algebra, denoted by $\mathcal{O}_\mathbf{X}$, which generalises the Cuntz--Pimsner algebra associated to a single Hilbert bimodules \cite{pimsner, katsura}. In general, Fowler's Cuntz--Pimsner algebra need not contain a faithful copy of $A$, and even when it does, examples in the appendix of \cite{raeburn--sims--yeend} show that a representation of $\mathcal{O}_X$ that is faithful on the copy of $A$ need not be faithful on the generalised fixed point algebra $\mathcal{O}_\mathbf{X}^\nu$ (where $\nu$ is the canonical gauge coaction of $G$ on $\mathcal{O}_X$). Furthermore, in contrast with the Cuntz--Pimsner algebras associated to Hilbert bimodules by both Pimsner \cite{pimsner} and Katsura \cite{katsura}, Fowler's Cuntz--Pimsner algebra need not be a quotient of $\mathcal{NT}_\mathbf{X}$. 

To overcome these issues, Sims and Yeend defined a $C^*$-algebra $\mathcal{NO}_{\mathbf{X}}$ generated by a universal Cuntz--Nica--Pimsner covariant representation of $\mathbf{X}$ \cite{sims--yeend}. They showed that $\mathcal{NO}_{\mathbf{X}}$ coincides with Katsura's Cuntz--Pimsner algebra when $P=\N$, and coincides with Fowler's $\mathcal{O}_{\mathbf{X}}$ whenever $A$ acts faithfully and compactly on each $\mathbf{X}_p$ and each pair of elements in $P$ has a common upper bound. Furthermore, Sims and Yeend showed that their Cuntz--Nica--Pimsner algebras generalise the Cuntz--Krieger algebras associated to finitely aligned higher-rank graphs \cite{raeburn--sims--yeend}, and the boundary quotients of Toeplitz algebras associated to quasi-lattice ordered groups studied by Crisp and Laca \cite{crisp--laca}. By viewing $\mathcal{NO}_\mathbf{X}$ as a co-universal algebra, Carlsen, Larsen, Sims, and Vittadello later proved a gauge-invariant uniqueness theorem for Cuntz--Nica--Pimsner algebras \cite{carlsen--larsen--sims--vittadello}.  

The object of this article is to study how the $C^*$-algebras $\mathcal{NT}_\mathbf{X}$ and $\mathcal{NO}_\mathbf{X}$ associated to a compactly aligned product system $\mathbf{X}$ over a quasi-lattice ordered group $(G,P)$ decompose given a decomposition of $(G,P)$. In particular, we show how a decomposition of the underlying quasi-lattice group as a semidirect product of quasi-lattice ordered groups gives a decomposition of the associated $C^*$-algebras. 

In \cite{deaconu}, Deaconu investigated what he called iterated Toeplitz and Cuntz--Pimsner algebras. Unfortunately, some of Deaconu's proofs lack detail, and it is not clear which of his various hypotheses are necessary to make the procedure work. Indeed, our original motivation for this article was to verify and better understand the arguments in \cite{deaconu}. The results in this article show that Deaconu's iterative procedure can be extended to quasi-lattice ordered groups that are more general than $(\Z^2,\N^2)$. Applying this iterative procedure to compactly aligned product systems over $\N^k$ immediately enables us to extend many of Katsura's results from \cite{katsura} to Nica--Toeplitz and Cuntz--Nica--Pimsner algebras. In the future, it would be interesting to see if this iterative procedure can be combined with the Pimsner--Voiculescu exact sequence \cite[Theorem~8.6]{katsura} to describe the $K$-theory of Cuntz--Nica--Pimsner algebras associated to product systems over $\N^k$ (we are particularly interested in the situation where the product system comes from a higher-rank graph).  

In Section~\ref{background} we review the necessary background material for product systems and their associated $C^*$-algebras. We begin by recapping the definitions of product systems, representations, compact alignment, and Nica covariance. Next we present the definition of Cuntz--Pimsner covariance developed by Sims and Yeend \cite{sims--yeend} and show how it relates to Fowler's notion of Cuntz--Pimsner covariance in \cite{fowler}. We also summarise the key results from the literature that we will make use of throughout the article. 

In Section~\ref{factorising NT alg} we examine the Nica--Toeplitz algebra associated to a compactly aligned product system $\mathbf{Z}$ over a quasi-lattice ordered group of the form $(G\rtimes_\alpha H, P\rtimes_\alpha Q)$, where $(G,P)$ and $(H,Q)$ are themselves quasi-lattice ordered groups.  If we let $\mathbf{X}$ be the product system corresponding to the fibres of $\mathbf{Z}$ associated to the semigroup $P$, then our main result (Theorem~\ref{isomorphisms of Nica--Toeplitz algebras}) shows that there exists a product system $\mathbf{Y}^\mathcal{NT}$ over $(H,Q)$, whose coefficient algebra is the Nica--Toeplitz algebra of $\mathbf{X}$, such that the Nica--Toeplitz algebras of $\mathbf{Y}^\mathcal{NT}$ and $\mathbf{Z}$ coincide. This result generalises \cite[Lemma~4.1]{deaconu}. 

In Section~\ref{factorising CNP alg}, we extend the results from Section~\ref{factorising NT alg} to Cuntz--Nica--Pimsner algebras associated to product systems. Our main result (Theorem~\ref{isomorphism of Cuntz Nica Pimsner algebras}) shows that there exists a product system $\mathbf{Y}^\mathcal{NO}$ over $(H,Q)$, whose coefficient algebra is $\mathcal{NO}_\mathbf{X}$, such that the Cuntz--Nica--Pimsner algebras of $\mathbf{Y}^\mathcal{NO}$ and $\mathbf{Z}$ coincide. This result generalises the second part of \cite[Lemma~4.2]{deaconu}. The main difficulty in establishing Theorem~\ref{isomorphism of Cuntz Nica Pimsner algebras} is finding sufficient conditions to ensure that $\mathcal{NO}_\mathbf{X}$ acts faithfully on each fibre of $\mathbf{Y}^\mathcal{NO}$, which makes the Cuntz--Pimsner covariance relation in $\mathcal{NO}_{\mathbf{Y}^\mathcal{NO}}$ tractable. We deal with this difficulty in Proposition~\ref{injectivity of induced actions}. At the end of Section~\ref{factorising CNP alg} we present three examples that illustrate our procedure and demonstrate the utility of the decompositions given by Theorem~\ref{isomorphisms of Nica--Toeplitz algebras} and Theorem~\ref{isomorphism of Cuntz Nica Pimsner algebras}. 

In Section~\ref{mixed cases} we finish by examining what we call relative Cuntz--Nica--Pimsner algebras. In particular, we consider the Cuntz--Nica--Pimsner algebra of the product system defined in Section~\ref{factorising NT alg}, as well as the Nica--Toeplitz algebra of the product system defined in Section~\ref{factorising CNP alg}. The main result of Section~\ref{mixed cases} (Theorem~\ref{isom of mixed case algebras}) generalises the first part of \cite[Lemma~4.2]{deaconu}. 

\section{Preliminaries}
\label{background}

The majority of background that we will require can be found in \cite{lance, fowler, sims--yeend}.

\subsection{Hilbert bimodules}

Let $A$ be a $C^*$-algebra. An inner-product $A$-module is a complex vector space $X$ equipped with a right action of $A$, and a map $\langle \cdot, \cdot \rangle_A:X\times X\rightarrow A$, (complex) linear in its second argument, satisfying the following conditions: 
\begin{enumerate}[label=\upshape(\roman*)]
\item $\langle x,y\rangle_A=\langle y,x\rangle_A^*$;
\item $\langle x,y\cdot a \rangle_A=\langle x,y\rangle_Aa$;
\item $\langle x,x\rangle_A\geq 0$ in $A$; and
\item $\langle x,x\rangle_A=0$ if and only if $x=0$;
\end{enumerate}
for any $x,y\in X$ and $a\in A$. The formula $\norm{x}_X:=\norm{\langle x,x\rangle_A}_A^{1/2}$ defines a norm on $X$ by \cite[Proposition~1.1]{lance}, and we say that $X$ is a Hilbert $A$-module if $X$ is complete with respect to this norm. 

Let $X$ be a Hilbert $A$-module. We say that a map $T:X\rightarrow X$ is adjointable if there exists a map $T^*:X\rightarrow X$ such that $\langle Tx, y\rangle_A=\langle x, T^*y\rangle_A$ for each $x,y\in X$. Adjointable operators are automatically linear and continuous, and have a unique adjoint. Equipping the collection of adjointable operators on $X$, which we denote by $\mathcal{L}_A(X)$, with the operator norm gives it the structure of a $C^*$-algebra. For each $x,y\in X$ there is an adjointable operator $\Theta_{x,y}$ defined by $\Theta_{x,y}(z):=x\cdot \langle y, z\rangle_A$. We call such operators generalised rank one operators. The closed subspace $\mathcal{K}_A(X):=\cspan\{\Theta_{x,y}:x,y\in X\}$, elements of which we call generalised compact operators, forms an essential ideal of $\mathcal{L}_A(X)$.

A Hilbert $B$--$A$-bimodule (also called a $B$--$A$ $C^*$-correspondence in the literature) is a Hilbert $A$-module $X$ equipped with a left action of $B$ by adjointable operators, i.e. there exists a homomorphism $\phi:B\rightarrow \mathcal{L}_A(X)$. To simplify notation, we will often write $b\cdot x$ for $\phi(b)(x)$. When $A=B$, we say that $X$ is a Hilbert $A$-bimodule. Since each $\phi(b)$ is by definition adjointable, and so $A$-linear, we have that $b\cdot (x\cdot a)=(b\cdot x)\cdot a$ for each $a\in A$, $b\in B$, and $x\in X$.

A particularly simple (and important) example occurs when $X=A$. Letting $A$ act on $X$ by left and right multiplication, and equipping $X$ with the $A$-valued inner product $\langle a,b\rangle_A:=a^*b$, gives a Hilbert $A$-bimodule, which we denote by ${}_A A_A$. The map $\Theta_{a,b}\mapsto ab^*$ extends to an isomorphism from $\mathcal{K}_A({}_A A_A)$ to $A$, whilst $\mathcal{L}_A({}_A A_A)$ is isomorphic to the multiplier algebra of $A$, which we denote by $\mathcal{M}(A)$. 

Every Hilbert $A$-module $X$ is nondegenerate in the sense that the span of the set $\{x\cdot a:x\in X, a\in A\}$ is dense in $X$. In particular, the Hewitt--Cohen--Blanchard factorisation theorem \cite[Proposition~2.31]{raeburn--williams} says that for each $x\in X$, there exists a unique $x'\in X$ such that $x=x'\cdot \langle x', x'\rangle_A$. In general, a Hilbert $A$-bimodule need not be left nondegenerate in the sense that $X=\cspan\{a\cdot x:x\in X,a\in A\}$. 

We form the balanced tensor product of two Hilbert $A$-bimodules $X$ and $Y$ as follows. We let $X\odot Y$ denote the algebraic tensor product of $X$ and $Y$ as complex vector spaces, and write $X\odot_A Y$ for the quotient by $\mathrm{span}\{x\cdot a \odot y- x\odot a\cdot y:x\in X, \, y\in Y, \, a\in A\}$ (we write $x\odot_A y$ for the coset containing $x\odot y$). The formula $\langle x\odot_A y, w\odot_A z\rangle_A:=\langle y, \langle x,w\rangle_A \cdot z\rangle_A$, determines a bounded $A$-valued sesquilinear form on $X\odot_A Y$. If we let $N$ be the subspace $\mathrm{span}\{n\in X\odot_A Y:\langle n,n\rangle_A=0\}$, then the formula  $\norm{z+N}:=\inf_{n\in N} \norm{\langle z+n,z+n\rangle_A}_A^{1/2}$ gives a norm on $(X\odot_A Y)/N$. We define the balanced tensor of $X$ and $Y$, which we denote by $X\otimes_A Y$, to be the completion of $(X\odot_A Y)/N$ with respect to this norm. 

Given two Hilbert $A$-bimodules $X$ and $Y$ and an adjointable operator $S\in \mathcal{L}_A(X)$, the formula $x\otimes_A y \mapsto (Sx)\otimes_A y$ extends to a map on all of $X\otimes_A Y$, which we denote by $S\otimes_A \mathrm{id}_Y$. It is straightforward to check that $S\otimes_A \mathrm{id}_Y$ is adjointable, with adjoint $S^*\otimes_A \mathrm{id}_Y$.

Given a Hilbert $B$--$A$-bimodule $X$, the theory of induced representations enables us to convert representations of $A$ into representations of $B$. If $\pi:A\rightarrow \mathcal{B}(\mathcal{H})$ is a nondegenerate representation of $A$ on a Hilbert space $\mathcal{H}$, then \cite[Proposition~2.66]{raeburn--williams} gives a representation $X\text{-}\mathrm{Ind}_A^B\pi:B\rightarrow \mathcal{B}\left(X\otimes_A \mathcal{H}\right)$ such that $\left(X\text{-}\mathrm{Ind}_A^B\pi\right)(b)(x\otimes_A h)=(b\cdot x)\otimes_A h$ for each $b\in B$, $x\in X$, and $h\in \mathcal{H}$. 

\subsection{Product systems of Hilbert bimodules and quasi-lattice ordered groups}

Let $P$ be a semigroup with identity $e$, and $A$ a $C^*$-algebra. A product system over $P$ with coefficient algebra $A$ is a semigroup $\mathbf{X}=\bigsqcup_{p\in P}\mathbf{X}_p$ such that:
\begin{enumerate}[label=\upshape(\roman*)]
\item for each $p\in P$, $\mathbf{X}_p\subseteq \mathbf{X}$ is a Hilbert $A$-bimodule;
\item $\mathbf{X}_e$ is equal to the Hilbert $A$-bimodule ${}_A A_A$;
\item for each $p,q\in P\setminus \{e\}$, there exists a Hilbert $A$-bimodule isomorphism $M_{p,q}^\mathbf{X}:\mathbf{X}_p\otimes_A \mathbf{X}_q\rightarrow \mathbf{X}_{pq}$ satisfying $M_{p,q}^\mathbf{X}(x\otimes_A y)=xy$ for each $x\in \mathbf{X}_p$ and $y\in \mathbf{X}_q$; and
\item[(iv)] multiplication in $\mathbf{X}$ by elements of $\mathbf{X}_e=A$ implements the left and right actions of $A$ on each $\mathbf{X}_p$; that is $xa=x\cdot a$ and $ax=a\cdot x$ for each $p\in P$, $a\in A$, and $x\in \mathbf{X}_p$. 
\end{enumerate}

We write $\phi_p:A\rightarrow \mathcal{L}_A(\mathbf{X}_p)$ for the homomorphism that implements the left action of $A$ on $\mathbf{X}_p$, i.e. $\phi_p(a)(x)=a\cdot x=ax$ for each $p\in P$, $a\in A$, and $x\in \mathbf{X}_p$. Multiplication in $\mathbf{X}$ is associative since $\mathbf{X}$ is a semigroup. Hence, $\phi_{pq}(a)(xy)=(\phi_p(a)x)y$ for all $p,q\in P$, $a\in A$, $x\in \mathbf{X}_p$, and $y\in \mathbf{X}_q$. We also write $\langle \cdot,\cdot\rangle_A^p$ for the $A$-valued inner-product on $\mathbf{X}_p$. 

By properties (ii) and (iv), for each $p\in P$ there exists $A$-linear inner-product preserving maps $M_{p,e}^\mathbf{X}:\mathbf{X}_p\otimes_A \mathbf{X}_e \rightarrow \mathbf{X}_p$ and $M_{e,p}^\mathbf{X}:\mathbf{X}_e\otimes_A \mathbf{X}_p\rightarrow \mathbf{X}_p$ such that $M_{p,e}^\mathbf{X}(x\otimes_A a)=xa=x\cdot a$ and $M_{e,p}^\mathbf{X}(a\otimes_A x)=ax=a\cdot x$ for each $a\in \mathbf{X}_e=A$ and $x\in \mathbf{X}_p$. Each $M_{p,e}^\mathbf{X}$ is surjective by the Hewitt--Cohen--Blanchard factorisation theorem, and so an $A$-bimodule isomorphism. On the other hand, the maps $M_{e,p}^\mathbf{X}$ need not be surjective, since we do not require that each $\mathbf{X}_p$ is (left) nondegenerate. 

For each $p\in P\setminus \{e\}$ and $q\in P$, we define a homomorphism $\iota_p^{pq}:\mathcal{L}_A\left(\mathbf{X}_p\right)\rightarrow \mathcal{L}_A\left(\mathbf{X}_{pq}\right)$ by 
\[
\iota_p^{pq}(S):=M_{p,q}^\mathbf{X}\circ (S\otimes_A \mathrm{id}_{\mathbf{X}_q})\circ (M_{p,q}^\mathbf{X})^{-1}
\]
for each $S\in \mathcal{L}_A\left(\mathbf{X}_p\right)$. Equivalently, $\iota_p^{pq}$ is characterised by the formula $\iota_p^{pq}(S)(xy)=(Sx)y$ for each $S\in \mathcal{L}_A\left(\mathbf{X}_p\right)$, $x\in \mathbf{X}_p$, $y\in \mathbf{X}_q$. 
Since $\mathbf{X}_e \otimes_A \mathbf{X}_q$ need not be isomorphic to $\mathbf{X}_q$, we cannot necessarily define a map from $\mathcal{L}_A\left(\mathbf{X}_e\right)$ to $\mathcal{L}_A\left(\mathbf{X}_q\right)$ using the above procedure. However, as $\mathcal{K}_A\left(\mathbf{X}_e\right)=\mathcal{K}_A\left({}_A A_A\right)\cong A$, we can define $\iota_e^q:\mathcal{K}_A\left(\mathbf{X}_e\right)\rightarrow \mathcal{L}_A\left(\mathbf{X}_q\right)$ by $\iota_e^q(a):=\phi_q(a)$. For notational purposes, we define $\iota_p^r: \mathcal{L}_A\left(\mathbf{X}_p\right)\rightarrow \mathcal{L}_A\left(\mathbf{X}_r\right)$ to be the zero map whenever $p,r\in P$ and $r\neq pq$ for all $q\in P$. It is also useful to know what happens when we compose these homomorphisms: routine calculations show that for any $p,q,r\in P$, we have $\iota_{pq}^{pqr}\circ \iota_p^{pq}=\iota_p^{pqr}$.

A quasi-lattice ordered group $(G,P)$ consists of a group $G$ and a subsemigroup $P$ of $G$ such that $P\cap P^{-1}=\{e\}$, and with respect to the partial order on $G$ induced by $p\leq q \Leftrightarrow p^{-1}q\in P$, any two elements $p,q\in G$ which have a common upper bound in $P$ have a least common upper bound in $P$. It is straightforward to show that if two elements in $G$ have a least common upper bound in $P$, then this least common upper bound is unique. If it exists, we write $p\vee q$ for the least common upper bound of $p,q\in G$. For $p,q\in G$, we write $p\vee q=\infty$ if $p$ and $q$ have no common upper bound in $P$, and $p\vee q<\infty$ otherwise. We say that $P$ is directed if $p\vee q<\infty$ for every $p,q\in P$. 

Let $(G,P)$ be a quasi-lattice ordered group and $\mathbf{X}$ a product system over $P$ with coefficient algebra $A$. We say that $\mathbf{X}$ is compactly aligned if $\iota_p^{p\vee q}(S)\iota_q^{p\vee q}(T)\in \mathcal{K}_A(\mathbf{X}_{p\vee q})$ whenever $S\in \mathcal{K}_A(\mathbf{X}_p)$ and $T\in \mathcal{K}_A(\mathbf{X}_q)$ for some $p,q\in P$ with $p\vee q<\infty$. It is important to note that this condition does not imply that either $\iota_p^{p\vee q}(S)$ or $\iota_q^{p\vee q}(T)$ is compact. 

\subsection{Representations of compactly aligned product systems, Nica covariance, and the Nica--Toeplitz algebra}

Let $(G,P)$ be a quasi-lattice ordered group and $\mathbf{X}$ a compactly aligned product system over $P$ with coefficient algebra $A$. A representation of $\mathbf{X}$ in a $C^*$-algebra $B$ is a map $\psi:\mathbf{X}\rightarrow B$ satisfying the following relations:
\begin{enumerate}
\item[(T1)] each $\psi_p:=\psi|_{\mathbf{X}_p}$ is a linear map, and $\psi_e$ is a homomorphism;
\item[(T2)] $\psi_p(x)\psi_q(y)=\psi_{pq}(xy)$ for all $p,q\in P$ and $x\in \mathbf{X}_p$, $y\in \mathbf{X}_q$; and
\item[(T3)] $\psi_p(x)^*\psi_p(y)=\psi_e(\langle x,y\rangle_A^p)$ for all $p\in P$ and $x,y\in \mathbf{X}_p$.
\end{enumerate}

Relations (T1) and (T3) show that a representation $\psi$ is always norm-decreasing, and isometric if and only $\psi_e$ is injective. Furthermore, by \cite[Proposition~8.11]{raeburn} there exists a homomorphism $\psi^{(p)}:\mathcal{K}_A\left(\mathbf{X}_p\right)\rightarrow B$ for each $p\in P$, such that $\psi^{(p)}\left(\Theta_{x,y}\right)=\psi_p(x)\psi_p(y)^*$ for all $x,y\in \mathbf{X}_p$. 

A representation $\psi:\mathbf{X}\rightarrow B$ is said to be Nica covariant if, for any $p,q\in P$ and $S\in \mathcal{K}_A(\mathbf{X}_p)$, $T\in \mathcal{K}_A(\mathbf{X}_q)$, we have
\begin{align*}
\psi^{(p)}(S)\psi^{(q)}(T)=
\begin{cases}
\psi^{(p\vee q)}\left(\iota_p^{p\vee q}(S)\iota_q^{p\vee q}(T)\right) & \text{if $p\vee q<\infty$}\\
0 & \text{otherwise.}
\end{cases}
\end{align*}
Making use of the Hewitt--Cohen--Blanchard factorisation theorem, it follows that 
$\psi_p(\mathbf{X}_p)^*\psi_q(\mathbf{X}_q)\subseteq\cspan\{\psi_{p^{-1}(p\vee q)}(\mathbf{X}_{p^{-1}(p\vee q)})\psi_{q^{-1}(p\vee q)}(\mathbf{X}_{q^{-1}(p\vee q)})^*\}$ if $p\vee q<\infty$, and is $\{0\}$ otherwise. 

Theorem~2.10 of \cite{loring} gives the existence of a $C^*$-algebra $\mathcal{NT}_\mathbf{X}$ and a Nica covariant representation $i_\mathbf{X}:\mathbf{X}\rightarrow \mathcal{NT}_\mathbf{X}$, that are universal in the following sense:
\begin{enumerate}[label=\upshape(\roman*)]
\item the image of $i_\mathbf{X}$ generates $\mathcal{NT}_\mathbf{X}$; and
\item given any other Nica covariant representation $\psi:\mathbf{X}\rightarrow B$, there exists a homomorphism $\psi_*:\mathcal{NT}_\mathbf{X}\rightarrow B$ such that $\psi_* \circ i_\mathbf{X}=\psi$. 
\end{enumerate}
We call $\mathcal{NT}_\mathbf{X}$ the Nica--Toeplitz algebra of $\mathbf{X}$. Since $i_\mathbf{X}$ generates $\mathcal{NT}_\mathbf{X}$, it follows that  $\mathcal{NT}_\mathbf{X}=\cspan\left\{i_\mathbf{X}(x)i_\mathbf{X}(y)^*:x,y\in \mathbf{X}\right\}$.

Theorem~3.2 of \cite{NTuniqueness} characterises the faithful representations of Nica--Toeplitz algebras. Suppose that $(G,P)$ is a quasi-lattice ordered group with $G$ amenable, $\mathbf{X}$ is a compactly aligned product system over $P$ with coefficient algebra $A$, and $\psi:\mathbf{X}\rightarrow \mathcal{B}(\mathcal{H})$ is a Nica covariant representation of $\mathbf{X}$ on a Hilbert space $\mathcal{H}$. Then the induced homomorphism $\psi_*:\mathcal{NT}_\mathbf{X}\rightarrow \mathcal{B}(\mathcal{H})$ is faithful provided for any finite set $K\subseteq P\setminus \{e\}$, the representation
\[
A\ni a \mapsto \psi_e(a)\prod_{t\in K}\big(\mathrm{id}_\mathcal{H}-\mathrm{proj}_{\overline{\psi_t(\mathbf{X}_t)\mathcal{H}}}\big)\in \mathcal{B}(\mathcal{H})
\]
is faithful.

Given a compactly aligned product system $\mathbf{X}$, we let $\mathcal{F}_\mathbf{X}:=\bigoplus_{p\in P} \mathbf{X}_p$ denote the space of sequences $(x_p)_{p\in P}$ such that $x_p\in \mathbf{X}_p$ for all $p\in P$ and $\sum_{p\in P}\langle x_p,x_p\rangle_A$ converges in $A$. Proposition~1.1 of  \cite{lance} shows that there exists an $A$-valued inner product on $\mathcal{F}_\mathbf{X}$ such that $\left\langle (x_p)_{p\in P}, (y_p)_{p\in P}\right\rangle_A=\sum_{p\in P}\langle x_p, y_p\rangle_A$, and that $\mathcal{F}_\mathbf{X}$ is complete with respect to the induced norm. Letting $A$ act pointwise from the left and right gives $\mathcal{F}_\mathbf{X}$ the structure of a Hilbert $A$-bimodule, which we call the Fock space of $\mathbf{X}$. Lemma~5.3 of \cite{fowler} then shows that there exists an isometric Nica covariant representation $l:\mathbf{X}\rightarrow \mathcal{L}_A(\mathcal{F}_\mathbf{X})$ such that $l_p(x)(y_q)_{q\in P}=\left(xy_q\right)_{q\in P}$ for each $p\in P$, $x\in \mathbf{X}_p$, and $(y_q)_{q\in P}\in \mathcal{F}_\mathbf{X}$. We call $l$ the Fock representation of $\mathbf{X}$.

\subsection{Cuntz--Pimsner covariance and the Cuntz--Nica--Pimsner algebra}

The notion of Cuntz--Pimsner covariance for representations of compactly aligned product systems introduced by Sims and Yeend in \cite{sims--yeend} is somewhat complicated and requires some preliminary definitions to formulate. We set $I_e:=A$ and
$I_p:=\bigcap_{e<q\leq p}\ker(\phi_q)$ for each $p\in P\setminus \{e\}$. For each $p\in P$ we define the Hilbert $A$-bimodule
\[
\widetilde{\mathbf{X}}_p:=\bigoplus_{q\leq p} \mathbf{X}_q\cdot I_{q^{-1}p}
,\]
where $ \mathbf{X}_q\cdot I_{q^{-1}p}=\{x\cdot a:x\in \mathbf{X}_q, \, a\in I_{p^{-1}q}\}$. 
For each $p\in P$, we write $\widetilde{\phi}_p:A\rightarrow \mathcal{L}_A\big(\widetilde{\mathbf{X}}_p\big)$ for the homomorphism defined by 
\[
\big(\widetilde{\phi}_p(a)(x)\big)_q:=\phi_q(a)(x_q)
\]
for each $a\in A$, $x\in \widetilde{\mathbf{X}}_p$, and $q\leq p$. 

Whilst the collection of bimodules $\big\{\widetilde{\mathbf{X}}_p:p\in P\big\}$ resembles a product system, we warn that in general $\widetilde{\mathbf{X}}_p\otimes_A \widetilde{\mathbf{X}}_q$ need not be isomorphic to
$\widetilde{\mathbf{X}}_{pq}$. However, for each $p,q\in P$ with $p\neq e$, we still have a homomorphism $\widetilde{\iota}_p^{\,q}:\mathcal{L}_A\big(\widetilde{\mathbf{X}}_p\big)\rightarrow \mathcal{L}_A\big(\widetilde{\mathbf{X}}_q\big)$ characterised by the formula
\[
\big(\widetilde{\iota}_p^{\,q}(S)(x)\big)_r=\iota_p^r(S)(x_r)
\]
for each $x\in \widetilde{\mathbf{X}}_p$ and $r\leq q$. Additionally, after identifying $\mathcal{K}_A\big(\widetilde{\mathbf{X}}_e\big)$ with $A$, we define $\widetilde{\iota}_e^{\,q}:\mathcal{K}_A\big(\widetilde{\mathbf{X}}_e\big)\rightarrow \mathcal{L}_A\big(\widetilde{\mathbf{X}}_q\big)$ to be $\widetilde{\phi}_q$ for each $q\in P$. 

To formulate the Cuntz--Pimsner covariance condition, we require another definition. Given a quasi-lattice ordered group $(G,P)$, we say that predicate statement $\mathcal{P}(s)$ (where $s\in P$) is true for large $s$ if, given any $p\in P$, there exists $q\geq p$, such that $\mathcal{P}(s)$ is true whenever $s\geq q$.

Finally, we are ready to present the definition of Cuntz--Pimsner covariance originally formulated by Sims and Yeend \cite[Definition~3.9]{sims--yeend}. We give a definition only in the situation that all of the homomorphisms $\widetilde{\phi}_p:A\rightarrow \mathcal{L}_A\big(\widetilde{\mathbf{X}}_p\big)$ are injective. We say that a representation $\psi:\mathbf{X}\rightarrow B$ is Cuntz--Pimsner covariant if, for any finite set $F\subseteq P$ and any choice of compact operators $\left\{T_p\in \mathcal{K}_A\left(\mathbf{X}_p\right):p\in F\right\}$, we have that 
\[
\sum_{p\in F} \widetilde{\iota}_p^{\,s} (T_p)=0\in \mathcal{L}_A\big(\widetilde{\mathbf{X}}_s\big) \quad \text{for large $s$} \quad \Rightarrow \quad \sum_{p\in F} \psi^{(p)}(T_p)=0.
\] 
We say that a representation is Cuntz--Nica--Pimsner covariant if it is both Nica covariant and Cuntz--Pimsner covariant. 

In Proposition~3.12 of \cite{sims--yeend} it is shown that for each compactly aligned product system $\mathbf{X}$ (with each $\widetilde{\phi}_p:A\rightarrow \mathcal{L}_A\big(\widetilde{\mathbf{X}}_p\big)$ injective), there exists a $C^*$-algebra $\mathcal{NO}_\mathbf{X}$, which we call the Cuntz--Nica--Pimsner algebra of $\mathbf{X}$, and a Cuntz--Nica--Pimsner covariant representation $j_\mathbf{X}:\mathbf{X}\rightarrow \mathcal{NO}_\mathbf{X}$, that are universal in the following sense
\begin{enumerate}[label=\upshape(\roman*)]
\item $\mathcal{NO}_\mathbf{X}$ is generated by the image of $j_\mathbf{X}$;
\item if $\psi:\mathbf{X}\rightarrow B$ is any other Cuntz--Nica--Pimsner covariant representation of $\mathbf{X}$, then there exists a homomorphism $\Pi\psi:\mathcal{NO}_\mathbf{X}\rightarrow B$ such that $\Pi\psi \circ j_\mathbf{X}=\psi$. 
\end{enumerate}
It follows that $\mathcal{NO}_\mathbf{X}$ is a quotient of $\mathcal{NT}_\mathbf{X}$. We write $q_\mathbf{X}:\mathcal{NT}_\mathbf{X}\rightarrow \mathcal{NO}_\mathbf{X}$ for the quotient homomorphism, which is characterised by $q_\mathbf{X}\circ i_\mathbf{X}=j_\mathbf{X}$. In \cite[Theorem~4.1]{sims--yeend}, it is shown that the universal Cuntz--Nica--Pimsner covariant representation $j_\mathbf{X}$ is always isometric. 

Lemma~3.15 of \cite{sims--yeend} shows that the requirement that each of the homomorphisms $\widetilde{\phi}_p$ is injective is automatic, provided either
\begin{enumerate}[label=\upshape(\roman*)]
\item each $\phi_p$ is injective; or
\item every nonempty bounded subset of $P$ has a maximal element in the following sense:
\begin{equation}
\label{bounded subsets of P with maximal element}
\text{\parbox{.85\textwidth}{
If $S\subseteq P$ is nonempty and there exists $q\in P$ such that $p\leq q$ for all $p\in S$, then there exists $p\in S$ such that $p\not\leq p'$ for all $p'\in S\setminus \{p\}$.
}}
\end{equation}
\end{enumerate}
In the situation where each $\phi_p$ is injective, the Hilbert $A$-bimodules $\widetilde{\mathbf{X}}_p$ and $\mathbf{X}_p$ are isomorphic. Moreover, this isomorphism intertwines $\phi_p$ and $\widetilde{\phi}_p$, as well as $\iota_p^r$ and $\widetilde{\iota}_p^{\,r}$. In this situation a  representation $\psi:\mathbf{X}\rightarrow B$ is Cuntz--Pimsner covariant provided, whenever $F\subseteq P$ is finite and $\left\{T_p\in \mathcal{K}_A\left(\mathbf{X}_p\right):p\in F\right\}$ is such that $\sum_{p\in F} \iota_p^s (T_p)=0\in \mathcal{L}_A\left(\mathbf{X}_s\right)$ for large $s$, we have $\sum_{p\in F} \psi^{(p)}(T_p)=0$.

In \cite{fowler}, Fowler defined a representation $\psi$ of a product system $\mathbf{X}$ over a semigroup $P$ (with each $\phi_p$ injective) to be Cuntz--Pimsner covariant if, for every $p\in P$, the Toeplitz representation $(\psi_p,\psi_e)$ of the Hilbert $A$-bimodule $\mathbf{X}_p$ is Cuntz--Pimsner covariant (in the sense that for each $a\in \phi_p^{-1}(\mathcal{K}_A(\mathbf{X}_p))$, we have $\psi^{(p)}(\phi_p(a))=\psi_e(a)$). Proposition~5.1 of \cite{sims--yeend} shows the relationship between Fowler's notion of Cuntz--Pimsner covariance and that defined by Sims and Yeend. If $\mathbf{X}$ is a compactly aligned product system with each $\phi_p$ injective, $P$ is directed, and $\psi:\mathbf{X}\rightarrow B$ is a representation, then
\begin{enumerate}[label=\upshape(\roman*)]
\item
If $\psi$ is Cuntz--Pimsner covariant in the sense of \cite[Definition~3.9]{sims--yeend}, then $(\psi_p,\psi_e)^{(1)}(\phi_p(a))=\psi_e(a)$ for $a\in \phi_p^{-1}(\mathcal{K}_A(\mathbf{X}_p))$ and $p\in P$;
\item
If $\phi_p(A)\subseteq \mathcal{K}_A(\mathbf{X}_p)$ and $\psi^{(p)}(\phi_p(a))=\psi_e(a)$ for each $a\in A$, $p\in P$, then $\psi$ is Cuntz--Pimsner covariant in the sense of \cite[Definition~3.9]{sims--yeend}. 
\end{enumerate}

For any discrete group $G$, the universal property of the group $C^*$-algebra $C^*(G)$ induces a homomorphism $\delta_G:C^*(G)\rightarrow C^*(G)\otimes C^*(G)$ such that $\delta_G(i_G(g))=i_G(g)\otimes i_G(g)$ for each $g\in G$ (we use an unadorned $\otimes$ to denote the minimal tensor product of $C^*$-algebras). Recall that a (full) coaction of $G$ on a $C^*$-algebra $A$, is an injective homomorphism $\delta:A\rightarrow A\otimes C^*(G)$, that satisfies the coaction identity $(\delta\otimes \mathrm{id}_{C^*(G)})\circ \delta=(\mathrm{id}_A\otimes \delta_G)\circ \delta$, and is nondegenerate in the sense that $A\otimes C^*(G)=\cspan\{\delta(A)(1_{\mathcal{M}(A)}\otimes C^*(G))\}$ (where $1_{\mathcal{M}(A)}$ is the identity of the multiplier algebra of $A$). For those readers interested in learning more about coactions, we suggest \cite[Appendix~A]{echterhoff--kaliszewski-quigg--raeburn}. 

As noted in \cite[Remark~4.5]{carlsen--larsen--sims--vittadello} there exists a gauge coaction $\nu_\mathbf{X}:\mathcal{NO}_\mathbf{X}\rightarrow \mathcal{NO}_\mathbf{X}\otimes C^*(G)$ such that $\nu_\mathbf{X}(j_{\mathbf{X}_p}(x))=j_{\mathbf{X}_p}(x)\otimes i_G(p)$ for each $p\in P$ and $x\in \mathbf{X}_p$.

Carlsen, Larsen, Sims, and Vittadello proved the following gauge-invariant uniqueness theorem for Cuntz--Nica--Pimsner algebras \cite[Corollary~4.12]{carlsen--larsen--sims--vittadello}. If $(G,P)$ is a quasi-lattice ordered group with $G$ amenable, and $\mathbf{X}$ is a compactly aligned product system over $P$ with either each $\phi_p$ injective, or each $\widetilde{\phi}_p$ injective and $P$ directed, then a surjective homomorphism $\phi:\mathcal{NO}_\mathbf{X}\rightarrow B$ is injective if and only if
\begin{enumerate}[label=\upshape(\roman*)]
\item $\phi|_{j_\mathbf{X}(A)}$ is injective; and
\item there exists a coaction $\beta:B\rightarrow B\otimes C^*(G)$ such that $\beta\circ \phi=\left(\phi \otimes \mathrm{id}_{C^*(G)}\right)\circ \nu_\mathbf{X}$. 
\end{enumerate}

\section{Combining quasi-lattice ordered groups}
\label{Combining quasi-lattice ordered groups}

The next result shows that the direct product of quasi-lattice ordered groups is quasi-lattice ordered. Furthermore, it provides sufficient conditions for a semidirect product of quasi-lattice ordered groups to be quasi-lattice ordered. 

\begin{lem}
\label{semidirect product of QLOG}
Let $(G,P)$ and $(H,Q)$ be quasi-lattice ordered groups. If $\alpha: H \rightarrow \mathrm{Aut}(G)$ is a group homomorphism with $\alpha_H(P)\subseteq P$, then the semidirect product $(G \rtimes_\alpha H, P\rtimes_\alpha Q)$ is a quasi-lattice ordered group.
\begin{proof}
Since $P$ and $Q$ are subsemigroups of $G$ and $H$ respectively, and $\alpha_Q(P)\subseteq P$, we see that $P\rtimes_\alpha Q$ is a subsemigroup of $G \rtimes_\alpha H$.

Next, we show that $(P\rtimes_\alpha Q) \cap (P\rtimes_\alpha Q)^{-1}=\{e_{G\rtimes_\alpha H}\}$. If $(g,h)\in  (P\rtimes_\alpha Q) \cap (P\rtimes_\alpha Q)^{-1}$, then 
$
(g,h)^{-1}=(\alpha_{h^{-1}}(g^{-1}), h^{-1})\in P\rtimes_\alpha Q.
$
Thus, $h, h^{-1}\in Q$, which forces $h=e_H$ since $(H,Q)$ is quasi-lattice ordered. Hence, $g, \alpha_{h^{-1}}(g^{-1})=g^{-1}\in P$, which forces $g=e_G$ since $(G,P)$ is also quasi-lattice ordered. Therefore, $(g,h)=(e_G,e_H)=e_{G\rtimes_\alpha H}$. 

Finally, we show that the order on the semidirect product $(G \rtimes_\alpha H, P\rtimes_\alpha Q)$ is the product order, i.e. $(g,h)\leq (g',h')$ in $(G \rtimes_\alpha H, P\rtimes_\alpha Q)$ if and only if $g\leq g'$ in $(G,P)$ and $h\leq h'$ in $(H,Q)$. Let $g,g'\in G$ and $h,h'\in H$. Suppose that $(g,h)\leq (g',h')$ in $(G \rtimes_\alpha H, P\rtimes_\alpha Q)$. Then
$
(g,h)^{-1}(g',h')=(\alpha_{h^{-1}}(g^{-1}g'),h^{-1}h')\in P\rtimes_\alpha Q.
$
Therefore, $h^{-1}h'\in Q$ (so that $h\leq h'$) and 
$
g^{-1}g'=\alpha_h(\alpha_{h^{-1}}(g^{-1}g'))\in \alpha_Q(P)\subseteq P
$
(so that $g\leq g'$). Conversely, if $g\leq g'$ and $h\leq h'$, then $h^{-1}h'\in Q$ and $\alpha_{h^{-1}}(g^{-1}g')\in \alpha_H(P)\subseteq P$, and so $(g,h)^{-1}(g',h')\in P\rtimes_\alpha Q$. Thus, $(g,h)\leq (g',h')$. We conclude that $(G \rtimes_\alpha H, P\rtimes_\alpha Q)$ is a quasi-lattice ordered group and for any $(g,h), (g',h')\in P\rtimes_\alpha Q$, we have
\begin{align*}
(g,h)\vee (g',h')&=
\begin{cases}
(g\vee g', h\vee h') & \text{if $g\vee g'< \infty$ and $h\vee h'< \infty$}\\
\infty & \text{otherwise.} 
\end{cases}
\qedhere
\end{align*}
\end{proof}
\end{lem}

Whilst the conditions in Lemma~\ref{semidirect product of QLOG} are sufficient for a semidirect product of quasi-lattice ordered groups to be quasi-lattice ordered, Proposition~2.2 of \cite{laca--raeburn2} shows that they are not necessary.

\section{Iterating the Nica--Toeplitz construction}
\label{factorising NT alg}

In this section we investigate product systems over semidirect products of quasi-lattice ordered groups of the sort appearing in Lemma~\ref{semidirect product of QLOG}. More precisely, given a product system $\mathbf{Z}$ (with coefficient algebra $A$) over a quasi-lattice ordered group of the form $(G\rtimes_\alpha H, P\rtimes_\alpha Q)$, with $(G,P)$ and $(H,Q)$ quasi-lattice ordered groups and $\alpha_H(P)\subseteq P$, we will show that there exists a product system $\mathbf{X}$ (also with coefficient algebra $A$) over $(G,P)$ sitting inside $\mathbf{Z}$, and a product system $\mathbf{Y}^\mathcal{NT}$ over $(H,Q)$ with coefficient algebra $\mathcal{NT}_\mathbf{X}$, such that the Nica--Toeplitz algebras of $\mathbf{Z}$ and $\mathbf{Y}^\mathcal{NT}$ are isomorphic. In passing from the product system $\mathbf{Z}$ to the product system $\mathbf{Y}^\mathcal{NT}$, we have in a sense decreased the size of the product system at the expense of increasing the size of the coefficient algebra, without losing any $C^*$-algebraic information in the process.

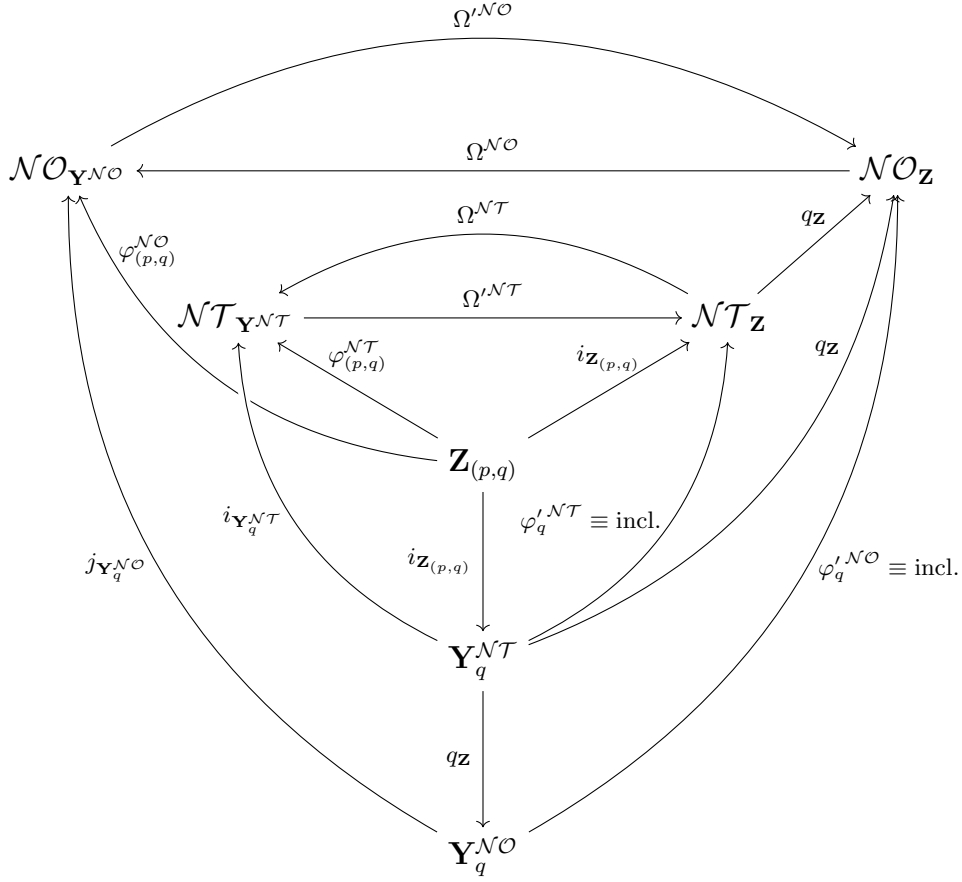
\begin{figure}[h]
\centering
\begin{tikzpicture}[scale=1.3]
\node (A) at (3,0) {$\mathbf{Y}_q^\mathcal{NT}$};
\node (B) at (3,2) {$\mathbf{Z}_{(p,q)}$};
\node (C) at (0.5,3.5) {$\mathcal{NT}_{\mathbf{Y}^\mathcal{NT}}$};
\node (D) at (5.5,3.5) {$\mathcal{NT}_\mathbf{Z}$};
\node (E) at (7.2,5) {$\mathcal{NO}_\mathbf{Z}$};
\node (F) at (3,-2) {$\mathbf{Y}_q^\mathcal{NO}$};
\node (G) at (-1.2,5) {$\mathcal{NO}_{\mathbf{Y}^\mathcal{NO}}$};
\path[->,font=\scriptsize]
(B) edge [bend left] node[right, very near end]{$\varphi_{(p,q)}^\mathcal{NO}$} (G)
(A) edge [bend left, line width=4pt,draw=white](C)
(A) edge [bend left] node[left]{$i_{\mathbf{Y}_q^\mathcal{NT}}$} (C)
(A) edge [bend right] node[left]{${\varphi_q'}^\mathcal{NT}\equiv \mathrm{incl.}$} (D)
(B) edge node[left]{$i_{\mathbf{Z}_{(p,q)}}$} (A)
(B) edge node[above=0.1cm]{$\varphi_{(p,q)}^\mathcal{NT}$} (C)
(C) edge node[above]{${\Omega'}^\mathcal{NT}$} (D)
(D) edge [bend right] node[above]{$\Omega^\mathcal{NT}$} (C)
(D) edge node[above=0.1cm]{$q_\mathbf{Z}$} (E)
(A) edge node[left]{$q_\mathbf{Z}$} (F)
(A) edge [bend right] node[left, near end]{$q_\mathbf{Z}$} (E)
(F) edge [bend right] node[right]{${\varphi_q'}^\mathcal{NO}\equiv \mathrm{incl.}$} (E)
(F) edge [bend left] node[left]{$j_{\mathbf{Y}_q^\mathcal{NO}}$} (G)
(G) edge [bend left] node[above]{${\Omega'}^\mathcal{NO}$} (E)
(E) edge node[above]{$\Omega^\mathcal{NO}$} (G)
(B) edge node[above=0.1cm]{$i_{\mathbf{Z}_{(p,q)}}$} (D);
\end{tikzpicture}
\caption{The homomorphisms $\Omega^\mathcal{NT}$,  $\Omega'^{\mathcal{NT}}$ and $\Omega^\mathcal{NO}$,  $\Omega'^{\mathcal{NO}}$}
\label{NO isomorphism}
\end{figure}

\begin{figure}[h]
\centering
\begin{tikzpicture}[scale=1.35]
\node (A) at (0,0) {$\mathcal{NT}_\mathbf{X}=\mathbf{Y}_{e_H}^\mathcal{NT}$};
\node (B) at (0,2) {$\mathbf{X}_p$};
\node (C) at (2,2) {$\mathcal{NT}_\mathbf{Z}$};
\node (D) at (0,-2) {$\mathcal{NO}_\mathbf{X}=\mathbf{Y}_{e_H}^\mathcal{NO}$};
\node (E) at (4,2) {$\mathcal{NO}_\mathbf{Z}$};
\path[->,font=\scriptsize]
(B) edge node[left]{$i_{\mathbf{X}_p}$} (A)
(B) edge node[above]{$i_{\mathbf{Z}_{\left(p,e_H\right)}}$} (C)
(A) edge node[right]{$\phi_\mathbf{X}^\mathcal{NT}$} (C)
(A) edge node[left]{$q_\mathbf{X}$} (D)
(C) edge node[above]{$q_\mathbf{Z}$} (E)
(D) edge node[right]{$\phi_\mathbf{X}^\mathcal{NO}$} (E)
(B) edge [bend right=45] node[left]{$j_{\mathbf{X}_p}$} (D)
(B) edge [bend left=45] node[above]{$j_{\mathbf{Z}_{\left(p,e_H\right)}}$} (E);
\end{tikzpicture}
\caption{The homomorphisms $\phi_{\mathbf{X}}^\mathcal{NT}$ and $\phi_{\mathbf{X}}^\mathcal{NO}$}
\label{NO homomorphism}
\end{figure}
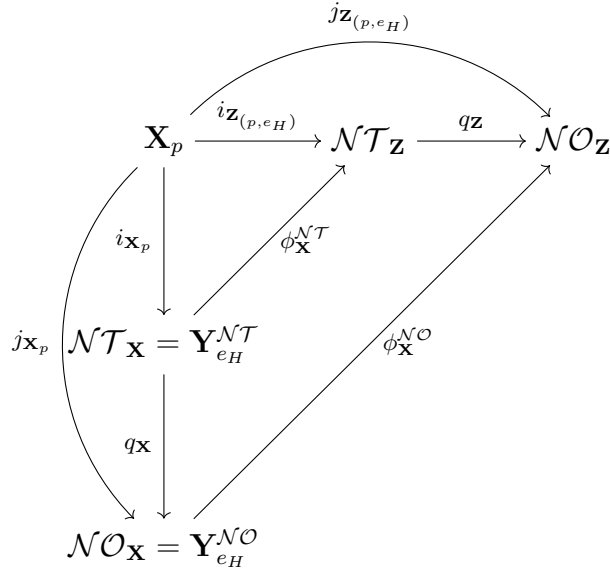

To help readers keep track of everything that is going on, we first provide a brief overview of the key results that we will prove and a pair of commutative diagrams in Figures~\ref{NO isomorphism} and~\ref{NO homomorphism} that summarise the various spaces, and the maps between them, that we are going to be working with. We hope that if readers happen to lose their way in Sections~\ref{Constructing the product system} and \ref{isomorphisms of Nica Toeplitz algebras} they will be able to return to this overview and the two diagrams for assistance. 

\begin{enumerate}
\item[(1)] In Proposition~\ref{inclusion of NT_X in NT_Z} we define a product system $\mathbf{X}\subseteq \mathbf{Z}$ over $(G,P)$, and show that the inclusion of $\mathbf{X}$ in $\mathbf{Z}$ induces a homomorphism $\phi_{\mathbf{X}}^\mathcal{NT}$ from $\mathcal{NT}_\mathbf{X}$ to $\mathcal{NT}_\mathbf{Z}$ such that $\phi_{\mathbf{X}}^\mathcal{NT}\circ i_\mathbf{X}=i_\mathbf{Z}$.
\item[(2)] In Proposition~\ref{injectivity of inclusion of NT_X in NT_Z} we argue that the homomorphism $\phi_{\mathbf{X}}^\mathcal{NT}$ is injective. 
\item[(3)] In Proposition~\ref{definition of Y_q}, we use $\phi_{\mathbf{X}}^\mathcal{NT}$ to construct a collection $\big\{\mathbf{Y}_q^\mathcal{NT}:q\in Q\big\}$ of Hilbert $\mathcal{NT}_\mathbf{X}$-modules inside $\mathcal{NT}_\mathbf{Z}$. 
\item[(4)] In Proposition~\ref{right actions on each Y_q}, we use $\phi_{\mathbf{X}}^\mathcal{NT}$ to show that each $\mathbf{Y}_q^\mathcal{NT}$ carries a left action of $\mathcal{NT}_\mathbf{X}$ by adjointable operators. 
\item[(5)] In Propositions~\ref{Y is a product system} and~\ref{Y is compactly aligned}, we show that $\mathbf{Y}^\mathcal{NT}:=\bigsqcup_{q\in Q} \mathbf{Y}_q^\mathcal{NT}$ is a compactly aligned product system over $(H,Q)$ with coefficient algebra $\mathcal{NT}_\mathbf{X}$. 
\end{enumerate}

Once we have the product system $\mathbf{Y}^\mathcal{NT}$, we prove that $\mathcal{NT}_{\mathbf{Y}^\mathcal{NT}}\cong \mathcal{NT}_\mathbf{Z}$.

\begin{enumerate}
\item[(6)] In Proposition~\ref{varphi is a representation} we use the universal Nica covariant representations of $\mathbf{Z}$ and $\mathbf{Y}^\mathcal{NT}$ to construct a representation $\varphi^{\mathcal{NT}}$ of $\mathbf{Z}$ in $\mathcal{NT}_{\mathbf{Y}^\mathcal{NT}}$.
\item[(7)] In Proposition~\ref{varphi is Nica covariant}, we prove that $\varphi^{\mathcal{NT}}$ is Nica covariant, and hence induces a homomorphism $\Omega^\mathcal{NT}:\mathcal{NT}_\mathbf{Z}\rightarrow \mathcal{NT}_{\mathbf{Y}^\mathcal{NT}}$ such that $\Omega^\mathcal{NT}\circ i_{\mathbf{Z}_{(p,q)}}=i_{\mathbf{Y}_q^\mathcal{NT}}$. 
\item[(8)] In Proposition~\ref{varphi' is Nica covariant} we show that the inclusion of $\mathbf{Y}^\mathcal{NT}$ in $\mathcal{NT}_\mathbf{Z}$ is a Nica covariant representation, and hence induces a homomorphism $\Omega'^{\mathcal{NT}}: \mathcal{NT}_{\mathbf{Y}^\mathcal{NT}}\rightarrow \mathcal{NT}_\mathbf{Z}$ such that $\Omega'^{\mathcal{NT}}\circ i_{\mathbf{Y}_q^\mathcal{NT}}$ is the inclusion of $\mathbf{Y}_q^\mathcal{NT}$ in $\mathcal{NT}_\mathbf{Z}$.
\item[(9)] In Theorem~\ref{isomorphisms of Nica--Toeplitz algebras} we prove that $\Omega^\mathcal{NT}$ and $\Omega'^{\mathcal{NT}}$ are mutually inverse isomorphisms. 
\end{enumerate}

In summary, we will show that for every $p\in P$ and $q\in Q$, the maps in the interior parts of Figures~\ref{NO isomorphism} and~\ref{NO homomorphism} exist and make the diagrams commutative (the exterior parts of the diagrams will be used in the setup for Section~\ref{factorising CNP alg}).

\subsection{Constructing the product systems}
\label{Constructing the product system}

\begin{hypoth}
We will assume that $(G,P)$ and $(H,Q)$ are quasi-lattice ordered groups and $\alpha: H \rightarrow \mathrm{Aut}(G)$ is a homomorphism such that $\alpha_H(P)\subseteq P$. By Lemma~\ref{semidirect product of QLOG}, $(G \rtimes_\alpha H, P\rtimes_\alpha Q)$ is quasi-lattice ordered. Also, $\mathbf{Z}$ will be a compactly aligned product system over $(G \rtimes_\alpha H, P\rtimes_\alpha Q)$ with coefficient algebra $A$. 
\end{hypoth}

Since $G\cong G\rtimes_\alpha \{e_H\}$ is a subgroup of $G\rtimes_\alpha H$ and $P\cong P\rtimes_\alpha \{e_H\}$ is a subsemigroup of $P\rtimes_\alpha Q$, we obtain a product system over the quasi-lattice ordered group $(G,P)$ by considering just those fibres of $\mathbf{Z}$ coming from  $P\rtimes_\alpha \{e_H\}$.

\begin{prop}
\label{inclusion of NT_X in NT_Z}
For each $p\in P$, let $\mathbf{X}_p:=\mathbf{Z}_{(p,e_H)}$. Then $\mathbf{X}:=\bigsqcup_{p\in P}\mathbf{X}_p$ is a compactly aligned product system over $(G,P)$ with coefficient algebra $A$. The inclusion of $\mathbf{X}$ in $\mathbf{Z}$ induces a homomorphism $\phi_{\mathbf{X}}^\mathcal{NT}:\mathcal{NT}_\mathbf{X}\rightarrow \mathcal{NT}_\mathbf{Z}$ such that $\phi_{\mathbf{X}}^\mathcal{NT}\circ i_\mathbf{X}=i_\mathbf{Z}$. 
\begin{proof}
For $p\in P$, define $\psi_p^{\mathcal{NT}}: \mathbf{X}_p \rightarrow \mathcal{NT}_\mathbf{Z}$ by $\psi_p^{\mathcal{NT}}:=i_{\mathbf{Z}_{(p,e_H)}}$. We claim that $\psi^{\mathcal{NT}}$ is a Nica covariant representation of $\mathbf{X}$ in $\mathcal{NT}_\mathbf{Z}$. Since $i_\mathbf{Z}$ is a representation, we know that each $\psi_p^{\mathcal{NT}}$ is linear and $\psi_{e_G}^{\mathcal{NT}}$ is a homomorphism. If $x\in \mathbf{X}_p$ and $z\in \mathbf{X}_r$, then 
\begin{align*}
\psi_p^{\mathcal{NT}}(x)\psi_r^{\mathcal{NT}}(z)
=i_{\mathbf{Z}_{\left(p,e_H\right)}}(x)i_{\mathbf{Z}_{\left(r,e_H\right)}}(z)
=i_{\mathbf{Z}_{\left(pr,e_H\right)}}(xz)=\psi_{pr}^{\mathcal{NT}}(xz).
\end{align*}
For $x,z\in \mathbf{X}_p$, we have
\begin{align*}
\psi_p^{\mathcal{NT}}(x)^*\psi_p^{\mathcal{NT}}(z)
&=i_{\mathbf{Z}_{\left(p,e_H\right)}}(x)^*i_{\mathbf{Z}_{\left(p,e_H\right)}}(z)\\
&=i_{\mathbf{Z}_{\left(e_G,e_H\right)}}\big(\langle x, z\rangle_A^{(p,e_H)}\big)
=\psi_{e_G}^{\mathcal{NT}}\big(\langle x, z\rangle_A^p\big).
\end{align*}
Thus, $\psi^{\mathcal{NT}}$ is a representation. It remains to check that $\psi^{\mathcal{NT}}$ is Nica covariant. Fix $S\in \mathcal{K}_A(\mathbf{X}_p)$ and $T\in \mathcal{K}_A(\mathbf{X}_r)$. Since ${\psi^{\mathcal{NT}}}^{(t)}=i_\mathbf{Z}^{((t,e_H))}$ for any $t\in P$, we see that
\begin{align*}
{\psi^{\mathcal{NT}}}^{(p)}(S){\psi^{\mathcal{NT}}}^{(r)}(T)&=i_\mathbf{Z}^{((p,e_H))}(S)i_\mathbf{Z}^{((r,e_H))}(T).
\end{align*}
If $(p, e_H)\vee (r,e_H)=\infty$ (which is precisely when $p\vee r =\infty$), Nica covariance of $i_\mathbf{Z}$ tells us that ${\psi^{\mathcal{NT}}}^{(p)}(S){\psi^{\mathcal{NT}}}^{(r)}(T)=0$. On the other hand, when $p\vee r <\infty$, we have $(p, e_H)\vee (r,e_H)=(p\vee r,e_H)<\infty$, and so
\begin{align*}
{\psi^{\mathcal{NT}}}^{(p)}(S){\psi^{\mathcal{NT}}}^{(r)}(T)
&=
i_\mathbf{Z}^{\left((p\vee r,e_H)\right)}\Big(\iota_{(p,e_H)}^{(p\vee r,e_H)}(S)\iota_{(r,e_H)}^{(p\vee r,e_H)}(T)\Big)\\
&=
{\psi^{\mathcal{NT}}}^{(p\vee r)}\big(\iota_p^{p\vee r}(S)\iota_r^{p\vee r}(T)\big).
\end{align*}
Thus,  $\psi^{\mathcal{NT}}$ is Nica covariant, and there exists a homomorphism $\phi_{\mathbf{X}}^\mathcal{NT}:\mathcal{NT}_\mathbf{X}\rightarrow \mathcal{NT}_\mathbf{Z}$ such that $\phi_{\mathbf{X}}^\mathcal{NT}(i_\mathbf{X}(x))=\psi^{\mathcal{NT}}(x)=i_\mathbf{Z}(x)$ for $x\in \mathbf{X}$. 
\end{proof}
\end{prop}

When $G$ is amenable the homomorphism $\phi_{\mathbf{X}}^\mathcal{NT}$ is injective. 

\begin{prop}
\label{injectivity of inclusion of NT_X in NT_Z}
Suppose that $G$ is an amenable group. Then $\phi_{\mathbf{X}}^\mathcal{NT}:\mathcal{NT}_\mathbf{X}\rightarrow \mathcal{NT}_\mathbf{Z}$ is injective.
\begin{proof}
Fix a faithful nondegenerate representation $\rho:A\rightarrow \mathcal{B}(\mathcal{H})$ of $A$ on a Hilbert space $\mathcal{H}$. Let $l_*:\mathcal{NT}_\mathbf{Z}\rightarrow \mathcal{L}_A(\mathcal{F}_\mathbf{Z})$ denote the Fock representation of $\mathcal{NT}_\mathbf{Z}$. To show that $\phi_{\mathbf{X}}^\mathcal{NT}$ is faithful, it suffices to show that the representation
\[
\mathcal{F}_\mathbf{Z}\text{-}\mathrm{Ind}_A^{\mathcal{NT}_\mathbf{X}}\rho=\left(l_*\circ \phi_{\mathbf{X}}^\mathcal{NT}\right)\otimes_A \mathrm{id}_\mathcal{H}:\mathcal{NT}_\mathbf{X}\rightarrow \mathcal{B}(\mathcal{F}_\mathbf{Z}\otimes_A \mathcal{H})
\]
is faithful. Denote by $\varrho:=\left(l_*\circ \phi_{\mathbf{X}}^\mathcal{NT}\circ i_\mathbf{X}\right)\otimes_A \mathrm{id}_\mathcal{H}$ the Nica covariant representation of $\mathbf{X}$ that induces the homomorphism $\left(l_*\circ \phi_{\mathbf{X}}^\mathcal{NT}\right)\otimes_A \mathrm{id}_\mathcal{H}$. Our aim is to prove that
\begin{enumerate}[label=\upshape(\roman*)]
\item
for each $p\in P\setminus \{e_G\}$,
$
A\otimes_A \mathcal{H}=\mathbf{X}_{e_G}\otimes_A \mathcal{H}\subseteq \left(\varrho_p(\mathbf{X}_p)(\mathcal{F}_\mathbf{Z}\otimes_A \mathcal{H})\right)^\perp;
$
\item
$A$ acts faithfully (via $\varrho_{e_G}$) on $A\otimes_A \mathcal{H}$.
\end{enumerate}
To see why this suffices, suppose that (i) and (ii) hold. For each $t\in P\setminus \{e_G\}$ let $P_t^\varrho:=\mathrm{proj}_{\overline{\varrho_t(\mathbf{X}_t)\left(\mathcal{F}_\mathbf{Z}\otimes_A \mathcal{H}\right)}}$. Then the representation 
\[
A\ni a \mapsto \varrho_{e_G}(a)\prod_{t\in K}\left(1-P_t^\varrho\right)\in \mathcal{B}\left(\mathcal{F}_\mathbf{Z}\otimes_A \mathcal{H}\right)
\]
is faithful for each finite subset $K\subseteq P\setminus \{e_G\}$. Since $G$ is amenable, \cite[Theorem~3.2]{NTuniqueness}  implies that $\left(l_*\circ \phi_{\mathbf{X}}^\mathcal{NT}\right)\otimes_A \mathrm{id}_\mathcal{H}$ is faithful as required. We now prove (i) and (ii).

If $p\in P$, then
\begin{align*}
\varrho_p(\mathbf{X}_p)(\mathcal{F}_\mathbf{Z}\otimes_A \mathcal{H})
&=\bigoplus_{(s,t)\in P\rtimes_\alpha Q}\mathbf{Z}_{\left(p,e_H\right)(s,t)}\otimes_A \mathcal{H}\\
&=\bigoplus_{(s,t)\in P\rtimes_\alpha Q:p\leq s}\mathbf{Z}_{(s,t)}\otimes_A \mathcal{H}.
\end{align*}
Now suppose that $p\in P\setminus \left\{e_G\right\}$. We suppose that $p\leq e_G$ and derive a contradiction. Then $p^{-1}=p^{-1}e_G\in P$, which forces $p=e_G$, since $P\cap P^{-1}=\{e_G\}$. Thus, $p\not\leq e_G$.  Hence for any $a\in \mathbf{X}_{e_G}=A$, $z\in \mathbf{Z}_{(s,t)}$ with $p\leq s$, and $h,g\in \mathcal{H}$, since $A=\mathbf{X}_{e_G}$ and $\mathbf{Z}_{(s,t)}$ are orthogonal in $\mathcal{F}_\mathbf{X}$, we see that
$
\left\langle a\otimes_A h, z\otimes_A g\right\rangle_\C
=\left\langle h, \langle a, z\rangle_A \cdot g\right\rangle_\C
=0.
$
Since inner-products are linear and continuous, we conclude that 
$
A\otimes_A \mathcal{H}\subseteq \left(\varrho_p(\mathbf{X}_p)(\mathcal{F}_\mathbf{Z}\otimes_A \mathcal{H})\right)^\perp,
$
which proves (i). 

It remains to show that $A$ acts faithfully on $A\otimes_A \mathcal{H}$. This follows from \cite[Corollary~2.74]{raeburn--williams} since $A$ acts faithfully on itself and the representation $\rho$ is faithful. 
\end{proof} 
\end{prop}

Using the injective homomorphism $\phi_{\mathbf{X}}^\mathcal{NT}:\mathcal{NT}_\mathbf{X}\rightarrow \mathcal{NT}_\mathbf{Z}$, we can construct a Hilbert $\mathcal{NT}_\mathbf{X}$-module $\mathbf{Y}_q^\mathcal{NT}$ for each $q\in Q$.

\begin{prop}
\label{definition of Y_q}
Suppose that $G$ is an amenable group so that the homomorphism $\phi_{\mathbf{X}}^\mathcal{NT}$ is injective by Proposition~\ref{injectivity of inclusion of NT_X in NT_Z}. For each $q\in Q\setminus \{e_H\}$, define
\begin{align*}
\mathbf{Y}_q^\mathcal{NT}:&=\cspan\big\{i_{\mathbf{Z}_{(e_G,q)}}(x)\phi_{\mathbf{X}}^\mathcal{NT}(b):x\in \mathbf{Z}_{(e_G,q)}, \ b\in \mathcal{NT}_X\big\}\subseteq \mathcal{NT}_\mathbf{Z}
\end{align*}
and let $\mathbf{Y}_{e_H}^\mathcal{NT}:={}_{\mathcal{NT}_\mathbf{X}}\left(\mathcal{NT}_\mathbf{X}\right)_{\mathcal{NT}_\mathbf{X}}$
Then for $q\in Q\setminus \{e_H\}$, $\mathbf{Y}_q^\mathcal{NT}$ carries a right action of $\mathcal{NT}_\mathbf{X}$ given by 
$
y\cdot b:=
y\phi_{\mathbf{X}}^\mathcal{NT}(b)
$
for each $y\in \mathbf{Y}_q^\mathcal{NT}$ and $b\in \mathcal{NT}_\mathbf{X}$.
There is an $\mathcal{NT}_\mathbf{X}$-valued inner-product on $\mathbf{Y}_q^\mathcal{NT}$ such that
$
\phi_\mathbf{X}^\mathcal{NT}\big(\langle y,w\rangle_{\mathcal{NT}_\mathbf{X}}^q\big)=y^*w
$
. With this structure, each $\mathbf{Y}_q^\mathcal{NT}$ is a Hilbert $\mathcal{NT}_\mathbf{X}$-module. 
\begin{proof}
Since multiplication in $\mathcal{NT}_\mathbf{Z}$ is continuous, the right action of $\mathcal{NT}_\mathbf{X}$ on $\mathbf{Y}_q^\mathcal{NT}$ is well-defined. Next, we check that for $y, w\in \mathbf{Y}_q^\mathcal{NT}$ we have $y^*w\in \phi_{\mathbf{X}}^\mathcal{NT}(\mathcal{NT}_X)$. Since $\phi_{\mathbf{X}}^\mathcal{NT}$ is a homomorphism, it suffices to check the case where $y=i_{\mathbf{Z}_{\left(e_G,q\right)}}(x)\phi_{\mathbf{X}}^\mathcal{NT}(b)$ and $w=i_{\mathbf{Z}_{\left(e_G,q\right)}}(z)\phi_{\mathbf{X}}^\mathcal{NT}(c)$ for some $x, z\in \mathbf{Z}_{\left(e_G,q\right)}$, and $b, c\in \mathcal{NT}_\mathbf{X}$. Since $i_\mathbf{Z}$ is a representation, we see that
\begin{align*}
y^*w
&=\phi_{\mathbf{X}}^\mathcal{NT}(b)^*i_{\mathbf{Z}_{\left(e_G,q\right)}}(x)^*i_{\mathbf{Z}_{\left(e_G,q\right)}}(z)\phi_{\mathbf{X}}^\mathcal{NT}(c)\\
&=\phi_{\mathbf{X}}^\mathcal{NT}(b^*)i_{\mathbf{Z}_{\left(e_G,e_H\right)}}\big(\langle x, z\rangle_A^{\left(e_G,q\right)}\big)\phi_{\mathbf{X}}^\mathcal{NT}(c)\\
&=\phi_{\mathbf{X}}^\mathcal{NT}\big(b^* i_{\mathbf{X}_{e_G}}\big(\langle x,z\rangle_A^{\left(e_G,q\right)}\big)c\big)
\in \phi_{\mathbf{X}}^\mathcal{NT}\left(\mathcal{NT}_\mathbf{X}\right).
\end{align*}
Thus, since $\phi_\mathbf{X}^\mathcal{NT}$ is injective, we may define $\langle \cdot, \cdot\rangle_{\mathcal{NT}_\mathbf{X}}^q: \mathbf{Y}_q^\mathcal{NT}\times  \mathbf{Y}_q^\mathcal{NT}\rightarrow \mathcal{NT}_\mathbf{X}$ by $\langle y,w\rangle_{\mathcal{NT}_\mathbf{X}}^q :=\left({\phi_{\mathbf{X}}^\mathcal{NT}}\right)^{-1}(y^*w)$ for each $y, w\in \mathbf{Y}_q^\mathcal{NT}$. Clearly $\langle \cdot, \cdot\rangle_{\mathcal{NT}_\mathbf{X}}^q$ is complex linear in its second argument. If $y, w\in \mathbf{Y}_q^\mathcal{NT}$ and $b\in \mathcal{NT}_\mathbf{X}$, then
\begin{align*}
\langle y,w\cdot b\rangle_{\mathcal{NT}_\mathbf{X}}^q
=\left({\phi_{\mathbf{X}}^\mathcal{NT}}\right)^{-1}\left(y^*w\phi_{\mathbf{X}}^\mathcal{NT}(b)\right)
=\left(\phi_{\mathbf{X}}^\mathcal{NT}\right)^{-1}\left(y^*w\right)b
=\langle y,w\rangle_{\mathcal{NT}_\mathbf{X}}^q b,
\end{align*}
and
\begin{align*}
\langle y,w\rangle_{\mathcal{NT}_\mathbf{X}}^q=\left({\phi_{\mathbf{X}}^\mathcal{NT}}\right)^{-1}\left(y^*w\right)=\left({\phi_{\mathbf{X}}^\mathcal{NT}}\right)^{-1}(w^*y)^*=\big(\langle w,y\rangle_{\mathcal{NT}_\mathbf{X}}^q\big)^*.
\end{align*}
Also, if $y\in \mathbf{Y}_q^\mathcal{NT}$, then
\[
\langle y,y\rangle_{\mathcal{NT}_\mathbf{X}}^q=\left({\phi_{\mathbf{X}}^\mathcal{NT}}\right)^{-1}(y^*y)\geq 0
\]
since $y^*y\geq 0$ in $\mathcal{NT}_\mathbf{Z}$. Moreover, if $\langle y,y\rangle_{\mathcal{NT}_\mathbf{X}}^q=0$, then $y^*y=0$ because $\phi_\mathbf{X}^\mathcal{NT}$ is injective, which forces $y=0$. Lastly, we show that the norm induced by $\langle \cdot, \cdot\rangle_{\mathcal{NT}_\mathbf{X}}^q$ is just the norm on $\mathcal{NT}_\mathbf{Z}$ restricted to $\mathbf{Y}_q^\mathcal{NT}$. Since $\phi_{\mathbf{X}}^\mathcal{NT}$ is isometric, we see that for any $y\in \mathbf{Y}_q^\mathcal{NT}$, 
\begin{align*}
\norm{y}_{\mathbf{Y}_q^\mathcal{NT}}^2
:=\big\|\langle y,y\rangle_{\mathcal{NT}_\mathbf{X}}^q\big\|_{\mathcal{NT}_\mathbf{X}}
&=\big\|\big(\phi_{\mathbf{X}}^\mathcal{NT}\big)^{-1}(y^*y)\big\|_{\mathcal{NT}_\mathbf{X}}\\
&=\norm{y^*y}_{\mathcal{NT}_\mathbf{Z}}
=\norm{y}^2_{\mathcal{NT}_\mathbf{Z}}.
\end{align*}
As $\mathbf{Y}_q^\mathcal{NT}$ is closed in $\mathcal{NT}_\mathbf{Z}$, we see that $\mathbf{Y}_q^\mathcal{NT}$ is complete with respect to the norm induced by $\langle \cdot, \cdot\rangle_{\mathcal{NT}_\mathbf{X}}^q$. We conclude that $\mathbf{Y}_q^\mathcal{NT}$ is a Hilbert $\mathcal{NT}_\mathbf{X}$-module.
\end{proof}
\end{prop}

\begin{rem}
\label{alternative description of Y_q}
Since $\alpha_H(P)\subseteq P$, we have that $\alpha_q(P)=P$ for each $q\in Q$. Hence, for any $q\in Q\setminus \{e_H\}$, it can be shown that
\begin{align*}
\mathbf{Y}_q^\mathcal{NT}
=\cspan\big\{i_{\mathbf{Z}_{\left(p,q\right)}}\big(\mathbf{Z}_{\left( p,q\right)}\big)i_{\mathbf{Z}_{\left(r,e_H\right)}}\big(\mathbf{Z}_{\left(r,e_H\right)}\big)^*:p,r\in P\big\}.
\end{align*}
Furthermore,
\begin{align*}
\mathbf{Y}_{e_H}^\mathcal{NT}\cong \phi_{\mathbf{X}}^\mathcal{NT}\left(\mathcal{NT}_\mathbf{X}\right)
=\cspan\big\{i_{\mathbf{Z}_{\left(p,e_H\right)}}\big(\mathbf{Z}_{\left( p,e_H\right)}\big)i_{\mathbf{Z}_{\left(r,e_H\right)}}\big(\mathbf{Z}_{\left(r,e_H\right)}\big)^*:p,r\in P\big\}.
\end{align*}
\end{rem}

\begin{rem}
\label{Y_q contains Z_{(e,q)}}
The module $\mathbf{Y}_q^\mathcal{NT}$ contains $i_{\mathbf{Z}_{\left(e_G,q\right)}}(\mathbf{Z}_{\left(e_G,q\right)})$: if $z\in \mathbf{Z}_{\left(e_G,q\right)}$, and $z'\in \mathbf{Z}_{\left(e_G,q\right)}$ is chosen so that $z=z'\cdot \langle z',z'\rangle_A^{\left(e_G,q\right)}$, then
\begin{align*}
i_{\mathbf{Z}_{\left(e_G,q\right)}}(z)
&=i_{\mathbf{Z}_{\left(e_G,q\right)}}(z')i_{\mathbf{Z}_{\left(e_G,e_H\right)}}\big(\langle z',z'\rangle_A^{\left(e_G,q\right)}\big)\\
&=i_{\mathbf{Z}_{\left(e_G,q\right)}}(z')\phi_{\mathbf{X}}^\mathcal{NT}\big(i_{\mathbf{X}_{e_G}}\big(\langle z',z'\rangle_A^{\left(e_G,q\right)}\big)\big)
\in \mathbf{Y}_q^\mathcal{NT}.
\end{align*}
\end{rem}

We now show that $\mathbf{Y}_q^\mathcal{NT}$ also carries a left action of $\mathcal{NT}_\mathbf{X}$ by adjointable operators for each $q\in Q\setminus \{e_H\}$, and hence is a Hilbert $\mathcal{NT}_\mathbf{X}$-bimodule.

\begin{prop}
\label{right actions on each Y_q}
Suppose that $G$ is an amenable group so that the Hilbert $\mathcal{NT}_\mathbf{X}$-module $\mathbf{Y}_q^\mathcal{NT}$ of Proposition~\ref{definition of Y_q} is defined. For each $q\in Q\setminus \{e_H\}$, there exists a homomorphism $\Phi_q^{\mathcal{NT}}: \mathcal{NT}_\mathbf{X}\rightarrow \mathcal{L}_{\mathcal{NT}_\mathbf{X}}\left(\mathbf{Y}_q^\mathcal{NT}\right)$ such that 
$
\Phi_q^{\mathcal{NT}}(b)(y)=\phi_{\mathbf{X}}^\mathcal{NT}(b)y 
$
for each $b\in \mathcal{NT}_\mathbf{X}$ and $y\in \mathbf{Y}_q^\mathcal{NT}$.
\begin{proof}
Firstly, we check that $\phi_{\mathbf{X}}^\mathcal{NT}\left(\mathcal{NT}_\mathbf{X}\right)\mathbf{Y}_q^\mathcal{NT}\subseteq \mathbf{Y}_q^\mathcal{NT}$. Since the semidirect product $(G\rtimes_\alpha H, P\rtimes_\alpha Q)$ has the product order, for any $s,t\in P$ we have
\begin{align*}
\left(s,e_H\right)\left(t,e_H\right)^{-1}\left(\left(t,e_H\right)\vee \left(e_G,q\right)\right)
=\left(e_G,q\right)\left(\alpha_{q^{-1}}(s),e_H\right)
\end{align*}
and
\begin{align*}
\left(e_G,q\right)^{-1}\left(\left(t,e_H\right)\vee \left(e_G,q\right)\right)=\left(\alpha_{q^{-1}}(t),e_H\right).
\end{align*}
Thus, for any $x\in \mathbf{X}_s$, $y\in \mathbf{X}_t$, $z\in \mathbf{Z}_{\left(e_G,q\right)}$, and $b\in \mathcal{NT}_\mathbf{X}$, using the Nica covariance of $i_\mathbf{Z}$, we have
\begin{align*}
\phi&_{\mathbf{X}}^\mathcal{NT}\left(i_{\mathbf{X}_s}(x)i_{\mathbf{X}_t}(y)^*\right)i_{\mathbf{Z}_{\left(e_G,q\right)}}(z)\phi_{\mathbf{X}}^\mathcal{NT}(b)\\
&=
i_{\mathbf{Z}_{\left(s,e_H\right)}}(x)i_{\mathbf{Z}_{\left(t,e_H\right)}}(y)^*i_{\mathbf{Z}_{\left(e_G,q\right)}}(z)\phi_{\mathbf{X}}^\mathcal{NT}(b)\\
&\in
\cspan\Big\{i_{\mathbf{Z}_{\left(e_G,q\right)}}(\mathbf{Z}_{\left(e_G,q\right)})\phi_{\mathbf{X}}^\mathcal{NT}\Big(i_{\mathbf{X}_{\alpha_{q^{-1}}(s)}}\big(\mathbf{X}_{\alpha_{q^{-1}}(s)}\big)i_{\mathbf{X}_{\alpha_{q^{-1}}(t)}}\big(\mathbf{X}_{\alpha_{q^{-1}}(t)}\big)^*b\Big)\Big\}.
\end{align*}
Since 
\[
\big\{i_{\mathbf{X}_s}(\mathbf{X}_s)i_{\mathbf{X}_t}(\mathbf{X}_t)^*:s,t\in P\big\} 
\quad \text{and} \quad
\big\{i_{\mathbf{Z}_{\left(e_G,q\right)}}(\mathbf{Z}_{\left(e_G,q\right)})\phi_{\mathbf{X}}^\mathcal{NT}\left(\mathcal{NT}_\mathbf{X}\right)\big\}
\]
span dense subspaces of $\mathcal{NT}_\mathbf{X}$ and $\mathbf{Y}_q^\mathcal{NT}$ respectively, whilst $\phi_{\mathbf{X}}^\mathcal{NT}$ is a homomorphism and multiplication in $\mathcal{NT}_\mathbf{Z}$ is bilinear and continuous, we conclude that $\phi_{\mathbf{X}}^\mathcal{NT}\left(\mathcal{NT}_\mathbf{X}\right)\mathbf{Y}_q^\mathcal{NT}$ is contained in $\mathbf{Y}_q^\mathcal{NT}$. Thus, for each $b\in \mathcal{NT}_\mathbf{X}$, we may define $\Phi_q^{\mathcal{NT}}(b):\mathbf{Y}_q^\mathcal{NT}\rightarrow\mathbf{Y}_q^\mathcal{NT}$ by $\Phi_q^{\mathcal{NT}}(b)(y):=\phi_{\mathbf{X}}^\mathcal{NT}(b)y$ for each $y\in \mathbf{Y}_q^\mathcal{NT}$. Next, we claim that $\Phi_q^{\mathcal{NT}}(b)$ is adjointable with $\Phi_q^{\mathcal{NT}}(b)^*=\Phi_q^{\mathcal{NT}}(b^*)$. To see this, observe that for any $y,w\in \mathbf{Y}_q^\mathcal{NT}$,
\begin{align*}
\left\langle \phi_{\mathbf{X}}^\mathcal{NT}(b)y,w\right\rangle_{\mathcal{NT}_\mathbf{X}}^q
=\left(\phi_{\mathbf{X}}^\mathcal{NT}\right)^{-1}\left(y^*\phi_{\mathbf{X}}^\mathcal{NT}(b^*)w\right)
=\left\langle y,\phi_{\mathbf{X}}^\mathcal{NT}(b^*)w\right\rangle_{\mathcal{NT}_\mathbf{X}}^q.
\end{align*}
Finally, since $\phi_{\mathbf{X}}^\mathcal{NT}$ is linear and multiplicative, the map $b\mapsto \Phi_q^{\mathcal{NT}}(b)$ is also linear and multiplicative. Thus, $b\mapsto \Phi_q^{\mathcal{NT}}(b)$ is a homomorphism from $\mathcal{NT}_\mathbf{X}$ to $\mathcal{L}_{\mathcal{NT}_\mathbf{X}}\left(\mathbf{Y}_q^\mathcal{NT}\right)$.
\end{proof}
\end{prop}

Next we show that $\mathbf{Y}^\mathcal{NT}:=\bigsqcup_{q\in Q} \mathbf{Y}_q^\mathcal{NT}$ can be viewed as a product system over the quasi-lattice ordered group $(H,Q)$ with coefficient algebra $\mathcal{NT}_\mathbf{X}$.

\begin{prop}
\label{Y is a product system}
Suppose that $G$ is an amenable group so that the Hilbert $\mathcal{NT}_\mathbf{X}$-bimodule $\mathbf{Y}_q^\mathcal{NT}$ from Propositions~\ref{definition of Y_q} and \ref{right actions on each Y_q} is defined. Let $\mathbf{Y}^\mathcal{NT}:=\bigsqcup_{q\in Q} \mathbf{Y}_q^\mathcal{NT}$. Then $\mathbf{Y}^\mathcal{NT}$ is a product system over $(H,Q)$ with coefficient algebra $\mathcal{NT}_\mathbf{X}$, and multiplication given by multiplication in $\mathcal{NT}_\mathbf{Z}$.
\begin{proof}
We already know from Propositions~\ref{definition of Y_q} and~\ref{right actions on each Y_q} that each $\mathbf{Y}_q^\mathcal{NT}$ is a Hilbert $\mathcal{NT}_\mathbf{X}$-bimodule, and $\mathbf{Y}_{e_H}^\mathcal{NT}={}_{\mathcal{NT}_\mathbf{X}}\left(\mathcal{NT}_\mathbf{X}\right)_{\mathcal{NT}_\mathbf{X}}$ by definition. If we equip $\mathbf{Y}^\mathcal{NT}$ with the associative multiplication from $\mathcal{NT}_\mathbf{Z}$ and identify $\mathcal{NT}_\mathbf{X}$ with $\phi_{\mathbf{X}}^\mathcal{NT}\left(\mathcal{NT}_\mathbf{X}\right)\subseteq \mathcal{NT}_\mathbf{Z}$, then $\mathbf{Y}^\mathcal{NT}$ becomes a semigroup. It is straightforward to check that multiplication in $\mathbf{Y}^\mathcal{NT}$ by elements of $\mathbf{Y}_{e_H}^\mathcal{NT}=\mathcal{NT}_\mathbf{X}$  implements the left and right actions of $\mathcal{NT}_\mathbf{X}$ on each $\mathbf{Y}_q^\mathcal{NT}$. 

For $\mathbf{Y}^\mathcal{NT}$ to be a product system, it remains to show that there exists a Hilbert $\mathcal{NT}_\mathbf{X}$-bimodule isomorphism $M_{q,t}^{\mathbf{Y}^\mathcal{NT}}:\mathbf{Y}_q^\mathcal{NT}\otimes_{\mathcal{NT}_\mathbf{X}}\mathbf{Y}_t^\mathcal{NT}\rightarrow \mathbf{Y}_{qt}^\mathcal{NT}$ for each $q, t\in Q\setminus \{e_H\}$ such that 
$
M_{q,t}^{\mathbf{Y}^\mathcal{NT}}(y\otimes_{\mathcal{NT}_\mathbf{X}} w)=yw
$
for each $y\in \mathbf{Y}_q^\mathcal{NT}$ and $w\in \mathbf{Y}_t^\mathcal{NT}$. We begin by checking that $\mathbf{Y}_q^\mathcal{NT} \mathbf{Y}_t^\mathcal{NT}\subseteq \mathbf{Y}_{qt}^\mathcal{NT}$ for each $q,t\in Q$. Making use of Proposition~\ref{right actions on each Y_q}, we see that
\begin{align*}
\mathbf{Y}&_q^\mathcal{NT} \mathbf{Y}_t^\mathcal{NT}\\
&=\cspan\big\{i_{\mathbf{Z}_{\left(e_G,q\right)}}\left(\mathbf{Z}_{\left(e_G,q\right)}\right)\phi_{\mathbf{X}}^\mathcal{NT}\left(\mathcal{NT}_\mathbf{X}\right)i_{\mathbf{Z}_{\left(e_G,t\right)}}(\mathbf{Z}_{\left(e_G,t\right)})\phi_{\mathbf{X}}^\mathcal{NT}\left(\mathcal{NT}_\mathbf{X}\right)\big\}\\
&\subseteq\cspan\big\{i_{\mathbf{Z}_{\left(e_G,qt\right)}}\left(\mathbf{Z}_{\left(e_G,qt\right)}\right)\phi_{\mathbf{X}}^\mathcal{NT}\left(\mathcal{NT}_\mathbf{X}\right)\big\}
=\mathbf{Y}_{qt}^\mathcal{NT}.
\end{align*}
Next, observe that if $y,u\in \mathbf{Y}_q^\mathcal{NT}$ and $w,v\in \mathbf{Y}_t^\mathcal{NT}$, then
\begin{align*}\
\langle y\otimes_{\mathcal{NT}_\mathbf{X}}w, u\otimes_{\mathcal{NT}_\mathbf{X}}v \rangle_{\mathcal{NT}_\mathbf{X}}
&=\big\langle w, \langle y,u \rangle_{\mathcal{NT}_\mathbf{X}}^q\cdot v \big\rangle_{\mathcal{NT}_\mathbf{X}}^t\\
&=\big\langle w, \phi\left(\left(\phi_{\mathbf{X}}^\mathcal{NT}\right)^{-1}\left(y^*u\right)\right)v \big\rangle_{\mathcal{NT}_\mathbf{X}}^t\\
&=\left(\phi_{\mathbf{X}}^\mathcal{NT}\right)^{-1}\left(w^*y^*uv\right)\\
&=\left\langle yw, uv\right\rangle_{\mathcal{NT}_\mathbf{X}}^{qt}.
\end{align*}
Thus, the map $y\otimes_{\mathcal{NT}_\mathbf{X}} w\mapsto yw$ extends by linearity and continuity to a well-defined inner-product preserving map from $\mathbf{Y}_q^\mathcal{NT}\otimes_{\mathcal{NT}_\mathbf{X}}\mathbf{Y}_t^\mathcal{NT}$ to $\mathbf{Y}_{qt}^\mathcal{NT}$, which we denote by $M_{q,t}^{\mathbf{Y}^\mathcal{NT}}$. Clearly, $M_{q,t}^{\mathbf{Y}^\mathcal{NT}}$ is both left and right $\mathcal{NT}_\mathbf{X}$-linear. Surjectivity of the map $M_{q,t}^{\mathbf{Y}^\mathcal{NT}}$ follows from the fact that $\mathbf{Z}_{\left(e_G,qt\right)}=\mathbf{Z}_{\left(e_G,q\right)\left(e_G,t\right)}\cong \mathbf{Z}_{\left(e_G,q\right)}\otimes_A \mathbf{Z}_{\left(e_G,t\right)}$ (since $q\neq e_H$) and $i_{\mathbf{Z}_{\left(e_G,q\right)}}\left(\mathbf{Z}_{\left(e_G,q\right)}\right)\subseteq \mathbf{Y}_q^\mathcal{NT}$:
\begin{align*}
\mathbf{Y}_{qt}^\mathcal{NT}
&=\cspan\big\{i_{\mathbf{Z}_{\left(e_G,qt\right)}}\left(\mathbf{Z}_{\left(e_G,qt\right)}\right)\phi_{\mathbf{X}}^\mathcal{NT}\left(\mathcal{NT}_\mathbf{X}\right)\big\}\\
&=\cspan\big\{i_{\mathbf{Z}_{\left(e_G,q\right)}}\left(\mathbf{Z}_{\left(e_G,q\right)}\right)i_{\mathbf{Z}_{\left(e_G,t\right)}}(\mathbf{Z}_{\left(e_G,t\right)})\phi_{\mathbf{X}}^\mathcal{NT}\left(\mathcal{NT}_\mathbf{X}\right)\big\}\\
&\subseteq \mathbf{Y}_q^\mathcal{NT} \mathbf{Y}_t^\mathcal{NT}\\
&=M_{q,t}^{\mathbf{Y}^\mathcal{NT}}\left(\mathbf{Y}_q^\mathcal{NT}\otimes_{\mathcal{NT}_\mathbf{X}} \mathbf{Y}_t^\mathcal{NT}\right).
\end{align*}
Putting all of this together, we see that $M_{q,t}^{\mathbf{Y}^\mathcal{NT}}$ is an $\mathcal{NT}_\mathbf{X}$-bimodule isomorphism.
\end{proof}
\end{prop}

We now prove that if $A$ acts faithfully on $\mathbf{Z}_{\left(e_G,q\right)}$, then $\mathcal{NT}_\mathbf{X}$ acts faithfully on $\mathbf{Y}_q^\mathcal{NT}$. We will make use of this result later in Subsection~\ref{mixed cases} when we consider the Cuntz--Nica--Pimsner algebra of the product system $\mathbf{Y}^\mathcal{NT}$.

\begin{prop}
\label{faithfulness of induced action - Nica Toeplitz case}
Suppose that $G$ is an amenable group so that the product system $\mathbf{Y}^\mathcal{NT}$ from Proposition~\ref{Y is a product system} is defined. If $q\in Q$ and $A$ acts faithfully on $\mathbf{Z}_{\left(e_G,q\right)}$, then the homomorphism $\Phi^\mathcal{NT}_q:\mathcal{NT}_\mathbf{X}\rightarrow \mathcal{L}_{\mathcal{NT}_\mathbf{X}}\left(\mathbf{Y}_q^{\mathcal{NT}}\right)$ is injective. 
\begin{proof}
When $q=e_H$, the map $\Phi^\mathcal{NT}_q$ is just left multiplication on $\mathbf{Y}_{e_H}^{\mathcal{NT}}=\mathcal{NT}_\mathbf{X}$ by elements of $\mathcal{NT}_\mathbf{X}$, which is obviously faithful. So suppose that $q\in Q\setminus \left\{e_H\right\}$. Let $\rho:A\rightarrow \mathcal{B}(\mathcal{H})$ be a faithful nondegenerate representation of $A$ on a Hilbert space $\mathcal{H}$. To prove that $\Phi^\mathcal{NT}_q$ is faithful, it suffices to show that the induced representation
\begin{align*}
\big(\mathbf{Y}_q^{\mathcal{NT}}&\otimes_{\mathcal{NT}_\mathbf{X}}\mathcal{F}_\mathbf{X}\big)
\text{-}\mathrm{Ind}_A^{\mathcal{NT}_\mathbf{X}}\rho\\
&=
\Phi^\mathcal{NT}_q\otimes_{\mathcal{NT}_\mathbf{X}}\mathrm{id}_{\mathcal{F}_\mathbf{X}}\otimes_A \mathrm{id}_\mathcal{H}
:\mathcal{NT}_\mathbf{X}\rightarrow \mathcal{B}\left(\mathbf{Y}_q\otimes_{\mathcal{NT}_\mathbf{X}}\mathcal{F}_\mathbf{X}\otimes_A \mathcal{H}\right)
\end{align*}
is faithful. Let $\varrho:=\left(\Phi^\mathcal{NT}_q\otimes_{\mathcal{NT}_\mathbf{X}}\mathrm{id}_{\mathcal{F}_\mathbf{X}}\otimes_A \mathrm{id}_\mathcal{H}\right)\circ i_\mathbf{X}$ denote the Nica covariant representation of $\mathbf{X}$ that induces $\Phi^\mathcal{NT}_q\otimes_{\mathcal{NT}_\mathbf{X}}\mathrm{id}_{\mathcal{F}_\mathbf{X}}\otimes_A \mathrm{id}_\mathcal{H}$. Our aim is to prove that
\begin{enumerate}[label=\upshape(\roman*)]
\item
for each $p\in P\setminus \{e_G\}$,
\[
i_{\mathbf{Z}_{\left(e_G,q\right)}}\left(\mathbf{Z}_{\left(e_G,q\right)}\right)\otimes_{\mathcal{NT}_\mathbf{X}}A\otimes_A \mathcal{H}\subseteq \left(\varrho_p(\mathbf{X}_p)\left(\mathbf{Y}_q^{\mathcal{NT}}\otimes_{\mathcal{NT}_\mathbf{X}}\mathcal{F}_\mathbf{X}\otimes_A \mathcal{H}\right)\right)^\perp;
\]
\item
$A$ acts faithfully (via $\varrho_{e_G}$) on $i_{\mathbf{Z}_{\left(e_G,q\right)}}(\mathbf{Z}_{\left(e_G,q\right)})\otimes_{\mathcal{NT}_\mathbf{X}}A\otimes_A \mathcal{H}$.
\end{enumerate}
To see that this suffices, suppose for a moment that (i) and (ii) hold. For each $t\in P\setminus \{e_G\}$, let 
\[
P_t^\varrho:=\mathrm{proj}_{\overline{\varrho_t(\mathbf{X}_t)\left(\mathbf{Y}_q^{\mathcal{NT}}\otimes_{\mathcal{NT}_\mathbf{X}}\mathcal{F}_\mathbf{X}\otimes_A \mathcal{H}\right)}}.
\]
Then the representation 
\[
A\ni a \mapsto \varrho_{e_G}(a)\prod_{t\in K}\left(1-P_t^\varrho\right)\in \mathcal{B}\left(\mathbf{Y}_q^{\mathcal{NT}}\otimes_{\mathcal{NT}_\mathbf{X}}\mathcal{F}_\mathbf{X}\otimes_A \mathcal{H}\right)
\]
is faithful for each finite subset $K\subseteq P\setminus \{e_G\}$ . Since $G$ is amenable, \cite[Theorem~3.2]{NTuniqueness} then implies that $\Phi^\mathcal{NT}_q\otimes_{\mathcal{NT}_\mathbf{X}}\mathrm{id}_{\mathcal{F}_\mathbf{X}}\otimes_A \mathrm{id}_\mathcal{H}$ is faithful as required. 

We prove (i) and (ii). Let $p\in P$. For $s\in P$, we have $(p,e_H)(s,q)=(e,q)(\alpha_{q^{-1}}(ps),e_H)$, and so using the description of $\mathbf{Y}_q^{\mathcal{NT}}$ given in Remark~4.4, we see that
\begin{align*}
\Phi^\mathcal{NT}_q&\left(i_{\mathbf{X}_p}(\mathbf{X}_p)\right)\left(\mathbf{Y}_q^{\mathcal{NT}}\right)\\
&=i_{\mathbf{Z}_{\left(e_G,q\right)}}(\mathbf{Z}_{\left(e_G,q\right)})\cdot\cspan\big\{
i_{\mathbf{X}_{\alpha_{q^{-1}}(ps)}}\big(\mathbf{X}_{\alpha_{q^{-1}}(ps)}\big)i_{\mathbf{X}_t}(\mathbf{X}_t)^*:s,t\in P 
\big\}.
\end{align*}
For any $s,r\in P$, since $\alpha_{q^{-1}}(p)^{-1}\alpha_{q^{-1}}(ps)r=\alpha_{q^{-1}}(s)r\in \alpha_H(P)P\subseteq P$, we have that $\alpha_{q^{-1}}(p)\leq \alpha_{q^{-1}}(ps)r$.
Thus, 
\[
\cspan\big\{i_{\mathbf{X}_{\alpha_{q^{-1}}(ps)}}\big(\mathbf{X}_{\alpha_{q^{-1}}(ps)}\big)i_{\mathbf{X}_t}(\mathbf{X}_t)^*:s,t\in P\big\}\cdot\mathcal{F}_\mathbf{X}
\subseteq
\bigoplus_{m\in P:\alpha_{q^{-1}}(p)\leq m}\mathbf{X}_m,
\]
and so
\begin{align*}
\varrho_p(\mathbf{X}_p)&\left(\mathbf{Y}_q^{\mathcal{NT}}\otimes_{\mathcal{NT}_\mathbf{X}}\mathcal{F}_\mathbf{X}\otimes_A \mathcal{H}\right)\\
&\subseteq
i_{\mathbf{Z}_{\left(e_G,q\right)}}(\mathbf{Z}_{\left(e_G,q\right)})\otimes_{\mathcal{NT}_\mathbf{X}}\bigoplus_{m\in P:\alpha_{q^{-1}}(p)\leq m}\mathbf{X}_m\otimes_A \mathcal{H}.
\end{align*}
Now suppose that $p\neq e_G$. Since $\alpha_H(P)\subseteq P$ and $P\cap P^{-1}=\{e_G\}$ it follows that $\alpha_{q^{-1}}(p)\not\leq e_G$. Hence for any $z,w \in \mathbf{Z}_{\left(e_G,q\right)}$, $a\in \mathbf{X}_{e_G}=A$, $x\in \mathbf{X}_m$ with $\alpha_{q^{-1}}(p)\leq m$, and $h,g\in \mathcal{H}$, we see that
\begin{align*}
\big\langle i_{\mathbf{Z}_{\left(e_G,q\right)}}(z)
\otimes_{\mathcal{NT}_\mathbf{X}}a\otimes_A h, &i_{\mathbf{Z}_{\left(e_G,q\right)}}(w)\otimes_{\mathcal{NT}_\mathbf{X}}x\otimes_A g\big\rangle_\C\\
&=\big\langle h, \big\langle a, i_{\mathbf{X}_{e_G}}\big(\langle z,w\rangle_A^{\left(e_G,q\right)}\big)\cdot x \big\rangle_A \cdot g\big\rangle_\C
\end{align*}
must be zero (since $i_{\mathbf{X}_{e_G}}\big(\langle z,w\rangle_A^{\left(e_G,q\right)}\big)\cdot x \in \mathbf{X}_m$ which is orthogonal to $A=\mathbf{X}_{e_G}$ in the Fock space $\mathcal{F}_\mathbf{X}$).
Since the inner-product on $\mathbf{Y}_q^{\mathcal{NT}}\otimes_{\mathcal{NT}_\mathbf{X}}\mathcal{F}_\mathbf{X}\otimes_A \mathcal{H}$ is linear and continuous, we conclude that 
\[
i_{\mathbf{Z}_{\left(e_G,q\right)}}\left(\mathbf{Z}_{\left(e_G,q\right)}\right)\otimes_{\mathcal{NT}_\mathbf{X}}A\otimes_A \mathcal{H}\subseteq \left(\varrho_p\left(\mathbf{X}_p\right)\left(\mathbf{Y}_q^{\mathcal{NT}}\otimes_{\mathcal{NT}_\mathbf{X}}\mathcal{F}_\mathbf{X}\otimes_A \mathcal{H}\right)\right)^\perp
.\]

Next we check that $A$ acts faithfully on $i_{\mathbf{Z}_{\left(e_G,q\right)}}\left(\mathbf{Z}_{\left(e_G,q\right)}\right)\otimes_{\mathcal{NT}_\mathbf{X}}A\otimes_A \mathcal{H}$ via the homomorphism $\varrho_{e_G}=\big(\Phi^\mathcal{NT}_q\circ i_{\mathbf{X}_{e_G}}\big)\otimes_{\mathcal{NT}_\mathbf{X}}\mathrm{id}_{\mathcal{F}_\mathbf{X}}\otimes_A \mathrm{id}_\mathcal{H}$. Fix $a\in A\setminus \{0\}$. Since $A$ acts faithfully on $\mathbf{Z}_{\left(e_G,q\right)}$, we can choose $z\in \mathbf{Z}_{\left(e_G,q\right)}$ such that $a\cdot z\neq 0$. Since $\rho$ is faithful, we can then choose $h\in \mathcal{H}$ such that $\big\langle h, \langle a\cdot z, a\cdot z\rangle_A^{\left(e_G,q\right)}\cdot h \big\rangle_\C\neq 0$. By the Hewitt--Cohen--Blanchard factorisation theorem, we can write $z=z'\cdot \langle z',z'\rangle_A^{\left(e_G,q\right)}$ for some $z'\in \mathbf{Z}_{\left(e_G,q\right)}$. 
Routine calculations show that
\begin{align*}
\Big\langle i_{\mathbf{Z}_{\left(e_G,q\right)}}(a\cdot z')\otimes_{\mathcal{NT}_\mathbf{X}}\langle z',z'\rangle_A^{\left(e_G,q\right)},&i_{\mathbf{Z}_{\left(e_G,q\right)}}(a\cdot z')\otimes_{\mathcal{NT}_\mathbf{X}}\langle z',z'\rangle_A^{\left(e_G,q\right)}\Big\rangle_A\\
&=
\langle a\cdot z, a\cdot z\rangle_A^{\left(e_G,q\right)}.
\end{align*}
Hence, 
\begin{align*}
\Big\|\varrho_e(a)&\Big(i_{\mathbf{Z}_{\left(e_G,q\right)}}(z')\otimes_{\mathcal{NT}_\mathbf{X}} \langle z',z'\rangle_A^{\left(e_G,q\right)}\otimes_A h\Big)\Big\|^2\\
&=
\big\langle h, \langle a\cdot z, a\cdot z\rangle_A^{\left(e_G,q\right)}\cdot h \big\rangle_\C
\neq 0,
\end{align*}
and so $A$ acts faithfully on $i_{\mathbf{Z}_{\left(e_G,q\right)}}\left(\mathbf{Z}_{\left(e_G,q\right)}\right)\otimes_{\mathcal{NT}_\mathbf{X}}A\otimes_A \mathcal{H}$.
\end{proof}
\end{prop}

We now work towards showing that the product system $\mathbf{Y}^\mathcal{NT}$ is compactly aligned. The next result characterises the compact operators on each fibre of $\mathbf{Y}^\mathcal{NT}$. We need some more notation: given a $C^*$-algebra $B$, for each $b\in B$, we write $M_b\in \mathcal{L}_B(B_B)$ for the map defined by $M_b(c):=bc$ for each $c\in B$. 

\begin{lem}
\label{useful property for multiplication operators}
Suppose that $G$ is an amenable group so that the Hilbert $\mathcal{NT}_\mathbf{X}$-bimodule $\mathbf{Y}_q^\mathcal{NT}$ from Propositions~\ref{definition of Y_q} and \ref{right actions on each Y_q} is defined. For each $q\in Q$, if $b\in \mathcal{NT}_\mathbf{Z}$ is such that $M_b\in \mathcal{L}_{\mathcal{NT}_\mathbf{X}}\left(\mathbf{Y}_q^\mathcal{NT}\right)$, then
\[
\norm{M_b}_{\mathcal{L}_{\mathcal{NT}_\mathbf{X}}\left(\mathbf{Y}_q^\mathcal{NT}\right)}\leq \norm{b}_{\mathcal{NT}_\mathbf{Z}}
.\]
\begin{proof}
Since the norm on $\mathbf{Y}_q^\mathcal{NT}\subseteq\mathcal{NT}_\mathbf{Z}$ is just the restriction of the norm on $\mathcal{NT}_\mathbf{Z}$, and the norm on any $C^*$-algebra is submultiplicative, we have
\begin{align*}
\norm{M_b}_{\mathcal{L}_{\mathcal{NT}_\mathbf{X}}\left(\mathbf{Y}_q^\mathcal{NT}\right)}
&=\sup_{\substack{y\in \mathbf{Y}_q^\mathcal{NT}:\\ \norm{y}_{\mathbf{Y}_q^\mathcal{NT}}\leq 1}}\norm{M_b(y)}_{\mathbf{Y}_q^\mathcal{NT}}\\
&=\sup_{\substack{y\in \mathbf{Y}_q^\mathcal{NT}:\\ \norm{y}_{\mathcal{NT}_\mathbf{Z}}\leq 1}}\norm{by}_{\mathcal{NT}_\mathbf{Z}}
\leq \norm{b}_{\mathcal{NT}_\mathbf{Z}}.
\qedhere
\end{align*}
\end{proof}
\end{lem}

\begin{lem}
\label{compact operators in Y}
Suppose that $G$ is an amenable group so that the product system $\mathbf{Y}^\mathcal{NT}$ from Proposition~\ref{Y is a product system} is defined. For each $q\in Q$, let $\mathcal{NT}_\mathbf{Z}^{q}$ denote the closed subspace $\cspan\{i_{\mathbf{Z}_{(p,q)}}(\mathbf{Z}_{(p,q)})i_{\mathbf{Z}_{(r,q)}}(\mathbf{Z}_{(r,q)})^*:p,r\in P\}$.

\begin{enumerate}[label=\upshape(\roman*)]
\item
Let $q, t\in Q$ and $T\in \mathcal{K}_{\mathcal{NT}_\mathbf{X}}\left(\mathbf{Y}_q^\mathcal{NT}\right)$. Then $\iota_q^{qt}(T)\in \mathcal{L}_{\mathcal{NT}_\mathbf{X}}\left(\mathbf{Y}_{qt}^\mathcal{NT}\right)$ is left multiplication by an element of $\mathcal{NT}_\mathbf{Z}^{q}$ on $\mathbf{Y}_{qt}^\mathcal{NT} \subseteq \mathcal{NT}_\mathbf{Z}$. In particular, if $y,w\in \mathbf{Y}_q^\mathcal{NT}$, then $yw^*\in \mathcal{NT}_\mathbf{Z}^{q}$, and the rank one operator $\Theta_{y,w}\in \mathcal{K}_{\mathcal{NT}_\mathbf{X}}\left(\mathbf{Y}_q^\mathcal{NT}\right)$ satisfies
\[
\iota_q^{qt}(\Theta_{y,w})=
\begin{cases}
M_{yw^*} & \text{if $q\neq e_H$}\\
M_{\phi_{\mathbf{X}}^\mathcal{NT}(yw^*)} & \text{if $q=e_H$.}
\end{cases}
\]
\item
Let $b\in \mathcal{NT}_\mathbf{Z}^{q}$ and $q\neq e_H$. Then $M_b\in \mathcal{K}_{\mathcal{NT}_\mathbf{X}}\left(\mathbf{Y}_q^\mathcal{NT}\right)$.
\end{enumerate}
\begin{proof}
We begin by proving part (i) of the result. Fix $q, t\in Q$ and a rank one operator $\Theta_{y,w}\in \mathcal{K}_{\mathcal{NT}_\mathbf{X}}\left(\mathbf{Y}_q^\mathcal{NT}\right)$. If $q=e_H$, then for any $z\in \mathbf{Y}_t^\mathcal{NT}$, we have
\begin{align*}
\iota_{e_H}^t(\Theta_{y,w})(z)=\Phi_t^{\mathcal{NT}}(yw^*)(z)=\phi_{\mathbf{X}}^\mathcal{NT}(yw^*)z=M_{\phi_{\mathbf{X}}^\mathcal{NT}(yw^*)}(z).
\end{align*}
Thus, $\iota_{e_H}^{t}(\Theta_{y,w})=M_{\phi_{\mathbf{X}}^\mathcal{NT}(yw^*)}$. Moreover, $\phi_{\mathbf{X}}^\mathcal{NT}(yw^*)\in \phi_{\mathbf{X}}^\mathcal{NT}\left(\mathcal{NT}_\mathbf{X}\right)=\mathcal{NT}_\mathbf{Z}^{e_H}$. 

Now suppose that $q\neq e_H$. For any $u\in \mathbf{Z}_{\left(e_G,q\right)}$, $v\in  \mathbf{Z}_{\left(e_G,t\right)}$, and $b\in \mathcal{NT}_\mathbf{X}$, 
\begin{align*}
\iota_q^{qt}(\Theta_{y,w})\big(i_{\mathbf{Z}_{\left(e_G,qt\right)}}(uv)\phi_{\mathbf{X}}^\mathcal{NT}(b)\big)
&=\Theta_{y,w}\big(i_{\mathbf{Z}_{\left(e_G,q\right)}}(u)\big)i_{\mathbf{Z}_{\left(e_G,t\right)}}(v)\phi_{\mathbf{X}}^\mathcal{NT}(b)\\
&=yw^*i_{\mathbf{Z}_{\left(e_G,q\right)}}(u)i_{\mathbf{Z}_{\left(e_G,t\right)}}(v)\phi_{\mathbf{X}}^\mathcal{NT}(b)\\
&=M_{yw^*}\big(i_{\mathbf{Z}_{\left(e_G,qt\right)}}(uv)\phi_{\mathbf{X}}^\mathcal{NT}(b)\big).
\end{align*}
Since $\mathbf{Y}_{qt}^\mathcal{NT}=\cspan\big\{i_{\mathbf{Z}_{\left(e_G,q\right)}}\big(\mathbf{Z}_{\left(e_G,q\right)}\big)i_{\mathbf{Z}_{\left(e_G,t\right)}}\big(\mathbf{Z}_{\left(e_G,t\right)}\big)\phi_{\mathbf{X}}^\mathcal{NT}\left(\mathcal{NT}_\mathbf{X}\right)\big\}$
whilst $\iota_q^{qt}(\Theta_{y,w})$ and $M_{yw^*}$ are linear and continuous, we conclude that $\iota_q^{qt}(\Theta_{y,w})=M_{yw^*}\in \mathcal{L}_{\mathcal{NT}_\mathbf{X}}\left(\mathbf{Y}_{qt}^\mathcal{NT}\right)$. 

It remains to check that if $y, w\in \mathbf{Y}_q^\mathcal{NT}$, then $yw^*\in \mathcal{NT}_\mathbf{Z}^{q}$. Since $\mathcal{NT}_\mathbf{Z}^{q}$ is a closed subspace, it suffices to consider when $y=i_{\mathbf{Z}_{(p,q)}}(x)i_{\mathbf{Z}_{\left(r,e_H\right)}}(z)^*$ and $w=i_{\mathbf{Z}_{(m,q)}}(u)i_{\mathbf{Z}_{\left(n,e_H\right)}}(v)^*$ for some $x\in \mathbf{Z}_{(p,q)}$, $z\in \mathbf{Z}_{\left(r,e_H\right)}$, $u\in \mathbf{Z}_{(m,q)}$, and $v\in \mathbf{Z}_{\left(n,e_H\right)}$. If $r\vee n=\infty$, then the Nica covariance of $i_\mathbf{Z}$ gives
\[
yw^*=i_{\mathbf{Z}_{(p,q)}}(x)i_{\mathbf{Z}_{\left(r,e_H\right)}}(z)^*i_{\mathbf{Z}_{\left(n,e_H\right)}}(v)i_{\mathbf{Z}_{(m,q)}}(u)^*=0,
\]
which is certainly in $\mathcal{NT}_\mathbf{Z}^{q}$. On the other hand, if $r\vee n< \infty$, then $yw^*$ belongs to the closed span of 
\begin{align*}
i_{\mathbf{Z}_{(p\alpha_q(r^{-1}(r\vee n)),q)}}\big(\mathbf{Z}_{(p\alpha_q(r^{-1}(r\vee n)),q)}\big)i_{\mathbf{Z}_{(m\alpha_q(n^{-1}(r\vee n)),q)}}\big(\mathbf{Z}_{(m\alpha_q(n^{-1}(r\vee n)),q)}\big)^*,
\end{align*}
and so is contained in $\mathcal{NT}_\mathbf{Z}^{q}$.

We now prove part (ii) of the result. Let $b:=i_{\mathbf{Z}_{(p,q)}}(z)i_{\mathbf{Z}_{(r,q)}}(w)^*\in \mathcal{NT}_\mathbf{Z}^{q}$ where $z\in \mathbf{Z}_{(p,q)}$ and $w\in \mathbf{Z}_{(r,q)}$. Then $i_{\mathbf{Z}_{(p,q)}}(z),  i_{\mathbf{Z}_{(r,q)}}(w)\in \mathbf{Y}_q^\mathcal{NT}$ and, so by (i) we see that
\begin{align*}
M_{i_{\mathbf{Z}_{(p,q)}}(z)i_{\mathbf{Z}_{(r,q)}}(w)^*}
&=\iota_q^q\big(\Theta_{i_{\mathbf{Z}_{(p,q)}}(z),i_{\mathbf{Z}_{(r,q)}}(w)}\big)\\
&=\Theta_{i_{\mathbf{Z}_{(p,q)}}(z),i_{\mathbf{Z}_{(r,q)}}(w)}\in \mathcal{K}_{\mathcal{NT}_\mathbf{X}}\left(\mathbf{Y}_q^\mathcal{NT}\right). 
\end{align*}
Since the map $b\mapsto M_b$ is linear, and $\norm{M_b}_{\mathcal{L}_{\mathcal{NT}_\mathbf{X}}\left(\mathbf{Y}_q^\mathcal{NT}\right)}\leq \norm{b}_{\mathcal{NT}_\mathbf{Z}}$ whenever $b\in \mathcal{NT}_\mathbf{Z}$ is such that $M_b\in \mathcal{L}_{\mathcal{NT}_\mathbf{X}}\left(\mathbf{Y}_q^\mathcal{NT}\right)$ (by Lemma~\ref{useful property for multiplication operators}), we conclude that $M_b\in \mathcal{K}_{\mathcal{NT}_\mathbf{X}}\left(\mathbf{Y}_q^\mathcal{NT}\right)$ whenever $b\in  \cspan\big\{i_{\mathbf{Z}_{(p,q)}}\big(\mathbf{Z}_{(p,q)}\big)i_{\mathbf{Z}_{(r,q)}}\big(\mathbf{Z}_{(r,q)}\big)^*:p,r\in P\big\}=\mathcal{NT}_\mathbf{Z}^{q}$. 
\end{proof}
\end{lem}

\begin{lem}
\label{Q indexed subalgebras of NO_Z}
For any $q, t\in Q$, 
\[
\mathcal{NT}_\mathbf{Z}^{q}\mathcal{NT}_\mathbf{Z}^{t}
\subseteq
\begin{cases}
\mathcal{NT}_\mathbf{Z}^{q\vee t} & \text{if $q\vee t<\infty$}\\
\{0\} & \text{otherwise}.
\end{cases}
\]
In particular, each $\mathcal{NT}_\mathbf{Z}^{q}$ is a subalgebra of $\mathcal{NT}_\mathbf{Z}$. 
\begin{proof}
Since multiplication in $\mathcal{NT}_\mathbf{Z}$ is bilinear and continuous, it suffices to show that 
\[
i_{\mathbf{Z}_{(p,q)}}(x)i_{\mathbf{Z}_{(r,q)}}(z)^*i_{\mathbf{Z}_{(m,t)}}(u)i_{\mathbf{Z}_{(n,t)}}(v)^*\in 
\begin{cases}
\mathcal{NT}_\mathbf{Z}^{q\vee t}, & \text{if $q\vee t<\infty$}\\
\{0\}, & \text{otherwise}
\end{cases}
\]
for each $x\in \mathbf{Z}_{(p,q)}$, $y\in \mathbf{Z}_{(r,q)}$, $u\in \mathbf{Z}_{(m,t)}$, and $v\in \mathbf{Z}_{(n,t)}$.
If either $r\vee m=\infty$ or $q \vee t=\infty$, then $(r,q)\vee (m,t)=\infty$, and so $i_{\mathbf{Z}_{(p,q)}}(x)i_{\mathbf{Z}_{(r,q)}}(z)^*i_{\mathbf{Z}_{(m,t)}}(u)i_{\mathbf{Z}_{(n,t)}}(v)^*=0.$
On the other hand, if $r\vee m<\infty$ and $q \vee t<\infty$, then $(r,q)\vee (m,t)<\infty$, and so
$
i_{\mathbf{Z}_{(p,q)}}(x)i_{\mathbf{Z}_{(r,q)}}(z)^*i_{\mathbf{Z}_{(m,t)}}(u)i_{\mathbf{Z}_{(n,t)}}(v)^*
$
can be approximated by linear combinations of elements from the set
\[
i_{\mathbf{Z}_{(pr^{-1}(r\vee m),q\vee t)}}\big(\mathbf{Z}_{(pr^{-1}(r\vee m),q\vee t)}\big)
i_{\mathbf{Z}_{(nm^{-1}(r\vee m),q\vee t)}}\big(\mathbf{Z}_{(nm^{-1}(r\vee m),q\vee t)}\big)^*,
\]
and so belongs to $\mathcal{NT}_\mathbf{Z}^{q\vee t}$.
\end{proof}
\end{lem}

We are finally ready to prove that the product system $\mathbf{Y}^\mathcal{NT}$ is compactly aligned. 

\begin{prop}
\label{Y is compactly aligned}
Suppose that $G$ is an amenable group so that the product system $\mathbf{Y}^\mathcal{NT}$ from Proposition~\ref{Y is a product system} is defined. Then $\mathbf{Y}^\mathcal{NT}$ is compactly aligned. 
\begin{proof}
Let $S\in \mathcal{K}_{\mathcal{NT}_\mathbf{X}}\left(\mathbf{Y}_q^\mathcal{NT}\right)$ and $T\in \mathcal{K}_{\mathcal{NT}_\mathbf{X}}\left(\mathbf{Y}_t^\mathcal{NT}\right)$ with $q\vee t<\infty$. If $q=t=e_H$, then
\[
\iota_{q}^{q\vee t}(S)\iota_{t}^{q\vee t}(T)=\iota_{e_H}^{e_H}(S)\iota_{e_H}^{e_H}(T)=ST\in \mathcal{K}_{\mathcal{NT}_\mathbf{X}}\left(\mathbf{Y}_{e_H}^\mathcal{NT}\right)=\mathcal{K}_{\mathcal{NT}_\mathbf{X}}\left(\mathbf{Y}_{q\vee t}^\mathcal{NT}\right).
\]
Now suppose that $q\neq e_H$ or $t\neq e_H$. Thus, $q\vee t\neq e_H$. By Lemma~\ref{compact operators in Y}, $\iota_{q}^{q\vee t}(S)=M_b$ and $\iota_{t}^{q\vee t}(T)=M_c$ for some $b\in \mathcal{NT}_\mathbf{Z}^{q}$ and $c\in \mathcal{NT}_\mathbf{Z}^{t}$. Since $bc\in \mathcal{NT}_\mathbf{Z}^{q\vee t}$ by Lemma~\ref{Q indexed subalgebras of NO_Z}, we can use Lemma~\ref{compact operators in Y} again to see that
\begin{align*}
\iota_{q}^{q\vee t}(S)\iota_{t}^{q\vee t}(T)&=M_b M_c=M_{bc}\in \mathcal{K}_{\mathcal{NT}_\mathbf{X}}\left(\mathbf{Y}_{q\vee t}^\mathcal{NT}\right).
\qedhere
\end{align*}
\end{proof}
\end{prop}


\subsection{Isomorphisms of Nica--Toeplitz algebras}
\label{isomorphisms of Nica Toeplitz algebras}

For the product system $\mathbf{Y}^\mathcal{NT}$ defined in Subsection~\ref{Constructing the product system}, we will show that $\mathcal{NT}_{\mathbf{Y}^\mathcal{NT}}\cong \mathcal{NT}_\mathbf{Z}$. To do this we will use the universal property of each $C^*$-algebra to induce a homomorphism to the other, and then check that these homomorphisms are mutually inverse. In summary, we will show that the maps in Figure~\ref{NO isomorphism} exist and make the diagram commutative. 

To make our arguments easier to write down, we will identify the coefficient algebra $\mathcal{NT}_\mathbf{X}$ of $\mathbf{Y}^{\mathcal{NT}}$ with $\phi_{\mathbf{X}}^\mathcal{NT}\left(\mathcal{NT}_\mathbf{X}\right)\subseteq \mathcal{NT}_\mathbf{Z}$. Thus, every fibre of $\mathbf{Y}^\mathcal{NT}$ can be viewed as sitting inside $\mathcal{NT}_\mathbf{Z}$, and the left and right actions of $\mathcal{NT}_\mathbf{X}\cong\phi_{\mathbf{X}}^\mathcal{NT}\left(\mathcal{NT}_\mathbf{X}\right)$ on each $\mathbf{Y}_q^\mathcal{NT}$ are just multiplication in $\mathcal{NT}_\mathbf{Z}$.

To begin we get a homomorphism from $\mathcal{NT}_\mathbf{Z}$ to $\mathcal{NT}_{\mathbf{Y}^\mathcal{NT}}$ by exhibiting a Nica covariant representation of $\mathbf{Z}$ in $\mathcal{NT}_{\mathbf{Y}^\mathcal{NT}}$.

\begin{prop}
\label{varphi is a representation}
Suppose that $G$ is an amenable group so that the compactly aligned product system $\mathbf{Y}^\mathcal{NT}$ from Proposition~\ref{Y is a product system} is defined. Define $\varphi^{\mathcal{NT}}:\mathbf{Z}\rightarrow \mathcal{NT}_{\mathbf{Y}^\mathcal{NT}}$ by
\[
\varphi_{(p,q)}^{\mathcal{NT}}:=
i_{\mathbf{Y}_q^\mathcal{NT}}\circ i_{\mathbf{Z}_{(p,q)}}
,\]
for each $(p,q)\in P\rtimes_\alpha Q$. Then $\varphi^{\mathcal{NT}}$ is a representation of $\mathbf{Z}$.
\begin{proof}
Firstly, $\varphi_{(e_G,e_H)}^{\mathcal{NT}}=i_{\mathbf{Y}_{e_H}^\mathcal{NT}}\circ i_{\mathbf{Z}_{(e_G,e_H)}}$ is the composition of homomorphisms and so is a homomorphism. Similarly, $\varphi_{(p,q)}^{\mathcal{NT}}=i_{\mathbf{Y}_q^\mathcal{NT}}\circ i_{\mathbf{Z}_{(p,q)}}$ is the composition of linear maps, and so is linear for any $p\in P$, $q\in Q$. For any $z \in \mathbf{Z}_{(p,q)}$ and $w \in \mathbf{Z}_{(s,t)}$, since $i_\mathbf{Y}$ is a representation, we see that
\begin{align*}
\varphi_{(p,q)}^{\mathcal{NT}}(z)\varphi_{(s,t)}^{\mathcal{NT}}(w)
&=i_{\mathbf{Y}_q^\mathcal{NT}}\big(i_{\mathbf{Z}_{(p,q)}}(z)\big)i_{\mathbf{Y}_t}\big(i_{\mathbf{Z}_{(s,t)}}(w)\big)\\
&=i_{\mathbf{Y}_{qt}^\mathcal{NT}}\big(i_{\mathbf{Z}_{(p,q)}}(z)i_{\mathbf{Z}_{(s,t)}}(w)\big).
\end{align*}
Since $i_\mathbf{Z}$ is also a representation, this must be equal to
\begin{align*}
i_{\mathbf{Y}_{qt}^\mathcal{NT}}\big(i_{\mathbf{Z}_{(p,q)(s,t)}}(zw)\big)
&=i_{\mathbf{Y}_{qt}^\mathcal{NT}}\big(i_{\mathbf{Z}_{\left(p\alpha_q(s),qt\right)}}(zw)\big)\\
&=\varphi_{\left(p\alpha_q(s),qt\right)}^{\mathcal{NT}}(zw)
=\varphi_{(p,q)(s,t)}^{\mathcal{NT}}(zw), 
\end{align*}
and we see that $\varphi^{\mathcal{NT}}$ satisfies (T2).
For $(p, q)\in P\rtimes_\alpha Q$ and $z, w \in \mathbf{Z}_{(p,q)}$, since $i_\mathbf{Y}$ is a representation, we have
\begin{align*}
\varphi_{(p,q)}^{\mathcal{NT}}(z)^*\varphi_{(p,q)}^{\mathcal{NT}}(w)
&=i_{\mathbf{Y}_q^\mathcal{NT}}\big(i_{\mathbf{Z}_{(p,q)}}(z)\big)^*i_{\mathbf{Y}_q}\big(i_{\mathbf{Z}_{(p,q)}}(w)\big)\\
&=i_{\mathbf{Y}_{e_H}^\mathcal{NT}}\big(\big\langle i_{\mathbf{Z}_{(p,q)}}(z), i_{\mathbf{Z}_{(p,q)}}(w)\big\rangle_{\mathcal{NT}_\mathbf{X}}^q\big).
\end{align*}
Since $i_\mathbf{Z}$ is also a representation, this must be equal to
\begin{align*}
i_{\mathbf{Y}_{e_H}^\mathcal{NT}}\big(i_{\mathbf{Z}_{(p,q)}}(z)^*i_{\mathbf{Z}_{(p,q)}}(w)\big)
&=i_{\mathbf{Y}_{e_H}^\mathcal{NT}}\Big(i_{\mathbf{Z}_{\left(e_G,e_H\right)}}\big(\langle z,w\rangle_A^{(p,q)}\big)\Big)\\
&=\varphi_{\left(e_G,e_H\right)}^{\mathcal{NT}}\big(\langle z,w\rangle_A^{(p,q)}\big).
\end{align*}
Hence, $\varphi^{\mathcal{NT}}$ satisfies (T3), and so is a representation of $\mathbf{Z}$ in $\mathcal{NT}_\mathbf{Y}$. 
\end{proof}
\end{prop}

To show that $\varphi^{\mathcal{NT}}$ is Nica covariant, we first need a lemma.

\begin{lem}
\label{rep on compacts}
Suppose that $G$ is an amenable group so that the compactly aligned product system $\mathbf{Y}^\mathcal{NT}$ from Proposition~\ref{Y is a product system} is defined. Let $(p,q)\in P\rtimes_\alpha Q$ and $T\in \mathcal{K}_A\left(\mathbf{Z}_{(p,q)}\right)$. Then
\[
{\varphi^{\mathcal{NT}}}^{\left(\left(p,q\right)\right)}(T)=i_{\mathbf{Y}^\mathcal{NT}}^{(q)}\Big(M_{i_\mathbf{Z}^{((p,q))}(T)}\Big).
\]
\begin{proof}
It suffices to prove the result when $T$ is a rank one operator. To this end, fix $z,w\in \mathbf{Z}_{(p,q)}$. Lemma~\ref{compact operators in Y}, says that $\Theta_{i_{\mathbf{Z}_{(p,q)}}(z),i_{\mathbf{Z}_{(p,q)}}(w)}=M_{i_{\mathbf{Z}_{(p,q)}}(z)i_{\mathbf{Z}_{(p,q)}}(w)^*}=M_{i_{\mathbf{Z}}^{((p,q))}(\Theta_{z,w})}$, and so 
\begin{align*}
\varphi^{((p,q))}(\Theta_{z,w})
&=i_{\mathbf{Y}_q^\mathcal{NT}}\big(i_{\mathbf{Z}_{(p,q)}}(z)\big)i_{\mathbf{Y}_q}\big(i_{\mathbf{Z}_{(p,q)}}(w)\big)^*\\
&=i_{\mathbf{Y}^\mathcal{NT}}^{(q)}\Big(\Theta_{i_{\mathbf{Z}_{(p,q)}}(z),i_{\mathbf{Z}_{(p,q)}}(w)}\Big)
=i_{\mathbf{Y}^\mathcal{NT}}^{(q)}\Big(M_{i_{\mathbf{Z}}^{((p,q))}(\Theta_{z,w})}\Big).
\qedhere
\end{align*}
\end{proof}
\end{lem}

\begin{prop}
\label{varphi is Nica covariant}
Suppose that $G$ is an amenable group so that the compactly aligned product system $\mathbf{Y}^\mathcal{NT}$ from Proposition~\ref{Y is a product system} is defined. Then the representation $\varphi^{\mathcal{NT}}$ is Nica covariant, and so there exists a homomorphism $\Omega^{\mathcal{NT}}:\mathcal{NT}_\mathbf{Z}\rightarrow \mathcal{NT}_\mathbf{Y}$ such that 
\[
\Omega^{\mathcal{NT}}\circ i_{\mathbf{Z}_{(p,q)}}=\varphi_{(p,q)}^{\mathcal{NT}}=i_{\mathbf{Y}_q^\mathcal{NT}}\circ i_{\mathbf{Z}_{(p,q)}}
\]
for each $(p,q)\in P\rtimes_\alpha Q$.
\begin{proof}
Fix $S\in \mathcal{K}_A\left(\mathbf{Z}_{(p,q)}\right)$ and $T\in \mathcal{K}_A\left(\mathbf{Z}_{(s,t)}\right)$. Using Lemma~\ref{rep on compacts} we see that
\begin{align}
\label{first calc for NC}
{\varphi^{\mathcal{NT}}}^{((p,q))}(S){\varphi^{\mathcal{NT}}}^{((s,t))}(T)=i_{\mathbf{Y}^\mathcal{NT}}^{(q)}\Big(M_{i_\mathbf{Z}^{((p,q))}(S)}\Big)i_{\mathbf{Y}^\mathcal{NT}}^{(t)}\Big(M_{i_\mathbf{Z}^{((s,t))}(T)}\Big).
\end{align}
Since $i_{\mathbf{Y}^\mathcal{NT}}$ is Nica covariant, this is zero if $q\vee t=\infty$. If $q\vee t<\infty$, then \eqref{first calc for NC} gives
\begin{equation}
\label{second calc for NC}
\begin{aligned}
&{\varphi^{\mathcal{NT}}}^{((p,q))}(S){\varphi^{\mathcal{NT}}}^{((s,t))}(T)\\
&\qquad=
i_{\mathbf{Y}^\mathcal{NT}}^{(q\vee t)}\left(\iota_q^{q\vee t}\left(M_{i_\mathbf{Z}^{((p,q))}(S)}\right)\iota_t^{q\vee t}\left(M_{i_\mathbf{Z}^{((s,t))}(T)}\right)\right)\\
&\qquad=
i_{\mathbf{Y}^\mathcal{NT}}^{(q\vee t)}\left(M_{i_\mathbf{Z}^{((p,q))}(S)}M_{i_\mathbf{Z}^{((s,t))}(T)}\right)\\
&\qquad=
i_{\mathbf{Y}^\mathcal{NT}}^{(q\vee t)}\left(M_{i_\mathbf{Z}^{((p,q))}(S)i_\mathbf{Z}^{((s,t))}(T)}\right).
\end{aligned}
\end{equation}
Since $i_\mathbf{Z}$ is also Nica covariant and $(G\rtimes_\alpha H, P\rtimes_\alpha Q)$ has the product order, if $p\vee s=\infty$, then the last line is zero. If $p\vee s<\infty$, then another application of Lemma~\ref{rep on compacts} shows that \eqref{second calc for NC} yields 
\begin{align*}
{\varphi^{\mathcal{NT}}}^{((p,q))}(S){\varphi^{\mathcal{NT}}}^{((s,t))}(T)
&=
i_{\mathbf{Y}^\mathcal{NT}}^{(q\vee t)}\Big(M_{i_\mathbf{Z}^{((p\vee s,q\vee t))}\left(\iota_{(p,q)}^{(p\vee s,q\vee t)}(S)\iota_{(s,t)}^{(p\vee s,q\vee t)}(T)\right)}\Big)\\
&=
{\varphi^{\mathcal{NT}}}^{((p\vee s,q\vee t))}\Big(\iota_{(p,q)}^{(p\vee s,q\vee t)}(S)\iota_{(s,t)}^{(p\vee s,q\vee t)}(T)\Big).
\end{align*}
Thus,
\begin{align*}
&{\varphi^{\mathcal{NT}}}^{((p,q))}(S)\varphi^{((s,t))}(T)\\
&\qquad=
\begin{cases}
{\varphi^{\mathcal{NT}}}^{((p\vee s,q\vee t))}\Big(\iota_{(p,q)}^{(p\vee s,q\vee t)}(S)\iota_{(s,t)}^{(p\vee s,q\vee t)}(T)\Big) & \text{if $(p,q)\vee (s,t)<\infty$}\\
0 & \text{otherwise,}
\end{cases}
\end{align*}
and so $\varphi$ is Nica covariant. 
\end{proof}
\end{prop}

Next, we construct a homomorphism from $\mathcal{NT}_{\mathbf{Y}^\mathcal{NT}}$ to $\mathcal{NT}_\mathbf{Z}$ by exhibiting a Nica covariant representation of $\mathbf{Y}^\mathcal{NT}$ in $\mathcal{NT}_\mathbf{Z}$.

\begin{prop}
\label{varphi' is Nica covariant}
Suppose that $G$ is an amenable group so that the compactly aligned product system $\mathbf{Y}^\mathcal{NT}$ from Proposition~\ref{Y is a product system} is defined. For each $q\in Q$, let $\varphi_q'^{\mathcal{NT}}$ be the inclusion of $\mathbf{Y}_q^\mathcal{NT}$ in $\mathcal{NT}_\mathbf{Z}$. Then $\varphi'^{\mathcal{NT}}$ is a Nica covariant representation of $\mathbf{Y}^\mathcal{NT}$. Hence, there exists a homomorphism $\Omega'^{\mathcal{NT}}:\mathcal{NT}_{\mathbf{Y}^\mathcal{NT}}\rightarrow \mathcal{NT}_\mathbf{Z}$ such that 
\[
\Omega'^{\mathcal{NT}}\circ i_{\mathbf{Y}_q^\mathcal{NT}}={\varphi'_q}^{\mathcal{NT}}
\]
for each $q\in Q$.
\begin{proof}
It is trivial to check that $\varphi'^{\mathcal{NT}}$ is a representation. We now show that $\varphi'^{\mathcal{NT}}$ is Nica covariant. If $q\in Q$ and $b\in \mathcal{NT}_\mathbf{Z}^q$, then Lemma~\ref{compact operators in Y} tells us that $M_b\in \mathcal{K}_{\mathcal{NT}_\mathbf{X}}\left(\mathbf{Y}_q^\mathcal{NT}\right)$. We claim that ${\varphi'^{\mathcal{NT}}}^{(q)}(M_b)=b$. To see this, observe that if $z\in \mathbf{Z}_{(p,q)}$ and $w\in\mathbf{Z}_{(r,q)}$, then
\begin{align*}
{\varphi'^{\mathcal{NT}}}^{(q)}\left(M_{i_{\mathbf{Z}_{(p,q)}}(z)i_{\mathbf{Z}_{(r,q)}}(w)^*}\right)
&={\varphi'^{\mathcal{NT}}}^{(q)}\left(\Theta_{i_{\mathbf{Z}_{(p,q)}}(z),i_{\mathbf{Z}_{(r,q)}}(w)}\right)\\
&={\varphi'_q}^{\mathcal{NT}}\big(i_{\mathbf{Z}_{(p,q)}}(z)\big){\varphi'_q}^{\mathcal{NT}}\big(i_{\mathbf{Z}_{(r,q)}}(w)\big)^*\\
&=i_{\mathbf{Z}_{(p,q)}}(z)i_{\mathbf{Z}_{(r,q)}}(w)^*.
\end{align*}
Now fix $\Theta_{y,w}\in \mathcal{K}_{\mathcal{NT}_\mathbf{X}}\left(\mathbf{Y}_q^\mathcal{NT}\right)$ and $\Theta_{u,v}\in \mathcal{K}_{\mathcal{NT}_\mathbf{X}}\left(\mathbf{Y}_t^\mathcal{NT}\right)$. Making use of Lemma~\ref{Q indexed subalgebras of NO_Z} and Proposition~\ref{Y is compactly aligned}, we see that
\begin{align*}
{\varphi'^{\mathcal{NT}}}^{(q)}(\Theta_{y,w}){\varphi'^{\mathcal{NT}}}^{(t)}(\Theta_{u,v})
=\varphi'^{(q)}\left(M_{yw^*}\right)\varphi'^{(t)}\left(M_{uv^*}\right)
=yw^*uv^*,
\end{align*}
which is zero if $q\vee t=\infty$, and equal to 
\begin{align*}
{\varphi'^{\mathcal{NT}}}^{(q\vee t)}\left(M_{yw^*uv^*}\right)
&={\varphi'^{\mathcal{NT}}}^{(q\vee t)}\left(M_{yw^*} M_{uv^*}\right) \\
&={\varphi'^{\mathcal{NT}}}^{(q\vee t)}\left(\iota_q^{q\vee t}(\Theta_{y,w})\iota_t^{q\vee t}(\Theta_{u,v})\right) 
\end{align*}
if $q\vee t<\infty$. Hence, $\varphi'$ is Nica covariant. 
\end{proof}
\end{prop}

Putting all of this together, we get the following. 

\begin{thm}
\label{isomorphisms of Nica--Toeplitz algebras}
Suppose that $G$ is amenable so that the compactly aligned product system $\mathbf{Y}^\mathcal{NT}$ from Proposition~\ref{Y is a product system} is defined. Then the homomorphisms $\Omega^{\mathcal{NT}}:\mathcal{NT}_\mathbf{Z}\rightarrow \mathcal{NT}_{\mathbf{Y}^\mathcal{NT}}$ and $\Omega'^{\mathcal{NT}}:\mathcal{NT}_{\mathbf{Y}^\mathcal{NT}}\rightarrow \mathcal{NT}_\mathbf{Z}$ are mutually inverse. Thus, $\mathcal{NT}_\mathbf{Z}\cong \mathcal{NT}_{\mathbf{Y}^\mathcal{NT}}$. 
\begin{proof}
We begin by showing that $\Omega'^{\mathcal{NT}}\circ \Omega^{\mathcal{NT}}=\mathrm{id}_{\mathcal{NT}_\mathbf{Z}}$. For $(p,q)\in P\rtimes_\alpha Q$, we have
\begin{align*}
\left(\Omega'^{\mathcal{NT}}\circ \Omega^{\mathcal{NT}}\right)\circ i_{\mathbf{Z}_{(p,q)}}
=
\Omega'^{\mathcal{NT}}\circ i_{\mathbf{Y}_q^\mathcal{NT}}\circ  i_{\mathbf{Z}_{(p,q)}}
=
i_{\mathbf{Z}_{(p,q)}}.
\end{align*}
Since $\mathcal{NT}_\mathbf{Z}$ is generated by the image of $i_\mathbf{Z}$ and $\Omega'^{\mathcal{NT}}\circ \Omega^{\mathcal{NT}}$ is a homomorphism, we conclude that $\Omega'^{\mathcal{NT}}\circ \Omega^{\mathcal{NT}}=\mathrm{id}_{\mathcal{NT}_\mathbf{Z}}$. 

Next we check that $\Omega^{\mathcal{NT}}\circ \Omega'^{\mathcal{NT}}=\mathrm{id}_{\mathcal{NT}_{\mathbf{Y}^\mathcal{NT}}}$. For any $p\in P$,
\[
\left(\Omega^{\mathcal{NT}}\circ \Omega'^{\mathcal{NT}}\right)\circ \big(i_{\mathbf{Y}_{e_H}^\mathcal{NT}}\circ i_{\mathbf{Z}_{\left(p,e_H\right)}}\big)=\Omega^{\mathcal{NT}} \circ i_{\mathbf{Z}_{\left(p,e_H\right)}}=i_{\mathbf{Y}_{e_H}^\mathcal{NT}}\circ i_{\mathbf{Z}_{\left(p,e_H\right)}}.
\]
Since $\mathbf{Y}_{e_H}^\mathcal{NT}=\mathcal{NT}_\mathbf{X}\cong \phi_{\mathbf{X}}^{\mathcal{NT}}\left(\mathcal{NT}_\mathbf{X}\right)$ is generated by the images of each $i_{\mathbf{Z}_{\left(p,e_H\right)}}$ we conclude that $\left(\Omega^{\mathcal{NT}}\circ \Omega'^{\mathcal{NT}}\right)\circ i_{\mathbf{Y}_{e_H}^\mathcal{NT}}=i_{\mathbf{Y}_{e_H}^\mathcal{NT}}$. Now let $q\in Q\setminus \{e_H\}$ and $z\in \mathbf{Z}_{(m,q)}$, $w\in  \mathbf{Z}_{\left(n,e_H\right)}$. Since $\Omega'^{\mathcal{NT}}\circ i_{\mathbf{Y}_q^\mathcal{NT}}$ is the inclusion map and $\Omega^{\mathcal{NT}}$ is multiplicative, we see that
\begin{align*}
\big(\big(\Omega^{\mathcal{NT}}\circ \Omega'^{\mathcal{NT}}\big)\circ i_{\mathbf{Y}_q^\mathcal{NT}}\big)&\big(i_{\mathbf{Z}_{(m,q)}}(z)i_{\mathbf{Z}_{\left(n,e_H\right)}}(w)^*\big)\\
&=\Omega^{\mathcal{NT}}\big(i_{\mathbf{Z}_{(m,q)}}(z)\big)\Omega^{\mathcal{NT}}\big(i_{\mathbf{Z}_{\left(n,e_H\right)}}(w)\big)^*\\
&=i_{\mathbf{Y}_q^\mathcal{NT}}\big(i_{\mathbf{Z}_{(m,q)}}(z)\big)i_{\mathbf{Y}_{e_H}^\mathcal{NT}}\big(i_{\mathbf{Z}_{\left(n,e_H\right)}}(w)^*\big)\\
&=i_{\mathbf{Y}_q^\mathcal{NT}}\big(i_{\mathbf{Z}_{(m,q)}}(z)i_{\mathbf{Z}_{\left(n,e_H\right)}}(w)^*\big).
\end{align*}
Since $\mathbf{Y}_q^\mathcal{NT}=\cspan\big\{i_{\mathbf{Z}_{(m,q)}}\big(\mathbf{Z}_{(m,q)}\big)i_{\mathbf{Z}_{\left(n,e_H\right)}}\big(\mathbf{Z}_{\left(n,e_H\right)}\big)^*:m,n\in P\big\}$, and $i_{\mathbf{Y}^\mathcal{NT}}$ is linear and isometric, we see that $\left(\Omega^{\mathcal{NT}}\circ \Omega'^{\mathcal{NT}}\right)\circ i_{\mathbf{Y}_q^\mathcal{NT}}=i_{\mathbf{Y}_q^\mathcal{NT}}$ for $q\in Q\setminus \{e_H\}$. Since $\mathcal{NT}_{\mathbf{Y}^\mathcal{NT}}$ is generated by the image of $i_{\mathbf{Y}^\mathcal{NT}}$ and $\Omega^{\mathcal{NT}}\circ \Omega'^{\mathcal{NT}}$ is a homomorphism, we conclude that $\Omega^{\mathcal{NT}}\circ \Omega'^{\mathcal{NT}}=\mathrm{id}_{\mathcal{NT}_{\mathbf{Y}^\mathcal{NT}}}$.
\end{proof}
\end{thm}

There are a number of results in the literature relating the coefficient algebra of a Hilbert bimodule to its Toeplitz algebra. By decomposing $\N^k$ as $\N\times \N^{k-1}$, Theorem~\ref{isomorphisms of Nica--Toeplitz algebras} enables us to view the Nica--Toeplitz algebra of a compactly aligned product system over $\N^k$ as a $k$-times iterated Toeplitz algebra, and so immediately extend these results to compactly aligned product systems over $\N^k$. 

\begin{cor}
\label{extending literature results for NT algebras}
Let $\mathbf{Z}$ be a compactly aligned product system over $\N^k$ with coefficient algebra $A$. Then
\begin{enumerate}[label=\upshape(\roman*)]
\item The homomorphism $i_{\mathbf{Z}_0}$ induces an isomorphism between $K_*(A)$ and $K_*(\mathcal{NT}_\mathbf{Z})$\/{\rm ;}
\item If $A$ is separable and the fibres $\{\mathbf{Z}_{e_i}:1\leq i\leq k\}$ are countably generated as Hilbert $A$-modules, then $i_{\mathbf{Z}_0}$ induces a $KK$-equivalence between $A$ and $\mathcal{NT}_\mathbf{Z}$\/{\rm ;}
\item $A$ is exact if and only if $\mathcal{NT}_\mathbf{Z}$ is exact\/{\rm ;} and
\item $A$ is nuclear if and only if $\mathcal{NT}_\mathbf{Z}$ is nuclear. 
\end{enumerate}
\begin{proof}
Part (i) follows from \cite[Proposition~8.2]{katsura}.
Part (ii) follows from \cite[Theorem~4.4]{pimsner} (the hypothesis that $A$ is separable and each $\mathbf{Z}_{e_i}$ is countably generated ensures that the Nica--Toeplitz algebra we get at each step of the iterative procedure is separable). 
Part (iii) follows from \cite[Theorem~7.1]{katsura} and part (iv) from \cite[Theorem~7.2]{katsura}.
\end{proof}
\end{cor}

\section{Iterating the Cuntz--Nica--Pimsner construction}
\label{factorising CNP alg}

We would like to be able to replicate our work from Section~\ref{factorising NT alg} using Cuntz--Nica--Pimsner algebras in place of Nica--Toeplitz algebras. Specifically, can we construct a product system $\mathbf{Y}^\mathcal{NO}$ over $(H,Q)$ with coefficient algebra $\mathcal{NO}_\mathbf{X}$ such that the Cuntz--Nica--Pimsner algebras of $\mathbf{Z}$ and $\mathbf{Y}^\mathcal{NO}$ are isomorphic? 

The basic idea is to extend the interior portions of the commutative diagrams in Figures~\ref{NO isomorphism} and \ref{NO homomorphism} by applying the canonical quotient homomorphisms $q_\mathbf{X}:\mathcal{NT}_\mathbf{X}\rightarrow \mathcal{NO}_\mathbf{X}$ and $q_\mathbf{Z}:\mathcal{NT}_\mathbf{Z}\rightarrow \mathcal{NO}_\mathbf{Z}$ at the relevant places. Since every Cuntz--Nica--Pimsner covariant representation is by definition a Nica covariant representation, many of the results from Section~\ref{factorising NT alg} that we require work exactly as before. Sometimes we will impose additional hypotheses to get things to work, and unsurprisingly, our proofs often become more complicated. We will add in hypotheses as and when needed, since we are not sure if they are necessary, in the hope that future work may be able to relax/remove them. 

We now list the results that we will prove in Section~\ref{factorising CNP alg}, and summarise the various spaces and maps that we will be working. As before, we hope that if readers lose track of what is going on, they will be able to return to this point for assistance.

\begin{enumerate}
\item[(1)] In Proposition~\ref{inclusion of NO_X in NO_Z} we show that the inclusion of $\mathbf{X}$ in $\mathbf{Z}$ induces a homomorphism $\phi_{\mathbf{X}}^\mathcal{NO}$ from $\mathcal{NO}_\mathbf{X}$ to $\mathcal{NO}_\mathbf{Z}$ such that $\phi_{\mathbf{X}}^\mathcal{NO}\circ q_\mathbf{X}=q_\mathbf{Z}\circ \phi_{\mathbf{X}}^\mathcal{NT}$.
\item[(2)] In Proposition~\ref{injectivity of inclusion of NO_X in NO_Z} we show that the homomorphism $\phi_{\mathbf{X}}^\mathcal{NO}$ is injective. 
\item[(3)] In Proposition~\ref{existence of fibres in NO product system in NO_Z}, we use $\phi_{\mathbf{X}}^\mathcal{NT}$ to construct a collection of Hilbert $\mathcal{NO}_\mathbf{X}$-bimodules $\big\{\mathbf{Y}_q^\mathcal{NO}:q\in Q\big\}$ inside $\mathcal{NO}_\mathbf{Z}$. 
\item[(4)] In Proposition~\ref{existence of NO product system in NO_Z}, we show that $\mathbf{Y}^\mathcal{NO}:=\bigsqcup_{q\in Q} \mathbf{Y}_q^\mathcal{NO}$ is a compactly aligned product system over $(H,Q)$ with coefficient algebra $\mathcal{NO}_\mathbf{X}$. 
\item[(5)] In Proposition~\ref{injectivity of induced actions}, we find sufficient conditions for $\mathcal{NO}_\mathbf{X}$ to act faithfully on each fibre of $\mathbf{Y}^\mathcal{NO}$.
\end{enumerate}

Once we have the product system $\mathbf{Y}^\mathcal{NO}$, we prove that $\mathcal{NO}_{\mathbf{Y}^\mathcal{NO}}\cong \mathcal{NO}_\mathbf{Z}$. 

\begin{enumerate}
\item[(6)] In Proposition~\ref{varphi is CP covariant} we use the universal Cuntz--Nica--Pimsner covariant representations of $\mathbf{Z}$ and $\mathbf{Y}^\mathcal{NO}$ to construct a Cuntz--Nica--Pimsner representation $\varphi^{\mathcal{NO}}$ of $\mathbf{Z}$ in $\mathcal{NO}_{\mathbf{Y}^\mathcal{NO}}$. Using the universal property of $\mathcal{NO}_\mathbf{Z}$ we get a homomorphism $\Omega^{\mathcal{NO}}:\mathcal{NO}_\mathbf{Z}\rightarrow \mathcal{NO}_{\mathbf{Y}^\mathcal{NO}}$ such that $\Omega^{\mathcal{NO}}\circ j_{\mathbf{Z}_{(p,q)}}=j_{\mathbf{Y}_q^\mathcal{NO}}\circ j_{\mathbf{Z}_{(p,q)}}$.
\item[(7)] In Proposition~\ref{varphi' is CNP covariant} we find sufficient conditions for the inclusion of $\mathbf{Y}^\mathcal{NO}$ in $\mathcal{NO}_\mathbf{Z}$ to be a Cuntz--Nica--Pimsner covariant representation. The universal property of $\mathcal{NO}_{\mathbf{Y}^\mathcal{NO}}$ gives a homomorphism $\Omega'^{\mathcal{NO}}: \mathcal{NO}_{\mathbf{Y}^\mathcal{NO}}\rightarrow \mathcal{NO}_\mathbf{Z}$ such that $\Omega'^{\mathcal{NO}}\circ j_{\mathbf{Y}_q^\mathcal{NO}}$ is the inclusion of $\mathbf{Y}_q^\mathcal{NO}$ in $\mathcal{NO}_\mathbf{Z}$.
\item[(8)] In Theorem~\ref{isomorphism of Cuntz Nica Pimsner algebras} we prove that $\Omega^\mathcal{NO}$ and $\Omega'^{\mathcal{NO}}$ are mutually inverse isomorphisms. 
\end{enumerate}

In summary, we will show that for every $p\in P$ and $q\in Q$, the maps in the outer portions of Figures~\ref{NO isomorphism} and~\ref{NO homomorphism} exist and make the diagrams commutative. 

\subsection{Constructing the product systems}
\label{Constructing the CNP product system}

The first step is to check that the inclusion of $\mathbf{X}$ in $\mathbf{Z}$ induces a homomorphism from $\mathcal{NO}_\mathbf{X}$ to $\mathcal{NO}_\mathbf{Z}$, analogous to Proposition~\ref{inclusion of NT_X in NT_Z}. To do this we need a couple of lemmas. 

\begin{lem}
\label{useful factorisation result}
Let $T\in \mathcal{K}_A(\mathbf{X}_p)$ and $(s,t)\in P\rtimes_\alpha Q$. If $s\neq e_G$, then 
\begin{equation}
\label{useful factorisation result equation}
\begin{aligned}
\iota_{(p,e)}^{(s,t)}(T)=M_{\left(s,e_H\right),\left(e_G,t\right)}^\mathbf{Z}\circ \big(\iota_p^s(T)\otimes_A \mathrm{id}_{\mathbf{Z}_{\left(e_G,t\right)}}\big)\circ \big(M_{\left(s,e_H\right),\left(e_G,t\right)}^\mathbf{Z}\big)^{-1}. 
\end{aligned}
\end{equation}
\begin{proof}
If $y\in\mathbf{Z}_{\left(s,e_H\right)}$ and $z\in \mathbf{Z}_{\left(e_G,t\right)}$, then
\begin{align*}
&\Big(M_{\left(s,e_H\right),\left(e_G,t\right)}^\mathbf{Z}\circ \big(\iota_{e_G}^s(T)\otimes_A \mathrm{id}_{\mathbf{Z}_{\left(e_G,t\right)}}\big)\circ \big(M_{\left(s,e_H\right),\left(e_G,t\right)}^\mathbf{Z}\big)^{-1}\Big)(yz)\\
&\qquad=\iota_{e_G}^s(T)(y)z
=\phi_s(T)(y)z
=\phi_{(s,t)}(T)(yz)
=\iota_{\left(e_G,e_H\right)}^{(s,t)}(T)(yz).
\end{align*}
Since $M_{\left(s,e_H\right),\left(e_G,t\right)}^\mathbf{Z}\circ \big(\iota_{e_G}^s(T)\otimes_A \mathrm{id}_{\mathbf{Z}_{\left(e_G,t\right)}}\big)\circ \big(M_{\left(s,e_H\right),\left(e_G,t\right)}^\mathbf{Z}\big)^{-1}$ and $\iota_{\left(e_G,e_H\right)}^{(s,t)}(T)$ are linear and continuous, and everything in  $\mathbf{Z}_{(s,t)}$ can be approximated by linear combinations of elements from $\{yz:y\in\mathbf{Z}_{\left(s,e_H\right)}, \ z\in \mathbf{Z}_{\left(e_G,t\right)}\}$ (as $s\neq e_G$), we see that Equation~\ref{useful factorisation result equation} holds when $p=e_G$. 

Now suppose that $p\neq e_G$. If $p\not\leq s$, then $\left(p,e_H\right)\not\leq (s,t)$ and so both sides of Equation~\ref{useful factorisation result equation} are zero. So suppose that $p\leq s$. For any $x\in \mathbf{Z}_{\left(p,e_H\right)}$, $y\in\mathbf{Z}_{\left(p^{-1}s,e_H\right)}$, and $z\in \mathbf{Z}_{\left(e_G,t\right)}$, we see that
\begin{align*}
&\Big(M_{\left(s,e_H\right),\left(e_G,t\right)}^\mathbf{Z}\circ \big(\iota_p^s(T)\otimes_A \mathrm{id}_{\mathbf{Z}_{\left(e_G,t\right)}}\big)\circ \big(M_{\left(s,e_H\right),\left(e_G,t\right)}^\mathbf{Z}\big)^{-1}\Big)(xyz)\\
&\qquad=\iota_p^s(T)(xy)z
=\left((Tx)y\right)z
=(Tx)(yz)
=\iota_{\left(p,e_H\right)}^{(s,t)}(T)(xyz).
\end{align*}
Since $\mathbf{Z}_{(s,t)}=\cspan\{xyz:x\in \mathbf{Z}_{\left(p,e_H\right)},y\in\mathbf{Z}_{\left(p^{-1}s,e_H\right)},z\in \mathbf{Z}_{\left(e_G,t\right)}\}$, we conclude that \eqref{useful factorisation result equation} holds when $p\neq e_G$ as well. 
\end{proof}
\end{lem}

\begin{lem}
\label{not sure what to call it}
Suppose that $A$ acts faithfully on each fibre of $\mathbf{X}$ and each $\widetilde{\phi}_{(p,q)}$ is injective. Fix a finite set $F\subseteq P$ and a collection of compact operators $\left\{T_p\in \mathcal{K}_A(\mathbf{X}_p):p\in F\right\}$ such that 
\[
\sum_{p\in F} \iota_p^s(T_p)=0\in \mathcal{L}_A(\mathbf{X}_s) \quad \text{for large $s\in P$}.\]
Then
\[
\sum_{p\in F} \widetilde{\iota}_{(p,e_H)}^{\,\,(s,t)}(T_p)=0\in \mathcal{L}_A\big(\widetilde{\mathbf{Z}}_{(s,t)}\big) \quad \text{for large $(s,t)\in P\rtimes_\alpha Q$}.
\]
\begin{proof}
Let $(u,v)\in P\rtimes_\alpha Q$. Hence, we can choose $r\geq u$ such that whenever $s\geq r$, we have 
\[
\sum_{p\in F} \iota_p^s(T_p)=0\in \mathcal{L}_A(\mathbf{X}_s).
\]
Since $(G\rtimes_\alpha H, P\rtimes_\alpha Q)$ has the product order, we know that $(r,v)\geq (u,v)$. We claim that for any $(s,t)\geq (r,v)$, we have  
\[
\sum_{p \in F} \widetilde{\iota}_{(p,e_H)}^{\,\,(s,t)}(T_p)=0\in \mathcal{L}_A\big(\widetilde{\mathbf{Z}}_{(s,t)}\big).
\]
Fix $(s,t)\geq (r,v)$ and let $z\in \widetilde{\mathbf{Z}}_{(s,t)}$. We need to show that 
\[
\Big(\sum_{p \in F} \widetilde{\iota}_{(p,e_H)}^{\,\,(s,t)}(T_p)\Big)(z)
\]
is zero in $\widetilde{\mathbf{Z}}_{(s,t)}$, which is equivalent to showing that 
\[
\Big(\Big(\sum_{p \in F} \widetilde{\iota}_{(p,e_H)}^{\,\,(s,t)}(T_p)\Big)(z)\Big)_{(m,n)}
\]
is zero in $\mathbf{Z}_{(m,n)}$ for every $(m,n)\leq (s,t)$. With this in mind, let $(m,n)\leq (s,t)$. If $m\neq s$, then $\alpha_{n^{-1}}(m^{-1}s)\neq e_G$, and so 
\[
\left(e_G,e_H\right)<(\alpha_{n^{-1}}(m^{-1}s),e_H)\leq \left(\alpha_{n^{-1}}(m^{-1}s),n^{-1}t\right).
\]
As $A$ acts faithfully on each fibre of $\mathbf{X}$, it acts faithfully on $\mathbf{X}_{\alpha_{n^{-1}}(m^{-1}s)}=\mathbf{Z}_{(\alpha_{n^{-1}}(m^{-1}s),e_H)}$ in particular, and so
\begin{align*}
I_{(m,n)^{-1}(s,t)}
=I_{\left(\alpha_{n^{-1}}(m^{-1}s),n^{-1}t\right)}
=\bigcap_{\substack{\left(e_G,e_H\right)<(x,y)\\\leq \left(\alpha_{n^{-1}}(m^{-1}s),n^{-1}t\right)}}\ker\left(\phi_{(x,y)}\right)=\{0\}.
\end{align*}
Thus, if $m\neq s$, then $z_{(m,n)}\in \mathbf{Z}_{(m,n)}\cdot I_{(m,n)^{-1}(s,t)}=\{0\}$. Hence, for any $(m,n)\leq (s,t)$, we have
\begin{align*}
\Big(\Big(\sum_{p \in F} \widetilde{\iota}_{(p,e_H)}^{\,\,(s,t)}(T_p)\Big)(z)\Big)_{(m,n)}
=&\Big(\sum_{p \in F} \iota_{(p,e_H)}^{(m,n)}(T_p)\Big)\left(z_{(m,n)}\right)\\
=&
\begin{cases}
\left(\sum_{p \in F} \iota_{(p,e_H)}^{(s,n)}(T_p)\right)\left(z_{(s,n)}\right) & \text{if $m=s$}\\
0 & \text{otherwise.}
\end{cases}
\end{align*}
Thus it remains to show that $\big(\sum_{p \in F} \iota_{(p,e_H)}^{(s,n)}(T_p)\big)(z_{(s,n)})=0$.

If $s\neq e_G$, then Lemma~\ref{useful factorisation result} combined with the fact that $\sum_{p\in F} \iota_p^s(T_p)=0$, shows that $\sum_{p \in F}\iota_{\left(p,e_H\right)}^{(s,n)}(T_p)=0$. On the other hand, if $s=e_G$, then 
\begin{align*}
0
=\sum_{p\in F} \iota_p^s(T_p)
=
\begin{cases}
\iota_{e_G}^{e_G}\left(T_{e_G}\right) &\text{if $e_G\in F$}\\
0 &\text{otherwise}
\end{cases}
\quad=
\begin{cases}
T_{e_G} &\text{if $e_G\in F$}\\
0 &\text{otherwise},
\end{cases}
\end{align*}
and so either $e\not\in F$ or $T_{e_G}=0$. Thus,
\begin{align*}
\Big(\sum_{p \in F} \iota_{(p,e_H)}^{(s,n)}(T_p)\Big)\left(z_{(s,n)}\right)
&=
\Big(\sum_{p \in F} \iota_{(p,e_H)}^{\left(e_G,n\right)}(T_p)\Big)\left(z_{\left(e_G,n\right)}\right)\\
&=
\begin{cases}
\left(\phi_{\left(e_G,n\right)}\left(T_{e_G}\right)\right)\left(z_{\left(e_G,n\right)}\right) & \text{if $e_G\in F$}\\
0 & \text{otherwise}
\end{cases}
\end{align*}
is zero as required. We conclude that $\sum_{p \in F} \widetilde{\iota}_{(p,e_H)}^{\,\,(s,t)}(T_p)=0$ as claimed. 
\end{proof}
\end{lem}

\begin{prop}
\label{inclusion of NO_X in NO_Z}
Suppose that $A$ acts faithfully on each fibre of $\mathbf{X}$ and each $\widetilde{\phi}_{(p,q)}$ is injective. Then the inclusion of $\mathbf{X}$ in $\mathbf{Z}$ induces a homomorphism $\phi_{\mathbf{X}}^\mathcal{NO}:\mathcal{NO}_\mathbf{X}\rightarrow \mathcal{NO}_\mathbf{Z}$ such that $\phi_{\mathbf{X}}^\mathcal{NO}\left(j_\mathbf{X}(x)\right)=j_\mathbf{Z}(x)$ for each $x\in \mathbf{X}$. Thus, $\phi_{\mathbf{X}}^\mathcal{NO}\circ q_\mathbf{X}=q_\mathbf{Z}\circ \phi_{\mathbf{X}}^\mathcal{NT}$.
\begin{proof}
We need to exhibit a Cuntz--Nica--Pimsner covariant representation of $\mathbf{X}$ in $\mathcal{NO}_\mathbf{Z}$. For each $p\in P$, define $\psi_p^{\mathcal{NO}}: \mathbf{X}_p \rightarrow \mathcal{NO}_\mathbf{Z}$ by $\psi_p^{\mathcal{NO}}:=j_{\mathbf{Z}_{(p,e_H)}}$. Since $\psi_p^{\mathcal{NO}}=q_\mathbf{Z}\circ \psi_p^{\mathcal{NT}}$, where  $\psi^{\mathcal{NT}}$ is the Nica covariant representation of $\mathbf{X}$ from Proposition~\ref{inclusion of NT_X in NT_Z}, we see that $\psi^{\mathcal{NO}}$ is a Nica covariant representation of $\mathbf{X}$.

We now check that $\psi^\mathcal{NO}$ is Cuntz--Pimsner covariant. Suppose that $F\subseteq P$ is finite and $\left\{T_p\in \mathcal{K}_A(\mathbf{X}_p):p\in F\right\}$ is such that 
$\sum_{p\in F} \iota_p^s(T_p)=0\in \mathcal{L}_A(\mathbf{X}_s)$ for large $s\in P$. We need to show that $\sum_{p\in F} {\psi^{\mathcal{NO}}}^{(p)}(T_p)=0$. Since ${\psi^{\mathcal{NO}}}^{(p)}=j_\mathbf{Z}^{((p,e_H))}$ for each $p\in P$ and $j_\mathbf{Z}$ is Cuntz--Pimsner covariant, it suffices to prove that $\sum_{p\in F} \widetilde{\iota}_{(p,e_H)}^{\,\,(s,t)}(T_p)=0\in \mathcal{L}_A\big(\widetilde{\mathbf{Z}}_{(s,t)}\big)$ for large $(s,t)\in P\rtimes_\alpha Q$, which follows from Lemma~\ref{not sure what to call it}. Thus, $\psi^{\mathcal{NO}}$ is Cuntz--Pimsner covariant, and there is a homomorphism $\phi_{\mathbf{X}}^\mathcal{NO}:\mathcal{NO}_\mathbf{X}\rightarrow \mathcal{NO}_\mathbf{Z}$ such that $\phi_{\mathbf{X}}^\mathcal{NO}\left(j_\mathbf{X}(x)\right)=\psi^{\mathcal{NO}}(x)=j_\mathbf{Z}(x)$ for each $x\in \mathbf{X}$. Finally, since the image of $i_\mathbf{X}$ generates $\mathcal{NT}_\mathbf{X}$, and
\[
\phi_{\mathbf{X}}^\mathcal{NO}\circ q_\mathbf{X}\circ i_\mathbf{X}=\phi_{\mathbf{X}}^\mathcal{NO}\circ {j_\mathbf{X}=j_\mathbf{Z}|}_\mathbf{X}
=q_\mathbf{Z}\circ {i_\mathbf{Z}|}_\mathbf{X}=q_\mathbf{Z}\circ \phi_{\mathbf{X}}^\mathcal{NT}\circ i_\mathbf{X},
\]
we conclude that $\phi_{\mathbf{X}}^\mathcal{NO}\circ q_\mathbf{X}=q_\mathbf{Z}\circ \phi_{\mathbf{X}}^\mathcal{NT}$.
\end{proof}
\end{prop}

It is not clear whether the hypotheses for Proposition~\ref{inclusion of NO_X in NO_Z} can be relaxed. The hypothesis that each $\widetilde{\phi}_{(p,q)}$ is injective is used only to ensure that the Cuntz--Nica--Pimsner algebra $\mathcal{NO}_\mathbf{Z}$ exists. The assumption that $A$ acts faithfully on each fibre of $\mathbf{X}$ allows us to make use of Lemma~\ref{not sure what to call it}. If we tried to prove Lemma~\ref{not sure what to call it} without assuming that $A$ acts faithfully on each fibre of $\mathbf{X}$, we would need to prove that whenever $F\subseteq P$ is finite and $\left\{T_p\in \mathcal{K}_A(\mathbf{X}_p):p\in F\right\}$ is a collection of compact operators such that
$
\sum_{p\in F} \tilde{\iota}_p^{\,s}(T_p)
$
is zero in $\mathcal{L}_A\big(\widetilde{\mathbf{X}}_s\big)$ for large $s\in P$, then
$
\sum_{p\in F} \tilde{\iota}_{(p,e_H)}^{\,\,(s,t)}(T_p)
$
is zero in $\mathcal{L}_A\big(\widetilde{\mathbf{Z}}_{(s,t)}\big)$ for large $(s,t)\in P\rtimes_\alpha Q$. The following example using higher-rank graphs, shows that this need not be true. For the relevant background on higher-rank graphs and their associated Cuntz--Krieger algebras, see \cite{raeburn--sims--yeend}. For the relationship between product systems and higher-rank graphs see \cite{raeburn--sims} and \cite[Section~5.3]{sims--yeend}.

\begin{exam}
\label{non locally convex 2-graph}
Let $\Gamma$  be the $2$-graph consisting of three vertices $u,v,w$, one edge $\lambda$ of degree $(1,0)$, one edge $\mu$ of degree $(0,1)$, with ranges and sources $s(\lambda)=u$, $s(\mu)=w$, and $r(\lambda)=r(\mu)=v$. The associated compactly aligned product system $\mathbf{Z}\subseteq C^*(\Gamma)$ over $(\mathbb{Z}^2,\mathbb{N}^2)$ has fibres given by
\begin{align*}
A:=\mathbf{Z}_{(0,0)}&:=\cspan\{t_u,t_v,t_w\}=C^*(\{t_u,t_v,t_w\}), \quad\mathbf{Z}_{(1,0)}:=\cspan\{t_\lambda\},\\
\mathbf{Z}_{(0,1)}&:=\cspan\{t_\mu\}, \quad\mathbf{Z}_{(m,n)}:=\{0\} \quad\text{for all other $(m,n)\in \mathbb{N}^2$}.
\end{align*}
From this, we define another compactly aligned product system $\mathbf{X}$ over $(\mathbb{Z},\mathbb{N})$, whose fibres are given by $\mathbf{X}_n:=\mathbf{Z}_{(n,0)}$ for each $n\in \mathbb{N}$. Example~3.3 of \cite{sims--yeend} shows that for any $n\geq 1$ we have
\[
\widetilde{\mathbf{X}}_n=\cspan\big\{t_\nu:\nu\in \Gamma^{\leq (n,0)}\big\}=\cspan\{t_u,t_w,t_\lambda\},
\]
whilst for any $(n,m)\geq (1,1)$ we have 
\[
\widetilde{\mathbf{Z}}_{(n,m)}=\cspan\big\{t_\nu:\nu\in \Gamma^{\leq (n,m)}\big\}=\cspan\{t_u,t_v,t_\lambda,t_\mu\}.
\]
Consider $t_v\in A \cong \mathcal{K}_A\left(\mathbf{X}_0\right)$ and $\Theta_{t_\lambda,t_\lambda}\in \mathcal{K}_A\left(\mathbf{X}_1\right)$. Then 
\[
\tilde{\iota}_0^{\,n}(t_v)-\tilde{\iota}_1^{\,n}(\Theta_{t_\lambda,t_\lambda})=0\in \mathcal{L}_A\big(\widetilde{\mathbf{X}}_n\big) \quad \text{for any $n\geq 1$.}
\]
For any $(n,m)\geq (1,1)$, since $r(\mu)=v$, we have
\[
\tilde{\iota}_{(0,0)}^{\,(n,m)}(t_v)(t_\mu)=t_\mu.
\]
Since $\lambda$ and $\mu$ have no minimal common extensions, it follows that 
\[
\tilde{\iota}_{(1,0)}^{\,(n,m)}(\Theta_{t_\lambda,t_\lambda})(t_\mu)=t_\lambda t_\lambda^* t_\mu=0.
\]
Thus,
\[
\big(\tilde{\iota}_{(0,0)}^{\,(n,m)}(t_v)-\tilde{\iota}_{(1,0)}^{\,(n,m)}(\Theta_{t_\lambda,t_\lambda})\big)(t_\mu)=t_\mu,
\] 
and so 
\[
\tilde{\iota}_{(0,0)}^{\,(n,m)}(t_v)-\tilde{\iota}_{(1,0)}^{\,(n,m)}(\Theta_{t_\lambda,t_\lambda})\neq 0\in \mathcal{L}_A\big(\widetilde{\mathbf{Z}}_{(n,m)}\big).
\]
\end{exam} 

We should not be particularly surprised by Example~\ref{non locally convex 2-graph} --- the product system $\mathbf{Z}$ associated to $\Gamma$ contains no two dimensional information (since the higher-rank graph $\Gamma$ contains no paths of degree $(1,1)$). In \cite[Proposition~4.6]{k-graph iterated} we show that if the higher-rank graph is locally convex (see \cite[Definition~3.9]{raeburn--sims--yeend2}), then problems of this variety do not occur. As such, it would be interesting to see if there is a notion of local convexity for arbitrary product systems over $\mathbb{N}^k$ (or even more general quasi-lattice ordered groups) that ensures the inclusion of $\mathbf{X}$ in $\mathbf{Z}$ induces a homomorphism from $\mathcal{NO}_\mathbf{X}$ to $\mathcal{NO}_\mathbf{Z}$. 

Similar to Subsection~\ref{Constructing the product system}, to construct a product system sitting inside $\mathcal{NO}_\mathbf{Z}$ with coefficient algebra $\mathcal{NO}_\mathbf{X}$, we need to check that the homomorphism $\phi_{\mathbf{X}}^\mathcal{NO}:\mathcal{NO}_\mathbf{X}\rightarrow \mathcal{NO}_\mathbf{Z}$ is injective. The idea is to show that by restricting the canonical coaction of $G\rtimes_\alpha H$ on $\mathcal{NO}_\mathbf{Z}$ to the image of $\phi_{\mathbf{X}}^\mathcal{NO}$, we get a coaction of $G$ on $\phi_{\mathbf{X}}^\mathcal{NO}\left(\mathcal{NO}_\mathbf{X}\right)$. As such, we need to know when the group $C^*$-algebra of $G\rtimes_\alpha H$ contains a faithful copy of the group $C^*$-algebra of $G$.

\begin{lem}
\label{group algebra isomorphisms}
Let $G$ and $H$ be groups, and suppose that $\alpha:H \rightarrow \mathrm{Aut}(G)$ is a group homomorphism. Then there exists a homomorphism $\iota: C^*(G)\rightarrow C^*(G\rtimes_\alpha H)$ such that $\iota(i_G(g))=i_{G\rtimes_\alpha H}(g,e_H)$ for each $g\in G$ \rm{(}where $i_G$ and $i_{G\rtimes_\alpha H}$ are the universal unitary representations of $G$ and $G\rtimes_\alpha H$ respectively\rm{)}. If $G$ is amenable, then $\iota$ is injective. Thus, $C^*(G)\cong \cspan\{i_{G\rtimes_\alpha H}(G \rtimes_\alpha \{e_H\})\}\subseteq C^*(G\rtimes_\alpha H)$. 
\begin{proof}
It is straightforward to check that $g\mapsto i_{G\rtimes_\alpha H}((g,e_H))$ is a unitary representation of $G$ in $C^*(G\rtimes_\alpha H)$. The universal property of $C^*(G)$ then provides us with the homomorphism $\iota$. 

We now assume that $G$ is amenable. Hence,
\[
C^*(G)\cong C^*(\{T_g:g\in G\})\subseteq \mathcal{B}\left(\ell^2(G)\right),
\]
where the operator $T_g \in \mathcal{B}(\ell^2(G))$ is defined by 
$
T_g(f)(h):=f\left(g^{-1}h\right)
$
for $g\in G$, $f\in \ell^2(G)$, and $h\in G$. Similarly, the map $(g,h)\mapsto S_{(g,h)}\in \mathcal{B}(\ell^2(G\rtimes_\alpha H))$ where 
\[
S_{(g,h)}(f)(k,l):=f\big((g,h)^{-1}(k,l)\big)=f\big(\alpha_{h^{-1}}(g^{-1}k),h^{-1}l\big)
\]
for each $f\in \ell^2\left(G\rtimes_\alpha H\right)$ and $(k,l)\in G\rtimes_\alpha H$, is a unitary representation of $G\rtimes_\alpha H$. Clearly, if $f\in \ell^2(G)$, then the map $\tilde{f}:G\rtimes_\alpha H\rightarrow \mathbb{C}$ defined by 
$
\tilde{f}(g,h):=\delta_{h,e_H}f(g)
$
belongs to $\ell^2(G\rtimes_\alpha H)$ and $\norm{f}_{\ell^2(G)}=\big\|\tilde{f}\big\|_{\ell^2(G\rtimes_\alpha H)}$. 
Now let $F\subseteq G$ be a finite set. For any $f\in \ell^2(G)$ we have
\begin{align*}
\bigg\|\sum_{g\in F}S_{(g,e_H)}\big(\tilde{f}\big)\bigg\|_{\ell^2(G\rtimes_\alpha H)}^2
&=
\sum_{(k,h)\in G\rtimes_\alpha H}\bigg|\sum_{g\in F}S_{(g,e_H)}\big(\tilde{f}\big)(k,h)\bigg|^2\\
&=
\sum_{(k,h)\in G\rtimes_\alpha H}\bigg|\sum_{g\in F}\tilde{f}\left(g^{-1}k,h\right)\bigg|^2\\
&=\sum_{k\in G}\bigg|\sum_{g\in F}f\left(g^{-1}k\right)\bigg|^2\\
&=\sum_{k\in G}\bigg|\sum_{g\in F}T_g(f)(k)\bigg|^2\\
&=\bigg\|\sum_{g\in F}T_g(f)\bigg\|_{ \ell^2(G)}^2.
\end{align*}
Thus, 
\[
\bigg\|\sum_{g\in F}T_g\bigg\|_{ \mathcal{B}(\ell^2(G))}\leq  \bigg\|\sum_{g\in F}S_{(g,e_H)}\bigg\|_{\mathcal{B}(\ell^2(G\rtimes_\alpha H))}.
\]
Since $\iota$ is norm-decreasing, we have
\begin{align*}
\bigg\|\sum_{g\in F}i_G(g)\bigg\|_{C^*(G)}
&\geq \bigg\|\iota\bigg(\sum_{g\in F}i_G(g)\bigg)\bigg\|_{C^*(G\rtimes_\alpha H)}\\
&=\bigg\|\sum_{g\in F}i_{G\rtimes_\alpha H}(g,e_H)\bigg\|_{C^*(G\rtimes_\alpha H)}.
\end{align*}
Since $i_{G\rtimes_\alpha H}$ is the universal unitary representation of $G\rtimes_\alpha H$, we deduce that
\begin{align*}
\bigg\|\sum_{g\in F}i_G(g)\bigg\|_{C^*(G)}
&\geq\bigg\|\sum_{g\in F}S_{(g,e_H)}\bigg\|_{\mathcal{B}(\ell^2(G\rtimes_\alpha H))}\\
&\geq \bigg\|\sum_{g\in F}T_g\bigg\|_{\mathcal{B}(\ell^2(G))}
=\bigg\|\sum_{g\in F}i_G(g)\bigg\|_{C^*(G)},
\end{align*}
where the last equality follows from the amenability of $G$. Hence $\iota$ is isometric on the dense subspace $\mathrm{span}\left\{i_G(g):g\in G\right\}$ of $C^*(G)$. We conclude that $\iota$ is an isometry.
\end{proof}
\end{lem}

\begin{prop}
\label{injectivity of inclusion of NO_X in NO_Z}
Suppose that $A$ acts faithfully on each fibre of $\mathbf{X}$, and each $\widetilde{\phi}_{(p,q)}$ is injective, so that the homomorphism $\phi_{\mathbf{X}}^\mathcal{NO}:\mathcal{NO}_\mathbf{X}\rightarrow \mathcal{NO}_\mathbf{Z}$ from  Proposition~\ref{inclusion of NO_X in NO_Z} exists. If $G$ is an amenable group, then $\phi_{\mathbf{X}}^\mathcal{NO}$ is injective. 
\begin{proof}
Since $G$ is amenable, we can use \cite[Corollary~4.12]{carlsen--larsen--sims--vittadello} to show that $\phi_{\mathbf{X}}^\mathcal{NO}$ is injective. Firstly, we need to check that ${\phi_{\mathbf{X}}^\mathcal{NO}}|_{j_{\mathbf{X}_{e_G}(A)}}$ is injective. Suppose that $a\in A$ is such that $\phi_{\mathbf{X}}^\mathcal{NO}\big(j_{\mathbf{X}_{e_G}}(a)\big)=0$. Then $j_{\mathbf{Z}_{(e_G,e_H)}}(a)=0$, which forces $a=0$ since $j_\mathbf{Z}$ is isometric by \cite[Theorem~4.1]{sims--yeend}. Thus, ${\phi_{\mathbf{X}}^\mathcal{NO}}|_{j_{\mathbf{X}_{e_G}(A)}}$ is injective. 

Next, let $\nu_\mathbf{Z}:\mathcal{NO}_{\mathbf{Z}}\rightarrow \mathcal{NO}_{\mathbf{Z}}\otimes C^*(G\rtimes_\alpha H)$ denote the canonical gauge coaction of $G\rtimes_\alpha H$ on $\mathcal{NO}_{\mathbf{Z}}$. 
Since 
\begin{align*}
\nu_\mathbf{Z}\big(\phi_{\mathbf{X}}^\mathcal{NO}(j_{\mathbf{X}_p}(x))\big)
=\nu_\mathbf{Z}\big(j_{\mathbf{Z}_{\left(p,e_H\right)}}(x)\big)
&=j_{\mathbf{Z}_{\left(p,e_H\right)}}(x)\otimes i_{G\rtimes_\alpha H}\left(p,e_H\right)\\
&=\phi_{\mathbf{X}}^\mathcal{NO}\big(j_{\mathbf{X}_p}(x)\big)\otimes \iota\left(i_G(p)\right), 
\end{align*}
for any $x\in \mathbf{X}_p$, we can define $\beta:\phi_{\mathbf{X}}^\mathcal{NO}\left(\mathcal{NO}_\mathbf{X}\right)\rightarrow \phi_{\mathbf{X}}^\mathcal{NO}\left(\mathcal{NO}_\mathbf{X}\right)\otimes C^*(G)$ by 
\[
\beta:=\big(\mathrm{id}_{\phi_{\mathbf{X}}^\mathcal{NO}\left(\mathcal{NO}_\mathbf{X}\right)}\otimes \iota^{-1}\big)\circ \nu_\mathbf{Z}|_{\phi(\mathcal{NO}_{\mathbf{X}})}.
\]
We claim that $\beta$ is a coaction of $G$ on $\phi_{\mathbf{X}}^\mathcal{NO}\left(\mathcal{NO}_\mathbf{X}\right)$. 

Since $\nu_\mathbf{Z}$ and $\mathrm{id}_{\phi_{\mathbf{X}}^\mathcal{NO}\left(\mathcal{NO}_\mathbf{X}\right)}\otimes \iota^{-1}$ are injective homomorphisms, so is $\beta$. If $x\in \mathbf{X}_p$, then
\begin{align*}
\big(\big(\beta\otimes \mathrm{id}_{C^*(G)}\big)\circ \beta\big)&\big(\phi_{\mathbf{X}}^\mathcal{NO}(j_{\mathbf{X}_p}(x))\big)\\
&=\big(\beta\otimes \mathrm{id}_{C^*(G)}\big)\big(j_{\mathbf{Z}_{\left(p,e_H\right)}}(x)\otimes i_G(p)\big)\\
&=j_{\mathbf{Z}_{\left(p,e_H\right)}}(x)\otimes i_G(p)\otimes i_G(p)\\
&=\big(\mathrm{id}_{\phi_{\mathbf{X}}^\mathcal{NO}\left(\mathcal{NO}_\mathbf{X}\right)}\otimes \delta_G\big)\big(j_{\mathbf{Z}_{\left(p,e_H\right)}}(x)\otimes i_G(p)\big)\\
&=\big(\big(\mathrm{id}_{\phi_{\mathbf{X}}^\mathcal{NO}\left(\mathcal{NO}_\mathbf{X}\right)}\otimes \delta_G\big)\circ \beta\big)\big(\phi_{\mathbf{X}}^\mathcal{NO}(j_{\mathbf{X}_p}(x))\big).
\end{align*}
Since $\phi_{\mathbf{X}}^\mathcal{NO}\left(\mathcal{NO}_\mathbf{X}\right)$ is generated by the image of the homomorphism $\phi_{\mathbf{X}}^\mathcal{NO}\circ j_\mathbf{X}$, and both $\big(\beta\otimes \mathrm{id}_{C^*(G)}\big)\circ \beta$ and $\big(\mathrm{id}_{\phi_{\mathbf{X}}^\mathcal{NO}\left(\mathcal{NO}_\mathbf{X}\right)}\otimes \delta_G\big)\circ \beta$ are homomorphisms, we conclude that $\beta$ satisfies the coaction identity. For any $x\in \mathbf{X}_p$, $w\in \mathbf{X}_r$, and $g\in G$, we have 
\begin{align*}
\phi_{\mathbf{X}}^\mathcal{NO}&\big(j_{\mathbf{X}_p}(x)j_{\mathbf{X}_r}(w)^*\big)\otimes i_G(g)\\
&=\Big(\phi_{\mathbf{X}}^\mathcal{NO}\big(j_{\mathbf{X}_p}(x)j_{\mathbf{X}_r}(w)^*\big)\otimes i_G\big(pr^{-1}\big)\Big)\big(1_{\mathcal{M}\left(\phi_{\mathbf{X}}^\mathcal{NO}\left(\mathcal{NO}_\mathbf{X}\right)\right)}\otimes i_G\big(rp^{-1}g\big)\big)\\
&=\beta\Big(\phi_{\mathbf{X}}^\mathcal{NO}\big(j_{\mathbf{X}_p}(x)j_{\mathbf{X}_r}(w)^*\big)\Big)\big(1_{\mathcal{M}\left(\phi_{\mathbf{X}}^\mathcal{NO}\left(\mathcal{NO}_\mathbf{X}\right)\right)}\otimes i_G\big(rp^{-1}g\big)\big).
\end{align*}
Thus, 
\[
\phi_{\mathbf{X}}^\mathcal{NO}\left(\mathcal{NO}_\mathbf{X}\right)\otimes C^*(G)
=\cspan\big\{\beta\big(\phi_{\mathbf{X}}^\mathcal{NO}\left(\mathcal{NO}_\mathbf{X}\right)\big)\big(1_{\mathcal{M}\left(\phi_{\mathbf{X}}^\mathcal{NO}\left(\mathcal{NO}_\mathbf{X}\right)\right)}\otimes C^*(G)\big)\big\},
\]
and so $\beta$ is coaction nondegenerate. Finally, for any $x\in \mathbf{X}_p$, we have
\begin{align*}
\big(\beta\circ \phi_{\mathbf{X}}^\mathcal{NO}\big)\big(j_{\mathbf{X}_p}(x)\big)
&=\big(\mathrm{id}_{\phi_{\mathbf{X}}^\mathcal{NO}\left(\mathcal{NO}_\mathbf{X}\right)}\otimes \iota^{-1}\big)\big(\nu_\mathbf{Z}\big(j_{\mathbf{Z}_{\left(p,e_H\right)}}(x)\big)\big)\\
&=\big(\mathrm{id}_{\phi_{\mathbf{X}}^\mathcal{NO}\left(\mathcal{NO}_\mathbf{X}\right)}\otimes \iota^{-1}\big)\big(j_{\mathbf{Z}_{\left(p,e_H\right)}}(x)\otimes i_{G\rtimes_\alpha H}(p,e_H)\big)\\
&=\big(\big(\phi_{\mathbf{X}}^\mathcal{NO}\otimes \mathrm{id}_{C^*(G)}\big)\circ \nu_\mathbf{X}\big)\big(j_{\mathbf{X}_p}(x)\big).
\end{align*}
Since $\beta\circ \phi_{\mathbf{X}}^\mathcal{NO}$ and $\big(\phi_{\mathbf{X}}^\mathcal{NO}\otimes \mathrm{id}_{C^*(G)}\big)\circ \nu_\mathbf{X}$ are homomorphisms, and $j_\mathbf{X}$ generates $\mathcal{NO}_\mathbf{X}$, we see that $\beta\circ \phi_{\mathbf{X}}^\mathcal{NO}=\big(\phi_{\mathbf{X}}^\mathcal{NO}\otimes \mathrm{id}_{C^*(G)}\big)\circ \nu_\mathbf{X}$. By \cite[Corollary~4.12]{carlsen--larsen--sims--vittadello}, $\phi_{\mathbf{X}}^\mathcal{NO}$ is injective.
\end{proof}
\end{prop}

Almost identically to Subsection~\ref{Constructing the product system}, we can use the injective homomorphism $\phi_{\mathbf{X}}^\mathcal{NO}$ to construct a product system. Firstly, we construct a collection of Hilbert $\mathcal{NO}_\mathbf{X}$-bimodules $\big\{\mathbf{Y}_q^\mathcal{NO}:q\in Q\setminus \{e_H\}\big\}$. The idea is to make use of the collection of Hilbert $\mathcal{NT}_\mathbf{X}$-bimodules $\big\{\mathbf{Y}_q^\mathcal{NT}:q\in Q\setminus \{e_H\}\big\}$ defined in Propositions~\ref{definition of Y_q} and \ref{right actions on each Y_q} and apply the quotient maps $q_\mathbf{X}$ and $q_\mathbf{Z}$ at the appropriate places. 

For each $q\in Q\setminus \{e_H\}$, we let $\mathbf{Y}_q^\mathcal{NO}=q_\mathbf{Z}\left(\mathbf{Y}_q^\mathcal{NT}\right)$. We then show that there exists a right action of $\mathcal{NO}_\mathbf{X}$ on $\mathbf{Y}_q^\mathcal{NO}$, and a $\mathcal{NO}_\mathbf{X}$-valued inner product $\langle \cdot ,\cdot \rangle_{\mathcal{NO}_\mathbf{X}}^q$ on $\mathbf{Y}_q^\mathcal{NO}$, making the diagrams in Figure~\ref{NO right action and IP} commutative. 

\begin{figure}[h]
\centering
\begin{tikzpicture}[scale=1.5]
\node (A) at (0,0) {$\mathbf{Y}_q^\mathcal{NO}\times \mathcal{NO}_\mathbf{X}$};
\node (B) at (0,2) {$\mathbf{Y}_q^\mathcal{NT}\times \mathcal{NO}_\mathbf{X}$};
\node (C) at (2.75,0) {$\mathbf{Y}_q^\mathcal{NO}$};
\node (D) at (2.75,2) {$\mathbf{Y}_q^\mathcal{NT}$};
\node (E) at (4.25,0) {$\mathbf{Y}_q^\mathcal{NO}\times \mathbf{Y}_q^\mathcal{NO}$};
\node (F) at (4.25,2) {$\mathbf{Y}_q^\mathcal{NT}\times \mathbf{Y}_q^\mathcal{NT}$};
\node (G) at (7,0) {$\mathcal{NO}_\mathbf{X}$};
\node (H) at (7,2) {$\mathcal{NT}_\mathbf{X}$};
\path[->,font=\scriptsize]
(B) edge node[left]{$q_\mathbf{Z}\times q_\mathbf{X}$} (A)
(D) edge node[right]{$q_\mathbf{Z}$} (C)
(B) edge node[above]{$(y,a)\mapsto y\cdot a$} (D)
(A) edge node[below]{$(y,a)\mapsto y\cdot a$} (C)
(F) edge node[left]{$q_\mathbf{Z}\times q_\mathbf{Z}$} (E)
(H) edge node[right]{$q_\mathbf{X}$} (G)
(F) edge node[above]{$\langle \cdot ,\cdot \rangle_{\mathcal{NT}_\mathbf{X}}^q$} (H)
(E) edge node[below]{$\langle \cdot ,\cdot \rangle_{\mathcal{NO}_\mathbf{X}}^q$} (G);
\end{tikzpicture}
\caption{The Hilbert $\mathcal{NO}_\mathbf{X}$-module $\mathbf{Y}_q^\mathcal{NO}$}
\label{NO right action and IP}
\end{figure}

We also show that there exists a left action of $\mathcal{NO}_\mathbf{X}$, implemented by a homomorphism $\Phi_q^{\mathcal{NO}}: \mathcal{NO}_\mathbf{X}\rightarrow \mathcal{L}_{\mathcal{NO}_\mathbf{X}}\left(\mathbf{Y}_q^\mathcal{NO}\right)$, such that the diagram in Figure~\ref{left action and multiplication isomorphisms} is commutative. 

\begin{figure}[h]
\centering
\begin{tikzpicture}[scale=1.5]
\node (A) at (0,0) {$\mathcal{NO}_\mathbf{X}\times \mathbf{Y}_q^\mathcal{NO}$};
\node (B) at (0,2) {$\mathcal{NT}_\mathbf{X}\times \mathbf{Y}_q^\mathcal{NT}$};
\node (C) at (3.5,0) {$\mathbf{Y}_q^\mathcal{NO}$};
\node (D) at (3.5,2) {$\mathbf{Y}_q^\mathcal{NT}$};
\path[->,font=\scriptsize]
(B) edge node[left]{$q_\mathbf{X}\times q_\mathbf{Z}$} (A)
(D) edge node[right]{$q_\mathbf{Z}$} (C)
(B) edge node[above]{$(a,y)\mapsto \Phi_q^\mathcal{NT}(a)(y)$} (D)
(A) edge node[below]{$(a,y)\mapsto \Phi_q^\mathcal{NO}(a)(y)$} (C);
\end{tikzpicture}
\caption{The homomorphism $\Phi_q^{\mathcal{NO}}$ }
\label{left action and multiplication isomorphisms}
\end{figure}

\begin{prop}
\label{existence of fibres in NO product system in NO_Z}
Suppose that $G$ is an amenable group, $A$ acts faithfully on each fibre of $\mathbf{X}$, and each $\widetilde{\phi}_{(p,q)}$ is injective, so that the homomorphism $\phi_{\mathbf{X}}^\mathcal{NO}$ of Proposition~\ref{inclusion of NO_X in NO_Z} exists and is injective. For each $q\in Q\setminus \{e_H\}$, define
\begin{align*}
\mathbf{Y}_q^\mathcal{NO}:&=\cspan\big\{j_{\mathbf{Z}_{\left(e_G,q\right)}}(x)\phi_{\mathbf{X}}^\mathcal{NO}(b):x\in \mathbf{Z}_{\left(e_G,q\right)}, \ b\in \mathcal{NO}_\mathbf{X}\big\}\subseteq \mathcal{NO}_\mathbf{Z}.
\end{align*}
Then $\mathbf{Y}_q^\mathcal{NO}$ carries a right action of $\mathcal{NO}_\mathbf{X}$ such that
\[
y\cdot b=
y\phi_{\mathbf{X}}^\mathcal{NO}(b) \quad \text{for each $y\in \mathbf{Y}_q^\mathcal{NO}$ and $b\in \mathcal{NO}_\mathbf{X}$.}
\] 
For each $y, w\in \mathbf{Y}_q^\mathcal{NO}$, we have $y^*w\in \phi_{\mathbf{X}}^\mathcal{NO}\left(\mathcal{NO}_\mathbf{X}\right)$, and there is an $\mathcal{NO}_\mathbf{X}$ valued inner-product $\langle \cdot, \cdot\rangle_{\mathcal{NO}_\mathbf{X}}^q: \mathbf{Y}_q^\mathcal{NO}\times  \mathbf{Y}_q^\mathcal{NO}\rightarrow \mathcal{NO}_\mathbf{X}$ such that
\[
\langle y,w\rangle_{\mathcal{NO}_\mathbf{X}}^q =
\left(\phi_{\mathbf{X}}^\mathcal{NO}\right)^{-1}(y^*w) 
\quad \text{for each $y, w\in \mathbf{Y}_q^\mathcal{NO}$.}
\] 
With this structure, $\mathbf{Y}_q^\mathcal{NO}$ becomes a Hilbert $\mathcal{NO}_\mathbf{X}$-module. Furthermore, there exists a homomorphism $\Phi_q^{\mathcal{NO}}: \mathcal{NO}_\mathbf{X}\rightarrow \mathcal{L}_{\mathcal{NO}_\mathbf{X}}\left(\mathbf{Y}_q^\mathcal{NO}\right)$ such that 
\[
\Phi_q^{\mathcal{NO}}(b)(y)=\phi_{\mathbf{X}}^\mathcal{NO}(b)y \quad \text{for each $b\in \mathcal{NO}_\mathbf{X}$ and $y\in \mathbf{Y}_q^\mathcal{NO}$.}
\]
With this additional structure, $\mathbf{Y}_q^\mathcal{NO}$ becomes a Hilbert $\mathcal{NO}_\mathbf{X}$-bimodule. 
\begin{proof}
We have effectively already completed all of the necessary calculations in Subsection~\ref{Constructing the product system} to prove the result. Rather than just rerunning the arguments of Propositions~\ref{definition of Y_q} and \ref{right actions on each Y_q} with $\mathbf{Y}_q^\mathcal{NO}$ in place of $\mathbf{Y}_q^\mathcal{NT}$, we will show how these two spaces are related via the quotient maps on $\mathcal{NT}_\mathbf{X}$ and $\mathcal{NT}_\mathbf{Z}$ and use this to prove the result. The key observation is that the quotient homomorphisms $q_\mathbf{X}$ and $q_\mathbf{Z}$ intertwine the homomorphisms $\phi_\mathbf{X}^\mathcal{NT}$ and $\phi_\mathbf{X}^\mathcal{NO}$. 

Since $q_\mathbf{Z}\circ i_\mathbf{Z}=j_\mathbf{Z}$ and $q_\mathbf{Z}\circ \phi_\mathbf{X}^\mathcal{NT}=\phi_\mathbf{X}^\mathcal{NO}\circ q_\mathbf{X}$, we have
\begin{align*}
\mathbf{Y}_q^\mathcal{NO}
&=\cspan\big\{j_{\mathbf{Z}_{\left(e_G,q\right)}}(x)\phi_{\mathbf{X}}^\mathcal{NO}(b):x\in \mathbf{Z}_{\left(e_G,q\right)}, \ b\in \mathcal{NO}_\mathbf{X}\big\}\\
&=q_\mathbf{Z}\big(\cspan\big\{i_{\mathbf{Z}_{\left(e_G,q\right)}}(x)\phi_{\mathbf{X}}^\mathcal{NT}(b):x\in \mathbf{Z}_{\left(e_G,q\right)}, \ b\in \mathcal{NT}_\mathbf{X}\big\}\big)\\
&=q_\mathbf{Z}\left(\mathbf{Y}_q^\mathcal{NT}\right).
\end{align*}

We now show how the right actions of $\mathcal{NO}_\mathbf{X}$ on $\mathbf{Y}_q^\mathcal{NO}$ can be obtained from the right action of $\mathcal{NT}_\mathbf{X}$ on $\mathbf{Y}_q^\mathcal{NT}$. For any $a\in \mathcal{NT}_\mathbf{X}$ and $y\in \mathbf{Y}_q^\mathcal{NT}$, we have
\begin{equation}
\label{right action property}
\begin{aligned}
q_\mathbf{Z}\left(y\cdot a\right)
=q_\mathbf{Z}\left(y\phi_\mathbf{X}^\mathcal{NT}(a)\right)
=q_\mathbf{Z}(y)q_\mathbf{Z}\left(\phi_\mathbf{X}^\mathcal{NT}(a)\right)
=q_\mathbf{Z}(y)\phi_\mathbf{X}^\mathcal{NO}\left(q_\mathbf{X}(a)\right).
\end{aligned}
\end{equation}
If  $a'\in \mathcal{NT}_\mathbf{X}$ and $y'\in \mathbf{Y}_q^\mathcal{NT}$ with $q_\mathbf{X}(a)=q_\mathbf{X}(a')$ and $q_\mathbf{Z}(y)=q_\mathbf{Z}(y')$, then \eqref{right action property} is equal to 
\begin{align*}
q_\mathbf{Z}(y')\phi_\mathbf{X}^\mathcal{NO}\left(q_\mathbf{X}(a')\right)
=q_\mathbf{Z}(y')q_\mathbf{Z}\left(\phi_\mathbf{X}^\mathcal{NT}(a')\right)
=q_\mathbf{Z}\left(y'\phi_\mathbf{X}^\mathcal{NT}(a')\right)
=q_\mathbf{Z}\left(y'\cdot a'\right)
\end{align*}
Thus, $\mathbf{Y}_q^\mathcal{NO}$ carries a right action of $\mathcal{NO}_\mathbf{X}$ defined by the formula 
\[
q_\mathbf{Z}(y)\cdot q_\mathbf{X}(a):=q_\mathbf{Z}(y\cdot a) 
\quad \text{for each $a\in \mathcal{NT}_\mathbf{X}$ and $y\in \mathbf{Y}_q^\mathcal{NT}$.}
\]
Moreover, \eqref{right action property} shows that 
\begin{align*}
y\cdot a=y\phi_\mathbf{X}^\mathcal{NO}\left(a\right) \quad \text{for any $a\in \mathcal{NO}_\mathbf{X}$ and $y\in \mathbf{Y}_q^\mathcal{NO}$.}
\end{align*}

Next we show how the $\mathcal{NO}_\mathbf{X}$ valued inner-product on $\mathbf{Y}_q^\mathcal{NO}$ can be obtained from the $\mathcal{NT}_\mathbf{X}$ valued inner-product on $\mathbf{Y}_q^\mathcal{NT}$. For any $y,w\in \mathbf{Y}_q^\mathcal{NT}$, we have
\begin{equation}
\label{inner product property}
\begin{aligned}
q_\mathbf{X}\big(\langle y,w\rangle_{\mathcal{NT}_\mathbf{X}}^q\big)
=q_\mathbf{X}\big(\left(\phi_\mathbf{X}^\mathcal{NT}\right)^{-1}(y^*w)\big)
=\left(\phi_\mathbf{X}^\mathcal{NO}\right)^{-1}\left(q_\mathbf{Z}(y)^*q_\mathbf{Z}(w)\right).
\end{aligned}
\end{equation}
If $y',w'\in\mathbf{Y}_q^\mathcal{NT}$ with $q_\mathbf{Z}(y)=q_\mathbf{Z}(y')$ and $q_\mathbf{Z}(w)=q_\mathbf{Z}(w')$, then this is equal to 
\begin{align*}
\left(\phi_\mathbf{X}^\mathcal{NO}\right)^{-1}\left(q_\mathbf{Z}(y')^*q_\mathbf{Z}(w')\right)
&=\left(\phi_\mathbf{X}^\mathcal{NO}\right)^{-1}\left(q_\mathbf{Z}(y'^*w')\right)\\
&=q_\mathbf{X}\big(\left(\phi_\mathbf{X}^\mathcal{NT}\right)^{-1}(y'^*w')\big)
=q_\mathbf{X}\big(\langle y',w'\rangle_{\mathcal{NT}_\mathbf{X}}^q\big).
\end{align*} 
Thus, we can define $\langle \cdot, \cdot\rangle_{\mathcal{NO}_\mathbf{X}}^q: \mathbf{Y}_q^\mathcal{NO}\times  \mathbf{Y}_q^\mathcal{NO}\rightarrow \mathcal{NO}_\mathbf{X}$ by
\[
\langle q_\mathbf{Z}(y), q_\mathbf{Z}(w)\rangle_{\mathcal{NO}_\mathbf{X}}^q:=q_\mathbf{X}\big(\langle y,w\rangle_{\mathcal{NT}_\mathbf{X}}^q\big)
\quad \text{for any $y,w\in\mathbf{Y}_q^\mathcal{NT}$.} 
\]
Moreover, \eqref{inner product property} shows that 
\begin{align}
\label{what does the IP look like?}
\langle y,w\rangle_{\mathcal{NO}_\mathbf{X}}^q=\left(\phi_\mathbf{X}^\mathcal{NO}\right)^{-1}\left(y^*w\right) \quad \text{for any $y,w\in\mathbf{Y}_q^\mathcal{NO}$.}
\end{align}
Using~\eqref{what does the IP look like?}, it is elementary to show that $\big(\mathbf{Y}_q^\mathcal{NO}, \langle \cdot, \cdot\rangle_{\mathcal{NO}_\mathbf{X}}^q \big)$ is a Hilbert $\mathcal{NO}_\mathbf{X}$-module. 

It remains to check that $\mathbf{Y}_q^\mathcal{NO}$ carries a left action of $\mathcal{NO}_\mathbf{X}$ by adjointable operators. We show how the homomorphism $\Phi_q^{\mathcal{NO}}: \mathcal{NO}_\mathbf{X}\rightarrow \mathcal{L}_{\mathcal{NO}_\mathbf{X}}\left(\mathbf{Y}_q^\mathcal{NO}\right)$ can be obtained from the homomorphism
$\Phi_q^{\mathcal{NT}}: \mathcal{NT}_\mathbf{X}\rightarrow \mathcal{L}_{\mathcal{NT}_\mathbf{X}}\left(\mathbf{Y}_q^\mathcal{NT}\right)$. For any $a\in \mathcal{NT}_\mathbf{X}$ and $y\in \mathbf{Y}_q^\mathcal{NT}$, we have
\begin{equation}
\label{left action property}
\begin{aligned}
q_\mathbf{Z}\left(\Phi_q^\mathcal{NT}(a)(y)\right)
=q_\mathbf{Z}\left(\phi_\mathbf{X}^\mathcal{NT}(a)y\right)
&=q_\mathbf{Z}\left(\phi_\mathbf{X}^\mathcal{NT}(a)\right)q_\mathbf{Z}(y)\\
&=\phi_\mathbf{X}^\mathcal{NO}\left(q_\mathbf{X}(a)\right)q_\mathbf{Z}(y).
\end{aligned}
\end{equation}
If  $a'\in \mathcal{NT}_\mathbf{X}$ and $y'\in \mathbf{Y}_q^\mathcal{NT}$ with $q_\mathbf{X}(a)=q_\mathbf{X}(a')$ and $q_\mathbf{Z}(y)=q_\mathbf{Z}(y')$, then (\ref{left action property}) is equal to 
\begin{align*}
\phi_\mathbf{X}^\mathcal{NO}\left(q_\mathbf{X}(a')\right)q_\mathbf{Z}(y')
=q_\mathbf{Z}\left(\phi_\mathbf{X}^\mathcal{NT}(a')\right)q_\mathbf{Z}(y')
&=q_\mathbf{Z}\left(\phi_\mathbf{X}^\mathcal{NT}(a')y'\right)\\
&=q_\mathbf{Z}\left(\Phi_q^\mathcal{NT}(a')(y')\right).
\end{align*}
Thus, for $a\in \mathcal{NT}_\mathbf{X}$, there exists a map $\Phi_q^{\mathcal{NO}}\left(q_\mathbf{X}(a)\right): \mathbf{Y}_q^\mathcal{NO}\rightarrow \mathbf{Y}_q^\mathcal{NO}$  given by
\[
\Phi_q^{\mathcal{NO}}\left(q_\mathbf{X}(a)\right)\left(q_\mathbf{Z}(y)\right):=q_\mathbf{Z}\left(\Phi_q^\mathcal{NT}(a)(y)\right)
\quad \text{for each $y\in \mathbf{Y}_q^\mathcal{NT}$.}
\]
Moreover, \eqref{left action property} shows that 
\begin{align}
\label{what does the left action look like?}
\Phi_q^{\mathcal{NO}}\left(a\right)\left(y\right)=\phi_\mathbf{X}^\mathcal{NO}\left(a\right)y 
\quad \text{for any $a\in \mathcal{NO}_\mathbf{X}$ and $y\in \mathbf{Y}_q^\mathcal{NO}$}
\end{align}
Using~\eqref{what does the IP look like?} and~\eqref{what does the left action look like?}, it is routine to check that $\Phi_q^{\mathcal{NO}}\left(q_\mathbf{X}(a)\right)$ is an adjointable map on $\mathbf{Y}_q^\mathcal{NO}$ and the map $q_\mathbf{X}(a)\mapsto \Phi_q^{\mathcal{NO}}\left(q_\mathbf{X}(a)\right)$ is a homomorphism from $\mathcal{NO}_\mathbf{X}$ to $\mathcal{L}_{\mathcal{NO}_\mathbf{X}}\left(\mathbf{Y}_q^\mathcal{NO}\right)$, which we denote by $\Phi_q^{\mathcal{NO}}$. 
\end{proof}
\end{prop}

The next result shows that if $\mathbf{Y}_{e_H}^\mathcal{NO}:={}_{\mathcal{NO}_\mathbf{X}}\left(\mathcal{NO}_\mathbf{X}\right)_{\mathcal{NO}_\mathbf{X}}$, then $\mathbf{Y}^\mathcal{NO}:=\bigsqcup_{q\in Q} \mathbf{Y}_q^\mathcal{NO}$ has the structure of a compactly aligned product system. 

We use the Hilbert $\mathcal{NT}_\mathbf{X}$-bimodule isomorphism $M_{q,t}^{\mathbf{Y}^\mathcal{NT}}:\mathbf{Y}_q^\mathcal{NT}\otimes_{\mathcal{NT}_\mathbf{X}}\mathbf{Y}_t^\mathcal{NT}\rightarrow \mathbf{Y}_{qt}^\mathcal{NT}$ to define a Hilbert $\mathcal{NO}_\mathbf{X}$-bimodule isomorphism $M_{q,t}^{\mathbf{Y}^\mathcal{NO}}:\mathbf{Y}_q^\mathcal{NO}\otimes_{\mathcal{NO}_\mathbf{X}}\mathbf{Y}_t^\mathcal{NO}\rightarrow \mathbf{Y}_{qt}^\mathcal{NO}$, such that the diagram in Figure~\ref{multiplication isomorphisms} is commutative.

\begin{figure}[h]
\centering
\begin{tikzpicture}[scale=1.5]
\node (E) at (0,0) {$\mathbf{Y}_q^\mathcal{NO}\otimes_{\mathcal{NO}_\mathbf{X}} \mathbf{Y}_t^\mathcal{NO}$};
\node (F) at (0,2) {$\mathbf{Y}_q^\mathcal{NT}\otimes_{\mathcal{NT}_\mathbf{X}} \mathbf{Y}_t^\mathcal{NT}$};
\node (G) at (3.5,0) {$\mathbf{Y}_{qt}^\mathcal{NO}$};
\node (H) at (3.5,2) {$\mathbf{Y}_{qt}^\mathcal{NT}$};
\path[->,font=\scriptsize]
(F) edge node[left]{$q_\mathbf{Z}\otimes_{\mathcal{NT}_\mathbf{X}} q_\mathbf{Z}$} (E)
(H) edge node[right]{$q_\mathbf{Z}$} (G)
(F) edge node[above]{$M_{q,t}^{\mathbf{Y}^\mathcal{NT}}$} (H)
(E) edge node[below]{$M_{q,t}^{\mathbf{Y}^\mathcal{NO}}$} (G);
\end{tikzpicture}
\caption{The Hilbert $\mathcal{NO}_\mathbf{X}$-bimodule isomorphism $M_{q,t}^{\mathbf{Y}^\mathcal{NO}}$}
\label{multiplication isomorphisms}
\end{figure}

\begin{prop}
\label{existence of NO product system in NO_Z}
Suppose $G$ is amenable, $A$ acts faithfully on each fibre of $\mathbf{X}$, and each $\widetilde{\phi}_{(p,q)}$ is injective, so that the collection of Hilbert $\mathcal{NO}_\mathbf{X}$-bimodules $\big\{\mathbf{Y}_q^\mathcal{NO}:q\in Q\setminus \{e_H\}\big\}$ from Proposition~\ref{existence of fibres in NO product system in NO_Z} exists.  

Let $\mathbf{Y}_{e_H}^\mathcal{NO}:={}_{\mathcal{NO}_\mathbf{X}}\left(\mathcal{NO}_\mathbf{X}\right)_{\mathcal{NO}_\mathbf{X}}$, and for each $q\in Q\setminus \{e_H\}$, define $\mathbf{Y}_q^\mathcal{NO}$ as in Proposition~\ref{existence of fibres in NO product system in NO_Z}. Then $\mathbf{Y}^\mathcal{NO}:=\bigsqcup_{q\in Q} \mathbf{Y}_q^\mathcal{NO}$ is a compactly aligned product system over $(H,Q)$ with coefficient algebra $\mathcal{NO}_\mathbf{X}$, with multiplication in $\mathbf{Y}^\mathcal{NO}$ given by multiplication in $\mathcal{NO}_\mathbf{Z}$.
\begin{proof}
We have already shown that each $\mathbf{Y}_q^\mathcal{NO}$ is a Hilbert $\mathcal{NO}_\mathbf{X}$-bimodule. To show that $\mathbf{Y}^\mathcal{NO}:=\bigsqcup_{q\in Q} \mathbf{Y}_q^\mathcal{NO}$ is a product system over $(H,Q)$ with coefficient algebra $\mathcal{NO}_\mathbf{X}$ and multiplication inherited from $\mathcal{NO}_{\mathbf{Z}}$, we need only check that there exists a Hilbert $\mathcal{NO}_\mathbf{X}$-bimodule isomorphism $M_{q,t}^{\mathbf{Y}^\mathcal{NO}}:\mathbf{Y}_q^\mathcal{NO}\otimes_{\mathcal{NO}_\mathbf{X}}\mathbf{Y}_t^\mathcal{NO}\rightarrow \mathbf{Y}_{qt}^\mathcal{NO}$ for each $q, t\in Q\setminus \{e_H\}$ such that 
\begin{align}
\label{bimodule isomorphism on CNP product system}
M_{q,t}^{\mathbf{Y}^\mathcal{NO}}(y\otimes_{\mathcal{NO}_\mathbf{X}} w)=yw
\quad \text{for each $y\in \mathbf{Y}_q^\mathcal{NO}$, $w\in \mathbf{Y}_t^\mathcal{NO}$.}
\end{align}
Rather than just rerun our argument from Proposition~\ref{Y is a product system} with $\mathbf{Y}_q^\mathcal{NO}$ in place of $\mathbf{Y}_q^\mathcal{NT}$, we will show how these isomorphisms can be obtained from the Hilbert $\mathcal{NT}_\mathbf{X}$-bimodule isomorphisms $M_{q,t}^{\mathbf{Y}^\mathcal{NT}}:\mathbf{Y}_q^\mathcal{NT}\otimes_{\mathcal{NT}_\mathbf{X}}\mathbf{Y}_t^\mathcal{NT}\rightarrow \mathbf{Y}_{qt}^\mathcal{NT}$ using the quotient map $q_\mathbf{Z}$. Observe that if $y,y'\in \mathbf{Y}_q^\mathcal{NT}$ with $q_\mathbf{Z}(y)=q_\mathbf{Z}(y')$ and $w,w'\in \mathbf{Y}_t^\mathcal{NT}$ with $q_\mathbf{Z}(w)=q_\mathbf{Z}(w')$, then
\begin{align*}
\label{bimodule isomorphism on CNP product system calculation}
q_\mathbf{Z}\big(M_{q,t}^{\mathbf{Y}^\mathcal{NT}}(y\otimes_{\mathcal{NT}_\mathbf{X}} w)\big)
&=q_\mathbf{Z}\left(yw\right)\\
&=q_\mathbf{Z}\left(y'w'\right)
=q_\mathbf{Z}\big(M_{q,t}^{\mathbf{Y}^\mathcal{NT}}(y'\otimes_{\mathcal{NT}_\mathbf{X}} w')\big). 
\end{align*}
Hence, there is a well-defined map $M_{q,t}^{\mathbf{Y}^\mathcal{NO}}:\mathbf{Y}_q^\mathcal{NO}\otimes_{\mathcal{NO}_\mathbf{X}}\mathbf{Y}_t^\mathcal{NO}\rightarrow \mathbf{Y}_{qt}^\mathcal{NO}$ given by
\begin{align*}
M_{q,t}^{\mathbf{Y}^\mathcal{NO}}\big((q_\mathbf{Z}\otimes_{\mathcal{NO}_\mathbf{X}} q_\mathbf{Z})(z)\big):&=q_\mathbf{Z}\big(M_{q,t}^{\mathbf{Y}^\mathcal{NT}}(z)\big) \\
&\quad \text{for each $z\in \mathbf{Y}_q^\mathcal{NT}\otimes_{\mathcal{NT}_\mathbf{X}}\mathbf{Y}_t^\mathcal{NT}$,}
\end{align*}
which satisfies \eqref{bimodule isomorphism on CNP product system}. Since $M_{q,t}^{\mathbf{Y}^\mathcal{NT}}$ is surjective and $q_\mathbf{Z}\left(\mathbf{Y}_{qt}^\mathcal{NT}\right)=\mathbf{Y}_{qt}^\mathcal{NO}$, the map $M_{q,t}^{\mathbf{Y}^\mathcal{NO}}$ is also surjective. Routine calculations using \eqref{bimodule isomorphism on CNP product system} show that $M_{q,t}^{\mathbf{Y}^\mathcal{NO}}$ is inner-product preserving and left $\mathcal{NO}_\mathbf{X}$-linear. Hence, $M_{q,t}^{\mathbf{Y}^\mathcal{NO}}$ is a Hilbert $\mathcal{NO}_\mathbf{X}$-bimodule isomorphism. We conclude that $\mathbf{Y}^\mathcal{NO}$ has the structure of a product system. 

It remains to show that $\mathbf{Y}^\mathcal{NO}$ is compactly aligned. Fix $q,t\in Q$ with $q\vee t<\infty$, and let $S\in \mathcal{K}_{\mathcal{NO}_\mathbf{X}}\left(\mathbf{Y}_q^\mathcal{NO}\right)$ and $T\in \mathcal{K}_{\mathcal{NO}_\mathbf{X}}\left(\mathbf{Y}_t^\mathcal{NO}\right)$. We show that 
$
\iota_q^{q\vee t}(S)\iota_t^{q\vee t}(T)\in \mathcal{K}_{\mathcal{NO}_\mathbf{X}}\left(\mathbf{Y}_{q\vee t}^\mathcal{NO}\right).
$
If $q=e_H$ or $t=e_H$, the result is trivial, so we may as well suppose that $q,t\neq e_H$. Since $\iota_q^{q\vee t}$ and $\iota_t^{q\vee t}$ are linear and continuous, as is multiplication in $\mathcal{L}_{\mathcal{NO}_\mathbf{X}}\left(\mathbf{Y}_{q\vee t}^\mathcal{NO}\right)$, we may as well assume that $S$ and $T$ are rank one operators. Hence, $S=\Theta_{q_\mathbf{Z}(x), q_\mathbf{Z}(y)}$ and $T=\Theta_{q_\mathbf{Z}(u), q_\mathbf{Z}(v)}$ for some $x,y\in \mathbf{Y}_q^\mathcal{NT}$ and $u,v\in \mathbf{Y}_t^\mathcal{NT}$. Then for any $z\in \mathbf{Y}_{q\vee t}^\mathcal{NT}$, we have
\begin{align*}
\Big(\iota_q^{q\vee t}\left(\Theta_{q_\mathbf{Z}(x), q_\mathbf{Z}(y)}\right)&\iota_t^{q\vee t}\left(\Theta_{q_\mathbf{Z}(u), q_\mathbf{Z}(v)}\right)\Big)(q_\mathbf{Z}(z))\\
&=\left(M_{q_\mathbf{Z}(x)q_\mathbf{Z}(y)^*}M_{q_\mathbf{Z}(u)q_\mathbf{Z}(v)^*}\right)(q_\mathbf{Z}(z))\\
&=q_\mathbf{Z}\left(M_{xy^*uv^*}(z)\right).
\end{align*}
Lemma~\ref{Q indexed subalgebras of NO_Z} and Proposition~\ref{Y is compactly aligned} show that $M_{xy^*uv^*}\in \mathcal{K}_{\mathcal{NT}_\mathbf{X}}\left(\mathbf{Y}_{q\vee t}^\mathcal{NT}\right)$. Thus, the previous line can be approximated by sums of the form 
\[
q_\mathbf{Z}\left(\Theta_{\alpha,\beta}(z)\right)=\Theta_{q_\mathbf{Z}(\alpha),q_\mathbf{Z}(\beta)}\left(q_\mathbf{Z}(z)\right),
\]
where $\alpha,\beta \in \mathbf{Y}_{q\vee t}^\mathcal{NT}$. Since $z\in \mathbf{Y}_{q\vee t}^\mathcal{NT}$ was arbitrary and $\mathbf{Y}_{r}^\mathcal{NO}=q_\mathbf{Z}\left(\mathbf{Y}_{r}^\mathcal{NT}\right)$ for each $r\in Q\setminus \{e_H\}$, we conclude that $\iota_q^{q\vee t}\left(\Theta_{q_\mathbf{Z}(x), q_\mathbf{Z}(y)}\right)\iota_t^{q\vee t}\left(\Theta_{q_\mathbf{Z}(u), q_\mathbf{Z}(v)}\right)$ can be approximated by sums of operators in $\mathcal{K}_{\mathcal{NO}_\mathbf{X}}\left(\mathbf{Y}_{q\vee t}^\mathcal{NO}\right)$.  Thus,
\[
\iota_q^{q\vee t}\left(\Theta_{q_\mathbf{Z}(x), q_\mathbf{Z}(y)}\right)\iota_t^{q\vee t}\left(\Theta_{q_\mathbf{Z}(u), q_\mathbf{Z}(v)}\right)
\in \mathcal{K}_{\mathcal{NO}_\mathbf{X}}\left(\mathbf{Y}_{q\vee t}^\mathcal{NO}\right)
\]
as required. 
\end{proof}
\end{prop}

Shortly, we will examine the Cuntz--Nica--Pimsner algebra of the product system $\mathbf{Y}^\mathcal{NO}$. To make the calculations tractable, we seek sufficient conditions for the homomorphisms $\Phi_q^{\mathcal{NO}}:\mathcal{NO}_\mathbf{X}\rightarrow \mathcal{L}_{\mathcal{NO}_\mathbf{X}}\left(\mathbf{Y}_q^\mathcal{NO}\right)$ that implement the left action of $\mathcal{NO}_\mathbf{X}$ on the fibres of $\mathbf{Y}^\mathcal{NO}$ to be injective. Again, we will make use of \cite[Corollary~4.12]{carlsen--larsen--sims--vittadello}. The main step is showing that $\Phi_q^{\mathcal{NO}}\left(\mathcal{NO}_\mathbf{X}\right)\subseteq \mathcal{L}_{\mathcal{NO}_\mathbf{X}}\left(\mathbf{Y}_q^\mathcal{NO}\right)$ carries a coaction of $G$ that intertwines $\Phi_q^\mathcal{NO}$ with the canonical coaction of $G$ on $\mathcal{NO}_\mathbf{X}$. We first need some preliminary results.

\begin{lem}
\label{interaction between coactions}
Suppose that $G$ is an amenable group, $A$ acts faithfully on each fibre of $\mathbf{X}$, and each $\widetilde{\phi}_{(p,q)}$ is injective, so that the product system $\mathbf{Y}^{\mathcal{NO}}$ from Proposition~\ref{existence of NO product system in NO_Z} exists. Then
\[
\nu_\mathbf{Z}\circ \phi_{\mathbf{X}}^\mathcal{NO}=\left(\phi_{\mathbf{X}}^\mathcal{NO}\otimes \iota\right)\circ \nu_\mathbf{X}.
\]
Furthermore, if $a\in \mathcal{NO}_\mathbf{X}$ and $y\in \mathbf{Y}_q^\mathcal{NO}$, then
\[
\nu_\mathbf{Z}\left(\Phi_q^{\mathcal{NO}}(a)(y)\right)=\big(\big(\Phi_q^{\mathcal{NO}}\otimes \iota\big)(\nu_\mathbf{X}(a))\big)\left(\nu_\mathbf{Z}(y)\right).
\]
\begin{proof}
Since both $\nu_\mathbf{Z}\circ \phi_{\mathbf{X}}^\mathcal{NO}$ and $\left(\phi_{\mathbf{X}}^\mathcal{NO}\otimes \iota\right)\circ\nu_\mathbf{X}$ are homomorphisms, it suffices to check that they agree on the generators of $\mathcal{NO}_\mathbf{X}$. If $x\in \mathbf{X}_p$, then
\begin{align*}
\nu_\mathbf{Z}\left(\phi_{\mathbf{X}}^\mathcal{NO}\left(j_{\mathbf{X}_p}(x)\right)\right)
=\nu_\mathbf{Z}\big(j_{\mathbf{Z}_{\left(p,e_H\right)}}(x)\big)
&=j_{\mathbf{Z}_{\left(p,e_H\right)}}(x)\otimes i_{G\rtimes H}\left(p,e_H\right)\\
&=\left(\phi_{\mathbf{X}}^\mathcal{NO}\otimes \iota\right)\left(j_{\mathbf{X}_p}(x)\otimes i_G(p)\right)\\
&=\left(\phi_{\mathbf{X}}^\mathcal{NO}\otimes \iota\right)\left(\nu_\mathbf{X}\left(j_{\mathbf{X}_p}(x)\right)\right),
\end{align*}
as required. Since $\nu_\mathbf{Z}$ is a homomorphism, we see that for any $a\in \mathcal{NO}_\mathbf{X}$ and $y\in \mathbf{Y}_q^\mathcal{NO}$,
\begin{align*}
\nu_\mathbf{Z}\left(\Phi_q^{\mathcal{NO}}(a)(y)\right)
=\nu_\mathbf{Z}\left(\phi_{\mathbf{X}}^\mathcal{NO}(a)y\right)
&=\nu_\mathbf{Z}\left(\phi_{\mathbf{X}}^\mathcal{NO}(a)\right)\nu_\mathbf{Z}(y)\\
&=\left(\phi_{\mathbf{X}}^\mathcal{NO}\otimes \iota\right)(\nu_\mathbf{X}(a))\nu_\mathbf{Z}(y)\\
&=\big(\big(\Phi_q^{\mathcal{NO}}\otimes \iota\big)(\nu_\mathbf{X}(a))\big)\left(\nu_\mathbf{Z}(y)\right).
\qedhere
\end{align*}
\end{proof}
\end{lem}

\begin{lem}
\label{range of coaction}
Suppose that $G$ is an amenable group, $A$ acts faithfully on each fibre of $\mathbf{X}$, and each $\widetilde{\phi}_{(p,q)}$ is injective, so that the product system $\mathbf{Y}^{\mathcal{NO}}$ from Proposition~\ref{existence of NO product system in NO_Z} exists. Then for each $q\in Q$,
\begin{align*}
&\big(\mathrm{id}_{\mathbf{Y}_q^\mathcal{NO}}\otimes \iota\big)\big(\mathbf{Y}_q^\mathcal{NO}\otimes C^*(G)\big)\\
&\qquad=\cspan\big\{\nu_\mathbf{Z}\big(\mathbf{Y}_q^\mathcal{NO}\big)\big(1_{\mathcal{M}(\mathcal{NO}_\mathbf{Z})}\otimes i_{G\rtimes_\alpha H}(G\times \{q^{-1}\})\big)\big\}. 
\end{align*}
\begin{proof}
We begin by showing $\subseteq$. Fix $g\in G$ and $z\in \mathbf{Z}_{(p,q)}$, $w\in \mathbf{Z}_{\left(r,e_H\right)}$. Since 
\[
\left(p,q\right)\left(r,e_H\right)^{-1}\left(r\alpha_{q^{-1}}\left(p^{-1}g\right),q^{-1}\right)=\left(g,e_H\right),
\]
we see that 
\begin{align*}
&\big(\mathrm{id}_{\mathbf{Y}_q^\mathcal{NO}}\otimes \iota\big)\big(j_{\mathbf{Z}_{(p,q)}}(z)j_{\mathbf{Z}_{\left(r,e_H\right)}}(w)^*\otimes i_G(g)\big)\\
&\qquad=\nu_\mathbf{Z}\big(j_{\mathbf{Z}_{(p,q)}}(z)j_{\mathbf{Z}_{\left(r,e_H\right)}}(w)^*\big)\big(1_{\mathcal{M}(\mathcal{NO}_\mathbf{Z})}\otimes i_{G\rtimes_\alpha H}(r\alpha_{q^{-1}}(p^{-1}g),q^{-1})\big)\\
&\qquad\in
\cspan\big\{\nu_\mathbf{Z}\big(\mathbf{Y}_q^\mathcal{NO}\big)\big(1_{\mathcal{M}(\mathcal{NO}_\mathbf{Z})}\otimes i_{G\rtimes_\alpha H}(G\times \{q^{-1}\})\big)\big\}. 
\end{align*}
Since $\mathbf{Y}_q^\mathcal{NO}=\cspan\big\{j_{\mathbf{Z}_{(p,q)}}(z)j_{\mathbf{Z}_{\left(r,e_H\right)}}(w)^*:p,r\in P, \, z\in \mathbf{Z}_{(p,q)}, \, w\in \mathbf{Z}_{\left(r,e_H\right)}\big\}$ and $\mathrm{span}\left\{i_G(g):g\in G\right\}$ is dense in $C^*(G)$, we conclude that $\subseteq$ holds. 

We now prove $\supseteq$. Fix $p,r\in P$, $g\in G$, $z\in \mathbf{Z}_{(p,q)}$, and $w\in \mathbf{Z}_{\left(r,e_H\right)}$. Since
$
(p,q)(r,e_H)^{-1}(g,q^{-1})=(p\alpha_q(r^{-1}g),e_H),
$
we have
\begin{align*}
&\nu_\mathbf{Z}\big(j_{\mathbf{Z}_{(p,q)}}(z)j_{\mathbf{Z}_{\left(r,e_H\right)}}(w)^*\big)\big(1_{\mathcal{M}(\mathcal{NO}_\mathbf{Z})}\otimes i_{G\rtimes_\alpha H}(g,q^{-1})\big)\\
&\qquad=
\big(\mathrm{id}_{\mathbf{Y}_q^\mathcal{NO}}\otimes \iota\big)
\big(j_{\mathbf{Z}_{(p,q)}}(z)j_{\mathbf{Z}_{\left(r,e_H\right)}}(w)^*\otimes i_G(p\alpha_q(r^{-1}g))\big)\\
&\qquad
\in\big(\mathrm{id}_{\mathbf{Y}_q^\mathcal{NO}}\otimes \iota\big)\big(\mathbf{Y}_q^\mathcal{NO}\otimes C^*(G)\big).
\end{align*}
Since $\mathbf{Y}_q^\mathcal{NO}=\cspan\big\{j_{\mathbf{Z}_{(p,q)}}(z)j_{\mathbf{Z}_{\left(r,e_H\right)}}(w)^*:p,r\in P, \, z\in \mathbf{Z}_{(p,q)}, \, w\in \mathbf{Z}_{\left(r,e_H\right)}\big\}$ and $g\in G$ was arbitrary, we conclude that $\supseteq$ holds.
\end{proof}
\end{lem}

\begin{lem}
\label{coaction is well defined}
Suppose that $G$ is an amenable group, $A$ acts faithfully on each fibre of $\mathbf{X}$, and each $\widetilde{\phi}_{(p,q)}$ is injective, so that the product system $\mathbf{Y}^{\mathcal{NO}}$ from Proposition~\ref{existence of NO product system in NO_Z} exists.  If $\Phi_q^{\mathcal{NO}}(a)=0\in \mathcal{L}_{\mathcal{NO}_\mathbf{X}}\left(\mathbf{Y}_q^\mathcal{NO}\right)$ for some $a\in \mathcal{NO}_\mathbf{X}$, then 
\[
\left(\Phi_q^{\mathcal{NO}}\otimes \mathrm{id}_{C^*(G)}\right)(\nu_\mathbf{X}(a))=0\in \mathcal{L}_{\mathcal{NO}_\mathbf{X}}\left(\mathbf{Y}_q^\mathcal{NO}\right)\otimes C^*(G).
\]
\begin{proof}
Suppose that $a\in \mathcal{NO}_\mathbf{X}$ is such that $\Phi_q^{\mathcal{NO}}(a)=0 \in\mathcal{L}_{\mathcal{NO}_\mathbf{X}}\left(\mathbf{Y}_q^\mathcal{NO}\right)$. We want to show that $\left(\Phi_q^{\mathcal{NO}}\otimes \mathrm{id}_{C^*(G)}\right)(\nu_\mathbf{X}(a))=0\in \mathcal{L}_{\mathcal{NO}_\mathbf{X}}\left(\mathbf{Y}_q^\mathcal{NO}\right)\otimes C^*(G)$. Since the external tensor product $\mathcal{L}_{\mathcal{NO}_\mathbf{X}}\left(\mathbf{Y}_q^\mathcal{NO}\right)\otimes C^*(G)$ is isomorphic to $\mathcal{L}_{\mathcal{NO}_\mathbf{X}}\left(\mathbf{Y}_q^\mathcal{NO}\right)\otimes \mathcal{L}_{C^*(G)}(C^*(G))$, which embeds isometrically in $\mathcal{L}_{\mathcal{NO}_\mathbf{X}\otimes C^*(G)}\left(\mathbf{Y}_q^\mathcal{NO}\otimes C^*(G)\right)$ (see \cite[Chapter~4]{lance} for the details regarding exterior tensor products of Hilbert modules), it suffices to show that $\left(\Phi_q^{\mathcal{NO}}\otimes \mathrm{id}_{C^*(G)}\right)(\nu_\mathbf{X}(a))$ acts as the zero operator on $\mathbf{Y}_q^\mathcal{NO}\otimes C^*(G)$. By Lemma~\ref{range of coaction}, it suffices to show that $\left(\Phi_q^{\mathcal{NO}}\otimes \mathrm{id}_{C^*(G)}\right)(\nu_\mathbf{X}(a))$ is the zero operator on
\[
\cspan\Big\{\big(\mathrm{id}_{\mathbf{Y}_q^\mathcal{NO}}\otimes \iota^{-1}\big)\Big(\nu_\mathbf{Z}\big(\mathbf{Y}_q^\mathcal{NO}\big)\big(1_{\mathcal{M}(\mathcal{NO}_\mathbf{Z})}\otimes i_{G\rtimes_\alpha H}(G\times \{q^{-1}\})\big)\Big)\Big\}.
\]
To do this we will use the fact that
\begin{align*}
&\left(\Phi_q^{\mathcal{NO}}\otimes \mathrm{id}_{C^*(G)}\right)(\nu_\mathbf{X}(a))\circ\big(\mathrm{id}_{\mathbf{Y}_q^\mathcal{NO}}\otimes \iota^{-1}\big)\\
&\qquad=\big(\mathrm{id}_{\mathbf{Y}_q^\mathcal{NO}}\otimes \iota^{-1}\big)\circ \left(\Phi_q^{\mathcal{NO}}\otimes \iota\right)(\nu_\mathbf{X}(a)),
\end{align*} 
and $\left(\Phi_q^{\mathcal{NO}}\otimes \iota\right)(\nu_\mathbf{X}(a))$ is left multiplication by $\left(\phi_{\mathbf{X}}^\mathcal{NO}\otimes \iota\right)(\nu_\mathbf{X}(a))$ on 
\begin{align*}
\mathbf{Y}_q^\mathcal{NO}\otimes C^*(G\rtimes_\alpha H)
&\subseteq \mathcal{NO}_\mathbf{Z}\otimes C^*(G\rtimes_\alpha H)\\
& \subseteq \mathcal{M}\left(\mathcal{NO}_\mathbf{Z}\right)\otimes C^*(G\rtimes_\alpha H).
\end{align*}
Since multiplication in $\mathcal{M}\left(\mathcal{NO}_\mathbf{Z}\right)\otimes C^*(G\rtimes_\alpha H)$ is associative, we see that for any $y\in \mathbf{Y}_q^\mathcal{NO}$ and $g\in G$,
\begin{equation}
\label{equation for coaction to be well defined}
\begin{aligned}
&\big(\Phi_q^{\mathcal{NO}}\!\otimes\! \mathrm{id}_{C^*(G)}\big)(\nu_\mathbf{X}(a))\Big(\big(\mathrm{id}_{\mathbf{Y}_q^\mathcal{NO}}\!\otimes\! \iota^{-1}\big)\big(\nu_\mathbf{Z}(y)\big(1_{\mathcal{M}(\mathcal{NO}_\mathbf{Z})}\!\otimes\! i_{G\rtimes_\alpha H}(g,q^{-1})\big)\big)\Big)\\
&\quad=\!
\big(\mathrm{id}_{\mathbf{Y}_q^\mathcal{NO}}\!\otimes\! \iota^{-1}\big)\Big(\!\big(\big(\Phi_q^{\mathcal{NO}}\!\otimes\! \iota\big)(\nu_\mathbf{X}(a))(\nu_\mathbf{Z}(y))\big)\big(1_{\mathcal{M}(\mathcal{NO}_\mathbf{Z})}\!\otimes\! i_{G\rtimes_\alpha H}(g,q^{-1})\big)\!\Big).
\end{aligned}
\end{equation}
By Lemma~\ref{interaction between coactions}, we have 
\[
\left(\Phi_q^{\mathcal{NO}}\otimes \iota\right)(\nu_\mathbf{X}(a))\left(\nu_\mathbf{Z}(y)\right)=\nu_\mathbf{Z}\left(\Phi_q^{\mathcal{NO}}(a)(y)\right)=0.
\]
Thus,~\eqref{equation for coaction to be well defined} is zero, and we conclude that $\left(\Phi_q^{\mathcal{NO}}\otimes \mathrm{id}_{C^*(G)}\right)(\nu_\mathbf{X}(a))$ is the zero operator on $\mathbf{Y}_q^\mathcal{NO}\otimes C^*(G)$ as required. 
\end{proof}
\end{lem}

We are now ready to show that the homomorphism $\Phi_q^{\mathcal{NO}}:\mathcal{NO}_\mathbf{X}\rightarrow \mathcal{L}_{\mathcal{NO}_\mathbf{X}}\left(\mathbf{Y}_q^{\mathcal{NO}}\right)$ is injective.

\begin{prop}
\label{injectivity of induced actions}
Suppose that $G$ is an amenable group, $A$ acts faithfully on each fibre of $\mathbf{X}$, and each $\widetilde{\phi}_{(p,q)}$ is injective, so that the product system $\mathbf{Y}^{\mathcal{NO}}$ from Proposition~\ref{existence of NO product system in NO_Z} exists. If $A$ acts faithfully on $\mathbf{Z}_{\left(e_G,q\right)}$, then $\Phi_q^{\mathcal{NO}}:\mathcal{NO}_\mathbf{X}\rightarrow \mathcal{L}_{\mathcal{NO}_\mathbf{X}}\left(\mathbf{Y}_q^\mathcal{NO}\right)$ is injective. 
\begin{proof}
Since $G$ is amenable we can use \cite[Corollary~4.12]{carlsen--larsen--sims--vittadello} to prove the result. We begin by checking that ${\Phi_q^{\mathcal{NO}}}|_{j_{\mathbf{X}_{e_G}}}$ is injective. Suppose that $a\in A$ is such that $\Phi_q^{\mathcal{NO}}\big(j_{\mathbf{X}_{e_G}}(a)\big)=0$. For any $z\in \mathbf{Z}_{\left(e_G,q\right)}$, we then have that
\[
0=\Phi_q^{\mathcal{NO}}\big(j_{\mathbf{X}_{e_G}}(a)\big)\big(j_{\mathbf{Z}_{\left(e_G,q\right)}}(z)\big)=j_{\mathbf{Z}_{\left(e_G,e_H\right)}}(a)j_{\mathbf{Z}_{\left(e_G,q\right)}}(z)=j_{\mathbf{Z}_{\left(e_G,q\right)}}(a\cdot z).
\]
Since $j_{\mathbf{Z}_{\left(e_G,q\right)}}$ is isometric and $A$ acts faithfully on $\mathbf{Z}_{\left(e_G,q\right)}$, we deduce that $a=0$, and so $j_{\mathbf{X}_{e_G}}(a)=0$. 

We now need to show that there is a coaction $\beta$ of $G$ on $\Phi_q^{\mathcal{NO}}\left(\mathcal{NO}_\mathbf{X}\right)$ that intertwines $\Phi_q^{\mathcal{NO}}$ with the canonical gauge coaction of $G$ on $\mathcal{NO}_\mathbf{X}$. Lemma~\ref{coaction is well defined} shows that there is a well-defined map $\beta:\Phi_q^{\mathcal{NO}}\left(\mathcal{NO}_\mathbf{X}\right)\rightarrow \Phi_q^{\mathcal{NO}}\left(\mathcal{NO}_\mathbf{X}\right)\otimes C^*(G)$ such that
\[
\beta\circ\Phi_q^{\mathcal{NO}}=\left(\Phi_q^{\mathcal{NO}}\otimes \mathrm{id}_{C^*(G)}\right)\circ\nu_\mathbf{X}
\]
Since $\Phi_q^{\mathcal{NO}}$, $\Phi_q^{\mathcal{NO}}\otimes \mathrm{id}_{C^*(G)}$, and $\nu_\mathbf{X}$ are all homomorphisms, we see that $\beta$ is a homomorphism. 

We now check that $\beta$ satisfies the coaction identity. We see that
\begin{align*}
\left(\beta\otimes \mathrm{id}_{C^*(G)}\right)&\circ \beta\circ\Phi_q^{\mathcal{NO}}\\
&=\left(\beta\otimes \mathrm{id}_{C^*(G)}\right)\circ \left(\Phi_q^{\mathcal{NO}}\otimes \mathrm{id}_{C^*(G)}\right)\circ \nu_\mathbf{X}\\
&=\left(\left( \beta\circ \Phi_q^{\mathcal{NO}}\right)\otimes \mathrm{id}_{C^*(G)}\right)\circ \nu_\mathbf{X}\\
&=\left(\Phi_q^{\mathcal{NO}}\otimes \mathrm{id}_{C^*(G)}\otimes \mathrm{id}_{C^*(G)}\right)\circ \left(\nu_\mathbf{X}\otimes \mathrm{id}_{C^*(G)}\right)\circ \nu_\mathbf{X}.
\end{align*}
Since $\nu_\mathbf{X}$ satisfies the coaction identity $\left(\nu_\mathbf{X}\otimes \mathrm{id}_{C^*(G)}\right)\circ \nu_\mathbf{X}=\left(\mathrm{id}_{\mathcal{NO}_\mathbf{X}}\otimes \delta_G\right)\circ\nu_\mathbf{X}$, the last line is equal to 
\begin{align*}
\big(\Phi_q^{\mathcal{NO}}&\otimes \mathrm{id}_{C^*(G)}\otimes \mathrm{id}_{C^*(G)}\big)\circ \left(\mathrm{id}_{\mathcal{NO}_\mathbf{X}}\otimes \delta_G\right)\circ \nu_\mathbf{X}\\
&=\big(\mathrm{id}_{\Phi_q^{\mathcal{NO}}\left(\mathcal{NO}_\mathbf{X}\right)}\otimes \delta_G\big)\circ \left(\Phi_q^{\mathcal{NO}}\otimes \mathrm{id}_{C^*(G)} \right)\circ \nu_\mathbf{X}\\
&=\big(\mathrm{id}_{\Phi_q^{\mathcal{NO}}\left(\mathcal{NO}_\mathbf{X}\right)}\otimes \delta_G\big)\circ \beta\circ\Phi_q^{\mathcal{NO}}. 
\end{align*}
Thus, $(\beta\otimes \mathrm{id}_{C^*(G)})\circ \beta=(\mathrm{id}_{\Phi_q(\mathcal{NO}_\mathbf{X})}\otimes \delta_G)\circ \beta$, and so $\beta$ satisfies the coaction identity. 

We also need to show that $\beta$ is coaction nondegenerate. We have
\begin{align*}
\cspan&\big\{\beta\big(\Phi_q^{\mathcal{NO}}\left(\mathcal{NO}_\mathbf{X}\right)\big)\big(1_{\mathcal{M}\left(\Phi_q^{\mathcal{NO}}\left(\mathcal{NO}_\mathbf{X}\right)\right)}\otimes C^*(G)\big)\big\}\\
&=\cspan\big\{\big(\big(\Phi_q^{\mathcal{NO}} \otimes \mathrm{id}_{C^*(G)}\big)\left(\nu_\mathbf{X}\left(\mathcal{NO}_\mathbf{X}\right)\right)\big)\big(1_{\mathcal{M}\left(\Phi_q^{\mathcal{NO}}\left(\mathcal{NO}_\mathbf{X}\right)\right)}\otimes C^*(G)\big)\big\}\\
&=\left(\Phi_q^{\mathcal{NO}} \otimes \mathrm{id}_{C^*(G)}\right)\cspan\left\{\nu_\mathbf{X}\left(\mathcal{NO}_\mathbf{X}\right)\left(1_{\mathcal{M}\left(\mathcal{NO}_\mathbf{X}\right)}\otimes C^*(G)\right)\right\}.
\end{align*}
Since $\nu_\mathbf{X}$ is coaction nondegenerate, this is equal to
\[
\left(\Phi_q^{\mathcal{NO}} \otimes \mathrm{id}_{C^*(G)}\right)\left(\mathcal{NO}_\mathbf{X}\otimes C^*(G)\right)=\Phi_q^{\mathcal{NO}}\left(\mathcal{NO}_\mathbf{X}\right)\otimes C^*(G),
\]
and we see that $\beta$ is coaction nondegenerate.

Finally, we check that $\beta$ is injective. Suppose that $\beta\left(\Phi_q^{\mathcal{NO}}(a)\right)=0$ for some $a\in \mathcal{NO}_\mathbf{X}$. We must show that $\Phi_q^{\mathcal{NO}}(a)=0$. We will make use of the fact that $\left(\Phi_q^{\mathcal{NO}}\otimes \iota\right)(\nu_\mathbf{X}(a))$ is left multiplication by $\left(\phi_\mathbf{X}^{\mathcal{NO}}\otimes \iota\right)(\nu_\mathbf{X}(a))$ on 
\[
\mathbf{Y}_q^\mathcal{NO}\otimes C^*(G\rtimes_\alpha H)\subseteq \mathcal{NO}_\mathbf{Z}\otimes C^*(G\rtimes_\alpha H) \subseteq \mathcal{M}\left(\mathcal{NO}_\mathbf{Z}\right)\otimes C^*(G\rtimes_\alpha H),
\] 
and multiplication in $\mathcal{M}\left(\mathcal{NO}_\mathbf{Z}\right)\otimes C^*(G\rtimes_\alpha H)$ is associative. For $y\in \mathbf{Y}_q^\mathcal{NO}$, Lemma~\ref{interaction between coactions} shows that
\begin{equation}
\label{equation to show coaction is injective}
\begin{aligned}
&\nu_\mathbf{Z}\left(\Phi_q^{\mathcal{NO}}(a)(y)\right)\\
&\qquad=\left(\Phi_q^{\mathcal{NO}}\otimes \iota\right)(\nu_\mathbf{X}(a))\left(\nu_\mathbf{Z}(y)\right)\\
&\qquad=\left(\Phi_q^{\mathcal{NO}}\otimes \iota\right)(\nu_\mathbf{X}(a))\big(\nu_\mathbf{Z}(y)\left(1_{\mathcal{M}(\mathcal{NO}_\mathbf{Z})}\otimes i_{G\rtimes_\alpha H}\left(e_G,q^{-1}\right)\right)\big)\\
&\qquad \qquad \qquad \times
\left(1_{\mathcal{M}(\mathcal{NO}_\mathbf{Z})}\otimes i_{G\rtimes_\alpha H}\left(e_G,q\right)\right).
\end{aligned}
\end{equation}
Since
\[
\nu_\mathbf{Z}(y)\left(1_{\mathcal{M}(\mathcal{NO}_\mathbf{Z})}\otimes i_{G\rtimes_\alpha H}\left(e_G,q^{-1}\right)\right)\in \mathbf{Y}_q^\mathcal{NO}\otimes \iota \left(C^*(G)\right)
\]
by Lemma~\ref{range of coaction}, and
\begin{align*}
&\big(\Phi_q^{\mathcal{NO}}\otimes\iota\big)(\nu_\mathbf{X}(a))\\
&\qquad=\big(\mathrm{id}_{\mathbf{Y}_q^\mathcal{NO}}\otimes \iota\big)\circ \big(\Phi_q^{\mathcal{NO}}\otimes\mathrm{id}_{C^*(G)}\big)(\nu_\mathbf{X}(a))\circ \big(\mathrm{id}_{\mathbf{Y}_q^\mathcal{NO}}\otimes \iota^{-1}\big)\\
&\qquad=\big(\mathrm{id}_{\mathbf{Y}_q^\mathcal{NO}}\otimes \iota\big)\circ \beta\big(\Phi_q^{\mathcal{NO}}(a)\big) \circ \big(\mathrm{id}_{\mathbf{Y}_q^\mathcal{NO}}\otimes \iota^{-1}\big)
=0,
\end{align*}
we conclude that \eqref{equation to show coaction is injective} is zero. As $\nu_\mathbf{Z}$ is injective, we have $\Phi_q^{\mathcal{NO}}(a)(y)=0$. As $y\in \mathbf{Y}_q^\mathcal{NO}$ was arbitrary, it follows that $\Phi_q^{\mathcal{NO}}(a)=0$. Hence, $\beta$ is injective. 

Putting all of this together and applying \cite[Corollary~4.12]{carlsen--larsen--sims--vittadello}, we conclude that $\Phi_q^{\mathcal{NO}}$ is injective. 
\end{proof}
\end{prop}

In the next subsection we will investigate the Cuntz--Nica--Pimsner algebra of $\mathbf{Y}^\mathcal{NO}$. To make our calculations tractable we will assume that $A$ acts faithfully on each fibre of $\mathbf{Z}$ --- by Proposition~\ref{injectivity of induced actions} this then implies that $\mathcal{NO}_\mathbf{X}$ acts faithfully on each fibre of $\mathbf{Y}^\mathcal{NO}$. Moreover, if $A$ acts faithfully on each fibre of $\mathbf{Z}$, then $A$ acts faithfully on each fibre of $\mathbf{X}$ and each $\widetilde{\phi}_{(p,q)}$ is injective, ensuring that all of the results from Subsection~\ref{Constructing the CNP product system} hold.

\subsection{Isomorphisms of Cuntz--Nica--Pimsner algebras}
\label{isomorphisms of CNP alg}

We now work towards showing that $\mathcal{NO}_{\mathbf{Y}^\mathcal{NO}}\cong \mathcal{NO}_\mathbf{Z}$. To do this we will use the universal property of each $C^*$-algebra to induce a homomorphism from one to the other, and then check that these homomorphisms are the inverses of each other. To make our arguments easier to write down, we will identify the coefficient algebra $\mathcal{NO}_\mathbf{X}$ of $\mathbf{Y}^{\mathcal{NO}}$ with $\phi_\mathbf{X}^{\mathcal{NO}}\left(\mathcal{NO}_\mathbf{X}\right)\subseteq \mathcal{NO}_\mathbf{Z}$. Thus, every fibre of $\mathbf{Y}^\mathcal{NO}$ can be viewed as sitting inside $\mathcal{NO}_\mathbf{Z}$, and the left and right actions of $\mathcal{NO}_\mathbf{X}\cong\phi_\mathbf{X}^{\mathcal{NO}}\left(\mathcal{NO}_\mathbf{X}\right)$ on each $\mathbf{Y}_q^\mathcal{NO}$ are multiplication in $\mathcal{NO}_\mathbf{Z}$.

We begin by getting a homomorphism from $\mathcal{NO}_\mathbf{Z}$ to $\mathcal{NO}_{\mathbf{Y}^\mathcal{NO}}$ by exhibiting a Cuntz--Nica--Pimsner covariant representation of $\mathbf{Z}$ in $\mathcal{NO}_{\mathbf{Y}^\mathcal{NO}}$. First, we need a couple of lemmas.

\begin{lem}
\label{tech result 1}
Suppose that $(K,R)$ is a quasi-lattice ordered group.  Let $F\subseteq R$ be finite. Then
\[
F=\left\{r\in F:r\leq m\right\}\sqcup \left\{r\in F:r\vee m=\infty\right\}
\]
for large $m\in R$.
\begin{proof}
Clearly, $\left\{r\in F:r\leq m\right\}\sqcup \left\{r\in F:r\vee m=\infty\right\}\subseteq F$  for any $m\in R$. Thus, it remains to show that $F\subseteq\left\{r\in F:r\leq m\right\}\sqcup \left\{r\in F:r\vee m=\infty\right\}$ for large $m$. Fix $p\in R$, and let $F'$ be a maximal element of the collection 
$
\left\{F'\subseteq F\cup \{p\}: \text{$F'$ contains $p$ and is bounded above}\right\}
$
(partially ordered by set inclusion), which exists since $F$ is finite whilst $\{p\}\subseteq F \cup \{p\}$ is bounded above (by $p$) and contains $p$. Let $n$ be an upper bound for $F'$. Since $p\in F'$, we have that $p\leq n$. Suppose that $m\in R$ with $n\leq m$. Let $r\in F$. If $r\in F'$ then $r\leq n\leq m$. Alternatively, if $r\in F\setminus F'$, then $r\vee m=\infty$ (if $r\vee m<\infty$, then $r\vee m$ is an upper bound for $F'\cup \{r\}$, which contradicts the maximality of $F'$). Thus, $F\subseteq\left\{r\in F:r\leq m\right\}\sqcup \left\{r\in F:r\vee m=\infty\right\}$. Since $p\in R$ was arbitrary, we conclude that $F\subseteq\left\{r\in F:r\leq m\right\}\sqcup \left\{r\in F:r\vee m=\infty\right\}$ for large $m$. 
\end{proof}
\end{lem}

\begin{lem}
\label{tech result 2}
Suppose that $G$ is an amenable group, and $A$ acts faithfully on each fibre of $\mathbf{Z}$, so that the product system $\mathbf{Y}^{\mathcal{NO}}$ from Proposition~\ref{existence of NO product system in NO_Z} exists. Let $(p,q)\in P\rtimes_\alpha Q$ and $T\in \mathcal{K}_A(\mathbf{Z}_{(p,q)})$. If $m\in Q$ and $q\vee m=\infty$, then
\[
M_{j_\mathbf{Z}^{((p,q))}(T)}=0\in \mathcal{K}_{\mathcal{NO}_\mathbf{X}}\left(\mathbf{Y}_m^\mathcal{NO}\right). 
\]
\begin{proof}
Fix $T\in \mathcal{K}_A(\mathbf{Z}_{(p,q)})$ and $m\in Q$ with $m\vee q=\infty$. Let $z\in \mathbf{Z}_{(r,m)}$ and $w\in \mathbf{Z}_{\left(n,e_H\right)}$. Since $(p,q)\vee (r,m)=\infty$ as $q\vee m=\infty$, the Nica covariance of $j_\mathbf{Z}$ gives
\begin{align*}
M_{j_\mathbf{Z}^{\left((p,q)\right)}(T)}&\big(j_{\mathbf{Z}_{(r,m)}}(z)j_{\mathbf{Z}_{\left(n,e_H\right)}}(w)^*\big)
=j_\mathbf{Z}^{((p,q))}(T)j_{\mathbf{Z}_{(r,m)}}(z)j_{\mathbf{Z}_{\left(n,e_H\right)}}(w)^*
=0.
\end{align*}
As $\mathbf{Y}_m^\mathcal{NO}=\cspan\big\{j_{\mathbf{Z}_{(r,m)}}(z)j_{\mathbf{Z}_{\left(n,e_H\right)}}(w)^*:r,n\in P, z\in \mathbf{Z}_{(r,m)}, w\in \mathbf{Z}_{\left(n,e_H\right)}\big\}$, we conclude that $M_{j_\mathbf{Z}^{((p,q))}(T)}=0\in \mathcal{K}_{\mathcal{NO}_\mathbf{X}}\left(\mathbf{Y}_m^\mathcal{NO}\right)$.
\end{proof}
\end{lem}

\begin{prop}
\label{varphi is CP covariant}
Suppose that $G$ is an amenable group, and $A$ acts faithfully on each fibre of $\mathbf{Z}$, so that the product system $\mathbf{Y}^{\mathcal{NO}}$ from Proposition~\ref{existence of NO product system in NO_Z} exists, and $\mathcal{NO}_\mathbf{X}$ acts faithfully on each fibre of $\mathbf{Y}^{\mathcal{NO}}$ by Proposition~\ref{injectivity of induced actions}.

Define $\varphi^\mathcal{NO}:\mathbf{Z}\rightarrow \mathcal{NO}_{\mathbf{Y}^\mathcal{NO}}$ by 
\[
\varphi_{(p,q)}^{\mathcal{NO}}:=j_{\mathbf{Y}_q^\mathcal{NO}}\circ j_{\mathbf{Z}_{(p,q)}}.
\]
Then $\varphi^{\mathcal{NO}}:\mathbf{Z}\rightarrow \mathcal{NO}_{\mathbf{Y}^\mathcal{NO}}$ is a Cuntz--Nica--Pimsner covariant representation, and so induces a homomorphism $\Omega^{\mathcal{NO}}:\mathcal{NO}_\mathbf{Z}\rightarrow \mathcal{NO}_{\mathbf{Y}^\mathcal{NO}}$ such that 
\[
\Omega^{\mathcal{NO}}\circ j_{\mathbf{Z}_{(p,q)}}=\varphi_{(p,q)}^{\mathcal{NO}}=j_{\mathbf{Y}_q^\mathcal{NO}}\circ j_{\mathbf{Z}_{(p,q)}}
\]
for each $(p,q)\in P\rtimes_\alpha Q$.
\begin{proof}
Since $j_{\mathbf{Y}^\mathcal{NO}}$ and $j_\mathbf{Z}$ are Nica covariant representations, the same reasoning as in the proofs of Proposition~\ref{varphi is a representation}, Lemma~\ref{rep on compacts}, and Proposition~\ref{varphi is Nica covariant} shows that $\varphi^{\mathcal{NO}}$ is a Nica covariant representation of $\mathbf{Z}$. We now show that it is Cuntz--Pimsner covariant. 

Let $F$ be a finite subset of $P\rtimes_\alpha Q$ and suppose that we have a set of compact operators $\big\{T_{(p,q)}\in \mathcal{K}_A\left(\mathbf{Z}_{(p,q)}\right):(p,q)\in F\big\}$ such that
\[
\sum_{(p,q)\in F}\iota_{(p,q)}^{(s,t)}\left(T_{(p,q)}\right)=0\in \mathcal{L}_A\left(\mathbf{Z}_{(s,t)}\right)
\]
for large $(s,t)\in P\rtimes_\alpha Q$. Since $j_\mathbf{Z}$ is Cuntz--Pimsner covariant, we know that
\[
\sum_{(p,q)\in F}j_\mathbf{Z}^{((p,q))}(T_{(p,q)})=0.
\]
To show that $\varphi^{\mathcal{NO}}$ is Cuntz--Pimsner covariant, we need to show that
\[
\sum_{(p,q)\in F}{\varphi^{\mathcal{NO}}}^{((p,q))}(T_{(p,q)})=0\in \mathcal{NO}_{\mathbf{Y}^\mathcal{NO}},
\]
which by Lemma~\ref{rep on compacts} is equivalent to showing that 
\[
\sum_{(p,q)\in F} j_{\mathbf{Y}^\mathcal{NO}}^{(q)}\Big(M_{j_\mathbf{Z}^{((p,q))}(T_{(p,q)})}\Big)=0\in \mathcal{NO}_{\mathbf{Y}^\mathcal{NO}}.
\]
Since $j_{\mathbf{Y}^\mathcal{NO}}$ is a Cuntz--Pimsner covariant representation of $\mathbf{Y}^\mathcal{NO}$, it suffices to show that
\[
\sum_{(p,q)\in F}\iota_q^m\Big(M_{j_\mathbf{Z}^{((p,q))}\left(T_{(p,q)}\right)}\Big)=0\in \mathcal{L}_{\mathcal{NO}_\mathbf{X}}\left(\mathbf{Y}_m^\mathcal{NO}\right) \quad \text{for large $m\in Q$,}
\]
which is equivalent to showing that 
\[
\sum_{(p,q)\in F:\,q\leq m}M_{j_\mathbf{Z}^{((p,q))}\left(T_{(p,q)}\right)}=0\in\mathcal{L}_{\mathcal{NO}_\mathbf{X}}\left(\mathbf{Y}_m^\mathcal{NO}\right)
\]
for large $m\in Q$. Making use of Lemma~\ref{tech result 1}, we see that for large $m\in Q$,
\begin{align*}
\sum_{\substack{(p,q)\in F:\\q\leq m}}M_{j_\mathbf{Z}^{((p,q))}(T_{(p,q)})}
&=\sum_{(p,q)\in F}M_{j_\mathbf{Z}^{((p,q))}(T_{(p,q)})}-\sum_{\substack{(p,q)\in F:\\q\vee m=\infty}}M_{j_\mathbf{Z}^{((p,q))}(T_{(p,q)})}\\
&=M_{\sum_{(p,q)\in F}j_\mathbf{Z}^{((p,q))}(T_{(p,q)})}-\sum_{\substack{(p,q)\in F:\\q\vee m=\infty}}M_{j_\mathbf{Z}^{((p,q))}(T_{(p,q)})}\\
&=-\sum_{\substack{(p,q)\in F:\\q\vee m=\infty}}M_{j_\mathbf{Z}^{((p,q))}(T_{(p,q)})}.
\end{align*}
It then follows from Lemma~\ref{tech result 2}, that the last line is zero, as required. Hence, $\varphi^{\mathcal{NO}}$ is Cuntz--Pimsner covariant. The final statement follows from the universal property of $\mathcal{NO}_\mathbf{Z}$. 
\end{proof}
\end{prop}

We now work towards getting a homomorphism from $\mathcal{NO}_{\mathbf{Y}^\mathcal{NO}}$ to $\mathcal{NO}_\mathbf{Z}$ inverse to the one from $\mathcal{NO}_\mathbf{Z}$ to $\mathcal{NO}_{\mathbf{Y}^\mathcal{NO}}$ just constructed. Our plan is to exhibit a Nica covariant representation of $\mathbf{Y}^\mathcal{NO}$ in $\mathcal{NO}_\mathbf{Z}$ and then use the universal property of $\mathcal{NO}_{\mathbf{Y}^\mathcal{NO}}$ to induce the required homomorphism. 

In Proposition~\ref{injectivity of induced actions} we found conditions on the group $G$ and the product system $\mathbf{Z}$ to ensure that the left actions of $\mathcal{NO}_\mathbf{X}$ on each $\mathbf{Y}_q^\mathcal{NO}$ are faithful. In the next result, we exhibit sufficient conditions for these actions to be by compact operators. 

\begin{lem}
\label{induced actions are by compacts}
Suppose that $G$ is an amenable group, and $A$ acts faithfully on each fibre of $\mathbf{Z}$, so that the product system $\mathbf{Y}^{\mathcal{NO}}$ from Proposition~\ref{existence of NO product system in NO_Z} exists. For any $q\in Q$ and $a\in \phi_{\left(e_G,q\right)}^{-1}\left(\mathcal{K}_A\left(\mathbf{Z}_{\left(e_G,q\right)}\right)\right)$, we have
\begin{equation}
\label{bumping up compacts}
\begin{aligned}
\Phi_q^\mathcal{NO}\big(j_{\mathbf{X}_{e_G}}(a)\big)
=M_{j_\mathbf{Z}^{\left(\left(e_G,q\right)\right)}(\phi_{\left(e_G,q\right)}(a))}
\in \mathcal{K}_{\mathcal{NO}_\mathbf{X}}\left(\mathbf{Y}_q^\mathcal{NO}\right).
\end{aligned}
\end{equation}
In particular, if $A$ acts compactly on $\mathbf{Z}_{\left(e_G,q\right)}$, then $\mathcal{NO}_\mathbf{X}$ acts compactly on $\mathbf{Y}_q^\mathcal{NO}$. 
\begin{proof}
Fix $q\in Q$ and $a\in \phi_{\left(e_G,q\right)}^{-1}\left(\mathcal{K}_A\left(\mathbf{Z}_{\left(e_G,q\right)}\right)\right)$. Hence, for some choice of $\mu_{k_i},\nu_{k_i}\in \mathbf{Z}_{\left(e_G,q\right)}$, we can write 
\[
\phi_{\left(e_G,q\right)}(a)=\lim_{i\rightarrow \infty}\sum_{k_i=1}^{n_i}\Theta_{\mu_{k_i},\nu_{k_i}}\in \mathcal{K}_A\left(\mathbf{Z}_{\left(e_G,q\right)}\right).
\]
For any $z\in \mathbf{Z}_{\left(e_G,q\right)}$ and $b\in \mathcal{NO}_\mathbf{X}$, we have that
\begin{align*}
&\Phi_q^\mathcal{NO}\big(j_{\mathbf{X}_{e_G}}(a)\big)\big(j_{\mathbf{Z}_{\left(e_G,q\right)}}(z)\phi_\mathbf{X}^\mathcal{NO}(b)\big)\\
&\qquad=j_{\mathbf{Z}_{\left(e_G,q\right)}}\left(\phi_{\left(e_G,q\right)}(a)(z)\right)\phi_\mathbf{X}^\mathcal{NO}(b)\\
&\qquad=j_{\mathbf{Z}_{\left(e_G,q\right)}}\bigg(\lim_{i\rightarrow \infty}\sum_{k_i=1}^{n_i}\mu_{k_i}\cdot \langle\nu_{k_i},z\rangle_A^{\left(e_G,q\right)}\bigg)\phi_\mathbf{X}^\mathcal{NO}(b).
\end{align*}
Since $j_{\mathbf{Z}_{\left(e_G,q\right)}}$ is linear and norm-decreasing, this is equal to 
\begin{align*}
&\lim_{i\rightarrow \infty}\sum_{k_i=1}^{n_i}j_{\mathbf{Z}_{\left(e_G,q\right)}}\big(\mu_{k_i}\cdot \langle\nu_{k_i},z\rangle_A^{\left(e_G,q\right)}\big)\phi_\mathbf{X}^\mathcal{NO}(b)\\
&\qquad=\lim_{i\rightarrow \infty}\sum_{k_i=1}^{n_i}M_{j_\mathbf{Z}^{\left(\left(e_G,q\right)\right)}(\Theta_{\mu_{k_i},\nu_{k_i}})}\big(j_{\mathbf{Z}_{\left(e_G,q\right)}}(z)\phi_\mathbf{X}^\mathcal{NO}(b)\big).
\end{align*}
Since the map $b\mapsto M_b$ is linear and $\norm{M_b}_{\mathcal{L}_{\mathcal{NO}_\mathbf{X}}\left(\mathbf{Y}_q^\mathcal{NO}\right)}\leq \norm{b}_{\mathcal{NO}_\mathbf{Z}}$ whenever $b\in \mathcal{NO}_\mathbf{Z}$ is such that $M_b\in \mathcal{L}_{\mathcal{NO}_\mathbf{X}}\left(\mathbf{Y}_q^\mathcal{NO}\right)$ (see Lemma~\ref{useful property for multiplication operators}), and $j_\mathbf{Z}^{\left(\left(e_G,q\right)\right)}$ is also linear and continuous, the previous line is equal to 
\begin{align*}
&M_{j_\mathbf{Z}^{\left(\left(e_G,q\right)\right)}\big(\lim_{i\rightarrow \infty}\sum_{k_i=1}^{n_i}\Theta_{\mu_{k_i},\nu_{k_i}}\big)}\big(j_{\mathbf{Z}_{\left(e_G,q\right)}}(z)\phi_\mathbf{X}^\mathcal{NO}(b)\big)\\
&\qquad=M_{j_\mathbf{Z}^{\left(\left(e_G,q\right)\right)}\big(\phi_{\left(e_G,q\right)}(a)\big)}\big(j_{\mathbf{Z}_{\left(e_G,q\right)}}(z)\phi_\mathbf{X}^\mathcal{NO}(b)\big).
\end{align*}
As $\Phi_q^\mathcal{NO}\big(j_{\mathbf{X}_{e_G}}(a)\big)$ and $M_{j_\mathbf{Z}^{\left(\left(e_G,q\right)\right)}(\phi_{\left(e_G,q\right)}(a))}$ are linear and continuous, whilst $\mathbf{Y}_q^\mathcal{NO}=\cspan\big\{j_{\mathbf{Z}_{\left(e_G,q\right)}}(z)\phi_{\mathbf{X}}^{\mathcal{NO}}(b):z\in \mathbf{Z}_{\left(e_G,q\right)},  b\in \mathcal{NO}_\mathbf{X}\big\}$, we conclude that 
$
\Phi_q^\mathcal{NO}\big(j_{\mathbf{X}_{e_G}}(a)\big)
$
is left multiplication by $j_\mathbf{Z}^{\left(\left(e_G,q\right)\right)}\left(\phi_{\left(e_G,q\right)}(a)\right)$. To establish \eqref{bumping up compacts}, it remains to show that multiplication by $j_\mathbf{Z}^{\left(\left(e_G,q\right)\right)}\left(\phi_{\left(e_G,q\right)}(a)\right)$ is a compact operator on$\mathbf{Y}_q^\mathcal{NO}$. Since $j_\mathbf{Z}^{\left(\left(e_G,q\right)\right)}\left(\phi_{\left(e_G,q\right)}(a)\right)\in \mathcal{NO}_\mathbf{Z}^q$, this follows from Lemma~\ref{compact operators in Y}. 

Now suppose that $A$ acts compactly on $\mathbf{Z}_{\left(e_G,q\right)}$. Fix $x\in \mathbf{X}_p$ and choose $x'\in  \mathbf{X}_p$ so that $x=x'\cdot \langle x',x'\rangle_A^p$. Since 
$\langle x',x'\rangle_A^p\in \phi_{\left(e_G,q\right)}^{-1}\left(\mathcal{K}_A\left(\mathbf{Z}_{\left(e_G,q\right)}\right)\right)$, we can apply the first part of the lemma to see that
\begin{align*}
\Phi_q^\mathcal{NO}\left(j_{\mathbf{X}_p}(x)\right)
=\Phi_q^\mathcal{NO}\left(j_{\mathbf{X}_p}(x')\right)\Phi_q^\mathcal{NO}\big(j_{\mathbf{X}_{e_G}}(\langle x',x'\rangle_A^p)\big)
\in \mathcal{K}_{\mathcal{NO}_\mathbf{X}}\left(\mathbf{Y}_q^\mathcal{NO}\right).
\end{align*}
Since $\Phi_q^{\mathcal{NO}}$ is a homomorphism and $\mathcal{NO}_\mathbf{X}$ is generated by the image of $j_\mathbf{X}$, we have that $\Phi_q^{\mathcal{NO}}\left(\mathcal{NO}_\mathbf{X}\right)\subseteq \mathcal{K}_{\mathcal{NO}_\mathbf{X}}\left(\mathbf{Y}_q^\mathcal{NO}\right)$, and so $\mathcal{NO}_\mathbf{X}$ acts compactly on $\mathbf{Y}_q^\mathcal{NO}$. 
\end{proof}
\end{lem}

We are almost ready to show that the inclusion of $\mathbf{Y}^\mathcal{NO}$ in $\mathcal{NO}_\mathbf{Z}$ is a Cuntz--Nica--Pimsner covariant representation. Before we do so, we need one last lemma.

\begin{lem}
\label{directed set and CP covariance}
Let $q\in Q$ and $a\in \phi_{\left(e_G,q\right)}^{-1}\left(\mathcal{K}_A\left(\mathbf{Z}_{\left(e_G,q\right)}\right)\right)$.
If $Q$ is directed, then 
\[
j_{\mathbf{Z}_{\left(e_G,e_H\right)}}\left( a\right)=j_\mathbf{Z}^{\left(\left(e_G,q\right)\right)}\left(\phi_{\left(e_G,q\right)}\left(a\right)\right).
\]
\begin{proof}
Since $j_\mathbf{Z}$ is Cuntz--Pimsner covariant, it suffices to show that 
\begin{equation}
\label{equation for varphi' is CNP covariant}
\begin{aligned}
\iota_{\left(e_G,e_H\right)}^{(s,t)}\left(a\right)-\iota_{\left(e_G,q\right)}^{(s,t)}\left(\phi_{\left(e_G,q\right)}\left(a\right)\right)
&=0 \in \mathcal{L}_A\left(\mathbf{Z}_{(s,t)}\right)\\
&\quad\text{for large $(s,t)\in P\rtimes_\alpha Q$.} 
\end{aligned}
\end{equation}
Fix $(m,n)\in P\rtimes_\alpha Q$. Since $Q$ is directed, $(m,n)\leq (m, n\vee q)$. Suppose that $(s,t)\in P\rtimes_\alpha Q$ with $(m,n\vee q)\leq (s,t)$. For any $z\in \mathbf{Z}_{\left(e_G,q\right)}$ and $w\in \mathbf{Z}_{\left(\alpha_{q^{-1}}(s),q^{-1}t\right)}$, we have
\begin{align*}
&\Big(\iota_{\left(e_G,e_H\right)}^{(s,t)}\left(a\right)-\iota_{\left(e_G,q\right)}^{(s,t)}\left(\phi_{\left(e_G,q\right)}\left(a\right)\right)\Big)\left(zw\right)\\
&\qquad=\phi_{(s,t)}\left(a\right)(zw)-\left(\phi_{\left(e_G,q\right)}\left(a\right)(z)\right)w\\
&\qquad=\left(\left(\phi_{\left(e_G,q\right)}\left(a\right)-\left(\phi_{\left(e_G,q\right)}\left(a\right)\right)\right)(z)\right)w
=0.
\end{align*}
Since $\iota_{\left(e_G,e_H\right)}^{(s,t)}\left(a\right)-\iota_{\left(e_G,q\right)}^{(s,t)}\left(\phi_{\left(e_G,q\right)}\left(a\right)\right)\in \mathcal{L}_A\left(\mathbf{Z}_{(s,t)}\right)$ is linear and continuous, and since $\mathbf{Z}_{(s,t)}=\cspan\{zw:z\in \mathbf{Z}_{\left(e_G,q\right)}, w\in \mathbf{Z}_{(\alpha_{q^{-1}}(s),q^{-1}t)}\}$, we conclude that \eqref{equation for varphi' is CNP covariant} holds. 
\end{proof}
\end{lem}

\begin{prop}
\label{varphi' is CNP covariant}
Suppose that $G$ is an amenable group, $A$ acts faithfully on each fibre of $\mathbf{Z}$, so that the product system $\mathbf{Y}^{\mathcal{NO}}$ from Proposition~\ref{existence of NO product system in NO_Z} exists, and $\mathcal{NO}_\mathbf{X}$ acts faithfully on each fibre of $\mathbf{Y}^{\mathcal{NO}}$ by Proposition~\ref{injectivity of induced actions}. Moreover, suppose that $A$ acts compactly on each $\mathbf{Z}_{\left(e_G,q\right)}$, so that $\mathcal{NO}_\mathbf{X}$ acts compactly on each fibre of $\mathbf{Y}^{\mathcal{NO}}$ by Lemma~\ref{induced actions are by compacts}.

For each $q\in Q$, let ${\varphi_q'}^{\mathcal{NO}}$ be the inclusion of $\mathbf{Y}_q^\mathcal{NO}$ in $\mathcal{NO}_\mathbf{Z}$. If $Q$ is directed, then $\varphi'^{\mathcal{NO}}$ is a Cuntz--Nica--Pimsner covariant representation of $\mathbf{Y}^\mathcal{NO}$. Hence, there exists a homomorphism $\Omega'^{\mathcal{NO}}:\mathcal{NO}_{\mathbf{Y}^\mathcal{NO}}\rightarrow \mathcal{NO}_\mathbf{Z}$ such that 
\[
\Omega'^{\mathcal{NO}}\circ j_{\mathbf{Y}_q^\mathcal{NO}}=\varphi'_q
\] 
for each $q\in Q$.
\begin{proof}
The same reasoning as in Proposition~\ref{varphi' is Nica covariant} shows that $\varphi'^{\mathcal{NO}}$ is a Nica covariant representation of $\mathbf{Y}^\mathcal{NO}$ in $\mathcal{NO}_\mathbf{Z}$. It remains to show that $\varphi'^{\mathcal{NO}}$ is Cuntz--Pimsner covariant. By Proposition~\ref{injectivity of induced actions} and Lemma~\ref{induced actions are by compacts}, $\mathcal{NO}_\mathbf{X}$ acts faithfully and compactly on each $\mathbf{Y}_q^\mathcal{NO}$. Hence, by \cite[Proposition~5.1]{sims--yeend}, it suffices to check that $\left(\varphi'^{\mathcal{NO}}\right)^{(q)}\circ \Phi_q^{\mathcal{NO}}=\varphi'^{\mathcal{NO}}_{e_H}$ for each $q\in Q$. This is clear when $q=e_H$, so we just need to worry about when $q\neq e_H$. As $\mathbf{Y}_{e_H}^\mathcal{NO}=\mathcal{NO}_\mathbf{X}\cong\phi_\mathbf{X}^{\mathcal{NO}}\left(\mathcal{NO}_\mathbf{X}\right)$ is generated by the image of ${j_\mathbf{Z}}|_{\mathbf{X}}$ it suffices to show that  
\begin{equation}
\label{sufficient for covariance}
\begin{aligned}
\big(\big(\varphi'^{\mathcal{NO}}\big)^{(q)}\circ \Phi_q^{\mathcal{NO}}\big)\left(j_\mathbf{Z}(x)\right)=\varphi'^{\mathcal{NO}}_{e_H}\left(j_\mathbf{Z}(x)\right) 
\quad \text{for all $x\in \mathbf{X}$.}
\end{aligned}
\end{equation}
If $x\in \mathbf{X}_p$ and $x'\in \mathbf{X}_p$ is chosen so that $x=x'\cdot \langle x',x'\rangle_A^p$, then
\[
\varphi'^{\mathcal{NO}}_{e_H}\big(j_{\mathbf{Z}_{\left(p,e_H\right)}}(x)\big)=j_{\mathbf{Z}_{\left(p,e_H\right)}}(x)=j_{\mathbf{Z}_{\left(p,e_H\right)}}(x')j_{\mathbf{Z}_{\left(e_G,e_H\right)}}(\langle x',x'\rangle_A^p).
\]
By the first part of Lemma~\ref{induced actions are by compacts}, we see that
\begin{align*}
\Phi_q^{\mathcal{NO}}\big(j_{\mathbf{Z}_{\left(p,e_H\right)}}(x)\big)
=M_{j_{\mathbf{Z}_{\left(p,e_H\right)}}(x')j_\mathbf{Z}^{\left(\left(e_G,q\right)\right)}\big(\phi_{\left(e_G,q\right)}\left(\langle x',x'\rangle_A^p\right)\big)}\in \mathcal{K}_{\mathcal{NO}_\mathbf{X}}\left(\mathbf{Y}_q^\mathcal{NO}\right).
\end{align*}
Since $j_{\mathbf{Z}_{\left(p,e_H\right)}}(x')j_\mathbf{Z}^{\left(\left(e_G,q\right)\right)}\left(\phi_{\left(e_G,q\right)}\left(\langle x',x'\rangle_A^p\right)\right)\in \mathcal{NO}_\mathbf{Z}^q$, it follows from the argument in Proposition~\ref{varphi' is Nica covariant} that 
\begin{align*}
&\big(\big(\varphi'^{\mathcal{NO}}\big)^{(q)}\circ \Phi_q^{\mathcal{NO}}\big)\big(j_{\mathbf{Z}_{\left(p,e_H\right)}}(x)\big)\\
&\qquad=\left(\varphi'^{\mathcal{NO}}\right)^{(q)}\bigg(M_{j_{\mathbf{Z}_{\left(p,e_H\right)}}(x')j_\mathbf{Z}^{\left(\left(e_G,q\right)\right)}\big(\phi_{\left(e_G,q\right)}\left(\langle x',x'\rangle_A^p\right)\big)}\bigg)\\
&\qquad=j_{\mathbf{Z}_{\left(p,e_H\right)}}(x')j_\mathbf{Z}^{\left(\left(e_G,q\right)\right)}\left(\phi_{\left(e_G,q\right)}(\langle x',x'\rangle_A^p)\right).
\end{align*}
Since 
\[
j_{\mathbf{Z}_{\left(e_G,e_H\right)}}(\langle x',x'\rangle_A^p)=j_\mathbf{Z}^{\left(\left(e_G,q\right)\right)}\left(\phi_{\left(e_G,q\right)}(\langle x',x'\rangle_A^p)\right)
,\]
by Lemma~\ref{directed set and CP covariance}, we conclude that \eqref{sufficient for covariance} holds.
\end{proof}
\end{prop}

\begin{rem}
\label{are all the hypoth necessary}
It is not clear if all of the hypotheses in Proposition~\ref{varphi' is CNP covariant} are necessary. Proposition~4.14 of \cite{k-graph iterated} suggests that the assumption of $A$ acting compactly on each $\mathbf{Z}_{(e_G,q)}$ is not necessary, at least in the situation where $\mathbf{Z}$ is a product system over $\N^k$ and we can make use of Katsura's work on gauge-invariant ideals of Cuntz--Pimsner algebras \cite[Theorem~8.6]{katsura2}.
\end{rem}

Putting all of this together, we can show that $\mathcal{NO}_\mathbf{Z}$ and $\mathcal{NO}_{\mathbf{Y}^\mathcal{NO}}$ are isomorphic. 

\begin{thm}
\label{isomorphism of Cuntz Nica Pimsner algebras}
Suppose that $G$ is an amenable group, $A$ acts faithfully on each fibre of $\mathbf{Z}$, so that the product system $\mathbf{Y}^{\mathcal{NO}}$ from Proposition~\ref{existence of NO product system in NO_Z} exists, and $\mathcal{NO}_\mathbf{X}$ acts faithfully on each fibre of $\mathbf{Y}^{\mathcal{NO}}$ by Proposition~\ref{injectivity of induced actions}. Moreover, suppose that $A$ acts compactly on each $\mathbf{Z}_{\left(e_G,q\right)}$, so that $\mathcal{NO}_\mathbf{X}$ acts compactly on each fibre of $\mathbf{Y}^{\mathcal{NO}}$ by Lemma~\ref{induced actions are by compacts}.

Suppose that $Q$ is directed, so that the homomorphism $\Omega'^{\mathcal{NO}}$ from Proposition~\ref{varphi' is CNP covariant} exists. Then the homomorphisms $\Omega^{\mathcal{NO}}:\mathcal{NO}_\mathbf{Z}\rightarrow \mathcal{NO}_{\mathbf{Y}^\mathcal{NO}}$ and $\Omega'^{\mathcal{NO}}:\mathcal{NO}_{\mathbf{Y}^\mathcal{NO}}\rightarrow \mathcal{NO}_\mathbf{Z}$ are mutually inverse. Thus, $\mathcal{NO}_\mathbf{Z}\cong \mathcal{NO}_{\mathbf{Y}^\mathcal{NO}}$. 
\begin{proof}
The same reasoning as in the proof of Theorem~\ref{isomorphisms of Nica--Toeplitz algebras} shows that $\Omega^{\mathcal{NO}}$ and $\Omega'^{\mathcal{NO}}$ are mutually inverse. 
\end{proof}
\end{thm}

Theorem~\ref{isomorphism of Cuntz Nica Pimsner algebras} enables us to view the Cuntz--Nica--Pimsner algebra of a compactly aligned product system over $\N^k$ (with faithful and compact left actions) as a $k$-times iterated Cuntz--Pimsner algebra. As in Corollary~\ref{extending literature results for NT algebras}, this immediately allows us to generalise results from the literature that relate the coefficient algebra of a Hilbert bimodule and its Cuntz--Pimsner algebra.

\begin{cor}
\label{extending literature results for CNP algebras}
Let $\mathbf{Z}$ be a product system over $\N^k$ with coefficient algebra $A$. Suppose $A$ acts faithfully and compactly on the fibres $\{\mathbf{Z}_{e_i}:1\leq i\leq k\}$. Then
\begin{enumerate}[label=\upshape(\roman*)]
\item $A$ is exact if and only if $\mathcal{NO}_\mathbf{Z}$ is exact\/{\rm ;}
\item If $A$ is nuclear, then $\mathcal{NO}_\mathbf{Z}$ is nuclear\/{\rm ;}
\item If each $\mathbf{Z}_{e_i}$ is countably generated as a Hilbert $A$-module and $A$ is separable, nuclear, and satisfies the universal coefficient theorem of \cite{rosenberg--schochet}, then $\mathcal{NO}_\mathbf{Z}$ satisfies the universal coefficient theorem.
\end{enumerate}
\begin{proof}
Firstly, we observe that since $A$ acts faithfully and compactly on each $\mathbf{Z}_{e_i}$, it follows from \cite[Proposition~4.7]{lance} that $A$ acts faithfully and compactly on every fibre of $\mathbf{Z}$. Since $A$ acts compactly on each fibre of $\mathbf{Z}$, \cite[Proposition~5.8]{fowler} says that $\mathbf{Z}$ is compactly aligned. Part (i) of the corollary then follows from \cite[Theorem~7.1]{katsura}, and part (ii) from \cite[Corollary~7.4]{katsura}.

If $A$ is separable and each $\mathbf{Z}_{e_i}$ is countably generated then the Cuntz--Nica--Pimsner algebra we get at each stage of the iterative procedure is separable. Furthermore, if $A$ is nuclear, then part (ii) ensures that the Cuntz--Nica--Pimsner algebra we get at each stage is also nuclear. Hence, part (iii) of the corollary follows from \cite[Proposition~8.8]{katsura}.
\end{proof}
\end{cor}

\begin{rem}
The Pimsner--Voiculescu exact sequence \cite[Theorem~8.6]{katsura} relates the $K$-theory of a Cuntz--Pimsner algebra associated to a Hilbert bimodule, to the $K$-theory of the bimodule's coefficient algebra. As such it may be possible to gain $K$-theoretic information about Cuntz--Nica--Pimsner algebras associated to compactly aligned product systems over $\N^k$ by viewing them as iterated Cuntz--Pimsner algebras. In particular we are interested to see if such an approach could be used to compute the $K$-theory of finitely aligned $k$-graph algebras when $k\geq 3$ (the $k=1$ case is covered by \cite[Theorem~6.1]{bates--hong--raeburn--szymanski} and the $k=2$ case by \cite[Proposition~3.16]{evans}). 
\end{rem}

\subsection{Examples}

We now present three examples that show what kind of product systems our two decomposition theorems (Theorem~\ref{isomorphisms of Nica--Toeplitz algebras} and Theorem~\ref{isomorphism of Cuntz Nica Pimsner algebras}) can be applied to, and what these results reveal about the structure of the associated Nica--Toeplitz and Cuntz--Nica--Pimsner algebras. 

Our first example/application looks at twisted crossed products by groups \cite{packer--raeburn}. We begin by recapping what we mean by twisted crossed products, and explain their relationship with product systems. Let $G$ be a (discrete) group, $A$ a $C^*$-algebra, $\beta:G\rightarrow \mathrm{Aut}(A)$ a homomorphism, and $\omega:G\times G\rightarrow \T$ a normalised $2$-cocycle, i.e. $\omega(e_G,e_G)=1$ and
\[
\omega(gh,k)\omega(g,h)=\omega(g,hk)\omega(h,k) \quad \text{for all $g,h,k\in G$.}
\]
The twisted crossed product $A\times_{\beta,\omega} G$ is by definition the universal $C^*$-algebra generated by a covariant representation of the twisted dynamical system $(A,G,\beta,\omega)$: that is, a homomorphism $i_A^{\beta,\omega}:A\rightarrow A\times_{\beta,\omega} G$ and multiplier unitaries $\{i_G^{\beta,\omega}(g):g\in G\}$ such that for $g,h\in G$ and $a\in A$, 
\[
i_G^{\beta,\omega}(g)i_G^{\beta,\omega}(h)=\omega(g,h)i_G^{\beta,\omega}(gh)
\]
\[
i_G^{\beta,\omega}(g)i_A^{\beta,\omega}(a)i_G^{\beta,\omega}(g)^*=i_A^{\beta,\omega}(\beta_g(a)).
\]

When $G$ contains a subsemigroup $P$ such that $(G,P)$ is quasi-lattice ordered, $P$ is directed, and $P$ generates $G$ as a group, \cite[\S~5]{carlsen--larsen--sims--vittadello} tells us how to realise $A\times_{\beta,\omega} G$ as a Cuntz--Nica--Pimsner algebra. 

We write $G^{\mathrm{op}}$ for the opposite group of $G$: as a set $G^{\mathrm{op}}$ is equal to $G$, with multiplication in $G^{\mathrm{op}}$ given by $g*_\mathrm{op}h:=hg$. Since the isomorphism $h\mapsto h^{-1}$ from $G$ to $G^{\mathrm{op}}$ maps $P$ to $P^{-1}$, we see that if $g,h\in G$, then $g\leq h$ in $(G,P)$ if and only if  $g^{-1}\leq h^{-1}$ in $(G^\mathrm{op},P^{-1})$. It follows that $(G^\mathrm{op},P^{-1})$ is quasi-lattice ordered: the least upper bound of $g,h\in G^\mathrm{op}$ in $(G^\mathrm{op},P^{-1})$ is equal to $(g^{-1}\vee h^{-1})^{-1}$ if the least upper bound $g^{-1}\vee h^{-1}$ of $g^{-1},h^{-1}$ in $(G,P)$ exists, and is $\infty$ otherwise. Observe that $P$ is directed (as a subset of $G$) if and only if $P^{-1}$ is directed (as a subset of $G^\mathrm{op}$).

We define a product system $\mathbf{Z}:=\bigsqcup_{p\in P^{-1}}\mathbf{Z}_p$ over $(G^\mathrm{op},P^{-1})$ as follows. For each $p\in P^{-1}$, $\mathbf{Z}_p$ is the Hilbert $A$-bimodule ${}_{\beta_p}A_A$, which is equal to $A$ as a vector space, with operations 
\[
\langle x,y\rangle_A^p:=x^*y \qquad a\cdot x:=\beta_p(a)x \qquad x\cdot a:=xa
\]
for all $x,y\in \mathbf{Z}_p$ and $a\in A$. For each $p,q\in P^{-1}$ there exists a Hilbert $A$-bimodule isomorphism $M^\mathbf{Z}_{p,q}:\mathbf{Z}_p\otimes_A \mathbf{Z}_q\rightarrow \mathbf{Z}_{qp}$ such that
\[
M^\mathbf{Z}_{p,q}(x\otimes_A y)=\overline{\omega(q,p)}\beta_q(x)y \quad \text{for all $x\in \mathbf{Z}_p$, $y\in \mathbf{Z}_q$.}
\]
It is not difficult to see that the homomorphism $\phi_p:A \rightarrow \mathcal{L}_A(\mathbf{Z}_p)$ implementing the left action of $A$ on $\mathbf{Z}_p$ is injective and takes values in $\mathcal{K}_A(\mathbf{Z}_p)$: if $\{u_\lambda\}_{\lambda\in \Lambda}$ is an approximate identity for $A$, then $\phi_p(a)=\lim_{\lambda\in \Lambda} \Theta_{\beta_p(a),u_\lambda}$. Thus, it follows from  \cite[Proposition~3.1.24]{fletcherphd} that $\mathbf{Z}$ is compactly aligned. Finally, \cite[Corollary~5.2]{carlsen--larsen--sims--vittadello} tells us that there is an isomorphism from $\mathcal{NO}_\mathbf{Z}$ to $A\times_{\beta,\omega} G$ that takes $j_{\mathbf{Z}_p}(x)$ to $i_G^{\beta,\omega}(p)^*i_A^{\beta,\omega}(x)$ for each $p\in P^{-1}$ and $x\in \mathbf{Z}_p$.

We now use our decomposition theorem for Cuntz--Nica--Pimsner algebras (Theorem~\ref{isomorphism of Cuntz Nica Pimsner algebras}) to show that the twisted crossed product of a $C^*$-algebra by a semidirect product $G\rtimes_\alpha H$ can be realised as an iterated twisted crossed product, first by $G$, and then by $H$. 

\begin{exam}
Let $G$ and $H$ be groups, and $\alpha:H\rightarrow \mathrm{Aut}(G)$ a homomorphism. Let $A$ be a $C^*$-algebra, $\beta:G\rtimes_\alpha H\rightarrow \mathrm{Aut}(A)$ a homomorphism, and $\omega:(G\rtimes_\alpha H)\times (G\rtimes_\alpha H)\rightarrow \T$ a normalised $2$-cocycle. Define normalised $2$-cocycles $\omega|_G:G\times G\rightarrow \T$ and $\omega|_H:H\times H\rightarrow \T$ by $\omega|_G(g_1,g_2):=\omega((g_1,e_H),(g_2,e_H))$ and $\omega|_H(h_1,h_2):=\omega((e_G,h_1),(e_G,h_2))$ respectively, and define $\beta|_G:G\rightarrow \mathrm{Aut}(A)$ by $(\beta|_G)_g:=\beta_{(g,e_H)}$. Then there exists a homomorphism $\widetilde{\beta}|_H:H\rightarrow \mathrm{Aut}(A\times_{\beta|_G,\omega|_G}G)$ such that 
\begin{equation}
\label{extended action}
\begin{aligned}
(\widetilde\beta|_H)_h\big(i_A^{\beta|_G,\omega|_G}(a)\big)&=i_A^{\beta|_G,\omega|_G}(\beta_{(e_G,h)}(a))\\
(\widetilde\beta|_H)_h\big(i_G^{\beta|_G,\omega|_G}(g)\big)&=\frac{\omega((e_G,h),(g,e_H))}{\omega((\alpha_h(g),e_H),(e_G,h))}i_G^{\beta|_G,\omega|_G}(\alpha_h(g))
\end{aligned}
\end{equation}
for each $a\in A$ and $g\in G$ (where in the second expression we are thinking of the extension of $(\widetilde\beta|_H)_h$ to the multiplier algebra of $A\times_{\beta|_G,\omega|_G}G$). 

Furthermore, suppose that $G$ contains a subsemigroup $P$, and $H$ contains a subsemigroup $Q$, such that $(G,P)$ and $(H,Q)$ are quasi-lattice ordered, $P$ and $Q$ are directed, and $P$ and $Q$ generate $G$ and $H$ respectively. If $G$ is amenable and $\alpha_H(P)\subseteq P$, then 
\[
(A\times_{\beta|_G,\omega|_G}G)\times_{\widetilde\beta|_H,\omega|_H}H \cong A\times_{\beta,\omega}(G\rtimes_\alpha H).
\]
Explicitly, this  isomorphism is given by
\begin{equation}
\label{explicit form}
\begin{aligned}
i_{H}^{\widetilde\beta|_H, \omega|_H}(h)^*
i_{A\times_{\beta|_G,\omega|_G} G}^{\widetilde\beta|_H, \omega|_H}&\big(i_G^{\beta|_G, \omega_G}(g)^*i_A^{\beta|_G, \omega_G}(a)\big)\\
&\mapsto
\overline{\omega((g,e_H), (e_G,h))}
i_{G\rtimes_\alpha H}^{\beta, \omega}(g,h)^*i_{A}^{\beta, \omega}(a)
\end{aligned}
\end{equation}
for $(g,h)\in G\rtimes_\alpha H$ and $a\in A$.
\begin{proof}
Firstly, we note that this is certainly not a novel result: \cite[Theorem~4.1]{packer--raeburn} deals with the more general situation where the topology on $G\rtimes_\alpha H$ need not be discrete. We also point out that the additional hypotheses on the groups $G$ and $H$ and the action $\alpha$ are to enable us to apply our decomposition theorem, as well as to allow us to view the three twisted crossed products $A\times_{\beta,\omega}(G\rtimes_\alpha H)$, $A\times_{\beta|_G,\omega|_G}G$, and $(A\times_{\beta|_G,\omega|_G}G)\times_{\widetilde\beta|_H,\omega|_H}H$ as Cuntz--Nica--Pimsner algebras, and are not required in \cite[Theorem~4.1]{packer--raeburn}.

For each $h\in H$, it is straightforward to check that the homomorphism $i_A^{\beta|_G,\omega|_G}\circ \beta_{(e_G,h)}$ and the multiplier unitaries 
\[
\bigg\{\frac{\omega((e_G,h),(g,e_H))}{\omega((\alpha_h(g),e_H),(e_G,h))}i_G^{\beta|_G,\omega|_G}(\alpha_h(g)):g\in G\bigg\}
\] 
give a covariant representation in $A\times_{\beta|_G,\omega|_G}G$ of the twisted dynamical system $(A,G,\beta|_G, \omega|_G)$. Hence, the universal property of $A\times_{\beta|_G,\omega|_G}G$ induces a homomorphism $(\widetilde\beta|_H)_h:A\times_{\beta|_G,\omega|_G}G\rightarrow A\times_{\beta|_G,\omega|_G}G$ satisfying \eqref{extended action}. One can check on generators that $(\widetilde\beta|_H)_{e_H}=\mathrm{id}_{A\times_{\beta|_G,\omega|_G}G}$ and $(\widetilde\beta|_H)_{h}\circ (\widetilde\beta|_H)_{k}=(\widetilde\beta|_H)_{hk}$ for each $h,k\in H$. Thus, the map $h\mapsto (\widetilde\beta|_H)_{h}$ gives a homomorphism from $H$ to $\mathrm{Aut}(A\times_{\beta|_G,\omega|_G}G)$, which we denote by $\widetilde\beta|_H$.

Now let $\mathbf{Z}:=\bigsqcup_{(p,q)\in (P\rtimes_\alpha Q)^{-1}} \mathbf{Z}_{(p,q)}:=\bigsqcup_{(p,q)\in (P\rtimes_\alpha Q)^{-1}} {}_{\beta_{(p,q)}}A_A$ be the compactly aligned product system over $((G\rtimes_\alpha H)^\mathrm{op}, (P\rtimes_\alpha Q)^{-1})$ associated to the twisted dynamical system $(A,G\rtimes_\alpha H, \beta, \omega)$. Thus, $ A\times _{\beta,\omega}(G\rtimes_\alpha H)\cong \mathcal{NO}_\mathbf{Z}$.

In order to apply our decomposition theorem to the product system $\mathbf{Z}$ we need to know that the quasi-lattice ordered group  $((G\rtimes_\alpha H)^\mathrm{op}, (P\rtimes_\alpha Q)^{-1})$ decomposes as the semidirect product of quasi-lattice ordered groups. It is routine to check that the map $\eta:(G\rtimes_\alpha H)^\mathrm{op}\rightarrow G^\mathrm{op}\rtimes_{\alpha\circ \mathrm{inv}_H}H^\mathrm{op}$ (where $\mathrm{inv}_H:H^\mathrm{op}\rightarrow H$ is the isomorphism that sends $h$ to $h^{-1}$) defined by $\eta(g,h):=(\alpha_{h^{-1}}(g),h)$ is an isomorphism. Moreover, since $\eta((P\rtimes_\alpha Q)^{-1})=P^{-1}\rtimes_{\alpha\circ \mathrm{inv}_H}Q^{-1}$, we have $(g,h)\leq (s,t)$ in $((G\rtimes_\alpha H)^\mathrm{op},(P\rtimes_\alpha Q)^{-1})$ if and only if $\eta(g,h)\leq \eta(s,t)$ in $(G^\mathrm{op}\rtimes_{\alpha\circ \mathrm{inv}_H}H^\mathrm{op}, P^{-1}\rtimes_{\alpha\circ \mathrm{inv}_H}Q^{-1})$. Hence, as in Proposition~\ref{inclusion of NT_X in NT_Z}, we may define a compactly aligned product system $\mathbf{X}:=\bigsqcup_{p\in P^{-1}}\mathbf{X}_p$ over $(G^\mathrm{op}, P^{-1})$ with coefficient algebra $A$ by setting $\mathbf{X}_p:=\mathbf{Z}_{\eta^{-1}(p,e_H)}=\mathbf{Z}_{(p,e_H)}={}_{(\beta|_G)_p}A_A$ for each $p\in P^{-1}$. Clearly, $\mathbf{X}$ is the product system associated to the twisted dynamical system $(A,G,\beta|_G, \omega|_G)$. By \cite[Corollary~5.2]{carlsen--larsen--sims--vittadello}, there is an isomorphism $\tau:\mathcal{NO}_\mathbf{X}\rightarrow A\times_{\beta|_G,\omega|_G} G$ such that
\begin{equation}
\label{formula for tau}
\tau(j_{\mathbf{X}_p}(x))=i_G^{\beta|_G,\omega|_G}(p)^*i_A^{\beta|_G,\omega|_G}(x) \quad \text{for $p\in P^{-1}$ and $x\in \mathbf{X}_p$.}
\end{equation} 

We now consider the compactly aligned product system $\mathbf{Y}^\mathcal{NO}$, as defined in Propositions~\ref{existence of fibres in NO product system in NO_Z} and \ref{existence of NO product system in NO_Z}, over $(H^\mathrm{op},Q^{-1})$ with coefficient algebra $\mathcal{NO}_\mathbf{X}$. For $q\in Q^{-1}$, we have $\eta^{-1}(e_G,q)=(e_G,q)$, and so the fibres of $\mathbf{Y}^\mathcal{NO}$ are given by
\[
\mathbf{Y}_q^\mathcal{NO}=\cspan\{j_{\mathbf{Z}_{(e_G,q)}}(\mathbf{Z}_{(e_G,q)})\phi_\mathbf{X}^\mathcal{NO}(\mathcal{NO}_\mathbf{X})\}.
\]
We claim that $\mathbf{Y}^\mathcal{NO}$ is isomorphic to the product system associated to the twisted dynamical system $(A \times_{\beta|_G,\omega|_G}G, H, \widetilde\beta|_H, \omega|_H)$, where we identify the respective coefficient algebras $\mathcal{NO}_\mathbf{X}$ and $A \times_{\beta|_G,\omega|_G}G$ via the isomorphism $\nu_{e_H}:=\tau$. For simplicity, we also identify $\mathcal{NO}_\mathbf{X}$ with its image in $\mathcal{NO}_\mathbf{Z}$ under the injective homomorphism $\phi_\mathbf{X}^\mathcal{NO}$ . 

Straightforward calculations show that for $q\in Q^{-1}$ there exists a linear inner-product preserving map $\nu_q:\mathbf{Y}_q^\mathcal{NO}\rightarrow {}_{(\widetilde\beta|_H)_q}(A \times_{\beta|_G,\omega|_G}G)_{(A \times_{\beta|_G,\omega|_G}G)}$ such that 
\begin{equation}
\label{definition of nu_q}
\nu_q\big(j_{\mathbf{Z}_{(e_G,q)}}(a)\phi_\mathbf{X}^\mathcal{NO}(s)\big)=i_A^{\beta|_G, \omega|_G}(a)\tau(s) \quad \text{for $a\in \mathbf{Z}_{(e_G,q)}$, $s\in \mathcal{NO}_\mathbf{X}$.}
\end{equation}
It is not difficult to see that each $\nu_q$ is surjective. Since every element of $A$ can be written as a product of two elements, it suffices to show that $i_G^{\beta|_G, \omega|_G}(g)i_A^{\beta|_G, \omega|_G}(bc)$ belongs to the range of $\nu_q$ for any $g\in G$ and $b,c\in A$. If we choose $s\in \mathcal{NO}_\mathbf{X}$ such that $\tau(s)=i_G^{\beta|_G, \omega|_G}(g)i_A^{\beta|_G, \omega|_G}(c)$, then 
\begin{align*}
\nu_q\big(j_{\mathbf{Z}_{(e_G,q)}}((\beta|_G)_g(b))\phi_\mathbf{X}^\mathcal{NO}(s)\big)
&=i_G^{\beta|_G, \omega|_G}(g)i_A^{\beta|_G, \omega|_G}(bc)
\end{align*}
as required. 

To see that the isomorphisms $\{\nu_q: q\in Q^{-1}\}$ implement an isomorphism of product systems between $\mathbf{Y}^\mathcal{NO}$ and the product system associated to $(A \times_{\beta|_G,\omega|_G}G, H, \widetilde\beta|_H, \omega|_H)$, it remains to show that 
\begin{equation}
\label{multiplication is compatible} 
\begin{aligned}
\overline{\omega|_H(t,q)}(\widetilde\beta|_H)_t&(\nu_q(y))\nu_t(w)=\nu_{tq}(yw)\\
&\text{for all $q,t\in Q^{-1}$ and $y\in \mathbf{Y}^\mathcal{NO}_q, w\in \mathbf{Y}^\mathcal{NO}_t$.} 
\end{aligned}
\end{equation}
By linearity and continuity, it suffices to verify \eqref{multiplication is compatible} in the situation where $y:=j_{\mathbf{Z}_{(e_G,q)}}(a)\phi_\mathbf{X}^\mathcal{NO}\big(j_{\mathbf{X}_p}(b)j_{\mathbf{X}_r}(c)^*\big)$ and $w:=j_{\mathbf{Z}_{(e_G,q)}}(d)\phi_\mathbf{X}^\mathcal{NO}(s)$ for arbitrary $q,t\in Q^{-1}$, $p,r\in P^{-1}$, $a,b,c,d\in A$, and $s\in \mathcal{NO}_\mathbf{X}$. The first thing to do is rewrite $yw$ in the form that allows us to apply the formula for $\nu_{tq}$ given by \eqref{definition of nu_q}. Observe that the least upper bound of $(r,e_H)$ and $(e_G,t)$ in the quasi-lattice ordered group $((G\rtimes_\alpha H)^\mathrm{op}, (P\rtimes_\alpha Q)^{-1})$ is 
\[
(\alpha_t(r),t)=(r,e_H)*_\mathrm{op} (e_G,t)=(e_G,t)*_\mathrm{op}(\alpha_t(r),e_H).
\] 
Using the Hilbert $A$-bimodule isomorphisms $\{M^\mathbf{Z}_{(m,n),(u,v)}:(m,n),(u,v)\in (P\rtimes_\alpha Q)^{-1}\}$ we can show that 
\begin{align*}
yw
&=
\frac{\omega((e_G,t),(r,e_H))\omega((\alpha_t(p),e_H),(e_G,t))}{\omega|_H(t,q)\omega((\alpha_t(r),e_H),(e_G,t))\omega((e_G,t),(p,e_H))}\\
&\qquad \times
j_{\mathbf{Z}_{(e_G,tq)}}\big(\beta_{(e_G,t)}(a)\beta_{(\alpha_t(p^{-1}),t)}(b)\big)\\
&\qquad\qquad \times
\phi_\mathbf{X}^\mathcal{NO}\big(j_{\mathbf{X}_{\alpha_t(p)}}(\beta_{(e_G,t)}(c^*))j_{\mathbf{X}_{\alpha_t(r)}}(\beta_{(\alpha_t(r),e_H)}(d^*))^*s\big).
\end{align*}
Applying the formulas for $\nu_{tq}$, $\nu_{t}$, and $\nu_{q}$ given by \eqref{definition of nu_q}, the formula for $\tau$ given by \eqref{formula for tau}, and the definition of $(\widetilde\beta|_H)_t$, we see that
\begin{align*}
\nu_{tq}(yw)
&=\frac{\omega((e_G,t),(r,e_H))\omega((\alpha_t(p),e_H),(e_G,t))}{\omega|_H(t,q)\omega((\alpha_t(r),e_H),(e_G,t))\omega((e_G,t),(p,e_H))}
\\
&\qquad \times
i_A^{\beta|_G, \omega|_G}\big(\beta_{(e_G,t)}(a)\beta_{(\alpha_t(p^{-1}),t)}(b)\big)i_G^{\beta|_G, \omega|_G}(\alpha_t(p))^*
\\
&\qquad\qquad \times
i_A^{\beta|_G, \omega|_G}\big(\beta_{(e_G,t)}(c^*)\beta_{(\alpha_t(r),e_H)}(d)\big)
i_G^{\beta|_G, \omega|_G}(\alpha_t(r))\tau(s)
\\
&=\frac{\omega((e_G,t),(r,e_H))\omega((\alpha_t(p),e_H),(e_G,t))}{\omega|_H(t,q)\omega((\alpha_t(r),e_H),(e_G,t))\omega((e_G,t),(p,e_H))}
\\
&\qquad \times
i_A^{\beta|_G, \omega|_G}\big(\beta_{(e_G,t)}(a)\big)i_G^{\beta|_G, \omega|_G}(\alpha_t(p))^*
\\
&
\qquad\qquad \times
i_A^{\beta|_G, \omega|_G}\big(\beta_{(e_G,t)}(bc^*)\big)i_G^{\beta|_G, \omega|_G}(\alpha_t(r))i_A^{\beta|_G, \omega|_G}(d)\tau(s)
\\
&=
\overline{\omega|_H(t,q)}(\widetilde\beta|_H)_t\big(i_A^{\beta|_G, \omega|_G}(a)i_G^{\beta|_G, \omega|_G}(p)^*i_A^{\beta|_G, \omega|_G}(bc^*)i_G^{\beta|_G, \omega|_G}(r)\big)
\\
&\qquad \times
i_A^{\beta|_G, \omega|_G}(d)\tau(s)
\\
&=
\overline{\omega|_H(t,q)}(\widetilde\beta|_H)_t(\nu_q(y))\nu_t(w)
\end{align*}
as required. Thus, $\mathbf{Y}^\mathcal{NO}$ is isomorphic to the product system associated to the twisted dynamical system $(A \times_{\beta|_G,\omega|_G}G, H, \widetilde\beta|_H, \omega|_H)$.

Hence, combining Theorem~\ref{isomorphism of Cuntz Nica Pimsner algebras} and \cite[Corollary~5.2]{carlsen--larsen--sims--vittadello}, we see that
\[
A\times _{\beta,\omega}(G\rtimes_\alpha H)\cong \mathcal{NO}_\mathbf{Z}\cong \mathcal{NO}_{\mathbf{Y}^\mathcal{NO}}
\cong (A\times_{\beta|_G,\omega|_G} G)\times_{\widetilde\beta|_H, \omega|_H}H.
\]
Chasing through the formulas for the isomorphisms given by Theorem~\ref{isomorphism of Cuntz Nica Pimsner algebras} and \cite[Corollary~5.2]{carlsen--larsen--sims--vittadello}, and the identification of $\mathbf{Y}_q^\mathcal{NO}$ with ${}_{(\widetilde\beta|_H)_q}(A \times_{\beta|_G,\omega|_G}G)_{(A \times_{\beta|_G,\omega|_G}G)}$ given by \eqref{definition of nu_q}, we see that the isomorphism from $(A\times_{\beta|_G,\omega|_G} G)\times_{\widetilde\beta|_H, \omega|_H}H$ to $A\times _{\beta,\omega}(G\rtimes_\alpha H)$ satisfies \eqref{explicit form}.
\end{proof}
\end{exam}

In our next example, we look at a result of Hao and Ng (\cite[Theorem~2.10]{hao--ng}) on crossed products of $C^*$-correspondences. Let $(\alpha,\gamma)$ be an action of a discrete group $G$ on a Hilbert $A$-bimodule $X$, i.e. $\alpha:G\rightarrow \mathrm{Aut}(A)$ and $\gamma:G\rightarrow \mathrm{Aut}(X)$ are homomorphisms such that for $g\in G$, $a\in A$, and $x,y\in X$, we have
\begin{enumerate}[label=\upshape(\roman*)]
\item $\langle \gamma_g(x), \gamma_g(y)\rangle_A=\alpha_g(\langle x,y\rangle_A)$;
\item $\gamma_g(x\cdot a)=\gamma_g(x)\cdot \alpha_g(a)$;
\item $\gamma_g(a\cdot x)=\alpha_g(a)\cdot \gamma_g(x)$.
\end{enumerate}
As in \cite[Section~3.7]{kasparov}, this data can be used to construct a Hilbert $A\rtimes_\alpha G$-bimodule (one can think of this process as extending the coefficient algebra of $X$ from $A$ to the (trivially twisted) crossed product $A\rtimes_\alpha G$). For $\xi,\eta\in C_c(G,X)$ and $u\in C_c(G,A)$ (which we view as a dense $*$-subalgebra of $A\rtimes_\alpha G$ via the map $u\mapsto \sum_{g\in G} i_A^\alpha(u(g))i_G^\alpha(g)$), we define $\xi\cdot u\in  C_c(G,X)$ and $\langle \xi,\eta\rangle_{A\rtimes_\alpha G}\in C_c(G,A)$ by 
\begin{align*}
(\xi\cdot u)(r)&:=\sum_{s\in G}\xi(s)\cdot \alpha_s(u(s^{-1}r))\\
\langle \xi,\eta\rangle_{A\rtimes_\alpha G}(r)&:=\sum_{s\in G}\alpha_{s^{-1}}(\langle \xi(s),\eta(sr)\rangle_A).
\end{align*}
These operations turn $C_c(G,X)$ into a pre-inner product $C_c(G,A)$-module, and we denote its completion by $X\rtimes_\gamma G$ (see \cite[Lemma~2.16]{raeburn--williams} for details on this process). The left action of $A\rtimes_\alpha G$ on $X\rtimes_\gamma G$ is given on functions by
\[
(u\cdot \xi)(r)=\sum_{s\in G}u(s)\cdot \gamma_s(\xi(s^{-1}r)).
\]
Our second example/application uses our decomposition theorem to show how the Cuntz--Pimsner algebras of $X$ and $X\rtimes_\gamma G$ are related. 

\begin{exam}
Let $(\alpha,\gamma)$ be an action of a discrete group $G$ on a Hilbert $A$-bimodule $X$. Then there exists a homomorphism $\beta:G\rightarrow \mathrm{Aut}(\mathcal{O}_X)$ such that 
\begin{equation}
\label{action on crossed product}
\beta_g(j_A(a))=j_A(\alpha_g(a)) \quad \text{and} \quad \beta_g(j_X(x))=j_X(\gamma_g(x))
\end{equation}
for each $g\in G$, $a\in A$, and $x\in X$. Suppose that $G$ is amenable and contains a subsemigroup $P$ that generates $G$, and such that $(G,P)$ is quasi-lattice ordered and $P$ is directed. Furthermore, suppose that the left action of $A$ on $X$ is faithful, nondegenerate, and by compacts. Then 
\[
\mathcal{O}_{X\rtimes_\gamma G}\cong \mathcal{O}_X\rtimes_\beta G. 
\]
Explicitly, this isomorphism is given by
\begin{equation}
\label{explicit form 2}
j_{X\rtimes_\gamma G}^{\otimes m}(x\delta_g)j_{X\rtimes_\gamma G}^{\otimes n}(y\delta_h)^*
\mapsto
i_{\mathcal{O}_X}^\beta\big(j_X^{\otimes m}(x)j_X^{\otimes n}(\gamma_{gh^{-1}}^{\otimes n}(y))^*\big)i_G^\beta(gh^{-1})
\end{equation}
for $m,n\geq 0$, $x\in X^{\otimes m}, y\in X^{\otimes n}$, $g,h\in G$. 
\begin{proof}
Before, we begin the proof, we point out that this result is not new. Indeed, Hao and Ng prove in \cite[Theorem~2.10]{hao--ng} that the Cuntz--Pimsner algebra of $X\rtimes_\gamma G$ and the crossed product associated to the dynamical system $(\mathcal{O}_X, \beta,G)$ coincide provided $G$ is amenable (their result also covers the situation where the topology on $G$ is not discrete). The additional constraints we have imposed on $G$ and $X$ are to enable us to realise the crossed product $\mathcal{O}_X\rtimes_\beta G$ as a Cuntz--Nica--Pimsner algebra and to apply our decomposition theorem. 

Routine calculations show for each $g\in G$, the maps $a\mapsto j_A(\alpha_g(a))$ and $x\mapsto j_X(\gamma_g(x))$ are a Cuntz--Pimsner covariant Toeplitz representation of $X$, and so induce a homomorphism $\beta_g:\mathcal{O}_X\rightarrow \mathcal{O}_X$ satisfying \eqref{action on crossed product}. Since $\alpha$ and $\gamma$ are automorphisms, it follows that the map $g\mapsto \beta_g$ is an action of $G$ on $\mathcal{O}_X$ by automorphisms. 

In order to prove that $\mathcal{O}_X\rtimes_\beta G\cong \mathcal{O}_{X\rtimes_\gamma G}$, we will first use the action $(\alpha,\gamma)$ to construct a compactly aligned product system $\mathbf{Z}$ over $(\Z\times G^\mathrm{op}, \N \times P^{-1})$.  We will then use two applications of our decomposition theorem to show that both $\mathcal{O}_X\rtimes_\beta G$ and $\mathcal{O}_{X\rtimes_\gamma G}$ are isomorphic to $\mathcal{NO}_\mathbf{Z}$. 

For each $(n,p)\in \N\times P^{-1}$, we let $\mathbf{Z}_{(n,p)}$ equal $X^{\otimes n}$ as a Hilbert $A$-module. The left action $\phi_{(n,p)}:A\rightarrow \mathcal{L}_A(\mathbf{Z}_{(n,p)})$ is then given on simple tensors by $\phi_{(n,p)}(a)(x_1\otimes_A \cdots \otimes_A x_n)=(\alpha_p(a)\cdot x_1)\otimes_A \cdots \otimes_A x_n$. Observe that $\mathbf{Z}_{(n,e_G)}$ is equal to $X^{\otimes n}$ as a Hilbert $A$-bimodule, whilst $\mathbf{Z}_{(0,p)}$ is equal to ${}_{\alpha_p}A_A$. For each $(n,p),(m,q)\in \N\times P^{-1}$, there exists a Hilbert $A$-bimodule isomorphism $M_{(n,p),(m,q)}^\mathbf{Z}:\mathbf{Z}_{(n,p)}\otimes_A \mathbf{Z}_{(m,q)}\rightarrow \mathbf{Z}_{(n+m,qp)}$ such that
\[
M_{(n,p),(m,q)}^\mathbf{Z}(x\otimes_A y)=\gamma_q^{\otimes n}(x)\otimes_A y
\]
for each $x\in \mathbf{Z}_{(n,p)}$ and $y\in \mathbf{Z}_{(m,q)}$ (where by convention $\gamma_q^{\otimes 0}:=\alpha_q$). These isomorphisms give $\mathbf{Z}:=\bigsqcup_{(n,p)\in \N\times P^{-1}}\mathbf{Z}_{(n,p)}$ the structure of a product system over $(\Z \times G^{\mathrm{op}},\N\times P^{-1})$. Since $A$ acts compactly on each ${}_{\alpha_p}A_A$ and on $X$, it follows from \cite[Proposition~4.7]{lance} that $A$ acts compactly on each $\mathbf{Z}_{(n,p)}$. Hence, by \cite[Proposition~3.1.24]{fletcherphd}, $\mathbf{Z}$ is compactly aligned. 

We now define $\mathbf{X}:=\bigsqcup_{n\in \N} \mathbf{Z}_{(n,e_G)}$. Then $\mathcal{NO}_\mathbf{X}\cong \mathcal{O}_X$. As in Proposition~\ref{inclusion of NO_X in NO_Z}, we let $\phi_\mathbf{X}^\mathcal{NO}:\mathcal{NO}_\mathbf{X}\rightarrow \mathcal{NO}_\mathbf{Z}$ denote the injective homomorphism induced by the inclusion of $\mathbf{X}$ in $\mathbf{Z}$. For each $p\in P^{-1}$, we define $\mathbf{Y}_p^\mathcal{NO}:=\cspan\{j_\mathbf{Z}(\mathbf{Z}_{(0,p)})\phi_\mathbf{X}^\mathcal{NO}(\mathcal{NO}_\mathbf{X})\}\subseteq \mathcal{NO}_\mathbf{Z}$, which has the structure of a Hilbert $\mathcal{NO}_\mathbf{X}$-bimodule, as in Proposition~\ref{existence of fibres in NO product system in NO_Z}. Equipping $\mathbf{Y}^\mathcal{NO}:=\bigsqcup_{p\in P^{-1}}\mathbf{Y}_p^\mathcal{NO}$ with the multiplication from $\mathcal{NO}_\mathbf{Z}$ gives $\mathbf{Y}^\mathcal{NO}$ the structure of a compactly aligned product system. We claim that $\mathbf{Y}^\mathcal{NO}$ is isomorphic to the product system associated to the dynamical system $(\mathcal{NO}_\mathbf{X}, G, \beta)$. Routine calculations show that the formula
\[
j_{\mathbf{Z}_{(0,p)}}(a)\phi_\mathbf{X}^\mathcal{NO}(t)\mapsto j_{\mathbf{X}_0}(a)t 
\]
extends by linearity and continuity to an inner product preserving map from $\mathbf{Y}_p^\mathcal{NO}$ to ${\mathcal{NO}_\mathbf{X}}_{\mathcal{NO}_\mathbf{X}}$, which we denote by $\tau_p$. The map $\tau_p$ is clearly right $\mathcal{NO}_\mathbf{X}$-linear, and is surjective because $X$ is nondegenerate. To see that the collection of maps $\{\tau_p:p\in P^{-1}\}$ give an isomorphism of product systems, it remains to show that $\beta_q(\tau_p(y))\tau_q(w)=\tau_{qp}(yw)$ for $p,q\in P^{-1}$, $y\in \mathbf{Y}_p^\mathcal{NO}, w\in \mathbf{Y}_q^\mathcal{NO}$. If $a\in\mathbf{Z}_{(0,p)}= {}_{\alpha_p}A_A, b\in \mathbf{Z}_{(0,q)}={}_{\alpha_q}A_A, s\in \mathcal{NO}_\mathbf{X}$, and $x\in \mathbf{X}_m, y\in \mathbf{X}_n$, then 
\begin{align*}
&j_{\mathbf{Z}_{(0,p)}}(a)
\phi_\mathbf{X}^\mathcal{NO}\big(j_{\mathbf{X}_m}(x)j_{\mathbf{X}_n}(y)^*\big)j_{\mathbf{Z}_{(0,q)}}(b)\phi_\mathbf{X}^\mathcal{NO}(s)\\
&\qquad=
j_{\mathbf{Z}_{(0,qp)}}(\alpha_q(a))\phi_\mathbf{X}^\mathcal{NO}\big(j_{\mathbf{X}_m}(\gamma_q^{\otimes m}(x))j_{\mathbf{X}_n}(b^*\cdot \gamma_q^{\otimes n}(y))^*s\big).
\end{align*}
Hence,
\begin{align*}
&\tau_{qp}\big(j_{\mathbf{Z}_{(0,p)}}(a)
\phi_\mathbf{X}^\mathcal{NO}(j_{\mathbf{X}_m}(x)j_{\mathbf{X}_n}(y)^*)
j_{\mathbf{Z}_{(0,p)}}(b)\phi_\mathbf{X}^\mathcal{NO}(s)\big)\\
&\qquad=
j_{\mathbf{X}_0}(\alpha_q(a))j_{\mathbf{X}_m}(\gamma_q^{\otimes m}(x))j_{\mathbf{X}_n}(b^*\cdot \gamma_q^{\otimes n}(y))^*s\\
&\qquad=\beta_q\big(j_{\mathbf{X}_0}(a)j_{\mathbf{X}_m}(x)j_{\mathbf{X}_n}(y)^*\big)j_{\mathbf{X}_0}(b)s\\
&\qquad=\beta_q\big(\tau_p\big(j_{\mathbf{Z}_{(0,p)}}(a)\phi_\mathbf{X}^\mathcal{NO}(j_{\mathbf{X}_m}(x)j_{\mathbf{X}_n}(x)^*)\big)\big)
\tau_q\big(j_{\mathbf{Z}_{(0,p)}}(b)\phi_\mathbf{X}^\mathcal{NO}(s)\big).
\end{align*}
By linearity and continuity it follows that $\beta_q(\tau_p(y))\tau_q(w)=\tau_{qp}(yw)$ for all $y\in \mathbf{Y}_p^\mathcal{NO}$ and $w\in \mathbf{Y}_q^\mathcal{NO}$ as claimed. Applying Theorem~\ref{isomorphism of Cuntz Nica Pimsner algebras}, we see that
\[
\mathcal{NO}_\mathbf{Z}\cong \mathcal{NO}_{\mathbf{Y}^\mathcal{NO}}\cong \mathcal{NO}_\mathbf{X}\rtimes_\beta G
\cong \mathcal{O}_X\rtimes_\beta G. 
\]

Next we define $\mathbf{X}':=\bigsqcup_{p\in P^{-1}} \mathbf{Z}_{(0,p)}$. Then $\mathcal{NO}_{\mathbf{X}'}$ is isomorphic to $A\rtimes_\alpha G$ via the map $j_{\mathbf{X}_p'}(a)\mapsto \alpha_{p^{-1}}(a)\delta_{p^{-1}}\in C_c(G,A)$ (which we denote by $\omega$). We let $\phi_{\mathbf{X}'}^\mathcal{NO}:\mathcal{NO}_{\mathbf{X}'}\rightarrow \mathcal{NO}_\mathbf{Z}$ denote the injective homomorphism induced by the inclusion of $\mathbf{X}'$ in $\mathbf{Z}$. We claim that ${\mathbf{Y}_1'}^\mathcal{NO}:=\cspan\{j_\mathbf{Z}(\mathbf{Z}_{(1,e_G)})\phi_{\mathbf{X}'}^\mathcal{NO}(\mathcal{NO}_{\mathbf{X}'})\}\subseteq \mathcal{NO}_\mathbf{Z}$ and $X\rtimes_\gamma G$ are isomorphic as Hilbert $\mathcal{NO}_{\mathbf{X}'}\cong A\rtimes_\alpha G$-bimodules. 

Calculations show that the formula
\[
j_{\mathbf{Z}_{(1,e_G)}}(x)\phi_{\mathbf{X}'}^\mathcal{NO}\big(j_{\mathbf{X}_p'}(a)j_{\mathbf{X}_q'}(b)^*\big)
\mapsto
\big(x\cdot \alpha_{p^{-1}}(ab^*)\big)\delta_{p^{-1}q} \in C_c(G,X)
\]
extends by linearity and continuity to an inner product preserving map from ${\mathbf{Y}_1'}^\mathcal{NO}$ to $X\rtimes_\gamma G$, which we denote by $\tau_1'$. 

For any $ a\in \mathbf{X}_p'={}_{\alpha_p}A_A, b\in \mathbf{X}_q'={}_{\alpha_q}A_A, c\in \mathbf{X}_s'={}_{\alpha_s}A_A, d\in \mathbf{X}_t'={}_{\alpha_t}A_A$, we have
\begin{align*}
&j_{\mathbf{X}_p'}(a)j_{\mathbf{X}_q'}(b)^*j_{\mathbf{X}_s'}(c)j_{\mathbf{X}_t'}(d)^*\\
&\qquad=j_{\mathbf{X}_{(q\vee s)q^{-1}p}'}(\alpha_{(q\vee s)q^{-1}}(ab^*))j_{\mathbf{X}_{(q\vee s)s^{-1}t}'}(\alpha_{(q\vee s)s^{-1}}(dc^*))^*.
\end{align*}
Hence if $x\in \mathbf{Z}_{(1,e_G)}=X$, then
\begin{align*}
&\tau_1'\big(j_{\mathbf{Z}_{(1,e_G)}}(x)\phi_{\mathbf{X}'}^\mathcal{NO}(j_{\mathbf{X}_p'}(a)j_{\mathbf{X}_q'}(b)^*)
\cdot (j_{\mathbf{X}_s'}(c)j_{\mathbf{X}_t'}(d)^*)\big)\\
&\qquad=\big(x\cdot \alpha_{p^{-1}}(ab^*\alpha_{qs^{-1}}(cd^*)\big)\delta_{p^{-1}qs^{-1}t}\\
&\qquad=\big(x\cdot \alpha_{p^{-1}}(ab^*)\delta_{p^{-1}q}\big)\cdot\big(\alpha_{s^{-1}}(cd^*)\delta_{s^{-1}t}\big)\\
&\qquad=\tau_1'\big(j_{\mathbf{Z}_{(1,e_G)}}(a)\phi_{\mathbf{X}'}^\mathcal{NO}(j_{\mathbf{X}_p'}(a)j_{\mathbf{X}_q'}(b)^*)\big)
\cdot 
\omega\big(j_{\mathbf{X}_s'}(c)j_{\mathbf{X}_t'}(d)^*\big),
\end{align*}
which shows that $\tau_1'$ is right $\mathcal{NO}_{\mathbf{X}'}\cong A\rtimes_\alpha G$-linear. Since $P$ generates $G$ as a group and $X$ is right nondegenerate as an $A$-module, the previous calculation also shows that $\tau_1'$ is surjective. 

Furthermore, we can show that if $x=x'\cdot \langle x', x'\rangle_A$, then
\begin{align*}
&j_{\mathbf{X}_p'}(a)j_{\mathbf{X}_q'}(b)^*\cdot j_{\mathbf{Z}_{(1,e_G)}}(x)\phi_{\mathbf{X}'}^\mathcal{NO}(j_{\mathbf{X}_s'}(c)j_{\mathbf{X}_t'}(d)^*)\\
&\qquad=
j_{\mathbf{Z}_{(1,e_G)}}\big(\gamma_{p^{-1}}(ab^*\cdot \gamma_q(x'))\big)\\
&\qquad\qquad\times
\phi_{\mathbf{X}'}^\mathcal{NO}\big(j_{\mathbf{X}_{(q\vee s)q^{-1}p}'}(\alpha_{q\vee s}(\langle x', x'\rangle_A))j_{\mathbf{X}_{(q\vee s)s^{-1}t}'}(\alpha_{(q\vee s)s^{-1}}(dc^*))^*\big).
\end{align*}
Hence, 
\begin{align*}
&\tau_1'\big(j_{\mathbf{X}_p'}(a)j_{\mathbf{X}_q'}(b)^*\cdot j_{\mathbf{Z}_{(1,e_G)}}(x)\phi_{\mathbf{X}'}^\mathcal{NO}(j_{\mathbf{X}_s'}(c)j_{\mathbf{X}_t'}(d)^*)\big)\\
&\qquad=
\big(\gamma_{p^{-1}}(ab^*\cdot \gamma_q(x'))\cdot \alpha_{p^{-1}q}(\langle x', x'\rangle_A\alpha_{s^{-1}}(cd^*))\big)
\delta_{p^{-1}qs^{-1}t}\\
&\qquad=
\big(\alpha_{p^{-1}}(ab^*)\cdot\gamma_{p^{-1}q}(x\cdot \alpha_{s^{-1}}(cd^*))\big)
\delta_{p^{-1}qs^{-1}t}\\
&\qquad=\big(\alpha_{p^{-1}}(ab^*)\delta_{p^{-1}q}\big)\cdot \big(x\cdot \alpha_{s^{-1}}(cd^*)\delta_{s^{-1}t}\big)\\
&\qquad=\omega\big(j_{\mathbf{X}_p'}(a)j_{\mathbf{X}_q'}(b)^*\big)\cdot \tau_1'\big(j_{\mathbf{Z}_{(1,e_G)}}(x)\phi_{\mathbf{X}'}^\mathcal{NO}(j_{\mathbf{X}_s'}(c)j_{\mathbf{X}_t'}(d)^*)\big),
\end{align*}
which shows that $\tau_1'$ is left $\mathcal{NO}_{\mathbf{X}'}\cong A\rtimes_\alpha G$-linear. Thus, ${\mathbf{Y}_1'}^\mathcal{NO}$ and $X\rtimes_\gamma G$ are isomorphic as Hilbert $\mathcal{NO}_{\mathbf{X}'}\cong A\rtimes_\alpha G$-bimodules, and so applying Theorem~\ref{isomorphism of Cuntz Nica Pimsner algebras}, we have
\[
\mathcal{NO}_\mathbf{Z}\cong \mathcal{NO}_{{\mathbf{Y}'}^\mathcal{NO}}
\cong \mathcal{O}_{{\mathbf{Y}_1'}^\mathcal{NO}}\cong \mathcal{O}_{X\rtimes_\gamma G}. 
\]
Hence,
\[
\mathcal{O}_{X\rtimes_\gamma G}\cong \mathcal{NO}_\mathbf{Z} \cong \mathcal{O}_X\rtimes_\beta G.
\]
Chasing through the isomorphisms given by Theorem~\ref{isomorphism of Cuntz Nica Pimsner algebras} as well as the identifications of ${\mathbf{Y}_1'}^\mathcal{NO}$ and ${\mathbf{Y}_1}^\mathcal{NO}$ with $X\rtimes_\gamma G$ and the product system associated to $(\mathcal{NO}_\mathbf{X}, G, \beta)$ respectively, shows that the isomorphism from $\mathcal{O}_{X\rtimes_\gamma G}$ to $\mathcal{O}_X\rtimes_\beta G$ satisfies \eqref{explicit form 2}.
\end{proof}
\end{exam}

In our final example we look at twisted $C^*$-algebras of topological higher-rank graphs \cite{armstrong--brownlowe}. Recall that for $k\in \N\setminus \{0\}$, a topological $k$-graph is a small category $\Lambda=(\mathrm{Obj}(\Lambda), \mathrm{Mor}(\Lambda),r,s,\circ)$ and a functor $d:\Lambda\rightarrow \N^k$ such that
\begin{enumerate}[label=\upshape(\roman*)]
\item
the sets of objects $\mathrm{Obj}(\Lambda)$ and morphisms $\mathrm{Mor}(\Lambda)$ in $\Lambda$ are second-countable, locally compact Hausdorff topological spaces;
\item
the range and source maps $r$ and $s$ are continuous, and $s$ is a local homeomorphism;
\item
if $\Lambda \times_c \Lambda:=\{(\lambda,\mu)\in \Lambda\times \Lambda: s(\lambda)=r(\lambda)\}$ is equipped with the subspace topology inherited from the product topology on $\Lambda\times \Lambda$, then composition $\circ: \Lambda \times_c \Lambda\rightarrow \Lambda$ is a continuous open map;
\item
the degree functor $d$ is continuous (where $\N^k$ has the discrete topology);
\item
$\Lambda$ satisfies the unique factorisation property: if $\lambda\in \Lambda$ and $d(\lambda)=m+n$ for some $m,n\in \N^k$, then there exist unique $\mu,\nu\in \Lambda$ with $d(\mu)=m$ and $d(\nu)=n$ such that $\lambda=\mu\nu$.
\end{enumerate}

For each $n\in \N^k$, $(\Lambda^0:=d^{-1}(0), \Lambda^n:=d^{-1}(n), r|_{\Lambda^n}, s|_{\Lambda^n})$ is a topological graph, and we let $\mathbf{Z}_n$ denote the associated topological graph correspondence (see \cite{katsura3} for the details of this construction). Lemma~3.18 of  \cite{armstrong--brownlowe} shows that if $\Lambda$ is source-free (in the sense that $\{\lambda\in \Lambda^{e_i}: r(\lambda)=v\}$ is nonempty for each $v\in\Lambda^0$ and each $i\in\{1,\ldots, k\}$) and proper (in the sense that the restriction $r|_{\Lambda^n}$ is a proper map for each $n\in \N^k$), then the left action of $C_0(\Lambda^0)$ on each $\mathbf{Z}_n$ is faithful and by compacts. 

Now fix a continuous $\T$-valued $2$-cocycle on $\Lambda$, i.e. a continuous map $c:\Lambda\times_c \Lambda \rightarrow \T$ satisfying
\begin{enumerate}[label=\upshape(\roman*)]
\item
$c(\lambda,\mu)c(\lambda\mu,\nu)=c(\lambda, \mu\nu)c(\mu,\nu)$ for all paths $\lambda,\mu,\nu\in \Lambda$ with $s(\lambda)=r(\mu)$ and $s(\mu)=r(\nu)$;
\item
$c(\lambda,s(\lambda))=c(r(\lambda),\lambda)=1$ for all $\lambda\in \Lambda$. 
\end{enumerate} 
We use this cocycle to define multiplication on $\mathbf{Z}:=\bigsqcup_{n\in \N^k}\mathbf{Z}_n$. For $f\in \mathbf{Z}_m, g\in \mathbf{Z}_n$ we use \cite[Proposition~4.1]{armstrong--brownlowe} to define $fg\in \mathbf{Z}_{m+n}$ by
\[
(fg)(\lambda\mu):=c(\lambda,\mu)f(\lambda)g(\mu)
\]
for $\lambda\in \Lambda^m, \mu\in \Lambda^n$ with $s(\lambda)=r(\mu)$. With this multiplication, $\mathbf{Z}$ has the structure of a compactly aligned product system over $(\Z^k,\N^k)$ with coefficient algebra $C_0(\Lambda^0)$. Following \cite[Definition~4.4]{armstrong--brownlowe}, the twisted Toeplitz algebra $\mathcal{T}C^*(\Lambda,c)$ and the twisted Cuntz--Krieger algebra $C^*(\Lambda,c)$ of $\Lambda$ are defined to be $\mathcal{NT}_\mathbf{Z}$ and $\mathcal{NO}_\mathbf{Z}$ respectively. 

\begin{exam}
Let $\Lambda$ be a source-free proper topological $k$-graph, and $c$ a continuous $\T$-valued $2$-cocycle on $\Lambda$. For each $i\in \{1,\ldots,k\}$, let $\Lambda^i:=\{\lambda\in \Lambda: d(\lambda)_i=0\}$ (which we view as a topological $(k-1)$-graph) and define $c^i:=c|_{\Lambda^i\times_c \Lambda^i}$. Then $\mathcal{T}C^*(\Lambda,c)$ and $C^*(\Lambda,c)$ can be realised as the Toeplitz and Cuntz--Pimsner algebras of Hilbert bimodules with coefficient algebras $\mathcal{T}C^*(\Lambda^i,c^i)$ and $C^*(\Lambda^i,c^i)$ respectively. Consequently, $\mathcal{T}C^*(\Lambda,c)$ and $C^*(\Lambda,c)$ can be viewed as iterated Toeplitz and iterated Cuntz--Pimsner algebras over $\mathcal{T}C^*(\Lambda^0)\cong C^*(\Lambda^0)\cong C_0(\Lambda^0)$. 

Furthermore, $\mathcal{T}C^*(\Lambda,c)$ and $C_0(\Lambda^0)$ are $KK$-equivalent and have the same $K$-theory. Both $\mathcal{T}C^*(\Lambda,c)$ and $C^*(\Lambda,c)$ are nuclear, exact, and satisfy the UCT. 

\begin{proof}
It is routine to check that the product system over $\N^{(k-1)}$ associated to $\Lambda^i$ and $c^i$, is just the restriction of the product system $\mathbf{Z}$ associated to $\Lambda$ and $c$ to the fibres $\{n\in \N^k:n_i=0\}\cong \N^{(k-1)}$, which we denote by $\mathbf{X}$. We let $\mathbf{Y}^{\mathcal{NT}}$ and $\mathbf{Y}^{\mathcal{NO}}$ be the product systems over $\N$, with coefficient algebras $\mathcal{NT}_\mathbf{X}=\mathcal{T}C^*(\Lambda^i,c^i)$ and $\mathcal{NO}_\mathbf{X}=C^*(\Lambda^i,c^i)$ respectively, as defined in Propositions~\ref{Y is a product system} and \ref{existence of NO product system in NO_Z}. By Theorems~\ref{isomorphisms of Nica--Toeplitz algebras} and \ref{isomorphism of Cuntz Nica Pimsner algebras},
\begin{align*}
\mathcal{T}C^*(\Lambda,c)&=\mathcal{NT}_\mathbf{Z}\cong \mathcal{NT}_{\mathbf{Y}^\mathcal{NT}}\cong \mathcal{T}_{\mathbf{Y}_1^\mathcal{NT}}\\
C^*(\Lambda,c)&=\mathcal{NO}_\mathbf{Z}\cong \mathcal{NO}_{\mathbf{Y}^\mathcal{NO}}\cong \mathcal{O}_{\mathbf{Y}_1^\mathcal{NO}}.
\end{align*}
Repeatedly applying this procedure shows that $\mathcal{T}C^*(\Lambda,c)$ and $C^*(\Lambda,c)$ can be realised as iterated Toeplitz and iterated Cuntz--Pimsner algebras over $C_0(\Lambda^0)$ respectively. 

Since each $\Lambda^n$ is locally compact, Hausdorff, and second-countable, each $\mathbf{Z}_{e_i}$ is countably generated as a $C_0(\Lambda^0)$-module \cite[Lemma~6.2]{katsura3} and $C_0(\Lambda^0)$ is separable \cite[Theorem~2.4]{chou}. Corollary~\ref{extending literature results for NT algebras}(ii) tells us that $\mathcal{T}C^*(\Lambda,c)$ and $C_0(\Lambda^0)$ are $KK$-equivalent (extending \cite[Lemma~6.5]{katsura3} from untwisted topological graphs to twisted higher-rank topological graphs). Consequently, $\mathcal{T}C^*(\Lambda,c)$ and $C_0(\Lambda^0)$ have the same $K$-theory and $\mathcal{T}C^*(\Lambda,c)$ satisfies the UCT \cite[Proposition~2.1(1)]{dadarlat}. Since $C_0(\Lambda^0)$ is commutative, it is nuclear (and hence exact), and so Corollary~\ref{extending literature results for NT algebras}(iv) and Corollary~\ref{extending literature results for CNP algebras}(ii) tell us that $\mathcal{T}C^*(\Lambda,c)$ and $C^*(\Lambda,c)$ are nuclear (and hence exact). Also, since $C_0(\Lambda^0)$ satisfies the UCT (it is commutative and separable), Corollary~\ref{extending literature results for CNP algebras}(iii) tells us that $C^*(\Lambda,c)$ satisfies the UCT (generalising \cite[Proposition~6.6]{katsura3}).
\end{proof}
\end{exam}

\begin{rem}
When $\Lambda$ is countable and discrete (i.e. $\Lambda$ is a higher-rank graph), \cite[Theorem~4.16]{k-graph iterated} shows that the hypothesis of $\Lambda$ being proper (i.e. $\Lambda$ is row-finite) can be dropped and the hypothesis that $\Lambda$ is source-free can be relaxed to asking that $\Lambda$ is locally-convex. 
\end{rem}

\section{Relative Cuntz--Nica--Pimsner algebras}
\label{mixed cases}

We now consider the Cuntz--Nica--Pimsner algebra of the product system $\mathbf{Y}^{\mathcal{NT}}$ defined in Section~\ref{factorising NT alg}. When the action $\alpha$ is trivial, we will show that there exists a product system $\mathbf{V}$ over $(H,Q)$, and a product system $\mathbf{W}^{\mathcal{NO}}$ over $(G,P)$ with coefficient algebra $\mathcal{NO}_\mathbf{V}$, such that $\mathcal{NO}_{\mathbf{Y}^{\mathcal{NT}}}$ and $\mathcal{NT}_{\mathbf{W}^{\mathcal{NO}}}$ are isomorphic. In a sense, the Cuntz--Pimsner covariance in $\mathcal{NO}_{\mathbf{Y}^{\mathcal{NT}}}$ can be moved into the coefficient algebra of $\mathbf{W}^{\mathcal{NO}}$. 

Since we are, in effect, only asking for Cuntz--Pimsner covariance in some of the fibres of $\mathbf{Z}$, we like to think of $\mathcal{NO}_{\mathbf{Y}^{\mathcal{NT}}}$ and $\mathcal{NT}_{\mathbf{W}^\mathcal{NO}}$ as relative Cuntz--Nica--Pimsner algebras. This nomenclature is motivated by the relative Cuntz--Krieger algebras of higher-rank graphs introduced by Sims \cite{sims}, and the relative Cuntz--Pimsner algebras introduced by Muhly and Solel \cite{muhly--solel} (which were studied further by Fowler, Muhly, and Raeburn \cite{fowler--muhly--raeburn}).

\begin{hypoth}
We assume that the action $\alpha:H \rightarrow \mathrm{Aut}(G)$ is trivial. Hence, $\mathbf{Z}$ is a compactly aligned product system over $(G\times H, P \times Q)\cong (H\times G, Q \times P)$. As such, we can swap the roles of $G$ and $H$, and $P$ and $Q$ in our results from Sections~\ref{factorising NT alg} and~\ref{factorising CNP alg}. To ensure that the $C^*$-algebras and product systems we want to work with actually exist, we assume that both $G$ and $H$ are amenable, $A$ acts faithfully on each $\mathbf{Z}_{\left(e_G,q\right)}$, and each $\widetilde{\phi}_{(p,q)}$ is injective. 
\end{hypoth}

We now summarise the setup, as well as fixing some notation. 

\begin{enumerate}
\item[(1)] For each $q\in Q$, we let $\mathbf{V}_q:=\mathbf{Z}_{\left(e_G,q\right)}$. Then $\mathbf{V}:=\bigsqcup_{q\in Q}\mathbf{V}_q$ is a compactly aligned product system over $(H,Q)$ with coefficient algebra $A$. 
\item[(2)] By Proposition~\ref{inclusion of NO_X in NO_Z}, since $A$ acts faithfully on each fibre of $\mathbf{V}$ and each homomorphism $\widetilde{\phi}_{(p,q)}:A\rightarrow \mathcal{L}_A\big(\widetilde{\mathbf{Z}}_{(p,q)}\big)$ is injective, there exists a homomorphism $\phi_\mathbf{V}^\mathcal{NO}:\mathcal{NO}_\mathbf{V}\rightarrow \mathcal{NO}_\mathbf{Z}$ such that $\phi_\mathbf{V}^\mathcal{NO}\circ j_\mathbf{V}=j_\mathbf{Z}$. Furthermore, since $H$ is an amenable group, Proposition~\ref{injectivity of inclusion of NO_X in NO_Z} says that $\phi_\mathbf{V}^\mathcal{NO}$ is injective.
\item[(3)] Since $H$ is amenable, $A$ acts faithfully on each fibre of $\mathbf{V}$, and each $\widetilde{\phi}_{(p,q)}$ is injective, Proposition~\ref{Y is a product system} and Proposition~\ref{Y is compactly aligned} give the existence of a compactly aligned product system $\mathbf{W}^\mathcal{NO}$ over $(G,P)$ with coefficient algebra $\mathcal{NO}_\mathbf{V}$, with fibres given by
\[
\mathbf{W}^\mathcal{NO}_p:=\cspan\big\{j_{\mathbf{Z}_{\left(p,e_H\right)}}(x)\phi_\mathbf{V}^\mathcal{NO}(b):x\in \mathbf{Z}_{\left(p,e_H\right)}, b\in \mathcal{NO}_\mathbf{V}\big\}
\]
for each $p\in P\setminus \{e_G\}$. For simplicity's sake, we will frequently identify the coefficient algebra of $\mathbf{W}^\mathcal{NO}$ with $\phi_\mathbf{V}^\mathcal{NO}\left(\mathcal{NO}_\mathbf{V}\right)\subseteq \mathcal{NO}_\mathbf{Z}$.  
\end{enumerate}

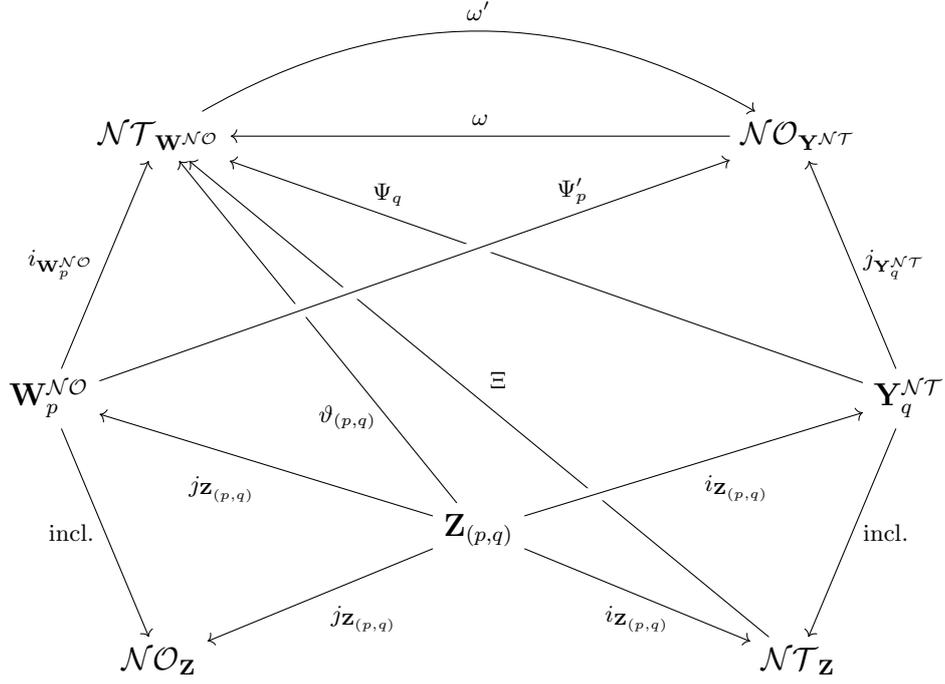
\begin{figure}[h]
\centering
\begin{tikzpicture}[scale=1.45]
\node (C) at (1,0) {$\mathcal{NO}_\mathbf{Z}$};
\node (D) at (3.9,1.2) {$\mathbf{Z}_{(p,q)}$};
\node (F) at (6.8,0) {$\mathcal{NT}_\mathbf{Z}$};
\node (H) at (0,2.4) {$\mathbf{W}_p^\mathcal{NO}$};
\node (I) at (1,4.8) {$\mathcal{NT}_{\mathbf{W}^\mathcal{NO}}$};
\node (J) at (7.8,2.4) {$\mathbf{Y}_q^\mathcal{NT}$};
\node (K) at (6.8,4.8) {$\mathcal{NO}_{\mathbf{Y}^\mathcal{NT}}$};
\path[->,font=\scriptsize]
(D) edge node[below right]{$j_{\mathbf{Z}_{(p,q)}}$} (C)
(D) edge node[below]{$i_{\mathbf{Z}_{(p,q)}}$} (F)
(H) edge node[left]{$\mathrm{incl.}$} (C)
(H) edge node[left]{$i_{\mathbf{W}_p^\mathcal{NO}}$} (I)
(J) edge node[right]{$\mathrm{incl.}$} (F)
(J) edge node[right]{$j_{\mathbf{Y}_q^\mathcal{NT}}$} (K)
(D) edge node[below left]{$j_{\mathbf{Z}_{(p,q)}}$} (H)
(F) edge node[above right]{$\Xi$} (I)
(D) edge [line width=6pt,draw=white](J)
(D) edge node[below right]{$i_{\mathbf{Z}_{(p,q)}}$} (J)
(D) edge node[left, near start]{$\vartheta_{(p,q)}$} (I)
(J) edge node[above, near end]{$\Psi_q$} (I)
(H) edge [line width=6pt,draw=white](K)
(H) edge node[above, near end]{$\Psi_p'$} (K)
(K) edge node[above]{$\omega$} (I)
(I) edge [bend left] node[above]{$\omega'$} (K);
\end{tikzpicture}
\caption{The homomorphisms $\omega$,  $\omega'$}
\label{relative isomorphisms}
\end{figure}

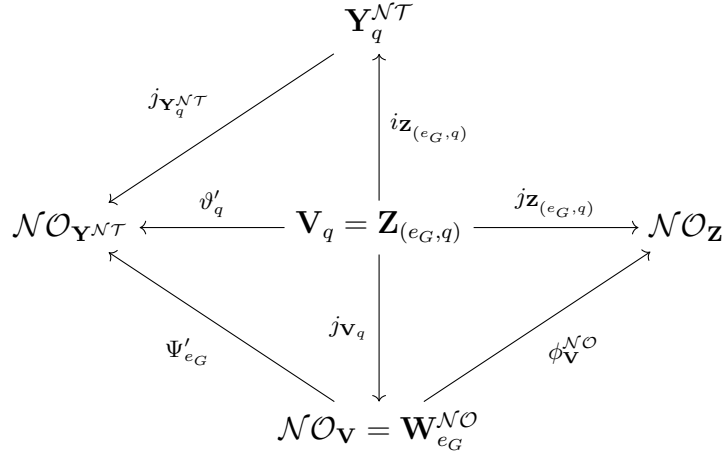
\begin{figure}[h]
\centering
\begin{tikzpicture}[scale=1.35]
\node (A) at (3,0) {$\mathcal{NO}_\mathbf{V}=\mathbf{W}_{e_G}^\mathcal{NO}$};
\node (B) at (0,2) {$\mathcal{NO}_{\mathbf{Y}^\mathcal{NT}}$};
\node (C) at (3,2) {$\mathbf{V}_q=\mathbf{Z}_{(e_G,q)}$};
\node (D) at (6,2) {$\mathcal{NO}_\mathbf{Z}$};
\node (E) at  (3,4) {$\mathbf{Y}_q^\mathcal{NT}$};
\path[->,font=\scriptsize]
(C) edge node[left]{$j_{\mathbf{V}_q}$} (A)
(A)edge node[below right]{$\phi_\mathbf{V}^{\mathcal{NO}}$} (D)
(C) edge node[above]{$j_{\mathbf{Z}_{\left(e_G,q\right)}}$} (D)
(C) edge node[right]{$i_{\mathbf{Z}_{\left(e_G,q\right)}}$} (E)
(E) edge node[above left]{$j_{\mathbf{Y}_q^\mathcal{NT}}$} (B)
(C) edge node[above]{$\vartheta'_q$} (B)
(A)edge node[below left]{$\Psi'_{e_G}$} (B)
;
\end{tikzpicture}
\caption{The homomorphisms $\phi_\mathbf{V}^{\mathcal{NO}}$,  $\Psi'_{e_G}$}
\label{relative homomorphisms}
\end{figure}

We exhibit homomorphisms $\omega:\mathcal{NO}_{\mathbf{Y}^\mathcal{NT}}\rightarrow \mathcal{NT}_{\mathbf{W}^\mathcal{NO}}$ and $\omega':\mathcal{NT}_{\mathbf{W}^\mathcal{NO}}\rightarrow \mathcal{NO}_{\mathbf{Y}^\mathcal{NT}}$, and show that they are inverses of each other. We now list the results that we will prove in Section~\ref{mixed cases}, and summarise the various spaces and maps that we will be working with in a pair of commutative diagrams (see Figures~\ref{relative isomorphisms} and~\ref{relative homomorphisms}). 

\begin{enumerate}
\item[(4)]
In Proposition~\ref{existence of xi hom and defn of Psi} we use the universal Cuntz--Nica--Pimsner covariant representation of $\mathbf{Z}$ along with the universal Nica covariant representation of $\mathbf{W}^\mathcal{NO}$ to define a Nica covariant representation $\vartheta$ of $\mathbf{Z}$ in $\mathcal{NT}_{\mathbf{W}^\mathcal{NO}}$. This gives the existence of a homomorphism $\Xi:\mathcal{NT}_\mathbf{Z}\rightarrow \mathcal{NT}_{\mathbf{W}^\mathcal{NO}}$ such that 
\[
\Xi\circ i_{\mathbf{Z}_{(p,q)}}=\vartheta_{(p,q)}=i_{\mathbf{W}^\mathcal{NO}_p}\circ j_{\mathbf{Z}_{(p,q)}} \text{ for $(p,q)\in P\times Q$.}
\]
We then show that restricting the homomorphism $\Xi$ to $\mathbf{Y}^\mathcal{NT}\subseteq \mathcal{NT}_\mathbf{Z}$ gives a Nica covariant representation of $\mathbf{Y}^\mathcal{NT}$, which we denote by $\Psi$. The idea is that $\Xi$ plays the same role as the inclusion map in Propositions~\ref{varphi' is Nica covariant} and~\ref{varphi' is CNP covariant}.
\item[(5)]
In Proposition~\ref{sufficient conditions for Psi' to be CNP cov} we find sufficient conditions for $\Psi$ to be Cuntz--Pimsner covariant, and use the universal property of $\mathcal{NO}_{\mathbf{Y}^{\mathcal{NT}}}$ to induce a homomorphism $\omega$ such that $\omega\circ j_{\mathbf{Y}^\mathcal{NT}}=\Psi$. 
\item[(6)]
In Proposition~\ref{defn of vartheta} we use the universal Nica covariant representation of $\mathbf{Z}$ and the universal Cuntz--Nica--Pimsner covariant representation of $\mathbf{Y}^\mathcal{NT}$ to define a Cuntz--Nica--Pimsner covariant representation $\vartheta'$ of $\mathbf{V}$ in $\mathcal{NO}_{\mathbf{Y}^{\mathcal{NT}}}$. The universal property of $\mathcal{NO}_\mathbf{V}$ then gives us a homomorphism $\Psi_{e_G}':\mathcal{NO}_\mathbf{V}\rightarrow \mathcal{NO}_{\mathbf{Y}^{\mathcal{NT}}}$ such that 
\[
\Psi_{e_G}'\circ  j_{\mathbf{V}_q}=j_{\mathbf{Y}^\mathcal{NT}_q}\circ i_{\mathbf{Z}_{\left(e_G,q\right)}} \text{ for $q\in Q$.}
\]
\item[(7)]
In Proposition~\ref{definition of Psi'} we use $\Psi_{e_G}'$ to construct a linear map 
$\Psi'_p:\mathbf{W}^\mathcal{NO}_p\rightarrow\mathcal{NO}_{\mathbf{Y}^\mathcal{NT}}$ such that 
\[
\Psi'_p\circ j_{\mathbf{Z}_{(p,q)}}=j_{\mathbf{Y}_q^\mathcal{NT}}\circ i_{\mathbf{Z}_{(p,q)}} \text{ for $(p,q)\in P\times Q$.}
\]
\item[(8)]
In Proposition~\ref{Psi' is repn of W} and \ref{Psi' is a Nica cov repn of W} we show that the collection of maps $\{\Psi'_p:p\in P\}$ gives a Nica covariant representation of the product system $\mathbf{W}^\mathcal{NO}$. The universal property of $\mathcal{NT}_{\mathbf{W}^\mathcal{NO}}$ then gives a homomorphism $\omega':\mathcal{NT}_{\mathbf{W}^\mathcal{NO}}\rightarrow \mathcal{NO}_{\mathbf{Y}^\mathcal{NT}}$ such that $\omega'\circ i_{\mathbf{W}^\mathcal{NO}}=\Psi'$.
\item[(9)]
In Theorem~\ref{isom of mixed case algebras} we prove that $\omega$ and $\omega'$ are mutually inverse isomorphisms. 
\end{enumerate}

In summary, we will show that for every $(p,q)\in P\times Q$, the maps in Figures~\ref{relative isomorphisms} and~\ref{relative homomorphisms} exist and make the diagrams commutative.

We begin by exhibiting a Cuntz--Nica--Pimsner covariant representation $\Psi$ of $\mathbf{Y}^\mathcal{NT}$ in $\mathcal{NT}_{\mathbf{W}^\mathcal{NO}}$. The idea is to produce a homomorphism from $\mathcal{NT}_\mathbf{Z}$ to $\mathcal{NT}_{\mathbf{W}^\mathcal{NO}}$, and then restrict this map to $\mathbf{Y}^\mathcal{NT}$. 

\begin{prop}
\label{existence of xi hom and defn of Psi}
Define $\vartheta:\mathbf{Z}\rightarrow \mathcal{NT}_{\mathbf{W}^\mathcal{NO}}$ by 
\[
\vartheta_{(p,q)}:=i_{\mathbf{W}^\mathcal{NO}_p}\circ j_{\mathbf{Z}_{(p,q)}}.
\]
Then $\vartheta$ is a Nica covariant representation of $\mathbf{Z}$, and so there exists a homomorphism $\Xi:\mathcal{NT}_\mathbf{Z}\rightarrow \mathcal{NT}_{\mathbf{W}^\mathcal{NO}}$ such that 
\[
\Xi\circ i_{\mathbf{Z}_{(p,q)}}=i_{\mathbf{W}^\mathcal{NO}_p}\circ j_{\mathbf{Z}_{(p,q)}}
\]
for each $(p,q)\in P\times Q$. If $\Psi: \mathbf{Y}^\mathcal{NT}\rightarrow  \mathcal{NT}_{\mathbf{W}^\mathcal{NO}}$ is defined by 
\[
\Psi_q:=\Xi|_{\mathbf{Y}^\mathcal{NT}_q},
\]
then $\Psi$ is a Nica covariant representation of $\mathbf{Y}^\mathcal{NT}$. 
\begin{proof}
Both $i_{\mathbf{W}^\mathcal{NO}_p}$ and $j_{\mathbf{Z}_{(p,q)}}$ are Nica covariant representations, and so the same reasoning as in Propositions~\ref{varphi is a representation} and \ref{varphi is Nica covariant} shows that $\vartheta$ is a Nica covariant representation. The universal property of $\mathcal{NT}_\mathbf{Z}$ then induces the homomorphism $\Xi$.

Since $\Xi$ is a homomorphism, it is elementary to check that $\Psi:=\Xi|_{\mathbf{Y}^\mathcal{NT}}$ is a representation. We claim that $\Psi^{(q)}(M_b)=\Xi(b)$ for each $b\in \mathcal{NT}_\mathbf{Z}^q$. If $z\in \mathbf{Z}_{(m,q)}$ and $w\in \mathbf{Z}_{(n,q)}$, then $i_{\mathbf{Z}_{(m,q)}}(z)i_{\mathbf{Z}_{(n,q)}}(w)^* \in \mathcal{NT}_\mathbf{Z}^q$ and 
\[
M_{i_{\mathbf{Z}_{(m,q)}}(z)i_{\mathbf{Z}_{(n,q)}}(w)^*}=\Theta_{i_{\mathbf{Z}_{(m,q)}}(z),i_{\mathbf{Z}_{(n,q)}}(w)}\in \mathcal{K}_{\mathcal{NT}_\mathbf{X}}\left(\mathbf{Y}^{\mathcal{NT}}_q\right).
\]
Hence,
\begin{align*}
\Psi^{(q)}\big(M_{i_{\mathbf{Z}_{(m,q)}}(z)i_{\mathbf{Z}_{(n,q)}}(w)^*}\big)
&=\Psi_q(i_{\mathbf{Z}_{(m,q)}}(z))\Psi_q(i_{\mathbf{Z}_{(n,q)}}(w))^*\\
&=\Xi\big(i_{\mathbf{Z}_{(m,q)}}(z)i_{\mathbf{Z}_{(n,q)}}(w)^* \big).
\end{align*}
By linearity, continuity, and an application of Lemma~\ref{useful property for multiplication operators}, we get our claim. 

We now prove that $\Psi$ is Nica covariant. Fix $b\in  \mathcal{NT}_\mathbf{Z}^q$ and $c\in  \mathcal{NT}_\mathbf{Z}^t$. By Lemma~\ref{Q indexed subalgebras of NO_Z}, $bc\in \mathcal{NT}_\mathbf{Z}^{(q\vee t)}$ if $q\vee t<\infty$, and is zero otherwise. Hence, $\Psi^{(q)}(M_b)\Psi^{(t)}(M_c)=\Xi(b)\Xi(c)=\Xi(bc)$ is zero if $q\vee t=\infty$, and equal to $\Psi^{(q\vee t)}(M_{bc})=\Psi^{(q\vee t)}(\iota_q^{q\vee t}(M_b) \iota_t^{q\vee t}(M_c))$ if $q\vee t<\infty$. By Lemma~\ref{compact operators in Y}, it follows that $\Psi$ is Nica covariant.
\end{proof}
\end{prop}

\begin{prop}
\label{sufficient conditions for Psi' to be CNP cov}
Suppose that $A$ acts compactly on each $\mathbf{Z}_{\left(e_G,q\right)}$ and $Q$ is directed. Then $\Psi$ is a Cuntz--Pimsner covariant representation of $\mathbf{Y}^\mathcal{NT}$, and so there exists a homomorphism $\omega:\mathcal{NO}_{\mathbf{Y}^\mathcal{NT}}\rightarrow\mathcal{NT}_{\mathbf{W}^\mathcal{NO}}$ such that 
\[
\omega \circ j_{\mathbf{Y}^\mathcal{NT}_q}=\Psi_q=\Xi|_{\mathbf{Y}^\mathcal{NT}_q}
\]
for each $q\in Q$.
\begin{proof}
Since $G$ is amenable and $A$ acts faithfully on each $\mathbf{Z}_{\left(e_G,q\right)}$, Proposition~\ref{faithfulness of induced action - Nica Toeplitz case} tells us that $\mathcal{NT}_X$ acts faithfully on each fibre of $\mathbf{Y}^\mathcal{NT}$. Additionally, since $A$ acts compactly on each $\mathbf{Z}_{\left(e_G,q\right)}$, the same reasoning as in the proof of Lemma~\ref{induced actions are by compacts} shows that $\mathcal{NT}_\mathbf{X}$ acts compactly on each fibre of $\mathbf{Y}^\mathcal{NT}$. Hence, to see that $\Psi$ is Cuntz--Pimsner covariant, it suffices by \cite[Proposition~5.1]{sims--yeend} to check that $\Psi^{(q)}\circ \Phi_q^{\mathcal{NT}}=\Psi_{e_H}$ for each $q\in Q$. As $\mathbf{Y}^\mathcal{NT}_{e_H}=\mathcal{NT}_X\cong \phi^\mathcal{NT}_\mathbf{X}(\mathcal{NT}_\mathbf{X})$ is generated by the image of $\phi^\mathcal{NT}_\mathbf{X}\circ i_\mathbf{X}=i_\mathbf{Z}|_\mathbf{X}$ it suffices to show that 
\begin{equation}
\label{sufficient for covariance 2}
\begin{aligned}
\big(\Psi^{(q)}\circ \Phi_q^{\mathcal{NT}}\big)(i_\mathbf{Z}(x))=\Psi_{e_H}\left(i_\mathbf{Z}(x)\right) \quad \text{for each $x\in \mathbf{X}$.} 
\end{aligned}
\end{equation}
To this end, fix $x\in \mathbf{X}_p$. Choosing $x'\in \mathbf{X}_p$ so that $x=x'\cdot \langle x',x'\rangle_A^p$, we have
\begin{align*}
\Psi_{e_H}\big(i_{\mathbf{Z}_{\left(p,e_H\right)}}(x)\big)
&=\Xi\big(i_{\mathbf{Z}_{\left(p,e_H\right)}}(x)\big)\\
&=i_{\mathbf{W}^\mathcal{NO}_p}\big( j_{\mathbf{Z}_{\left(p,e_H\right)}}(x)\big)\\
&=i_{\mathbf{W}^\mathcal{NO}_p}\big(j_{\mathbf{Z}_{\left(p,e_H\right)}}(x')\big)i_{\mathbf{W}^\mathcal{NO}_{e_G}}\big(j_{\mathbf{Z}_{\left(e_G,e_H\right)}}(\langle x',x'\rangle_A^p)\big).
\end{align*}
By the first part of Lemma~\ref{induced actions are by compacts}, we see that
\begin{align*}
\Phi_q^{\mathcal{NT}}\big(i_{\mathbf{Z}_{\left(p,e_H\right)}}(x)\big)
=M_{i_{\mathbf{Z}_{\left(p,e_H\right)}}(x')i_\mathbf{Z}^{\left(\left(e_G,q\right)\right)}\big(\phi_{\left(e_G,q\right)}(\langle x',x'\rangle_A^p)\big)}\in \mathcal{K}_{\mathcal{NT}_\mathbf{X}}\left(\mathbf{Y}_q^\mathcal{NT}\right)
\end{align*}
and $i_{\mathbf{Z}_{\left(p,e_H\right)}}(x')i_\mathbf{Z}^{\left(\left(e_G,q\right)\right)}\big(\phi_{\left(e_G,q\right)}(\langle x',x'\rangle_A^p)\big)\in \mathcal{NO}_\mathbf{Z}^q$. Thus, 
\begin{align*}
&\Psi^{(q)}\Big(\Phi^\mathcal{NT}_q\big(i_{\mathbf{Z}_{\left(p,e_H\right)}}(x)\big)\Big)\\
&\qquad=\Xi\Big(i_{\mathbf{Z}_{\left(p,e_H\right)}}(x')i_\mathbf{Z}^{\left(\left(e_G,q\right)\right)}\big(\phi_{\left(e_G,q\right)}(\langle x',x'\rangle_A^p)\big)\Big)\\
&\qquad=i_{\mathbf{W}^\mathcal{NO}_p}\big(j_{\mathbf{Z}_{(p,q)}}(x')\big)
i_{\mathbf{W}^\mathcal{NO}_{e_G}}\Big(j_\mathbf{Z}^{\left(\left(e_G,q\right)\right)}\big(\phi_{\left(e_G,q\right)}(\langle x',x'\rangle_A^p)\big)\Big).
\end{align*}
Thus, for \eqref{sufficient for covariance 2} to hold, we need
\[
j_{\mathbf{Z}_{\left(e_G,e_H\right)}}(\langle x',x'\rangle_A^p)=j_\mathbf{Z}^{\left(\left(e_G,q\right)\right)}\big(\phi_{\left(e_G,q\right)}(\langle x',x'\rangle_A^p)\big)
,
\]
which follows from Lemma~\ref{directed set and CP covariance}.
\end{proof}
\end{prop}

\begin{rem}
As in Remark~\ref{are all the hypoth necessary}, it is not clear if all of the hypotheses in Proposition~\ref{sufficient conditions for Psi' to be CNP cov} are necessary. We would like to be able to rerun the argument used in the proof of Proposition~\ref{varphi is CP covariant} (where we did not need $A$ to act compactly on each $\mathbf{Z}_{\left(e_G,q\right)}$, nor for $Q$ to be directed), but $i_{\mathbf{W}^\mathcal{NO}_p}$ need not be Cuntz--Pimsner covariant in general.
\end{rem}

It remains to exhibit a Nica covariant representation of $\mathbf{W}^\mathcal{NO}$ in $\mathcal{NO}_{\mathbf{Y}^\mathcal{NT}}$ that induces a homomorphism $\omega':\mathcal{NT}_{\mathbf{W}^\mathcal{NO}}\rightarrow\mathcal{NO}_{\mathbf{Y}^\mathcal{NT}}$ inverse to $\omega$. Unfortunately, defining this representation is more difficult than when we defined the representation $\Psi$ of $\mathbf{Y}^\mathcal{NT}$ in $\mathcal{NT}_{\mathbf{W}^\mathcal{NO}}$. Whilst each fibre of $\mathbf{W}^\mathcal{NO}$ sits inside the $C^*$-algebra $\mathcal{NO}_\mathbf{Z}$, in general there need not exist a homomorphism from $\mathcal{NO}_\mathbf{Z}$ to $\mathcal{NT}_{\mathbf{W}^\mathcal{NO}}$ which we can just restrict to $\mathbf{W}^\mathcal{NO}$. We get around this difficulty as follows. Firstly, we produce a homomorphism $\Psi'_{e_G}$ from $\mathbf{W}^\mathcal{NO}_{e_G}=\mathcal{NO}_\mathbf{V}$ to $\mathcal{NO}_{\mathbf{Y}^\mathcal{NT}}$ by exhibiting a Cuntz--Nica--Pimsner covariant representation of $\mathbf{V}$ in $\mathcal{NO}_{\mathbf{Y}^\mathcal{NT}}$. Secondly, we use the homomorphism $\Psi'_{e_G}$ to construct a collection of linear maps $\Psi'_p$ from $\mathbf{W}^\mathcal{NO}_p$ to $\mathcal{NO}_{\mathbf{Y}^\mathcal{NT}}$ for each $p\in P\setminus \{e_G\}$, and then argue that this collection forms a Nica covariant representation.

\begin{prop}
\label{defn of vartheta} 
Define $\vartheta':\mathbf{V}\rightarrow \mathcal{NO}_{\mathbf{Y}^\mathcal{NT}}$ by 
\[
\vartheta'_q:=j_{\mathbf{Y}^\mathcal{NT}_q}\circ i_{\mathbf{Z}_{\left(e_G,q\right)}}
\]
for each $q\in Q$. Then $\vartheta'$ is a Cuntz--Nica--Pimsner covariant representation of $\mathbf{V}$, and so there exists a homomorphism $\Psi'_{e_G}:\phi_\mathbf{V}^\mathcal{NO}(\mathcal{NO}_\mathbf{V})\rightarrow\mathcal{NO}_{\mathbf{Y}^\mathcal{NT}}$ such that 
\[
\Psi'_{e_G}\circ j_{\mathbf{Z}_{\left(e_G,q\right)}}=j_{\mathbf{Y}^\mathcal{NT}_q}\circ i_{\mathbf{Z}_{\left(e_G,q\right)}}.
\] 
\begin{proof}
The same reasoning as in the proof of Propositions~\ref{varphi is a representation} and~\ref{varphi is Nica covariant} shows that $\vartheta'$ is a Nica covariant representation of $\mathbf{V}$. We show that $\vartheta'$ is Cuntz--Pimsner covariant. Suppose $F\subseteq Q$ is finite and $\left\{T_q\in \mathcal{K}_A\left(\mathbf{V}_q\right):q\in F\right\}$ is such that $\sum_{q\in F}\iota_q^t(T_q)=0\in \mathcal{L}_A\left(\mathbf{V}_t\right)$ for large $t\in Q$. We need to show that $\sum_{q\in F}\vartheta'^{(q)}(T_q)=0\in\mathcal{NO}_{\mathbf{Y}^\mathcal{NT}}$. We claim that 
\begin{equation}
\label{CNP hypothesis for showing vartheta' is CNP cov}
\begin{aligned}
\sum_{q\in F}\iota_q^t\Big(M_{i_\mathbf{Z}^{\left(e_G,q\right)}(T_q)}\Big)=0\in \mathcal{L}_{\mathcal{NT}_\mathbf{X}}\left(\mathbf{Y}_t^\mathcal{NT}\right) \quad \text{for large $t$}. 
\end{aligned}
\end{equation}
Observe that for any compact operator $T\in \mathcal{K}_A\left(\mathbf{V}_q\right)$, and $z\in \mathbf{Z}_{\left(e_G,q\right)}$, $b\in \mathcal{NT}_\mathbf{X}$, we have
\[
\iota_q^t\Big(M_{i_\mathbf{Z}^{\left(e_G,q\right)}(T)}\Big)\big(i_{\mathbf{Z}_{\left(e_G,t\right)}}(z)\phi_\mathbf{X}^\mathcal{NT}(b)\big)=i_{\mathbf{Z}_{\left(e_G,t\right)}}\big(\iota_q^t(T)(z)\big)\phi_\mathbf{X}^\mathcal{NT}(b).
\]
Given $r\in Q$, fix $s\geq r$, such that $\sum_{q\in F}\iota_q^t(T_q)=0\in \mathcal{L}_A\left(\mathbf{V}_t\right)=\mathcal{L}_A\left(\mathbf{Z}_{\left(e_G,t\right)}\right)$ for all $t\geq s$. Then for any $t\geq s$, we see that
\begin{align*}
&\bigg(\sum_{q\in F}\iota_q^t\Big(M_{i_\mathbf{Z}^{\left(e_G,q\right)}(T_q)}\Big)\bigg)\left(\mathbf{Y}_t^\mathcal{NT}\right)\\
&\qquad=\cspan\bigg\{\bigg(\sum_{q\in F}\iota_q^t\Big(M_{i_\mathbf{Z}^{\left(e_G,q\right)}(T_q)}\Big)\bigg)
\Big(i_{\mathbf{Z}_{\left(e_G,t\right)}}(\mathbf{Z}_{\left(e_G,t\right)})\phi_\mathbf{X}^\mathcal{NT}(\mathcal{NT}_\mathbf{X})\Big)\bigg\}\\
&\qquad=\cspan\bigg\{i_{\mathbf{Z}_{\left(e_G,t\right)}}\Big(\Big(\sum_{q\in F}\iota_q^t\left(T_q\right)\Big)\big(\mathbf{Z}_{\left(e_G,t\right)}\big)\Big)\phi_\mathbf{X}^\mathcal{NT}(\mathcal{NT}_\mathbf{X})\bigg\}
=\{0\},
\end{align*}
and so \eqref{CNP hypothesis for showing vartheta' is CNP cov} holds. Thus, 
\[
\sum_{q\in F}\vartheta'^{(q)}(T_q)=\sum_{q\in F}j_{\mathbf{Y}^\mathcal{NT}}^{(q)}\Big(M_{i_\mathbf{Z}^{(e,q)}(T_q)}\Big)=0,
\]
since $j_{\mathbf{Y}^\mathcal{NT}}$ is Cuntz--Pimsner covariant, and so we conclude that $\vartheta'$ is Cuntz--Pimsner covariant.
\end{proof}
\end{prop}

\begin{prop}
\label{definition of Psi'}
For each $p\in P\setminus \{e_G\}$, there exists a norm-decreasing linear map $\Psi'_p:\mathbf{W}^\mathcal{NO}_p\rightarrow \mathcal{NO}_{\mathbf{Y}^\mathcal{NT}}$ such that 
\begin{align}
\label{first form of Psi'}
\Psi'_p\big(j_{\mathbf{Z}_{\left(p,e_H\right)}}(x)\phi_\mathbf{V}^\mathcal{NO}(b)\big)=j_{\mathbf{Y}_{e_H}^\mathcal{NT}}\big(i_{\mathbf{Z}_{\left(p,e_H\right)}}(x)\big)\Psi'_{e_G}\big(\phi_\mathbf{V}^\mathcal{NO}(b)\big)
\end{align}
for $x\in \mathbf{Z}_{\left(p,e_H\right)}$ and $b\in \mathcal{NO}_\mathbf{V}$. In particular, if $z\in \mathbf{Z}_{(p,q)}$ and $w\in \mathbf{Z}_{\left(e_G,t\right)}$, then 
\begin{align}
\label{useful form of Psi' for calcs}
\Psi'_p\big(j_{\mathbf{Z}_{(p,q)}}(z)j_{\mathbf{Z}_{\left(e_G,t\right)}}(w)^*\big)
=j_{\mathbf{Y}_q^\mathcal{NT}}\big(i_{\mathbf{Z}_{(p,q)}}(z)\big)j_{\mathbf{Y}_t^\mathcal{NT}}\big(i_{\mathbf{Z}_{\left(e_G,t\right)}}(w)\big)^*.
\end{align}
\begin{proof}
We claim that for any finite set $F\subseteq \mathbf{Z}_{\left(p,e_H\right)}\times \mathcal{NO}_{\mathbf{V}}$, 
\begin{equation}
\label{norm estimate for Psi'}
\begin{aligned}
&\Bigg\|\sum_{(x,b)\in F} j_{\mathbf{Y}_{e_H}^\mathcal{NT}}\big(i_{\mathbf{Z}_{\left(p,e_H\right)}}(x)\big)\Psi'_{e_G}\big(\phi_\mathbf{V}^\mathcal{NO}(b)\big)\Bigg\|_{\mathcal{NO}_{\mathbf{Y}^\mathcal{NT}}}\\
&\qquad
\leq \Bigg\|\sum_{(x,b)\in F} j_{\mathbf{Z}_{\left(p,e_H\right)}}(x)\phi_\mathbf{V}^\mathcal{NO}(b)\Bigg\|_{\mathbf{W}^\mathcal{NO}_p}. 
\end{aligned}
\end{equation}
Since $j_{\mathbf{Y}_{e_G}^\mathcal{NT}}$ is a homomorphism and $i_\mathbf{Z}$ is a representation, the $C^*$-identity can be used to show that
\begin{align*}
&\Bigg\|\sum_{(x,b)\in F} j_{\mathbf{Y}_{e_G}^\mathcal{NT}}\big(i_{\mathbf{Z}_{\left(p,e_H\right)}}(x)\big)\Psi'_{e_G}\big(\phi_\mathbf{V}^\mathcal{NO}(b)\big)\Bigg\|_{\mathcal{NO}_{\mathbf{Y}^\mathcal{NT}}}^2\\
&=\Bigg\|\!\sum_{\substack{(x,b),\\(y,c)\in F}} \!\Psi'_{e_G}\big(\phi_\mathbf{V}^\mathcal{NO}(b^*)\big)j_{\mathbf{Y}_{e_H}^\mathcal{NT}}\big(i_{\mathbf{Z}_{\left(e_G,e_H\right)}}\!\big(\langle x,y\rangle_A^{\left(p,e_H\right)}\big)\big)\Psi'_{e_G}\big(\phi_\mathbf{V}^\mathcal{NO}(c)\big)\Bigg\|_{\mathcal{NO}_{\mathbf{Y}^\mathcal{NT}}}\\
&=\Bigg\|\sum_{\substack{(x,b),\\(y,c)\in F}} \Psi'_{e_G}\Big(\phi_\mathbf{V}^\mathcal{NO}(b^*)j_{\mathbf{Z}_{\left(e_G,e_H\right)}}\big(\langle x,y\rangle_A^{\left(p,e_H\right)}\big)\phi_\mathbf{V}^\mathcal{NO}(c)\Big)\Bigg\|_{\mathcal{NO}_{\mathbf{Y}^\mathcal{NT}}},
\end{align*}
where the second equality comes from the fact that $\Psi'_{e_G}\circ j_{\mathbf{Z}_{\left(e_G,e_H\right)}}=\vartheta'_{e_H}=j_{\mathbf{Y}_{e_G}^\mathcal{NT}}\circ i_{\mathbf{Z}_{\left(e_G,e_H\right)}}$. Since $\Psi'_{e_G}$ is a homomorphism, this is no greater than
\begin{align*}
&\Bigg\|\sum_{\substack{(x,b),\\(y,c)\in F}} \phi_\mathbf{V}^\mathcal{NO}(b^*)j_{\mathbf{Z}_{\left(e_G,e_H\right)}}\big(\langle x,y\rangle_A^{\left(p,e_H\right)}\big)\phi_\mathbf{V}^\mathcal{NO}(c)\Bigg\|_{\mathcal{NO}_\mathbf{Z}}\\
&\qquad=\Bigg\|\sum_{\substack{(x,b),\\(y,c)\in F}} \phi_\mathbf{V}^\mathcal{NO}(b^*)j_{\mathbf{Z}_{\left(p,e_H\right)}}(x)^*j_{\mathbf{Z}_{\left(p,e_H\right)}}(y)\phi_\mathbf{V}^\mathcal{NO}(c)\Bigg\|_{\mathcal{NO}_\mathbf{Z}}\\
&\qquad=\Bigg\|\Bigg\langle \sum_{(x,b)\in F}j_{\mathbf{Z}_{\left(p,e_H\right)}}(x)\phi_\mathbf{V}^\mathcal{NO}(b),\sum_{(y,c)\in F}j_{\mathbf{Z}_{\left(p,e_H\right)}}(y)\phi_\mathbf{V}^\mathcal{NO}(c)\Bigg\rangle_{\mathcal{NO}_\mathbf{Z}}^p\Bigg\|_{\mathcal{NO}_\mathbf{Z}}\\
&\qquad=\Bigg\|\sum_{(x,b)\in F}j_{\mathbf{Z}_{\left(p,e_H\right)}}(x)\phi_{\mathbf{V}}^\mathcal{NO}(b)\Bigg\|^2_{{\mathbf{W}}^\mathcal{NO}_p}.
\end{align*}
Thus, \eqref{norm estimate for Psi'} holds. It follows that 
\[
j_{\mathbf{Z}_{\left(p,e_H\right)}}(x)\phi_\mathbf{V}^\mathcal{NO}(b)\mapsto
j_{\mathbf{Y}_{e_H}^\mathcal{NT}}\big(i_{\mathbf{Z}_{\left(p,e_H\right)}}(x)\big)\Psi'_{e_G}\big(\phi_\mathbf{V}^\mathcal{NO}(b)\big)
\]
determines a well-defined norm-decreasing linear map $\Psi'_p$ on $\mathbf{W}^\mathcal{NO}_p$. 
\end{proof}
\end{prop}

For our proofs of Lemma~\ref{Psi' satisfies T2} and Proposition~\ref{Psi' is a Nica cov repn of W} we will make use of the next remark.

\begin{rem}
\label{consequence of Nica cov explicit form}
Suppose that $(K,R)$ is a quasi-lattice ordered group and $\varsigma$ is a Nica covariant representation of a compactly aligned product system $\mathbf{U}$ over $R$ with coefficient algebra $B$. Let $x\in \mathbf{U}_p$ and $y\in \mathbf{U}_q$ for some $p,q\in R$ with $p,q\neq e_K$ and $p\vee q <\infty$, and choose $x'\in \mathbf{U}_p$ and $y'\in \mathbf{U}_q$ so that $x=x'\cdot \langle x',x'\rangle_B^p$ and $y=y'\cdot \langle y',y'\rangle_B^q$. Since $\mathbf{U}$ is compactly aligned and 
\begin{align*}
\mathbf{U}_{p\vee q}
&=\cspan\big\{\sigma\tau:\sigma\in \mathbf{U}_p, \, \tau\in \mathbf{U}_{p^{-1}(p\vee q)}\big\}\\
&=\cspan\big\{\eta\rho:\eta\in \mathbf{U}_q, \, \rho\in \mathbf{U}_{q^{-1}(p\vee q)}\big\},
\end{align*}
there exist $\sigma_{j_i}\in \mathbf{U}_p$, $\tau_{j_i}\in \mathbf{U}_{p^{-1}(p\vee q)}$, $\eta_{j_i}\in \mathbf{U}_q$, $\rho_{j_i}\in \mathbf{U}_{q^{-1}(p\vee q)}$ such that
\[
\iota_p^{p\vee q}\left(\Theta_{x',x'}\right)\iota_q^{p\vee q}\left(\Theta_{y',y'}\right)
=\lim_{i\rightarrow \infty}\sum_{j_i=1}^{k_i}\Theta_{\sigma_{j_i}\tau_{j_i}, \eta_{j_i}\rho_{j_i}}
\in \mathcal{K}_B\left(\mathbf{U}_{p\vee q}\right).
\]
Using relations (T1)--(T3) and the Nica covariance of $\varsigma$, it can be shown that
\[
\varsigma_p(x)^*\varsigma_q(y)=\lim_{i\rightarrow \infty}\sum_{j_i=1}^{k_i} 
\varsigma_{p^{-1}(p\vee q)}\left(\left\langle x',\sigma_{j_i}\right\rangle_B^p\tau_{j_i}\right) 
\varsigma_{q^{-1}(p\vee q)}\left(\left\langle y',\eta_{j_i}\right\rangle_B^q\rho_{j_i}\right)^*.
\]
Moreover, if $z\in \mathbf{U}_r$ and $w\in \mathbf{U}_s$, then
\begin{align*}
&\varsigma_r(z)\varsigma_p(x)^*\varsigma_q(y)\varsigma_s(w)^*\\
&\qquad=\lim_{i\rightarrow \infty}\sum_{j_i=1}^{k_i} 
\varsigma_{rp^{-1}(p\vee q)}\left(z\left\langle x',\sigma_{j_i}\right\rangle_B^p\tau_{j_i}\right) 
\varsigma_{sq^{-1}(p\vee q)}\left(w\left\langle y',\eta_{j_i}\right\rangle_B^q\rho_{j_i}\right)^*.
\end{align*}
\end{rem}

We now prove that the collection of maps $\{\Psi'_p:p\in P\}$ defined in Propositions~\ref{defn of vartheta} and~\ref{definition of Psi'} gives a representation of the product system $\mathbf{W}^\mathcal{NO}$. The proof that $\Psi'$ satisfies relation (T2)  is particularly onerous, so we present it first as a lemma. 

\begin{lem}
\label{Psi' satisfies T2}
The map $\Psi':\mathbf{W}^\mathcal{NO}\rightarrow \mathcal{NO}_{\mathbf{Y}^\mathcal{NT}}$ satisfies relation $\rm{(T2)}$. 
\begin{proof}
We need to show that 
\begin{align}
\label{Psi' satisfies T2 calculation}
\Psi_p'(x)\Psi_r'(y)=\Psi_{pr}'(xy) \quad \text{for any $p,r\in P$ and $x\in \mathbf{W}_p^\mathcal{NO}$, $y\in \mathbf{W}_r^\mathcal{NO}$}.
\end{align}

We begin by showing that \eqref{Psi' satisfies T2 calculation} holds when $r=e_G$. If $z\in \mathbf{Z}_{\left(p,e_H\right)}$ and $b,c\in \mathcal{NO}_\mathbf{V}$, using \eqref{first form of Psi'} and the multiplicativity of $\Psi'_{e_G}$ and $\phi_\mathbf{V}^\mathcal{NO}$, we see that
\begin{align*}
&\Psi'_p\big(j_{\mathbf{Z}_{\left(p,e_H\right)}}(z)\phi_\mathbf{V}^\mathcal{NO}(b)\big)\Psi'_{e_G}\big(\phi_\mathbf{V}^\mathcal{NO}(c)\big)\\
&\qquad=j_{\mathbf{Y}_{e_H}^\mathcal{NT}}\big(i_{\mathbf{Z}_{\left(p,e_H\right)}}(x)\big)\Psi'_{e_G}\big(\phi_\mathbf{V}^\mathcal{NO}(b)\big)\Psi'_{e_G}\big(\phi_\mathbf{V}^\mathcal{NO}(c)\big)\\
&\qquad=\Psi'_p\big(j_{\mathbf{Z}_{\left(p,e_H\right)}}(z)\phi_\mathbf{V}^\mathcal{NO}(b)\phi_\mathbf{V}^\mathcal{NO}(c)\big).
\end{align*}
Since $\mathbf{W}_p^\mathcal{NO}=\cspan\{j_{\mathbf{Z}_{\left(p,e_H\right)}}(z)\phi_\mathbf{V}^\mathcal{NO}(b):z\in \mathbf{Z}_{\left(p,e_H\right)}, b\in \mathcal{NO}_\mathbf{V}\}$ and $\mathbf{W}_{e_G}^\mathcal{NO}$ is isomorphic to $\phi_\mathbf{V}^\mathcal{NO}\left(\mathcal{NO}_\mathbf{V}\right)$, we conclude that \eqref{Psi' satisfies T2 calculation} holds when $r=e_G$.

We now move on to the case where $r\neq e_G$. Since $\Psi_p'$ and $\Psi_r'$ are linear and norm-decreasing, and multiplication in $\mathcal{NO}_{\mathbf{Y}^\mathcal{NT}}$ is linear and continuous, it suffices to prove that \eqref{Psi' satisfies T2 calculation} holds when $x=j_{\mathbf{Z}_{(p,q)}}(z)j_{\mathbf{Z}_{\left(e_G,t\right)}}(w)^*$ and $y=j_{\mathbf{Z}_{(r,m)}}(u)j_{\mathbf{Z}_{\left(e_G,n\right)}}(v)^*$ for $z\in \mathbf{Z}_{(p,q)}$, $w\in \mathbf{Z}_{\left(e_G,t\right)}$, $u\in \mathbf{Z}_{(r,m)}$, and $v\in \mathbf{Z}_{\left(e_G,n\right)}$.  We deal with the cases where $t=e_H$ and $t\neq e_H$ separately. 

If $t=e_H$, then $w\in \mathbf{Z}_{\left(e_G,e_H\right)}=A$, and so an application of \eqref{useful form of Psi' for calcs} gives
\begin{align*}
&\Psi_{pr}'\big(j_{\mathbf{Z}_{(p,q)}}(z)j_{\mathbf{Z}_{\left(e_G,e_H\right)}}(w)^*j_{\mathbf{Z}_{(r,m)}}(u)j_{\mathbf{Z}_{\left(e_G,n\right)}}(v)^*\big)\\
&\qquad=j_{\mathbf{Y}_{qm}^\mathcal{NT}}\big(i_{\mathbf{Z}_{(pr,qm)}}(zw^*u)\big)j_{\mathbf{Y}_n^\mathcal{NT}}\big(i_{\mathbf{Z}_{\left(e_G,n\right)}}(v)\big)^*.
\end{align*}
Another application of \eqref{useful form of Psi' for calcs}, shows that the previous line is equal to
\begin{align*}
&j_{\mathbf{Y}_{qm}^\mathcal{NT}}\big(i_{\mathbf{Z}_{(p,q)}}(z)i_{\mathbf{Z}_{\left(e_G,e_H\right)}}(w)^*i_{\mathbf{Z}_{(r,m)}}(u)\big)j_{\mathbf{Y}_n^\mathcal{NT}}\big(i_{\mathbf{Z}_{\left(e_G,n\right)}}(v)\big)^*\\
&\qquad=j_{\mathbf{Y}_{q}^\mathcal{NT}}\big(i_{\mathbf{Z}_{(p,q)}}(z)\big)
j_{\mathbf{Y}_{e_H}^\mathcal{NT}}\big(i_{\mathbf{Z}_{\left(e_G,e_H\right)}}(w)\big)^*\\
&\qquad\qquad\times
j_{\mathbf{Y}_{m}^\mathcal{NT}}\big(i_{\mathbf{Z}_{(r,m)}}(u)\big)
j_{\mathbf{Y}_n^\mathcal{NT}}\big(i_{\mathbf{Z}_{\left(e_G,n\right)}}(v)\big)^*\\
&\qquad=\Psi'_p\big(j_{\mathbf{Z}_{(p,q)}}(z)j_{\mathbf{Z}_{\left(e_G,e_H\right)}}(w)^*\big)
\Psi'_r\big(j_{\mathbf{Z}_{(r,m)}}(u)j_{\mathbf{Z}_{\left(e_G,n\right)}}(v)^*\big),
\end{align*}
as required.

It remains to deal with the situation where $t\neq e_H$. We make use of Remark~\ref{consequence of Nica cov explicit form} to rewrite the product $j_{\mathbf{Z}_{(p,q)}}(z)j_{\mathbf{Z}_{\left(e_G,t\right)}}(w)^*j_{\mathbf{Z}_{(r,m)}}(u)j_{\mathbf{Z}_{\left(e_G,n\right)}}(v)^*$ in the form required to apply \eqref{useful form of Psi' for calcs}. If $t\vee m=\infty$, then $\left(e_G,t\right)\vee (r,m)=\infty$, and so using the Nica covariance of $j_\mathbf{Z}$ and $j_{\mathbf{Y}^\mathcal{NT}}$ we see that both
\[
\Psi'_{pr}\big(j_{\mathbf{Z}_{(p,q)}}(z)j_{\mathbf{Z}_{\left(e_G,t\right)}}(w)^*j_{\mathbf{Z}_{(r,m)}}(u)j_{\mathbf{Z}_{\left(e_G,n\right)}}(v)^*\big)
\]
and
\begin{align*}
&\Psi'_p\big(j_{\mathbf{Z}_{(p,q)}}(z)j_{\mathbf{Z}_{\left(e_G,t\right)}}(w)^*\big)\Psi'_r\big(j_{\mathbf{Z}_{(r,m)}}(u)j_{\mathbf{Z}_{\left(e_G,n\right)}}(v)^*\big)\\
&\quad=
j_{\mathbf{Y}_q^\mathcal{NT}}\big(i_{\mathbf{Z}_{(p,q)}}(z)\big)j_{\mathbf{Y}_t^\mathcal{NT}}\big(i_{\mathbf{Z}_{\left(e_G,t\right)}}(w)\big)^*
j_{\mathbf{Y}_m^\mathcal{NT}}\big(i_{\mathbf{Z}_{(r,m)}}(u)\big)j_{\mathbf{Y}_n^\mathcal{NT}}\big(i_{\mathbf{Z}_{\left(e_G,n\right)}}(v)\big)^*
\end{align*}
are zero. Thus, we may as well suppose that $t\vee m<\infty$. Choose $w'\in \mathbf{Z}_{\left(e_G,t\right)}$ and $u'\in \mathbf{Z}_{\left(r,m\right)}$ so that $w=w'\cdot \langle w',w'\rangle_A^{\left(e_G,t\right)}$ and $u=u'\cdot \langle u',u'\rangle_A^{\left(r,m\right)}$. Since $\mathbf{Z}$ is compactly aligned, and
\begin{align*}
\mathbf{Z}_{\left(r,t\vee m\right)}
&=\cspan\big\{\sigma\tau:\sigma\in \mathbf{Z}_{\left(e_G,t\right)}, \tau\in \mathbf{Z}_{\left(r,t^{-1}(t\vee m)\right)}\big\}\\
&=\cspan\big\{\eta\rho:\eta\in \mathbf{Z}_{\left(r,m\right)}, \rho\in \mathbf{Z}_{\left(e_G,m^{-1}(t\vee m)\right)}\big\},
\end{align*}
we can write
\begin{align}
\label{compact alignment in Z explicit}
\iota_{\left(e_G,t\right)}^{\left(r,t\vee m\right)}\big(\Theta_{w',w'}\big)\iota_{\left(r,m\right)}^{\left(r,t\vee m\right)}\big(\Theta_{u',u'}\big)=\lim_{i\rightarrow \infty}\sum_{j_i=1}^{k_i}\Theta_{\sigma_{j_i}\tau_{j_i},\eta_{j_i}\rho_{j_i}}\in \mathcal{K}_A\left(\mathbf{Z}_{\left(r,t\vee m\right)}\right),
\end{align}
for some 
$\sigma_{j_i}\in \mathbf{Z}_{\left(e_G,t\right)}$, $\tau_{j_i}\in \mathbf{Z}_{\left(r,t^{-1}(t\vee m)\right)}$, $\eta_{j_i}\in \mathbf{Z}_{\left(r,m\right)}$, $\rho_{j_i}\in \mathbf{Z}_{\left(e_G,m^{-1}(t\vee m)\right)}$.
Combining \eqref{useful form of Psi' for calcs} and Remark~\ref{consequence of Nica cov explicit form}, we have
\begin{equation}
\label{another calculation for showing that Psi' satisfies T2}
\begin{aligned}
&\Psi'_{pr}\big(j_{\mathbf{Z}_{(p,q)}}(z)j_{\mathbf{Z}_{\left(e_G,t\right)}}(w)^*j_{\mathbf{Z}_{(r,m)}}(u)j_{\mathbf{Z}_{\left(e_G,n\right)}}(v)^*\big)\\
&\qquad=
\lim_{i\rightarrow \infty}\sum_{j_i=1}^{k_i}
j_{\mathbf{Y}_{qt^{-1}(t\vee m)}^\mathcal{NT}}\big(i_{\mathbf{Z}_{\left(pr,qt^{-1}(t\vee m)\right)}}\big(z\langle w', \sigma_{j_i}\rangle_A^{\left(e_G,t\right)}\tau_{j_i}\big)\big)\\
&\qquad\qquad\qquad\times
j_{\mathbf{Y}_{nm^{-1}(t\vee m)}^\mathcal{NT}}\big(i_{\mathbf{Z}_{\left(e_G,nm^{-1}(t\vee m)\right)}}\big(v\langle u', \eta_{j_i}\rangle_A^{\left(r,m\right)}\rho_{j_i}\big)\big)^*.
\end{aligned}
\end{equation}
We now calculate $\Psi'_p\big(j_{\mathbf{Z}_{(p,q)}}(z)j_{\mathbf{Z}_{\left(e_G,t\right)}}(w)^*\big)\Psi'_r\big(j_{\mathbf{Z}_{(r,m)}}(u)j_{\mathbf{Z}_{\left(e_G,n\right)}}(v)^*\big)$, which requires applying Remark~\ref{consequence of Nica cov explicit form} to the Nica covariant representation $j_{\mathbf{Y}^\mathcal{NT}}$. 
Using the fact that $w=w'\cdot \langle w',w'\rangle_A^{\left(e_G,t\right)}$ and $u=u'\cdot \langle u',u'\rangle_A^{\left(r,m\right)}$, we see that
$
i_{\mathbf{Z}_{\left(e_G,t\right)}}(w)=i_{\mathbf{Z}_{\left(e_G,t\right)}}(w')\cdot \langle i_{\mathbf{Z}_{\left(e_G,t\right)}}(w'), i_{\mathbf{Z}_{\left(e_G,t\right)}}(w')\rangle_{\mathcal{NT}_\mathbf{X}}^t
$ 
in $\mathbf{Y}_t^\mathcal{NT}$ and
$i_{\mathbf{Z}_{\left(r,m\right)}}(u)=i_{\mathbf{Z}_{\left(r,m\right)}}(u')\cdot \langle i_{\mathbf{Z}_{\left(r,m\right)}}(u'), i_{\mathbf{Z}_{\left(r,m\right)}}(u')\rangle_{\mathcal{NT}_\mathbf{X}}^m
$
in $\mathbf{Y}_m^\mathcal{NT}$. Furthermore, \eqref{compact alignment in Z explicit} implies that
\begin{align*}
&\iota_t^{t\vee m}\big(\Theta_{i_{\mathbf{Z}_{\left(e_G,t\right)}}(w'),i_{\mathbf{Z}_{\left(e_G,t\right)}}(w')}\big)
\iota_m^{t\vee m}\big(\Theta_{i_{\mathbf{Z}_{\left(r,m\right)}}(u'),i_{\mathbf{Z}_{\left(r,m\right)}}(u')}\big)\\
&\qquad=\lim_{i\rightarrow \infty}\sum_{j_i=1}^{k_i}\Theta_{i_{\mathbf{Z}_{\left(e_G,t\right)}}(\sigma_{j_i})i_{\mathbf{Z}_{\left(r,t^{-1}(t\vee m)\right)}}(\tau_{j_i}), i_{\mathbf{Z}_{\left(r,m\right)}}(\eta_{j_i})i_{\mathbf{Z}_{\left(e_G,m^{-1}(t\vee m)\right)}}(\rho_{j_i})}\\
&\qquad\in \mathcal{K}_{\mathcal{NT}_\mathbf{X}}\left(\mathbf{Y}_{t\vee m}^\mathcal{NT}\right),
\end{align*}
whilst 
$i_{\mathbf{Z}_{\left(e_G,t\right)}}(\sigma_{j_i})\in \mathbf{Y}_t^\mathcal{NT}$, 
$i_{\mathbf{Z}_{\left(r,t^{-1}(t\vee m)\right)}}(\tau_{j_i})\in \mathbf{Y}_{t^{-1}(t\vee m)}^\mathcal{NT}$, 
$i_{\mathbf{Z}_{\left(r,m\right)}}(\eta_{j_i})\in \mathbf{Y}_m^\mathcal{NT}$, and
$i_{\mathbf{Z}_{\left(e_G,m^{-1}(t\vee m)\right)}}(\rho_{j_i})\in \mathbf{Y}_{m^{-1}(t\vee m)}^\mathcal{NT}$.
Thus, \eqref{useful form of Psi' for calcs} and Remark~\ref{consequence of Nica cov explicit form} imply that
\begin{align*}
&\Psi'_p\big(j_{\mathbf{Z}_{(p,q)}}(z)j_{\mathbf{Z}_{\left(e_G,t\right)}}(w)^*\big)\Psi'_r\big(j_{\mathbf{Z}_{(r,m)}}(u)j_{\mathbf{Z}_{\left(e_G,n\right)}}(v)^*\big)\\
&\hspace{0.5em}=
j_{\mathbf{Y}_q^\mathcal{NT}}\big(i_{\mathbf{Z}_{(p,q)}}(z)\big)j_{\mathbf{Y}_t^\mathcal{NT}}\big(i_{\mathbf{Z}_{\left(e_G,t\right)}}(w)\big)^*
j_{\mathbf{Y}_m^\mathcal{NT}}\big(i_{\mathbf{Z}_{(r,m)}}(u)\big)j_{\mathbf{Y}_n^\mathcal{NT}}\big(i_{\mathbf{Z}_{\left(e_G,n\right)}}(v)\big)^*\\
&\hspace{0.5em}=
\lim_{i\rightarrow \infty}\sum_{j_i=1}^{k_i}
j_{\mathbf{Y}_{qt^{-1}(t\vee m)}^\mathcal{NT}}\hspace{-0.3em}\big(\hspace{-0.1em}i_{\mathbf{Z}_{(p,q)}}\hspace{-0.2em}(z)\hspace{-0.1em}\langle i_{\mathbf{Z}_{(e_G,t)}}\hspace{-0.2em}(w'),\hspace{-0.1em} i_{\mathbf{Z}_{(e_G,t)}}\hspace{-0.2em}(\sigma_{j_i})\rangle_{\mathcal{NT}_\mathbf{X}}^t \hspace{-0.1em} i_{\mathbf{Z}_{(r,t^{-1}(t\vee m))}}\hspace{-0.2em}(\tau_{j_i})\big)\\
&\qquad\times
j_{\mathbf{Y}_{nm^{-1}(t\vee m)}^\mathcal{NT}}\hspace{-0.3em}\big(\hspace{-0.1em}i_{\mathbf{Z}_{\left(e_G,n\right)}}\hspace{-0.2em}(v)\hspace{-0.1em}\langle i_{\mathbf{Z}_{\left(r,m\right)}}\hspace{-0.2em}(u'), \hspace{-0.1em} i_{\mathbf{Z}_{\left(r,m\right)}}\hspace{-0.2em}(\eta_{j_i})\rangle_{\mathcal{NT}_\mathbf{X}}^m \hspace{-0.1em}i_{\mathbf{Z}_{\left(e_G,m^{-1}(t\vee m)\right)}}\hspace{-0.2em}(\rho_{j_i})\big)^{\hspace{-0.2em}*}\\
&\hspace{0.5em}=
\lim_{i\rightarrow \infty}\sum_{j_i=1}^{k_i}
j_{\mathbf{Y}_{qt^{-1}(t\vee m)}^\mathcal{NT}}\hspace{-0.3em}\big(i_{\mathbf{Z}_{\left(pr,qt^{-1}(t\vee m)\right)}}\big(z\langle w', \sigma_{j_i}\rangle_A^{\left(e_G,t\right)}\tau_{j_i}\big)\big)\\
&\qquad\qquad\qquad\times
j_{\mathbf{Y}_{nm^{-1}(t\vee m)}^\mathcal{NT}}\hspace{-0.3em}\big(i_{\mathbf{Z}_{\left(e_G,nm^{-1}(t\vee m)\right)}}\big(v\langle u', \eta_{j_i}\rangle_A^{\left(r,m\right)}\rho_{j_i}\big)\big)^*,
\end{align*}
which is \eqref{another calculation for showing that Psi' satisfies T2}. Thus, 
\begin{align*}
&\Psi'_{pr}\big(j_{\mathbf{Z}_{(p,q)}}(z)j_{\mathbf{Z}_{\left(e_G,t\right)}}(w)^*j_{\mathbf{Z}_{(r,m)}}(u)j_{\mathbf{Z}_{\left(e_G,n\right)}}(v)^*\big)\\
&\qquad=\Psi'_p\big(j_{\mathbf{Z}_{(p,q)}}(z)j_{\mathbf{Z}_{\left(e_G,t\right)}}(w)^*\big)\Psi'_r\big(j_{\mathbf{Z}_{(r,m)}}(u)j_{\mathbf{Z}_{\left(e_G,n\right)}}(v)^*\big)
\end{align*}
when $t\neq e_H$. We conclude that the map $\Psi':\mathbf{W}^\mathcal{NO}\rightarrow \mathcal{NO}_{\mathbf{Y}^\mathcal{NT}}$ satisfies relation (T2).
\end{proof}
\end{lem}

\begin{prop}
\label{Psi' is repn of W}
The map $\Psi':\mathbf{W}^\mathcal{NO}\rightarrow \mathcal{NO}_{\mathbf{Y}^\mathcal{NT}}$ is a representation of $\mathbf{W}^\mathcal{NO}$. 
\begin{proof}
By construction, each $\Psi_p'$ is linear and $\Psi_{e_G}'$ is a homomorphism. Hence, $\Psi'$ satisfies (T1). We already showed that $\Psi'$ satisfies (T2) in Lemma~\ref{Psi' satisfies T2}. It remains to show that $\Psi'$ satisfies (T3).

Fix $z,w\in \mathbf{Z}_{\left(p,e_H\right)}$ and $b,c\in \mathcal{NO}_\mathbf{V}$. Making use of \eqref{first form of Psi'}, we see that
\begin{align*}
\Psi'_p
&
\big(j_{\mathbf{Z}_{\left(p,e_H\right)}}(z)\phi_\mathbf{V}^\mathcal{NO}(b)\big)^*\Psi'_p\big(j_{\mathbf{Z}_{\left(p,e_H\right)}}(w)\phi_\mathbf{V}^\mathcal{NO}(c)\big)\\
&
=\Big(j_{\mathbf{Y}_{e_H}^\mathcal{NT}}\big(i_{\mathbf{Z}_{\left(p,e_H\right)}}(z)\big)\Psi'_{e_G}\big(\phi_\mathbf{V}^\mathcal{NO}(b)\big)\Big)^*
\Big(j_{\mathbf{Y}_{e_H}^\mathcal{NT}}\big(i_{\mathbf{Z}_{\left(p,e_H\right)}}(w)\big)\Psi'_{e_G}\big(\phi_\mathbf{V}^\mathcal{NO}(c)\big)\Big)\\
&
=\Psi'_{e_G}\big(\phi_\mathbf{V}^\mathcal{NO}(b)^*\big)j_{\mathbf{Y}_{e_H}^\mathcal{NT}}\big(i_{\mathbf{Z}_{\left(e_G,e_H\right)}}\big(\langle z, w\rangle_A^{\left(p,e_H\right)}\big)\big)\Psi'_{e_G}\big(\phi_\mathbf{V}^\mathcal{NO}(c)\big).
\end{align*}
Since $j_\mathbf{Z}$ is a representation of $\mathbf{Z}$ and $\Psi'_{e_G}$ is a homomorphism satisfying $j_{\mathbf{Y}_{e_H}^\mathcal{NT}}\circ i_{\mathbf{Z}_{\left(e_G,e_H\right)}}=\Psi'_{e_G}\circ j_{\mathbf{Z}_{\left(e_G,e_H\right)}}$, the previous expression is equal to 
\begin{align*}
&
\Psi'_{e_G}\big(\phi_\mathbf{V}^\mathcal{NO}(b)^*j_{\mathbf{Z}_{\left(e_G,e_H\right)}}\big(\langle z, w\rangle_A^{\left(p,e_H\right)}\big)\phi_\mathbf{V}^\mathcal{NO}(c)\big)\\
&
\qquad=\Psi'_{e_G}\big(\phi_\mathbf{V}^\mathcal{NO}(b)^*j_{\mathbf{Z}_{\left(p,e_H\right)}}(z)^*j_{\mathbf{Z}_{\left(p,e_H\right)}}(w)\phi_\mathbf{V}^\mathcal{NO}(c)\big)\\
&
\qquad=\Psi'_{e_G}\big(\big\langle j_{\mathbf{Z}_{\left(p,e_H\right)}}(z)\phi_\mathbf{V}^\mathcal{NO}(b),j_{\mathbf{Z}_{\left(p,e_H\right)}}(w)\phi_\mathbf{V}^\mathcal{NO}(c)\big\rangle_{\mathcal{NO}_\mathbf{V}}^p\big).
\end{align*}
Since $\mathbf{W}^\mathcal{NO}_p=\cspan\big\{j_{\mathbf{Z}_{\left(p,e_H\right)}}\left(\mathbf{Z}_{\left(p,e_H\right)}\right)\phi_\mathbf{V}^\mathcal{NO}\left( \mathcal{NO}_\mathbf{V}\right)\big\}$ for each $p\in P$, $\Psi_p'$ is linear and norm-decreasing, and multiplication in $\mathcal{NO}_{\mathbf{Y}^\mathcal{NT}}$ is linear and continuous, we conclude that 
$\Psi'_p(x)^*\Psi'_p(y)=\Psi_{e_G}'\big(\langle x,y\rangle_{\mathcal{NO}_\mathbf{V}}^p\big)$ for each $x,y\in \mathbf{W}^\mathcal{NO}_p$. Thus, $\Psi'$ satisfies (T3), and we conclude that $\Psi'$ is a representation of $\mathbf{W}^\mathcal{NO}$
\end{proof}
\end{prop}

We can also show that the representation $\Psi'$ is Nica covariant. 

\begin{prop}
\label{Psi' is a Nica cov repn of W}
The representation $\Psi':\mathbf{W}^\mathcal{NO}\rightarrow \mathcal{NO}_{\mathbf{Y}^\mathcal{NT}}$ is Nica covariant, and so there exists a homomorphism $\omega':\mathcal{NT}_{\mathbf{W}^\mathcal{NO}}\rightarrow \mathcal{NO}_{\mathbf{Y}^\mathcal{NT}}$ such that 
\[
\omega'\circ i_{\mathbf{W}^\mathcal{NO}_p}=\Psi'_p
\]
for each $p\in P$. 
\begin{proof}
We need to show that if $S\in \mathcal{K}_{\mathcal{NO}_\mathbf{V}}(\mathbf{W}^\mathcal{NO}_p)$ and $T\in \mathcal{K}_{\mathcal{NO}_\mathbf{V}}(\mathbf{W}^\mathcal{NO}_r)$, then 
$\Psi'^{(p)}(S)\Psi'^{(r)}(T)$ is equal to $\Psi'^{(p\vee r)}(\iota_p^{p\vee r}(S)\iota_r^{p\vee r}(T))$ when $p\vee r<\infty$ and is zero otherwise. If $p=e_G$ or $r=e_G$, the result is trivial, so we suppose that $p,r\neq e_G$. Hence, by Lemma~\ref{compact operators in Y}, it suffices to show that
\begin{equation}
\label{what we need to show for Nica covariance}
\begin{aligned}
&\Psi'^{(p)}
\big(\Theta_{j_{\mathbf{Z}_{(p,q)}}(z),j_{\mathbf{Z}_{(p,t)}}(w)}\big)\Psi'^{(r)}\big(\Theta_{j_{\mathbf{Z}_{(r,m)}}(u),j_{\mathbf{Z}_{(r,n)}}(v)}\big)\\
&\qquad=
\begin{cases}
\Psi'^{(p\vee r)}\big(M_{j_{\mathbf{Z}_{(p,q)}}(z)j_{\mathbf{Z}_{(p,t)}}(w)^*j_{\mathbf{Z}_{(r,m)}}(u)j_{\mathbf{Z}_{(r,n)}}(v)^*}\big) & \text{if $p\vee r<\infty$}\\
0 &\text{otherwise.}
\end{cases}
\end{aligned}
\end{equation}
whenever $z\in \mathbf{Z}_{(p,q)}$, $w\in \mathbf{Z}_{(p,t)}$, $u\in \mathbf{Z}_{(r,m)}$, $v\in \mathbf{Z}_{(r,n)}$. 

We begin by showing that \eqref{what we need to show for Nica covariance} holds when $p\vee r=\infty$ or $t\vee m=\infty$. Since $j_{\mathbf{Z}_{(p,q)}}(z)j_{\mathbf{Z}_{(p,t)}}(w)^*j_{\mathbf{Z}_{(r,m)}}(u)j_{\mathbf{Z}_{(r,n)}}(v)^*=0$ if $t\vee m=\infty$, it suffices to show that
\begin{equation}
\label{first nica cov claim}
\begin{aligned}
\Psi'^{(p)}\big(\Theta_{j_{\mathbf{Z}_{(p,q)}}(z),j_{\mathbf{Z}_{(p,t)}}(w)}\big)\Psi'^{(r)}&\big(\Theta_{j_{\mathbf{Z}_{(r,m)}}(u),j_{\mathbf{Z}_{(r,n)}}(v)}\big)=0\\
&\text{if $p\vee r=\infty$ or $t\vee m=\infty$.}
\end{aligned}
\end{equation}
Observe that
\begin{equation}
\label{Nica cov calc 1}
\begin{aligned}
&\Psi'^{(p)}\big(\Theta_{j_{\mathbf{Z}_{(p,q)}}(z),j_{\mathbf{Z}_{(p,t)}}(w)}\big)\Psi'^{(r)}\big(\Theta_{j_{\mathbf{Z}_{(r,m)}}(u),j_{\mathbf{Z}_{(r,n)}}(v)}\big)\\
&\qquad=j_{\mathbf{Y}_q^\mathcal{NT}}\big(i_{\mathbf{Z}_{(p,q)}}(z)\big)j_{\mathbf{Y}_t^\mathcal{NT}}\big(i_{\mathbf{Z}_{(p,t)}}(w)\big)^*\\
&\qquad\qquad\times
j_{\mathbf{Y}_m^\mathcal{NT}}\big(i_{\mathbf{Z}_{(r,m)}}(u)\big)j_{\mathbf{Y}_n^\mathcal{NT}}\big(i_{\mathbf{Z}_{(r,n)}}(v)\big)^*.
\end{aligned}
\end{equation}
If $t\vee m=\infty$ then \eqref{Nica cov calc 1} is zero since $j_{\mathbf{Y}^\mathcal{NT}}$ is Nica covariant. Hence, we need to consider the situation where $t\vee m<\infty$ and $p\vee r=\infty$. Choose $w'\in \mathbf{Z}_{(p,t)}$ and $u'\in \mathbf{Z}_{(r,m)}$ so that $w=w'\cdot\langle w',w'\rangle_A^{(p,t)}$ and $u=u'\cdot\langle u',u'\rangle_A^{(r,m)}$. The Nica covariance of $j_{\mathbf{Y}^\mathcal{NT}}$ gives
\begin{align*}
&j_{\mathbf{Y}_t^\mathcal{NT}}\big(i_{\mathbf{Z}_{(p,t)}}(w)\big)^*j_{\mathbf{Y}_m^\mathcal{NT}}\big(i_{\mathbf{Z}_{(r,m)}}(u)\big)\\
&\qquad=j_{\mathbf{Y}_t^\mathcal{NT}}\big(i_{\mathbf{Z}_{(p,t)}}(w')\big)^*
j_{\mathbf{Y}^\mathcal{NT}}^{(t\vee m)}
\Big(
M_{i_\mathbf{Z}^{((p,t))}(\Theta_{w',w'})i_\mathbf{Z}^{((r,m))}(\Theta_{u',u'})}
\Big)\\
&\qquad \qquad\times
j_{\mathbf{Y}_m^\mathcal{NT}}\big(i_{\mathbf{Z}_{(r,m)}}(u')\big)^*,
\end{align*}
which is zero since $i_\mathbf{Z}$ is Nica covariant and $(p,t)\vee (r,m)=\infty$ (as $p\vee r=\infty$). Thus, \eqref{Nica cov calc 1} is zero, which completes the proof of \eqref{first nica cov claim}. 

It remains to prove that \eqref{what we need to show for Nica covariance} holds whenever $p\vee r<\infty$ and $t\vee m<\infty$. As in Lemma~\ref{Psi' satisfies T2}, we will need to make use of Remark~\ref{consequence of Nica cov explicit form} to rewrite things in a form that allows us to apply the description of $\Psi'$ given by \eqref{useful form of Psi' for calcs}. Again, choose $w'\in \mathbf{Z}_{(p,t)}$ and $u'\in \mathbf{Z}_{(r,m)}$ so that $w=w'\cdot\langle w',w'\rangle_A^{(p,t)}$ and $u=u'\cdot\langle u',u'\rangle_A^{(r,m)}$. Writing
\begin{align*}
&\iota_{(p,t)}^{(p\vee r,t\vee m)}\big(\Theta_{w',w'}\big)\iota_{(r,m)}^{(p\vee r,t\vee m)}\big(\Theta_{u',u'}\big)\\
&\qquad
=\lim_{i\rightarrow \infty}\sum_{j_i=1}^{k_i}\Theta_{\sigma_{j_i}\tau_{j_i},\eta_{j_i}\rho_{j_i}}\in \mathcal{K}_A\big(\mathbf{Z}_{(p\vee r,t\vee m)}\big),
\end{align*}
for some choice of
$\sigma_{j_i}\in \mathbf{Z}_{\left(p,t\right)}$, 
$\tau_{j_i}\in \mathbf{Z}_{\left(p^{-1}(p\vee r), t^{-1}(t\vee m)\right)}$, 
$\eta_{j_i}\in \mathbf{Z}_{\left(r,m\right)}$, and 
$\rho_{j_i}\in \mathbf{Z}_{\left(r^{-1}(p\vee r),m^{-1}(t\vee m)\right)}$,
Remark~\ref{consequence of Nica cov explicit form} tells us that
\begin{align*}
&j_{\mathbf{Z}_{(p,q)}}(z)j_{\mathbf{Z}_{(p,t)}}(w)^*j_{\mathbf{Z}_{(r,m)}}(u)j_{\mathbf{Z}_{(r,n)}}(v)^*\\
&\qquad=\lim_{i\rightarrow \infty}\sum_{j_i=1}^{k_i}
j_{\mathbf{Z}_{\left(p\vee r,qt^{-1}(t\vee m)\right)}}\big(z\langle w', \sigma_{j_i}\rangle_A^{\left(p,t\right)}\tau_{j_i}\big)\\
&\qquad\qquad\qquad\times
j_{\mathbf{Z}_{\left(p\vee r,nm^{-1}(t\vee m)\right)}}\big(v\langle u', \eta_{j_i}\rangle_A^{\left(r,m\right)}\rho_{j_i}\big)^*.
\end{align*}
Thus,
\begin{equation}
\label{nearly done}
\begin{aligned}
&\Psi'^{(p\vee r)}\big(M_{j_{\mathbf{Z}_{(p,q)}}(z)j_{\mathbf{Z}_{(p,t)}}(w)^*j_{\mathbf{Z}_{(r,m)}}(u)j_{\mathbf{Z}_{(r,n)}}(v)^*}\big)\\
&\qquad=
\lim_{i\rightarrow \infty}\sum_{j_i=1}^{k_i}
\Psi'^{(p\vee r)}\Big(
\Theta_{j_{\mathbf{Z}}(z\langle w', \sigma_{j_i}\rangle_A^{(p,t)}\tau_{j_i}), j_{\mathbf{Z}}(v\langle u', \eta_{j_i}\rangle_A^{(r,m)}\rho_{j_i})}
\Big)\\
&\qquad=
\lim_{i\rightarrow \infty}\sum_{j_i=1}^{k_i}
j_{\mathbf{Y}_{qt^{-1}(t\vee m)}^\mathcal{NT}}\big(i_{\mathbf{Z}_{\left(p\vee r,qt^{-1}(t\vee m)\right)}}\big(z\langle w', \sigma_{j_i}\rangle_A^{\left(p,t\right)}\tau_{j_i}\big)\big)\\
&\quad\qquad\qquad\qquad\times
j_{\mathbf{Y}_{nm^{-1}(t\vee m)}^\mathcal{NT}}\big(i_{\mathbf{Z}_{\left(p\vee r,nm^{-1}(t\vee m)\right)}}\big(v\langle u', \eta_{j_i}\rangle_A^{\left(r,m\right)}\rho_{j_i}\big)\big)^*.
\end{aligned}
\end{equation}
To complete the proof we will apply Remark~\ref{consequence of Nica cov explicit form} to the product system $\mathbf{Y}^\mathcal{NT}$ in order to show that \eqref{nearly done} is equal to the left hand side of \eqref{what we need to show for Nica covariance}. Since $w=w'\cdot\langle w',w'\rangle_A^{(p,t)}$, the Hewitt--Cohen--Blanchard factorisation of $i_{\mathbf{Z}_{\left(p,t\right)}}(w)$ in $\mathbf{Y}_t^\mathcal{NT}$ is 
$
i_{\mathbf{Z}_{\left(p,t\right)}}(w')\cdot \langle i_{\mathbf{Z}_{\left(p,t\right)}}(w'), i_{\mathbf{Z}_{\left(p,t\right)}}(w')\rangle_{\mathcal{NT}_\mathbf{X}}^t. 
$
Similarly, $i_{\mathbf{Z}_{\left(r,m\right)}}(u)$ is equal to $i_{\mathbf{Z}_{\left(r,m\right)}}(u')\cdot \langle i_{\mathbf{Z}_{\left(r,m\right)}}(u'), i_{\mathbf{Z}_{\left(r,m\right)}}(u')\rangle_{\mathcal{NT}_\mathbf{X}}^m$. 
Moreover, 
\begin{align*}
\iota_t^{t\vee m}&\big(\Theta_{i_{\mathbf{Z}_{\left(p,t\right)}}(w'),i_{\mathbf{Z}_{\left(p,t\right)}}(w')}\big)
\iota_m^{t\vee m}\big(\Theta_{i_{\mathbf{Z}_{\left(r,m\right)}}(u'),i_{\mathbf{Z}_{\left(r,m\right)}}(u')}\big)\\
&=
\lim_{i\rightarrow \infty}\sum_{j_i=1}^{k_i}\Theta_{i_{\mathbf{Z}}\left(\sigma_{j_i}\right)i_{\mathbf{Z}}\left(\tau_{j_i}\right), i_{\mathbf{Z}}\left(\eta_{j_i}\right)i_{\mathbf{Z}}\left(\rho_{j_i}\right)},
\end{align*}
whilst 
$i_{\mathbf{Z}_{\left(p,t\right)}}\left(\sigma_{j_i}\right)\in \mathbf{Y}_t^\mathcal{NT}$, 
$i_{\mathbf{Z}_{\left(p^{-1}(p\vee r),t^{-1}(t\vee m)\right)}}\left(\tau_{j_i}\right)\in \mathbf{Y}_{t^{-1}(t\vee m)}^\mathcal{NT}$, 
$i_{\mathbf{Z}_{\left(r,m\right)}}\left(\eta_{j_i}\right)\in \mathbf{Y}_m^\mathcal{NT}$, and
$i_{\mathbf{Z}_{\left(r^{-1}(p\vee r),m^{-1}(t\vee m)\right)}}\left(\rho_{j_i}\right)\in \mathbf{Y}_{m^{-1}(t\vee m)}^\mathcal{NT}$.
Hence, 
\begin{align*}
&\Psi'^{(p)}\Big(\Theta_{j_{\mathbf{Z}_{(p,q)}}(z),j_{\mathbf{Z}_{(p,t)}}(w)}\Big)\Psi'^{(r)}\Big(\Theta_{j_{\mathbf{Z}_{(r,m)}}(u),j_{\mathbf{Z}_{(r,n)}}(v)}\Big)\\
&\qquad=
j_{\mathbf{Y}_q^\mathcal{NT}}\big(i_{\mathbf{Z}_{(p,q)}}(z)\big)j_{\mathbf{Y}_t^\mathcal{NT}}\big(i_{\mathbf{Z}_{(p,t)}}(w)\big)^*
j_{\mathbf{Y}_m^\mathcal{NT}}\big(i_{\mathbf{Z}_{(r,m)}}(u)\big)j_{\mathbf{Y}_n^\mathcal{NT}}\big(i_{\mathbf{Z}_{(r,n)}}(v)\big)^*\\
&\qquad=
\lim_{i\rightarrow \infty}\sum_{j_i=1}^{k_i}
j_{\mathbf{Y}_{qt^{-1}(t\vee m)}^\mathcal{NT}}\big(i_{\mathbf{Z}_{\left(p\vee r,qt^{-1}(t\vee m)\right)}}\big(z\langle w', \sigma_{j_i}\rangle_A^{\left(p,t\right)}\tau_{j_i}\big)\big)\\
&\quad\qquad\qquad\qquad\times
j_{\mathbf{Y}_{nm^{-1}(t\vee m)}^\mathcal{NT}}\big(i_{\mathbf{Z}_{\left(p \vee r,nm^{-1}(t\vee m)\right)}}\big(v\langle u', \eta_{j_i}\rangle_A^{\left(r,m\right)}\rho_{j_i}\big)\big)^*,
\end{align*}
which is precisely \eqref{nearly done}. This completes the proof that $\Psi'$ is Nica covariant. 
\end{proof}
\end{prop}

Finally, we are ready to prove that  $\mathcal{NO}_{\mathbf{Y}^\mathcal{NT}}$ and $\mathcal{NT}_{\mathbf{W}^\mathcal{NO}}$ are isomorphic. 

\begin{thm}
\label{isom of mixed case algebras}
Suppose that $A$ acts compactly on each $\mathbf{Z}_{\left(e_G,q\right)}$ and $Q$ is directed, so that the homomorphism $\omega$ of Proposition~\ref{sufficient conditions for Psi' to be CNP cov} exists. Then $\omega:\mathcal{NO}_{\mathbf{Y}^\mathcal{NT}}\rightarrow\mathcal{NT}_{\mathbf{W}^\mathcal{NO}}$ and $\omega':\mathcal{NT}_{\mathbf{W}^\mathcal{NO}}\rightarrow \mathcal{NO}_{\mathbf{Y}^\mathcal{NT}}$ are mutually inverse isomorphisms. Thus, $\mathcal{NO}_{\mathbf{Y}^\mathcal{NT}}\cong\mathcal{NT}_{\mathbf{W}^\mathcal{NO}}$. 
\begin{proof}
Firstly, we show that $\omega\circ \omega'=\mathrm{id}_{\mathcal{NT}_{\mathbf{W}^\mathcal{NO}}}$. As $\mathcal{NT}_{\mathbf{W}^\mathcal{NO}}$ is generated by $i_{\mathbf{W}^\mathcal{NO}}$, it suffices to show that $\omega\circ \omega'\circ i_{\mathbf{W}^\mathcal{NO}}=i_{\mathbf{W}^\mathcal{NO}}$. If $z\in \mathbf{Z}_{(p,q)}$ and $w\in \mathbf{Z}_{\left(e_G,t\right)}$, then
\begingroup
\allowdisplaybreaks
\begin{align*}
\big(\omega\circ \omega'\big)&\Big(i_{\mathbf{W}_p^\mathcal{NO}}\big(j_{\mathbf{Z}_{(p,q)}}(z)j_{\mathbf{Z}_{\left(e_G,t\right)}}(w)^*\big)\Big)\\
&=\big(\omega\circ \Psi_p'\big)\big(j_{\mathbf{Z}_{(p,q)}}(z)j_{\mathbf{Z}_{\left(e_G,t\right)}}(w)^*\big)\\
&=\omega\Big(j_{\mathbf{Y}_q^\mathcal{NT}}\big(i_{\mathbf{Z}_{(p,q)}}(z)\big)j_{\mathbf{Y}_t^\mathcal{NT}}\big(i_{\mathbf{Z}_{\left(e_G,t\right)}}(w)^*\big)\Big)\\
&=\Xi\big(i_{\mathbf{Z}_{(p,q)}}(z)i_{\mathbf{Z}_{\left(e_G,t\right)}}(w)^*\big)\\
&= i_{\mathbf{W}_p^\mathcal{NO}}\big(j_{\mathbf{Z}_{(p,q)}}(z)\big) i_{\mathbf{W}_{e_G}^\mathcal{NO}}\big(j_{\mathbf{Z}_{\left(e_G,t\right)}}(w)\big)^*\\
&= i_{\mathbf{W}_p^\mathcal{NO}}\big(j_{\mathbf{Z}_{(p,q)}}(z)j_{\mathbf{Z}_{\left(e_G,t\right)}}(w)^*\big).
\end{align*}
\endgroup
Since $\mathbf{W}^\mathcal{NO}_p=\cspan\big\{j_{\mathbf{Z}_{(p,q)}}\left(\mathbf{Z}_{(p,q)}\right)j_{\mathbf{Z}_{\left(e_G,t\right)}}\left(\mathbf{Z}_{\left(e_G,t\right)}\right)^*:q,t\in Q\big\}$ for each $p\in P$, whilst both of the maps $\omega\circ \omega'\circ i_{\mathbf{W}^\mathcal{NO}}$ and $i_{\mathbf{W}^\mathcal{NO}}$ are linear and continuous, we conclude that $\omega\circ \omega' \circ i_{\mathbf{W}^\mathcal{NO}_p}= i_{\mathbf{W}^\mathcal{NO}_p}$ for each $p\in P$. Thus, $\omega\circ \omega'=\mathrm{id}_{\mathcal{NT}_{\mathbf{W}^\mathcal{NO}}}$.

Secondly, we check that $\omega'\circ \omega=\mathrm{id}_{\mathcal{NO}_{\mathbf{Y}^\mathcal{NT}}}$. As $\mathcal{NO}_{\mathbf{Y}^\mathcal{NT}}$ is generated by $j_{\mathbf{Y}^\mathcal{NT}}$, it suffices to check that $\omega'\circ \omega\circ j_{\mathbf{Y}^\mathcal{NT}}=j_{\mathbf{Y}^\mathcal{NT}}$. For any $z\in \mathbf{Z}_{(p,q)}$ and $w\in \mathbf{Z}_{\left(r,e_H\right)}$, we have
\begin{align*}
\big(\omega'\circ \omega\big)&\Big(j_{\mathbf{Y}_q^\mathcal{NT}}\big(i_{\mathbf{Z}_{(p,q)}}(z)i_{\mathbf{Z}_{\left(r,e_H\right)}}(w)^*\big)\Big)\\
&=\big(\omega'\circ \Xi'\big)\big(i_{\mathbf{Z}_{(p,q)}}(z)i_{\mathbf{Z}_{\left(r,e_H\right)}}(w)^*\big)\\
&=\omega'\Big(i_{\mathbf{W}_p^\mathcal{NO}}\big(j_{\mathbf{Z}_{(p,q)}}(z)\big)i_{\mathbf{W}_r^\mathcal{NO}}\big(j_{\mathbf{Z}_{\left(r,e_H\right)}}(w)\big)^*\Big)\\
&=\Psi_p'\big(j_{\mathbf{Z}_{(p,q)}}(z)\big)\Psi_r'\big(j_{\mathbf{Z}_{\left(r,e_H\right)}}(w)\big)^*\\
&=j_{\mathbf{Y}_q^\mathcal{NT}}\big(i_{\mathbf{Z}_{(p,q)}}(z)\big)j_{\mathbf{Y}_{e_H}^\mathcal{NT}}\big(i_{\mathbf{Z}_{\left(r,e_H\right)}}(w)\big)^*\\
&=j_{\mathbf{Y}_q^\mathcal{NT}}\big(i_{\mathbf{Z}_{(p,q)}}(z)i_{\mathbf{Z}_{\left(r,e_H\right)}}(w)\big)^*.
\end{align*}
Since $\mathbf{Y}^\mathcal{NT}_q=\cspan\big\{i_{\mathbf{Z}_{(p,q)}}\left(\mathbf{Z}_{(p,q)}\right)i_{\mathbf{Z}_{\left(r,e_H\right)}}\left(\mathbf{Z}_{\left(r,e_H\right)}\right)^*:p,q\in P\big\}$ for each $q\in Q$, whilst both of the maps $\omega'\circ \omega\circ j_{\mathbf{Y}^\mathcal{NT}}$ and $j_{\mathbf{Y}^\mathcal{NT}}$ are linear and continuous, we conclude that $\omega'\circ \omega \circ j_{\mathbf{Y}^\mathcal{NT}_q}= j_{\mathbf{Y}^\mathcal{NT}_q}$ for each $q\in Q$. Thus, $\omega'\circ \omega=\mathrm{id}_{\mathcal{NO}_{\mathbf{Y}^\mathcal{NT}}}$.
\end{proof}
\end{thm}

\section*{Acknowledgements}
The results in this article are from my PhD thesis. Thank you to my supervisors Adam Rennie and Aidan Sims at the University of Wollongong for their advice and encouragement during my PhD and during the writing of this article.

\end{document}